\documentclass[12pt,a4paper]{article}
\usepackage{amsmath,mathrsfs,amssymb,amsthm,amscd}
\usepackage[]{fontenc}
\usepackage{xy} 
\usepackage{enumerate}
\xyoption{all}

\numberwithin{equation}{section}
\newcommand{\nnpar}[1]{\bigbreak\noindent {\scshape #1}}\nopagebreak%
\theoremstyle{plain}
\newtheorem{theorem}[equation]{Theorem}
\newtheorem{corollary}[equation]{Corollary}
\newtheorem{proposition}[equation]{Proposition}
\newtheorem{lemma}[equation]{Lemma}

\theoremstyle{definition}
\newtheorem{definition}[equation]{Definition}
\newtheorem{notation}[equation]{Notation}
\newtheorem{example}[equation]{Example}

\newtheorem{remark}[equation]{Remark}

\DeclareMathOperator{\Gr}{Gr}
\DeclareMathOperator{\PV}{PV}

\DeclareMathOperator{\dv}{div}
\DeclareMathOperator{\diva}{\widehat{div}}
\DeclareMathOperator{\dega}{\widehat{deg}}
\DeclareMathOperator{\cha}{\widehat{CH}}
\DeclareMathOperator{\ca}{\widehat{c}}
\DeclareMathOperator{\CH}{CH}
\DeclareMathOperator{\Cl}{Cl}
\DeclareMathOperator{\za}{\widehat{Z}}
\DeclareMathOperator{\rata}{\widehat{Rat}}
\DeclareMathOperator{\alg}{alg}
\DeclareMathOperator{\pica}{\widehat{Pic}}

\DeclareMathOperator{\dd}{d}
\DeclareMathOperator{\Img}{Im}

\DeclareMathOperator{\zar}{zar}
\DeclareMathOperator{\ZAR}{ZAR}
\DeclareMathOperator{\ch}{ch}
\DeclareMathOperator{\cl}{cl}
\DeclareMathOperator{\rk}{rk}
\DeclareMathOperator{\sk}{sk}

\DeclareMathOperator{\BGL}{BGL}

\DeclareMathOperator{\Ker}{Ker}
\DeclareMathOperator{\Hom}{Hom}
\DeclareMathOperator{\Ho}{Ho}

\DeclareMathOperator{\Id}{id}
\DeclareMathOperator{\Pic}{Pic}

\DeclareMathOperator{\Spec}{Spec}

\DeclareMathOperator{\supp}{supp}
\DeclareMathOperator{\codim}{codim}
\DeclareMathOperator{\amap}{a}
\DeclareMathOperator{\cone}{cone}
\DeclareMathOperator{\Res}{Res}
\DeclareMathOperator{\bmap}{b}

\newcommand{\an}{\text{{\rm an}}}

\newcommand{\hol}{\text{{\rm hol}}}
\newcommand{\wlg}{\text{{\rm pre}}}
\newcommand{\w}{\text{{\rm pre}}}

\newcommand{\verti}{\text{{\rm vert}}}

\newcommand{\fin}{\text{{\rm fin}}}

\newcommand{\pure}{\text{{\rm pure}}}

\newcommand{\Zar}{\text{{\rm Zar}}}

\newcommand{\QQ}{{\mathbb Q}}

\newcommand{\KK}{{\mathbf K}}
\newcommand{\Gi}{{\mathcal {G}}}

\newcommand{\D}{\text{{\rm cur}}}
\newcommand{\cc}{{\mathcal{C}}}

\def\?{\ ???\ \immediate\write16{}%
\immediate\write16{Warning: There was still a question mark . . . }%
\immediate\write16{}}

\begin{document}
\setcounter{tocdepth}{2}
\setcounter{section}{-1}
\begin{titlepage}

\begin{center}

\vspace*{4cm}
{\Huge\bf Cohomological arithmetic \\ Chow rings}

\vspace*{2cm}
{\LARGE J. I. Burgos Gil{\footnotemark
\footnotetext{Partially supported by Grants DGI BFM2000-0799-C02-01
  and BFM2003-02914}
}, J. Kramer, U. K\"uhn}

\end{center}

\end{titlepage}

\thispagestyle{empty}

\newpage
 \phantom{aaa}
\begin{abstract}
We develop a theory of abstract arithmetic Chow rings, where the 
role of the fibers at infinity is played by a complex of abelian 
groups that computes a suitable cohomology theory. As particular 
cases of this formalism we recover the original arithmetic 
intersection theory of H.~Gillet and C.~Soul\'e for projective 
varieties. We introduce a theory of arithmetic Chow groups, which 
are covariant with respect to arbitrary proper morphisms, and we 
develop a theory of arithmetic Chow rings using a complex of 
differential forms with log-log singularities along a fixed normal
crossing divisor. This last theory is suitable for the study of
automorphic line bundles. In particular, we generalize the classical 
Faltings height with respect to logarithmically singular hermitian 
line bundles to higher dimensional cycles. As an application we 
compute the Faltings height of Hecke correspondences on a product 
of modular curves. 
\end{abstract}

\tableofcontents
\section{Introduction}

\nnpar{Background.} Arakelov Geometry was initiated by S.J.~Arakelov 
in his paper \cite{Arakelov:itdas}, where he shows that one can
\emph{compactify} a curve defined over the (spectrum of the) ring 
of integers of a number field by considering Green functions on the 
associated complex curve. Subsequently, Arakelov's ideas have been 
successfully refined by P.~Deligne \cite{Deligne:dc} and G.~Faltings 
\cite{Faltings:cas}, and generalized to higher dimensional arithmetic 
varieties by H.~Gillet and C.~Soul\'e who introduced arithmetic 
Chow groups for higher codimensional arithmetic cycles in \cite
{GilletSoule:ait}. An important application of this formalism was 
P.~Vojta's proof of Fal\-tings's theorem, formerly known as Mordell's 
conjecture (see \cite{Vojta:stcc}).

The use of analytical objects such as Green functions allows various
modifications and extensions of the formalism of arithmetic Chow
groups, which is adapted to different situations. Among others, we
mention S.~Zhang's work \cite{Zhang:_small} on admissible metrized 
line bundles and V.~Maillot's work \cite{Maillot:GAdvt} taking into 
account hermitian line bundles whose metrics are no longer smooth 
but still continuous; J.-B.~Bost's work \cite{Bost:Lfg} on $L^{2}_{1}
$-Green functions and A.~Moriwaki's work \cite{Moriwaki:ipacdGc} on 
degenerate Green currents; U.~K\"uhn's work \cite{Kuehn:gainc} 
treating hermitian line bundles on arithmetic surfaces having
logarithmically singular metrics, e.g., the line bundle of modular
forms equipped with the Petersson metric on a modular curve, and
J.I.~Burgos's work \cite{Burgos:acr}, \cite{Burgos:CDB} on arithmetic
Chow groups and Deligne-Beilinson cohomology.

\nnpar{Motivation.}
The main motivation of these notes is to extend U.~K\"uhn's 
generalized arithmetic intersection pairing on arithmetic surfaces 
to higher dimensional arithmetic varieties. While U.~K\"uhn's 
intersection pairing, which was motivated by preliminary 
results of J.~Kramer (see \cite{Kramer:atsjf}), is given in terms 
of an explicit formula for the arithmetic intersection number of 
two divisors (in the spirit of Arakelov), the development of a 
corresponding higher dimensional theory needs to be approached 
in a more abstract way. Moreover, we will extend this theory, 
not only to line bundles with singular metrics, but also to higher 
rank vector bundles such as the Hodge bundle that appear when 
considering fibrations of semi-abelian varieties. The study of 
higher rank vector bundles will be the subject of the forthcoming 
paper \cite{BurgosKramerKuehn:accavb}.  

\nnpar{Arithmetic intersection theory.} An arithmetic ring $(A,
\Sigma,F_{\infty})$ is a triple consisting of an excellent regular 
noetherian integral domain $A$, a finite non-empty set $\Sigma$ of 
monomorphisms $\sigma:A\longrightarrow\mathbb{C}$, and an antilinear 
involution $F_{\infty}:\mathbb{C}^{\Sigma}\longrightarrow\mathbb{C}^
{\Sigma}$ of $\mathbb{C}$-algebras. For simplicity, we will forget
about the antilinear involution $F_{\infty}$ in this introduction. 
Then, an arithmetic variety is a flat, regular scheme $X$ over $S=
\Spec(A)$ together with a complex analytic space $X_{\infty}$ 
obtained from $X$ by means of $\Sigma$; $S$ will be called the base
scheme. Intuitively, the elements of $\Sigma$ are the points at
infinity of $S$ providing a ``compactification'' of $S$, and the 
analytic space $X_{\infty}$ is the fiber at these points at infinity 
or archimedean fiber of $X$. An arithmetic intersection theory will 
involve three main ingredients: first, a geometric intersection 
theory over the scheme $X$, the \emph{geometric part}, second, a
``refined'' intersection theory over $X_{\infty}$, the \emph{analytic 
part}, and finally an interface relating the geometric and the 
analytic part. The main theme of this paper is to study the second 
and third of these ingredients, and we will rely on existing geometric
intersection theories.  

\nnpar{Geometric intersection theory.} It might be useful to review,
at this point, the geometric intersection theories we have at our 
disposal.

The first of these theories is based on the moving lemma to reduce 
the intersection product of two algebraic cycles to the case of 
proper intersection. In order to be able to apply the moving lemma 
we need the scheme $X$ under consideration to be quasi-projective 
and the arithmetic ring $A$ to be a field. Since we are interested 
in more general arithmetic rings, we cannot use this method.

The second approach is the deformation to the normal cone technique 
due to W.~Fulton and R.~MacPherson (see \cite {Fulton:IT}). This method 
is much more general, since the scheme $X$ need not be quasi-projective
and, moreover, it can also be applied to the case in which the base scheme
is the spectrum of a Dedekind domain. But in this case $X$ needs to be
not only regular, but smooth over $S$. Since most interesting arithmetic 
varieties are not smooth over the base scheme, this method is not general 
enough for our purposes. Note however that, in the case in which $X$ is 
smooth over $S$, this method provides an intersection pairing which 
is defined without tensoring with $\mathbb{Q}$. We should also note
that, in contrast to the case of an intersection product, this method 
can be used to define the inverse image morphism for a morphism between 
regular schemes of finite type over the spectrum of a regular, noetherian
ring without the assumption of smoothness. The lack of a general theorem 
of resolution of singularities prevents us from obtaining an intersection 
product from this inverse image morphism. 

However, as a third approach, one can use J.~de Jong's theorem on 
alterations \cite{deJong:alt} to reduce the intersection of algebraic 
cycles to the inverse image between regular schemes, and then apply 
the deformation to the normal cone technique. Nevertheless, these
alterations are finite morphisms whose degree is larger that one, in
general. Therefore, this method yields an intersection product with 
rational coefficients. 

The last general approach that we will mention (and the one introduced
originally by H.~Gillet and C.~Soul\'e) is to use the isomorphism 
between $K$-theory and Chow groups to transfer the ring structure of 
$K$-theory to the Chow groups. This method is valid for any regular,
noetherian scheme $X$ of finite Krull dimension. The main drawback of 
this method is that the isomorphism between $K$-theory and Chow groups 
is only true in general after tensoring with $\mathbb{Q}$. Therefore, 
one also obtains a rational valued intersection product. 

Finally, we note that the intersection product with a divisor on a 
regular scheme can be defined directly using a simple version of the 
moving lemma for divisors (see \cite{GilletSoule:ait}, 4.2.3.2), which 
holds in complete generality.

The $K$-theoretical method and the alteration method are the most 
general of the above methods. Since the $K$-theoretical approach was 
the one used by H.~Gillet and C.~Soul\'e in \cite{GilletSoule:ait} 
and, conceptually, fits very well with the cohomological approach 
we will be using for the analytical part, we also use it 
as the geometric part of our arithmetic intersection theory. But, 
since the geometric and the analytic part of our arithmetic intersection 
theory are isolated and related by a clear interface, we hope that 
the reader will have no difficulty in using any other of these geometric 
methods when applicable.    

\nnpar{Green currents and Green forms.}
We now discuss the refined intersection theories which have been used 
as the analytic part of a higher dimensional arithmetic intersection 
theory. 

The first one is the notion of a Green current introduced in \cite
{GilletSoule:ait}. Let $y$ be a $p$-codimensional cycle on $X$. A 
Green current for $y$ is a class of currents     
\begin{equation*}
g_{y}\in\widetilde D^{p-1,p-1}_{X_{\infty}}=D^{p-1,p-1}_{X_{\infty}}/
(\partial D^{p-2,p-1}_{X_{\infty}}+\bar\partial D^{p-1,p-2}_{X_{
\infty}}) 
\end{equation*}
satisfying the residue equation 
\begin{equation}
\label{eq:48}
\dd\dd^{c}g_{y}+\delta_{y}=\omega_{y}, 
\end{equation}
where $\delta_{y}$ is the current given by integration along the 
cycle $y_{\infty}$ and $\omega_{y}$ is a smooth form uniquely 
determined by \eqref{eq:48}.

If $y$ and $z$ are two cycles intersecting properly (i.e., with the
appropriate codimension) in $X_{\infty}$, the product of two Green 
currents, called the star-product or $\ast$-product, is given by the 
formula 
\begin{displaymath}
g_{y}\ast g_{z}=g_{y}\land\delta_{z}+\omega_{y}\land g_{z}.
\end{displaymath}
It is a Green current for the intersection cycle $y\cdot z$. This 
approach has many analytical difficulties. For instance, some care
has to be taken to define $g_{y}\land\delta_{z}$. Moreover, this
product is not obviously associative and commutative, and the proof
that this is indeed the case, is not trivial.

As we have already mentioned this approach has been generalized in
several directions. Typically, these generalizations allow the 
presence of certain singularities for the differential form $\omega_
{y}$. But usually only the case when $X_{\infty}$ is of dimension
one is treated (see \cite{Bost:Lfg}, \cite{Kuehn:gainc}), or one
does not obtain a full ring structure for the generalized arithmetic 
Chow rings (see \cite{Maillot:GAdvt}). 

There are other proposals for the definition of the product of
Green currents. For instance, B.~Harris and B.~Wang \cite
{HarrisWang:_arith} have given a definition of the star-product 
of two Green currents for non-properly intersecting cycles that 
depends on a deformation of one of the cycles, and N.~Dan \cite
{Dan:_prolon_green} has given a definition of the  star-product 
using meromorphic continuation of certain zeta functions.  

In \cite{Burgos:CDB} J.I.~Burgos introduced a new definition 
of Green forms along the following lines. To every complex algebraic
manifold $X$ (not necessarily compact), there is associated a graded
complex $\mathcal{D}^{\ast}_{\log}(X,\ast)$, which consists of smooth 
forms on $X$ with logarithmic singularities at infinity. For instance, 
if $X$ is proper, then
\begin{align*}
\mathcal{D}^{2p}_{\log}(X,p)=E^{p,p}(X)\cap (2\pi i)^{p}E^{2p}_
{\mathbb{R}}(X), \\
\mathcal{D}^{2p-1}_{\log}(X,p)=E^{p-1,p-1}(X)\cap (2\pi i)^{p-1}
E^{2p-2}_{\mathbb{R}}(X),
\end{align*}
where $E^{p,p}(X)$ is the space of smooth complex valued differential 
forms of type $(p,p)$ and $E_{\mathbb{R}}^{2p}(X)$ is the space of 
smooth real valued differential forms. The boundary morphism
\begin{displaymath}
\dd_{\mathcal{D}}:\mathcal{D}^{2p-1}_{\log}(X,p)\longrightarrow 
\mathcal{D}^{2p}_{\log}(X,p)
\end{displaymath}
is given by $\dd_{\mathcal{D}}\eta=-2\partial\bar\partial\eta$. 
Observe that, up to a normalization factor, this is the same
differential operator as the one that appears in the residue equation  
\eqref{eq:48}.

The complex $\mathcal{D}^{\ast}_{\log}(X,\ast)$ computes the real 
Deligne-Beilinson cohomology of $X$, which is denoted by $H^{\ast}_
{\mathcal{D}}(X,\mathbb{R}(p))$. If $Y$ is a closed subset of $X$,
then the real Deligne-Beilinson cohomology of $X$ with supports on 
$Y$, denoted by  $H^{\ast}_{\mathcal{D},Y}(X,\mathbb{R}(p))$, is the 
cohomology of the simple complex associated to the morphism of 
complexes     
\begin{displaymath}
\mathcal{D}^{\ast}_{\log}(X,\ast)\longrightarrow\mathcal{D}^{\ast}_
{\log}(X\setminus Y,\ast).
\end{displaymath}
Every $p$-codimensional algebraic cycle $y$ with support $Y$ defines 
a cohomology class $\cl(y)\in H^{2p}_{\mathcal{D},Y}(X,\mathbb{R}(p))$.
Moreover, if $W\subseteq X$ is a closed subvariety of codimension 
$p-1$ and $f\in K(W)^{\ast}$, then $f$ defines a class 
\begin{displaymath}
\cl(f)\in H^{2p-1}_{\mathcal{D}}(X\setminus\dv(f),\mathbb{R}(p)).
\end{displaymath}
These classes satisfy the compatibility condition 
\begin{displaymath}
\delta(\cl(f))=\cl(\dv(f)),
\end{displaymath}
where
\begin{displaymath}
\delta:H^{2p-1}_{\mathcal{D}}(X\setminus\dv(f),\mathbb{R}(p))
\longrightarrow H^{2p}_{\mathcal{D},\dv(f)}(X,\mathbb{R}(p))
\end{displaymath}
is the connecting homomorphism.

In this setting a Green form for a $p$-codimensional algebraic cycle
$y$ is a representative of the class $\cl(y)$. More explicitly, 
we write  
\begin{align*}
{\rm Z}\mathcal{D}^{\ast}_{\log}(X,\ast)&=\Ker(\dd_{\mathcal{D}}), \\
\widetilde{\mathcal{D}}^{\ast}_{\log}(X,\ast)&={\mathcal{D}}^{\ast}_
{\log}(X,\ast)/\Img(\dd_{\mathcal{D}}).
\end{align*}
The space of Green forms associated to a $p$-codimensional cycle
$y$ with support $Y$ is then given by
\begin{displaymath}
GE(y)=\left\{(\omega_{y},\widetilde g_{y})\in{\rm Z}\mathcal{D}^{2p}_
{\log}(X,p)\oplus\widetilde{\mathcal{D}}^{2p-1}_{\log}(X\setminus Y,p)
\big|\begin{aligned}\dd_{\mathcal{D}}\widetilde{g}_{y}&=\omega_{y}\\
\cl((\omega_{y},\widetilde{g}_{y}))&=\cl(y)\end{aligned}\right\}.
\end{displaymath}
The star-product of Green forms is now simply the cup product in
cohomology with supports. With this approach, the proof of the 
associativity and commutativity of the star-product is straightforward 
and completely formal. 

In \cite{Burgos:Gftp} and \cite{Burgos:CDB} it is proven that, when
$X$ is projective, the arithmetic Chow groups obtained by this method
agree with the ones obtained by the method of Gillet and Soul\'e. It
is interesting to note that all the analytical complexities appearing
in the proof of the associativity and commutativity of the star-product
in \cite{GilletSoule:ait} are needed to prove the compatibility of 
the two definitions.

In contrast, in the quasi-projective case, the groups obtained by 
this new method have better Hodge theoretical properties. For instance,
they possess a certain homotopy invariance with respect to vector 
bundles. Another advantage of this new definition is that it is very 
easy to make variants adapted to new problems just by changing the 
complex $\mathcal{D}_{\log}$. 

We should stress here that a Green form associated to a cycle 
is a representative of the cohomology class of the cycle with support 
in the same cycle. In order for the star-product of two Green forms 
to be a Green form for the intersection cycle, we need the cycles to 
intersect properly in $X_{\infty}$. Therefore, as was the case in 
\cite{GilletSoule:ait}, the arithmetic intersection product relies on 
the moving lemma for complex varieties. In particular, this implies 
that our varieties should be at least geometrically quasi-projective.

In the preprint \cite{Hu:sogf}, J.~Hu uses the homotopy invariance 
and the flexibility of the definition of Green forms to give a new
definition of the inverse image morphism of arithmetic Chow groups 
for a regular closed immersion by means of Fulton's deformation to 
the normal cone technique. The main result of this paper is the
construction of a specialization morphism for Green forms. Combining
this technique with J.~de Jong's result on alterations mentioned 
above, one can define the arithmetic intersection pairing without 
using any moving lemma, thereby removing the hypothesis of 
quasi-projectivity.

\nnpar{The interface between the geometric and the analytic part.} 
In the definition of arithmetic Chow groups by means of Green 
currents, this interface is implemented by very concrete objects. 
Namely, to any $p$-codimensional algebraic cycle $y$, one assigns 
the current ``integration along the cycle $y_{\infty}$'', denoted 
by $\delta_{y}$, and to every rational function $f$, one associates 
the current $-\log|f_{\infty}|^{2}$. Compared with to the approach of 
Green currents, this interface is more abstract. It is implemented by 
the theory of characteristic classes for cycles and rational functions. 
The two approaches are compatible because in the appropriate complexes, 
the currents $\delta_{y}$ and $-\log|f_{\infty}|^{2}$ represent the 
cohomology class of $y_{\infty}$ and $f_{\infty}$, respectively.  

\nnpar{Abstract arithmetic Chow groups.}
Recall that our main motivation is to extend U.~K\"uhn's generalized 
arithmetic intersection pairing on arithmetic surfaces to higher 
dimensional arithmetic varieties. In order to accomplish this goal 
we will use the flexibility of the Green form approach, changing 
the complex $\mathcal{D}_{\log}$ by a complex of differential forms 
with certain mild singularities along a fixed subvariety. Nevertheless, 
the arithmetic Chow groups that we will define in this way, and their
properties, depend strongly on the actual complex used. And there is  
no one best choice for this complex of singular differential 
forms. For instance, in this paper we introduce the complex of 
pre-log-log differential forms which, although it does not have all
the cohomological properties one would expect, is enough to define
an arithmetic intersection pairing and, in particular, the height 
with respect to log singular hermitian line bundles. On the other 
hand, in the paper \cite{BurgosKramerKuehn:accavb}, we introduce
the complex of log-log singular differential forms. This complex has 
the expected cohomological properties, but is slightly more difficult 
to handle because one has to bound all the derivatives of the functions
involved. Moreover, as we mentioned previously, in the literature
there are several other variants of arithmetic Chow groups with
singular differential forms. 

In addition, if one uses completely different kinds of complexes, one 
can obtain arithmetic Chow groups with new properties. For instance, 
in his PhD thesis \cite{Burgos:acr}, the first author constructed a 
fully covariant version of the arithmetic Chow groups by choosing a 
complex of currents instead of the complex of smooth differential 
forms. Similar arithmetic Chow groups were introduced independently 
by A.~Moriwaki in \cite{Moriwaki:ipacdGc}.     

In another direction, one may consider the following example. Let
$p_i$, $i=1,\dots,4$ be four different points in $\mathbb{P}^{1}_
{\mathbb{Q}}$. Then, the archimedean component of the height pairing 
of $p_{1}-p_{2}$ and $p_{3}-p_{4}$ is essentially given by the logarithm 
of the norm of the cross ratio of the four points. Thus, this height 
pairing has a canonical lifting from $\mathbb{R}$ to $\mathbb{C}^
{\ast}$, namely the cross ratio of the four points. This suggests 
that one can define a finer version of arithmetic Chow groups, where 
the role of real Deligne-Beilinson cohomology is played by integral 
Deligne-Beilinson cohomology.  

Even more, one can think of an adelic version of the arithmetic Chow 
groups, similar to the one introduced in \cite{BlochGilletSoule:naat}, 
but where each geometric fiber is substituted by a suitable complex
that computes a certain cohomology theory, or a theory, where only
certain geometric fibers are substituted by a cohomological complex. 

All these considerations lead us to the conclusion that it is
worthwhile to develop a formalism of arithmetic Chow groups with
respect to an abstract complex and to explore how the properties 
of the complex are reflected by the properties of the arithmetic 
Chow groups. In this way, different variants of arithmetic Chow 
groups can be obtained as particular cases. 

\nnpar{Results.} The main achievement of this paper is the systematic
development of the formalism of abstract arithmetic Chow groups in
arbitrary dimensions depending on a suitable cohomological complex.
Every choice of such a complex gives rise to new types of arithmetic
Chow groups with different properties reflected by the complexes under
consideration. As special cases we recover Burgos' version of the
the arithmetic intersection theory developed by H.~Gillet and
C.~Soul\'e in \cite{GilletSoule:ait}, we introduce a theory of
arithmetic Chow groups which are covariant with respect to arbitrary
proper morphisms, and we develop a theory of arithmetic Chow rings
using a complex of differential forms with pre-log-log singularities 
along a fixed normal crossings divisor. This latter theory is suitable 
for the study of automorphic line bundles. In particular, we generalize
the classical Faltings height with respect to a logarithmically
singular hermitian line bundle to higher dimensional cycles. As an
application we compute the Faltings height of Hecke correspondences 
on a product of modular curves.
 
This formalism of arithmetic Chow groups is an abstraction of 
\cite{Burgos:CDB}. Note however that the passage from the concrete
example of \cite{Burgos:CDB} to the abstract version presented here is
not completely straightforward. Although some constructions such as the
definition of truncated cohomology classes and their product are already 
(at least implicitly) present in \cite{Burgos:CDB}, others, like the 
notions of $\Gi$-complex and of covariant $f$-morphism of complexes, 
are new.

The basic idea of this paper is that the role of the complex $\mathcal
{D}_{\log}$ can be played by any graded complex of sheaves in the Zariski 
topology $\cc$. We only require two properties for this complex. The
first is  
that the hypercohomology of this complex of sheaves always agrees with 
the cohomology of the complex of global sections. In this way we can 
represent cohomology classes by concrete elements of this complex. A 
convenient way to ensure this is to ask the complex to satisfy a 
Mayer-Vietoris condition. The second property we require is that the
complex receives characteristic classes from cycles and rational 
functions (i.e., from $K_{1}$-chains). Typically, in order to ensure
the existence of characteristic classes, one imposes a series of axioms 
to the cohomology. Nevertheless, in many applications it is not convenient 
to impose too many conditions to the cohomology. Therefore, we use an 
auxiliary cohomology, given by a graded complex of sheaves $\Gi$ 
satisfying the axioms of \cite{Gillet:RRhK}; we call it a Gillet 
cohomology. The characteristic classes will be induced by a morphism 
$\mathfrak{c}_{\cc}:\Gi\longrightarrow\cc$ in the derived category.

In this abstract setting, instead of Green forms, we will define Green
objects that live in a space that we call ``truncated cohomology
classes''. These truncated cohomology classes are something between
relative cohomology classes and representatives of relative cohomology
classes. 

Of course, very little can be done with this minimal set of properties. 
Therefore, we undertake a complete study of how the properties of the 
complex $\cc$ are reflected by the properties of the Green objects and 
the arithmetic Chow groups. For instance, in order to have an intersection
product in the arithmetic Chow groups, we only need the existence of 
a cup product in the complex $\cc$, which is compatible with the cup 
product of the complex $\Gi$. This implies that the cup product in 
$\cc$ is compatible with the intersection product of cycles.   

We emphasize here that this abstract approach also simplifies many
difficulties that appear when working with Green currents. We have
already mentioned the proof of the associativity and commutativity of
the star-product, but our approach also provides a new and logically 
independent proof of the well-definedness of the arithmetic intersection 
product due to H.~Gillet and C.~Soul\'e. This proof does not rely on 
the $K_{1}$-chain moving lemma. We point out that in J.I.~Burgos's 
preceding work the arguments for proving the well-definedness of the 
arithmetic intersection product relied on the corresponding arguments 
in \cite{GilletSoule:ait}. We will give a more detailed discussion on 
how we avoid the $K_{1}$-chain moving lemma in remark \ref{rem:7}.
We also emphasize that some problems  with the $K_{1}$-chain moving
lemma have been discussed and successfully solved using completely
different  
techniques by W.~Gubler in \cite{Gubler:mlfkc}. 

Once the abstract theory is developed, in the subsequent chapters, we
study particular cases of this construction. For instance, we recover 
the original arithmetic intersection theory of H.~Gillet and C.~Soul\'e 
for projective varieties. We point out that we will show in \cite
{BurgosKramerKuehn:accavb} how to recover these arithmetic Chow groups 
for quasi-projective varieties as a particular case of our construction.
This example agrees with the theory developed in \cite{Burgos:CDB}.

As a second example we introduce a theory of arithmetic Chow groups 
which are covariant with respect to arbitrary proper morphisms. This 
construction was first introduced in \cite{Burgos:acr}, and a similar
construction can be found in \cite{Moriwaki:ipacdGc}.  

Furthermore, by choosing for $\cc$ a complex of forms satisfying certain 
growth conditions of log- and log-log-type, we obtain a theory which is compatible with the theories developed by J.-B.~Bost and U.~K\"uhn in 
the $1$-dimensional setting. This latter theory is specifically suited 
for the study of automorphic line bundles and allows to generalize the 
classical Faltings height with respect to a logarithmically singular 
hermitian line bundle to higher dimensional cycles. As an application, 
we compute the Faltings height of certain Hecke correspondences. Note 
that the same result has been obtained recently by P.~Autissier \cite
{Autissier:hhc} using the $1$-dimensional theory.

We point out that Bost's theory can also be seen as a particular case 
of our abstract setting. Nevertheless, the definition of the corresponding complex $\cc$ involves a mixture of $L^{2}$-forms, $L^{2}_{1}$-forms,
$L^{2}_{-1}$-currents, forms with logarithmic singularities, and currents; 
we will not write it  explicitly.

It would be interesting to extend the abstract setting of this paper
to cover also Hu's deformation to the normal cone technique. This would 
involve incorporating the specialization functor to the axiomatic 
system and asking for the existence of a specialization morphism at 
the level of complexes. 

\nnpar{Applications.} The theory developed in this paper is extensively 
used in the forthcoming paper \cite{BruinierBurgosKuehn:bpaihs} by 
J.~Bruinier, J.I.~Burgos and U.~K\"uhn, where explicit calculations 
for the arithmetic self-intersection number of the line bundle of 
modular forms and the Faltings height of Hirzebruch-Zagier cycles on 
Hilbert modular surfaces are carried out. Further calculations in 
this direction for other naturally metrized automorphic line bundles 
have been made in \cite{BruinierKuehn:iaghd}. In his forthcoming thesis 
G.~Freixas-Montplet will prove finiteness results for the height with 
respect to such naturally metrized line bundles, which generalize 
G.~Faltings' original result for points to cycles of higher dimensions. 
In the sequel \cite{BurgosKramerKuehn:accavb} of this paper, we will 
show that our abstract arithmetic Chow groups attached to the complex 
of forms having suitable growth conditions of log- and log-log-type 
combined with the work \cite{BurgosWang:hBC} allow us to construct 
arithmetic characteristic classes for vector bundles equipped with 
hermitian metrics, which are logarithmically singular along a divisor 
with normal crossings. These arithmetic characteristic classes give 
rise to operations on the arithmetic Chow groups, even for non-regular
arithmetic varieties. In addition, we show in \cite{BurgosKramerKuehn:accavb} 
that automorphic vector bundles equipped with the natural invariant 
metric (Petersson metric) on Shimura varieties of non-compact type 
are hermitian vector bundles of the type considered above.

The framework of our arithmetic Chow groups attached to forms 
having certain log- and log-log-type singularities is one of 
the key ingredients in order to formulate various conjectures: 
in this context, we mention a conjecture of K.~K\"ohler on 
arithmetic intersection numbers on the moduli space of principally 
polarized abelian varieties (see \cite{Koehler:hppat}); secondly, 
we mention a conjecture of V.~Maillot and D.~Roessler on arithmetic 
Chern numbers associated to fibrations of motives with complex 
multiplication (see \cite{MaillotRoessler:cdl}); finally, we 
emphasize the conjectures of S.~Kudla on Faltings heights and 
generating series for special cycles on Shimura varieties of 
orthogonal type (see \cite{Kudla:msri}, \cite{Kudla:desgf}, 

\nnpar{Outline of the paper.} 
Let us now give a more detailed outline of the contents of each 
chapter.

\nnpar{Chapter 1.}
In the first chapter we review various results relating $K$-theory, 
Chow groups and cohomology theories satisfying Gillet's axioms
\cite{Gillet:RRhK}; these facts will be needed in the sequel. For 
our purposes, the main interest in $K$-theory is that it provides 
a method to define an intersection product and an inverse image 
for algebraic cycles on regular, separated schemes of finite type 
over a base scheme, which is regular, separated, noetherian, and 
of finite Krull dimension. Most of the results needed from $K
$-theory deal only with the groups $K_{0}$ and $K_{1}$, and can be
found in the first chapter of \cite{Soule:lag}. The main 
exception to this is proposition \ref{prop:2}; its translation 
into Chow theory and cohomology theory, namely theorem \ref{thm:3}
and corollary \ref{cor:2}, is of crucial importance in chapter 4
in the course of the proof of the well-definedness of the intersection 
product of arithmetic cycles. There is very little new in this chapter
and its main purpose is to gather together all the needed results,
some of which are difficult to find explicitly in the literature.

\nnpar{Chapter 2.}
The second chapter is devoted to a systematic study of relative
cohomology groups $H^{\ast}(A,B)$ attached to a morphism $f:A
\longrightarrow B$ of abstract complexes of abelian groups and 
their product structure based on a product structure of the 
complexes under consideration. Furthermore, we study truncated 
relative cohomology groups $\widehat{H}^{\ast}(A,B)$ associated 
to the above data together with their product structure, which 
is the basis of the definition of the $*$-product in the third 
chapter. We have also included a discussion of the signs appearing
when considering multidimensional complexes, complexes of
complexes, and products between them. This chapter is an extended
and much more detailed version of the corresponding chapter of
\cite{Burgos:CDB}.  

\nnpar{Chapter 3.}
The aim of the third chapter is to develop an abstract theory of 
Green objects as elements of a suitable truncated cohomology theory. 
The main property for such a cohomology theory is that it receives
characteristic classes from $K$-theory, at least from $K_{0}$ and 
$K_{1}$, and that it satisfies some additional natural properties. 
More precisely, we fix a Gillet complex $\Gi=\Gi^{\ast}(\ast)$ over 
the site of regular schemes $X$ of finite type over a field $k$. A
graded complex $\cc=\cc^{\ast}(\ast)$ of sheaves of abelian groups 
together with a morphism $\mathfrak{c}_{\cc}:\Gi\longrightarrow\cc$ 
in the derived category will be called a $\Gi$-complex over $X$. A 
Green object for a $p$-codimensional cycle $y$ on $X$ with values 
in $\cc$ is then given by an element 
\begin{displaymath}
\mathfrak{g}_{y}\in\widehat{H}^{2p}(\cc(X,p),\cc(U,p))
\end{displaymath}
such that the class of $\mathfrak{g}_{y}$ equals the class of the 
cycle $y$ in the relative cohomology group $H^{2p}(\cc(X,p),\cc(U,
p))$; here $U$ denotes the complement of the support of $y$ in $X$. 
After studying the basic properties of such Green objects, we define 
a $*$-product for two Green objects using the techniques developed 
in the second chapter. We end this chapter with a proof of the 
associativity and commutativity of the $*$-product under suitable 
assumptions on the $\Gi$-complexes under consideration. The material
in this and the next chapter, although a generalization of the 
results of \cite{Burgos:CDB}, is new.

\nnpar{Chapter 4.}
In the fourth chapter we introduce generalized arithmetic Chow
groups for arithmetic varieties $X$ over an arithmetic ring. The 
idea behind this definition is that the arithmetic variety $X$ 
can be ``compactified'' by adding the associated complex manifold
$X_{\infty}$ or, more precisely, a truncated cohomology theory
on $X_{\infty}$ to the picture. More specifically, we proceed as
follows: After fixing a Gillet complex $\Gi$ over the site of 
schemes over the real numbers, we call the pair $(X,\cc)$ a $\Gi
$-arithmetic variety, when $X$ is an arithmetic variety and $\cc$
a $\Gi$-complex over the associated real variety $X_{\mathbb{R}}$.
The group of $p$-codimensional arithmetic cycles of $(X,\cc)$ is
then given by the set of pairs $(y,\mathfrak{g}_{y})$, where $y$
is a $p$-codimensional cycle on $X$ and $\mathfrak{g}_{y}$ is a 
Green object for (the class of) the cycle induced by $y$ on $X_
{\mathbb{R}}$. The $\cc$-arithmetic Chow group $\cha^{p}(X,\cc)$
is now obtained from the group of $p$-codimensional arithmetic 
cycles of $(X,\cc)$ by factoring out by a suitable rational
equivalence relation. We prove various properties for these 
generalized arithmetic Chow groups, emphasizing how the properties 
of the $\Gi$-complexes involved are reflected in the properties 
of the arithmetic Chow groups; a typical example is the ring 
structure of the direct sum $\bigoplus_{p\ge 0}\cha^{p}(X,\cc)_
{\mathbb{Q}}$.

\nnpar{Chapter 5.}
In all the examples of generalized arithmetic Chow groups 
presented in this paper, the underlying Gillet cohomology 
will be the Deligne-Beilinson cohomology. Therefore, we 
recall in the fifth chapter the basic definitions and facts 
of Deligne-Beilinson cohomology and homology, which will be 
needed in the sequel. In particular, we use the fact that real
Deligne-Beilinson cohomology can be computed as the sheaf 
cohomology of the Deligne algebra associated to the Dolbeault 
algebra of differential forms with logarithmic singularities 
at infinity. We denote this graded complex of sheaves by
$\mathcal{D}_{\log}$. Towards the end of this chapter, we 
give explicit representatives for the classes of cycles, 
in particular for the classes of divisors of sections of 
(hermitian) line bundles, in real Deligne-Beilinson cohomology 
in terms of the underlying singular differential forms. Most 
of the material in this chapter is well known; we include it 
for the convenience of the reader.

\nnpar{Chapter 6.}
In the sixth chapter we use our abstract theory of arithmetic
Chow groups to define contravariant and covariant arithmetic 
Chow groups starting with the Gillet complex $\Gi=\mathcal{D}_
{\log}$. The contravariant Chow groups, which were introduced 
in \cite{Burgos:CDB}, are obtained by considering $\mathcal{D}_
{\log}$ itself as a $\Gi$-complex. In this way, we obtain the 
arithmetic Chow groups $\cha^{\ast}(X,\mathcal{D}_{\log})$. By 
means of the properties of the Deligne algebra $\mathcal{D}_
{\log}$, we find that the direct sum $\bigoplus_{p\ge 0}\cha^
{p}(X,\mathcal{D}_{\log})_{\mathbb{Q}}$ has the structure of 
a commutative and associative ring. Furthermore, this ring 
coincides with the arithmetic Chow ring defined by Gillet and 
Soul\'e in \cite{GilletSoule:ait} for arithmetic varieties with 
projective generic fiber. The covariant Chow groups, which were 
introduced in \cite{Burgos:acr}, are defined using as a $\mathcal
{D}_{\log}$-complex a complex $\mathcal{D}_{\D}$, which is made 
out of certain currents and computes real Deligne-Beilinson 
homology. The properties of the $\mathcal{D}_{\log}$-complex 
$\mathcal{D}_{\D}$ show that the arithmetic Chow groups $\cha^
{\ast}(X,\mathcal{D}_{\D})$ are covariant for arbitrary proper 
morphisms and have the structure of a module over the 
contravariant Chow ring. We end this chapter with a reformulation 
of the definition of the height of a cycle, originally given in 
\cite{BostGilletSoule:HpvpGf}, in the framework of contravariant 
and covariant arithmetic Chow groups.

\nnpar{Chapter 7.}
In the seventh chapter we again fix the Deligne algebra $\mathcal
{D}_{\log}$ as our Gillet complex $\Gi$. We then construct a 
$\mathcal{D}_{\log}$-complex by means of differential forms 
satisfying (together with their derivatives with respect to 
$\partial$, $\bar{\partial}$, and $\partial\bar{\partial}$) 
certain growth conditions of log- and log-log-type. We call 
these differential forms pre-log and pre-log-log forms, 
respectively, and denote the corresponding $\mathcal{D}_{\log}
$-complex by $\mathcal{D}_{\wlg}$. The notation is justified 
in order to distinguish these differential forms from the forms 
satisfying the same type of growth conditions together with all 
their derivatives; the corresponding $\mathcal{D}_{\log}$-complex
will be introduced and studied in \cite{BurgosKramerKuehn:accavb}. 
By means of the properties of the complex $\mathcal{D}_{\wlg}$, 
the direct sum $\bigoplus_{p\ge 0}\cha^{p}(X,\mathcal{D}_{\wlg})_
{\mathbb{Q}}$ has the structure of a commutative and associative 
ring. This provides the desired higher dimensional extension of 
the generalized arithmetic intersection pairing introduced in 
\cite{Kuehn:gain}. A useful application of this formalism is
the extension of the definition of the height of a cycle with
respect to line bundles equipped with a hermitian metric, which
becomes logarithmically singular along a divisor with normal
crossings. As an illustration, we compute the arithmetic 
self-intersection number of the line bundle of modular forms
on the product of two modular curves equipped with the Petersson
metric. Furthermore, we determine the Faltings height of Hecke 
correspondences on this product of modular curves with respect 
to the line bundle of modular forms. The same result has been 
obtained recently by P.~Autissier \cite{Autissier:hhc} using a 
different approach. Related but more elaborate results in the 
case of Hilbert modular surfaces are contained in the paper
\cite{BruinierBurgosKuehn:bpaihs} mentioned above. The general 
theory of arithmetic characteristic classes of automorphic
vector bundles of arbitrary rank will be developed in \cite
{BurgosKramerKuehn:accavb}.

\nnpar{Acknowledgements.} In the course of preparing this manuscript,
we had many stimulating discussions with many colleagues. We would 
like to thank them all. In particular, we would like to express our 
gratitude to J.-B.~Bost, W.~Gubler, M.~Harris, S.~Kudla, V.~Maillot, 
D.~Roessler, C.~Soul\'e. We would also like to thank G.~Freixas for 
his careful proof reading of the manuscript. Furthermore, we would 
like to thank to EAGER, Arithmetic Geometry Network, Newton Institute, 
and the Institut Henri Poincar\'e for partial support of our work. 

\newpage
\section{Background results on $K$-theory}
\label{sec:rKt}

In this section we will review some facts relating to $K$-theory,
Chow groups and  cohomology theory that will be needed in the 
sequel. For our purposes, the main interest of $K$-theory is 
that it provides a method to define the intersection product 
and the inverse images of algebraic cycles for regular schemes. 
Most of the results we need are only concerned with the groups 
$K_{0}$; they are explained in chapter 1 of \cite{Soule:lag}. 
The main exception to this is proposition \ref{prop:2}. Its 
translation into Chow theory, theorem \ref{thm:3}, and cohomology  
theory will be used in the proof of the fact that the intersection  
product of arithmetic cycles is well defined. Note however that
proposition  
\ref{prop:2} is proven in \cite{Gillet:RRhK}, and is used in
\cite{GilletSoule:ait}.

In this section all schemes will be noetherian, separated and of
finite Krull dimension. Given an abelian group $A$ we will write
$A_{\mathbb{Q}}=A\otimes\mathbb{Q}$.

\subsection{$K$-theory}

\nnpar{$K$-theory of schemes.} 
For any exact category $\mathscr{A}$, Quillen has introduced a 
simplicial space $\Omega BQ\mathscr{A}$, and has defined the 
$K$-groups of the category $\mathscr{A}$ as the homotopy groups 
of this simplicial space (see \cite{Quillen:haKt}), i.e.,
\begin{displaymath}
K_{m}(\mathscr{A})=\pi _{m}\Omega BQ\mathscr{A}.
\end{displaymath}

Let $X$ be a scheme (recall that this means a noetherian, 
separated scheme of finite Krull dimension). Then, we will 
denote by $\mathscr{M}(X)$ the exact category of coherent 
sheaves on $X$, and by $\mathscr{P}(X)$ the exact category 
of locally free coherent sheaves on $X$. We write 
\begin{align*}
K'_{m}(X)&=K_{m}(\mathscr{M}(X)), \\
K_{m}(X)&=K_{m}(\mathscr{P}(X)).
\end{align*}
In particular, $K_{0}(X)$ is the Grothendieck group of locally 
free coherent sheaves on $X$, and $K'_{0}(X)$ is the Grothendieck 
group of coherent sheaves on $X$. 

We will use many standard facts about the $K$-theory of schemes. 
For more details the reader is referred, for instance, to 
\cite{Quillen:haKt}, or \cite{Srinivas:akt}. Let us quote 
some of these facts.

\nnpar{Functoriality of $K$-theory.} 
The first property we want to quote is the functoriality of 
$K$-theory. Let $F:\mathscr{A}\longrightarrow\mathscr{B}$ be 
an exact functor between exact categories. Then, there is an 
induced morphism 
\begin{displaymath}
\Omega BQ\mathscr{A}\longrightarrow\Omega BQ\mathscr{B},
\end{displaymath}
and hence morphisms
\begin{displaymath}
K_{m}(\mathscr{A})\longrightarrow K_{m}(\mathscr{B}).
\end{displaymath}
From this one can derive:

\begin{proposition}
\label{prop:4}
If $f:X\longrightarrow X'$ is a morphism of schemes, then the 
inverse image of locally free sheaves induces morphisms  
\begin{displaymath}
f^{\ast}:K_{m}(X')\longrightarrow K_{m}(X).
\end{displaymath}
With these morphisms, $K_{\ast}$ is a contravariant functor 
between the category of schemes to the category of graded 
abelian groups.
\hfill $\square$
\end{proposition}

\begin{proposition}
\label{prop:5}
If $f:Y\longrightarrow Y'$ is a proper morphism of schemes, then
the higher direct images of coherent sheaves induce morphisms
\begin{displaymath}
f_{\ast}:K'_{m}(Y)\longrightarrow K'_{m}(Y').
\end{displaymath}
With these morphisms, $K'_{\ast}$ is a covariant functor between 
the category of schemes with proper morphisms and the category 
of graded abelian groups. 

\hfill $\square$
\end{proposition}

\nnpar{Localization.} 
Given exact functors $\mathscr{A}\longrightarrow\mathscr{B}
\longrightarrow\mathscr{C}$ such that the composition $\mathscr{A}
\longrightarrow\mathscr{C}$ is zero and the induced maps $\Omega 
BQ(\mathscr{A})\longrightarrow\Omega BQ(\mathscr{B})\longrightarrow
\Omega BQ(\mathscr{C})$ form a fibration up to homotopy, then there 
is a long exact sequence 
\begin{displaymath}
\dots\longrightarrow K_{m}(\mathscr{A})\longrightarrow K_{m}
(\mathscr{B})\longrightarrow K_{m}(\mathscr{C})\overset{\delta}
{\longrightarrow}K_{m-1}(\mathscr{A})\longrightarrow\dots
\end{displaymath}
This applies to the case when $\mathscr{A}$ is a Serre subcategory 
of $\mathscr{B}$, and $\mathscr{C}=\mathscr{B}/\mathscr{A}$. 
In particular, this can be used to derive the following result
(\cite{Quillen:haKt}, \S7, 3.2):

\begin{proposition}[Localization]
\label{prop:eskts} 
Let $X$ be a scheme, $Y\subset X$ a closed subscheme, and 
$U=X\setminus Y$. Then, there is a long exact sequence
\begin{displaymath}
\dots\longrightarrow K_{m}'(Y)\longrightarrow K_{m}'(X)\longrightarrow
K_{m}'(U)\overset{\delta}{\longrightarrow}K_{m-1}'(Y)\longrightarrow
\dots 
\end{displaymath}
\hfill $\square$
\end{proposition}

\nnpar{$K$ and $K'$ of regular schemes.} 
Another important result is the comparison between the groups $K$ 
and $K'$ (\cite{Quillen:haKt}, \S7).

\begin{proposition}
\label{prop:6}
If $X$ is a regular scheme, then the natural morphism
\begin{displaymath}
K_{m}(X)\longrightarrow K_{m}'(X)
\end{displaymath}
is an isomorphism.
\hfill $\square$
\end{proposition}

\nnpar{$K$-theory with support.} 
Let $X$ be a scheme, $Y\subset X$ a closed subscheme, and 
$U=X\setminus Y$. Then, the \emph{$K$-theory groups of $X$ with 
support in $Y$} are defined by
\begin{displaymath}
K^{Y}_{m}(X)=\pi_{m+1}(\text{homotopy fiber}(BQ\mathscr{P}(X)
\longrightarrow BQ\mathscr{P}(U))).
\end{displaymath}
By definition, there is a long exact sequence
\begin{displaymath}
\dots\longrightarrow K_{m}^{Y}(X)\longrightarrow K_{m}(X)\longrightarrow
K_{m}(U)\overset{\delta}{\longrightarrow}K_{m-1}^{Y}(X)\longrightarrow 
\dots 
\end{displaymath}
Combining proposition \ref{prop:6}, the above exact sequence and
proposition \ref{prop:eskts}, we obtain 

\begin{proposition}[Purity]
\label{prop:1} 
If $X$ is a regular scheme and $Y\subset X$ a closed subscheme, 
then there is a natural isomorphism
\begin{displaymath}
K^{Y}_{m}(X)\longrightarrow K'_{m}(Y).
\end{displaymath}
\hfill $\square$
\end{proposition} 

\noindent
The purity property clearly implies the following excision property:

\begin{proposition}
\label{prop:3}
Let $X$ be a regular scheme, $Y$ a closed subscheme of $X$, and $V$ 
an open subscheme of $X$ satisfying $Y\subset V$. Then, the restriction 
map
\begin{displaymath}
K^{Y}_{m}(X)\longrightarrow K^{Y}_{m}(V)
\end{displaymath}
is an isomorphism.
\hfill $\square$
\end{proposition}

\subsection{$K$-theory as generalized sheaf cohomology}
\label{sec:ktgc}

The $K$-theory groups can also be realized as generalized cohomology
groups \cite{BrownGersten:haKt}. This allows us to treat $K$-theory
and cohomology in a uniform way, and is the basis of the construction
by Gillet \cite{Gillet:RRhK} of characteristic classes for higher
$K$-theory. In this section we will briefly recall the definition 
of generalized sheaf cohomology and how to see higher $K$-theory as
generalized sheaf cohomology. We will follow \cite{GilletSoule:fKt}
and appendix B of \cite{HuberWildeshaus:cmp}. There is also a short
account in \cite{Jeu:Zc}. The reader may also consult the book
\cite{Jardine:Gect} for a slightly different point of view. 

\nnpar{Homotopy theory of spaces.} 
Let us fix a regular noetherian scheme $S$ of finite Krull dimension.
Let $C$ be one of the following sites: $\zar(S)$ the small Zariski
site over $S$, $\ZAR(S)$ the big Zariski site of all schemes of finite 
type over $S$, or $\ZAR$ the big Zariski site of all noetherian schemes 
of bounded Krull dimension. Let $T=T(C)$ be the topos of sheaves over 
the site $C$ as in \cite{SGA4}. Let us denote by $sT$ the category of 
pointed simplicial objects in $T$. The elements of $sT$ will be called
spaces. Given a scheme $X$ in the site $C$, we can consider $X$ as a 
sheaf by
\begin{displaymath}   
U\longmapsto \Hom_{C}(U,X).
\end{displaymath}
We will also denote by $X$ the corresponding constant simplicial
object pointed by a disjoint base point.

\begin{definition}
\label{def:7}
A space is said \emph{to be constructed from schemes} if all 
components are representable by a scheme in the site plus a 
disjoint base point. If $P$ is a property of schemes, we will 
say that a space $X$ constructed from schemes satisfies $P$, if 
all the schematic parts of $X$ satisfy $P$. 
\end{definition}

The following result, due to Joyal, is fundamental for passing to the
homotopy 
category of spaces. A published proof in the non-pointed case can 
be found in \cite{Jardine:sp} (see \cite{GilletSoule:fKt}, and
\cite{HuberWildeshaus:cmp}).  

\begin{proposition}
The category $sT$ is a pointed closed model category in the sense 
of Quillen \cite{Quillen:ha}.
\hfill $\square$ 
\end{proposition}

A part of the definition of a closed model category comprises the
concepts  
of weak equivalence, fibrations and cofibrations, and, in particular, 
the concept of fibrant objects and cofibrant objects. We refer to
\cite{HuberWildeshaus:cmp}, appendix B, for the definition of these
concepts in the present setting.
 
Let us denote by $\Ho(sT)$ the homotopy category associated to the
closed model category $sT$. If $X$ and $Y$ are spaces, we denote by
$[X,Y]$ the morphism in $\Ho(sT)$. If $Y$ is fibrant, it is just the
homotopy classes of morphisms. If $Y$ is not fibrant, then $[X,Y]=
[X,\widetilde Y]$, where $\widetilde Y$ is a fibrant space weakly 
equivalent to $Y$. We will denote by $SX$ and $\Omega X$ the suspension 
and the loop space of a space $X$, respectively. The loop space functor 
$\Omega$ is the right adjoint functor of the suspension functor $S$
(see \cite{Quillen:ha}).

\nnpar{Generalized cohomology.}
\begin{definition}
Let $A$ be a space and let $X\in C$. Then, the \emph{cohomology 
space with coefficients in $A$} is defined for $m\ge 0$ by
\begin{displaymath}
H^{-m}_{sT}(X,A)=[S^{m}\wedge X,A],
\end{displaymath}
where $S^{m}\wedge X$ is the smash product between the pointed
$m$-dimensional sphere $S^{m}$ and $X$. Note that $S^{m}\wedge 
X$ is canonically isomorphic to the $m$-fold suspension $S^{m}X=
S\overset{m}{\dots}SX$. 
\end{definition}

\begin{remark}
If $X$ is a scheme in the site $\ZAR$ and $A$ is an element of
$sT(\ZAR)$, we will denote also by $A$ the restriction of $A$ to
$\ZAR(X)$ or to $\zar(X)$. By the argument at the beginning of 
the proof of \cite{GilletSoule:fKt}, proposition 5, the cohomology
groups $H^{-m}(X,A)$ are the same in the three sites.  
\end{remark}

\nnpar{Pseudo-flasque presheaves.} 
If $A$ is a fibrant space, the generalized cohomology can be computed 
as homotopy groups 
\begin{displaymath}
H^{-m}_{sT}(X,A)=[S^{m}\wedge X,A]=\pi_{m}(A(X))=\pi_{m}(\Hom(X,A)).
\end{displaymath}
More generally, the above equality is a property of pseudo-flasque
presheaves. 

\begin{definition}
A presheaf $A$ is \emph{pseudo-flasque}, if $A(U)$ is a fibrant
simplicial set for all $U$, $A(\emptyset)$ is contractible, and 
for each pair of open sets the diagram     
\begin{displaymath}
\begin{CD}
A(U\cup V)@>>>A(U) \\
@VVV @VVV \\
A(V)@>>>A(U\cap V)
\end{CD}
\end{displaymath}
is homotopically cartesian.
\end{definition}

\begin{proposition}[Brown-Gersten \cite{BrownGersten:haKt},
Gillet-Soul\'e \cite{GilletSoule:fKt}]
\label{prop:13}  
Let $A$ be a pseudo-flasque presheaf, and let $A'$ be the 
associated sheaf. Then, we have for any scheme $X$ in $\ZAR$
\begin{displaymath}
H^{-m}(X,A')=\pi_{m}(A(X)).
\end{displaymath}
\hfill $\square$
\end{proposition}

\nnpar{$K$-theory as generalized sheaf cohomology.} 
Let us denote by $\KK$ the pointed simplicial sheaf $\mathbb{Z}
\times\mathbb{Z}_{\infty}\BGL$ (see \cite{GilletSoule:fKt}, 3.1), 
and by $\KK^{N}$ the pointed simplicial sheaf $\mathbb{Z}\times
\mathbb{Z}_{\infty}\BGL_{N}$. If $X$ is a space, we define for 
$m\ge 0$ its \emph{$m$-th $K$-theory group} by $H^{-m}(X,\KK)$. 
We also define the \emph{unstable $m$-th $K$-theory group} by 
$H^{-m}(X,\KK^{N})$. In general, $H^{0}(X,\KK^{N})$ are only 
pointed spaces, but $ H^{-m}(X,\KK^{N})$ are abelian groups for 
$m\ge 1$. Following \cite{GilletSoule:fKt}, we define

\begin{definition}
We say that a space $X$ is \emph{$K$-coherent}, if the natural maps  
\begin{displaymath}
\lim_{\substack{\longrightarrow\\N}}H^{-m}(X,\KK^{N})
\longrightarrow H^{-m}(X,\KK)
\end{displaymath}
and 
\begin{displaymath}
\lim_{\substack{\longrightarrow\\N}}H^{m}(X,\pi_{-n}\KK^{N})
\longrightarrow H^{m}(X,\pi_{-n}\KK)
\end{displaymath}
are isomorphisms for all $m,n\ge 0$.
\end{definition}

The fact that the groups $H^{-m}(X,\KK)$ agree with the $K$-groups 
defined by Quillen is proved in \cite{GilletSoule:fKt}. They use the 
argument of \cite{BrownGersten:haKt} together with some modifications 
due to the different definitions. The main ingredients are the ``$+=Q$''
result of Quillen \cite{Quillen:haKt}, proposition \ref{prop:13}, and 
the fact that the homotopy category used in \cite{BrownGersten:haKt} 
agrees with the one considered in \cite{GilletSoule:fKt}.

\begin{proposition}
Let $X$ be a scheme in the site $C$. Then, there are canonical 
isomorphisms
\begin{displaymath}
K_{m}(X)\longrightarrow H^{-m}(X,\KK).
\end{displaymath}
Moreover, $X$ is $K$-coherent.
\hfill $\square$
\end{proposition}

From the above proposition and a little homological algebra 
(see {\cite{GilletSoule:fKt}, 3.2.3, \cite{Jeu:Zc}, 2.1 and
\cite{HuberWildeshaus:cmp}, B.2.3}) one can derive:

\begin{proposition}
Let $X$ be a noetherian space constructed from schemes of finite
Krull dimension (see definition \ref{def:7}) that is degenerate 
above some simplicial degree, i.e., there exists $N$ such that
$X=\sk_{N}X$, where $\sk_{N}X$ is the $N$-th skeleton of $X$. 
Then, $X$ is $K$-coherent.
\hfill $\square$
\end{proposition}

\nnpar{The mapping cone.} 
Let $f:X\longrightarrow Y$ be a morphism of spaces. We can define 
\emph{the mapping cone $C(X,Y)=C(f)$ of $f$} as the space 
\begin{displaymath}
C(X,Y)=\left(Y\amalg X\times I)\right/\sim\,\,,
\end{displaymath}
where $I$ is the simplicial unit interval and $\sim$ is the
equivalence relation generated by
\begin{align*}
(x,1)&\sim f(x), \\
(x,0)&\sim\ast, \\
(\ast,t)&\sim\ast,
\end{align*}
with $\ast$ denoting the distinguished point in all spaces.
Then, we write
\begin{displaymath}
H^{-m}(Y,X,\KK)=H^{-m}(C(f),\KK).
\end{displaymath}
There is a long exact sequence 
\begin{displaymath}
\dots \longrightarrow H^{-m}(Y,\KK)\longrightarrow H^{-m}(X,\KK)
\overset{\delta }{\longrightarrow}H^{-m+1}(Y,X,\KK)\longrightarrow 
\dots 
\end{displaymath}
The connection morphism of this exact sequence is defined as
follows. Let  
\begin{displaymath}
\alpha:S^{m}\wedge X\longrightarrow\KK
\end{displaymath}
be an element of $H^{-m}(X,\KK)=[S^{m}\wedge X,\KK]$. Then, $\delta
(\alpha)$ is given as the following composition of maps
\begin{align*}
S^{m-1}\wedge C(X,Y)&\cong C(S^{m-1}\wedge X,S^{m-1}\wedge Y)
\longrightarrow \\ 
&C(S^{m-1}\wedge X,\ast)\cong S^{1}\wedge S^{m-1}\wedge X=
S^{m}\wedge X\overset{\alpha}{\longrightarrow }\KK,
\end{align*}
where $\ast$ is the simplicial point; here we also use the fact that
the smash 
product and the mapping cone commute, and that $C(X,\ast)\cong S^{1}
\wedge X$ is the suspension. 

\nnpar{$K$-theory with support.} 
Let $X$ be a regular noetherian scheme of finite Krull dimension, 
$f:U\longrightarrow X$ an open immersion, and set $Y=X\setminus U$. 
The mapping cone of $f$ is a space $C(f)\in s\ZAR$. It is clear 
that there are canonical isomorphisms 
$$
K_{m}^{Y}(X)\longrightarrow K_{m}(C(f)),\ m\ge 0.
$$

\nnpar{Loday product in $K$-theory.} 
Let us discuss the existence of a product in $K$-theory and its 
basic properties. The choice of any bijection $\mathbb{N}\times
\mathbb{N}\longrightarrow\mathbb{N}$ induces a map \cite{Loday:Ktrg} 
(see \cite{HuberWildeshaus:cmp}, p. 103, for details in this setting)
\begin{displaymath}
(\mathbb{Z}\times\mathbb{Z}_{\infty}\BGL_{N})\wedge
(\mathbb{Z}\times\mathbb{Z}_{\infty}\BGL_{N})\longrightarrow 
\mathbb{Z}\times\mathbb{Z}_{\infty}\BGL.
\end{displaymath}
Hence, if $X$ and $Y$ are $K$-coherent spaces, there is a well-defined
external product  
\begin{displaymath}
H^{-m}(X,\KK)\times H^{-n}(Y,\KK)\longrightarrow 
H^{-m-n}(X\wedge Y,\KK)
\end{displaymath}
via
\begin{displaymath}
S^{m+n}\wedge X\wedge Y\cong S^{m}\wedge X\wedge S^{n}\wedge Y
\longrightarrow\KK^{N}\wedge\KK^{N}\longrightarrow\KK.
\end{displaymath}
When $X=Y$, composing with the pull-back of the diagonal map
$X\longrightarrow X\wedge X$, we obtain an associative and 
graded commutative ring structure (possibly without unit) in 
$H^{-\ast}(X,\KK)$. Clearly the product is functorial.

Let $X$, $Y$ and $Z$ be $K$-coherent spaces, and let $f:X
\longrightarrow Y$ be a morphism of spaces. Having defined 
the product for spaces and using the fact that the mapping cone 
commutes with the smash product, we observe that there is 
a product
\begin{displaymath}
H^{-m}(Y,X,\KK)\times H^{-n}(Z,\KK)\longrightarrow 
H^{-m-n}(Y\wedge Z,X\wedge Z,\KK).
\end{displaymath}
The following result is a version of \cite{Gillet:RRhK}, 
corollary 7.14. 

\begin{proposition}
\label{prop:2}
Let $X$, $Y$ and $Z$ be $K$-coherent spaces, and let $f:
X\longrightarrow Y$ be a morphism of spaces. If $\alpha 
\in H^{-m}(X,\KK)$ and $\beta\in H^{-n}(Z,\KK)$, then the
Loday product is compatible with long exact sequences, i.e.,
we have
\begin{align*}
\delta(\alpha\cdot\beta)&=\delta(\alpha)\cdot\beta\in 
H^{-m-n+1}(Y\wedge Z,X\wedge Z,\KK), \\  
\delta(\beta\cdot\alpha)&=(-1)^{n}\beta\cdot\delta(\alpha)\in 
H^{-m-n+1}(Z\wedge Y,Z\wedge X,\KK). 
\end{align*}
\end{proposition}
\begin{proof}
The first equality follows from the commutative diagram 
\begin{displaymath}
\xymatrix{S^{m+n}\wedge C(X,Y)\wedge Z\ar[d]\ar[r]&
S^{m}\wedge C(X,Y)\wedge S^{n}\wedge Z\ar[d] \\
C(S^{m+n}\wedge X\wedge Z,S^{m+n}\wedge Y\wedge Z)\ar[d]&
C(S^{m}\wedge X,S^{m}\wedge Y)\wedge S^{n}\wedge Z\ar[d] \\
C(S^{m+n}\wedge X\wedge Z,\ast)\ar[d]&
C(S^{m}\wedge X,\ast)\wedge S^{n}\wedge Z\ar[d] \\
S^{1}\wedge S^{m+n}\wedge X\wedge Z\ar[r]&
S^{1}\wedge S^{m}\wedge X\wedge S^{n}\wedge Z.}
\end{displaymath}
The second statement follows from a similar diagram that involves 
the isomorphism 
\begin{displaymath}
S^{n}\wedge S^{1}\wedge S^{m}\longrightarrow 
S^{1}\wedge S^{n}\wedge S^{m}.
\end{displaymath}
But identifying both sides with $S^{m+n+1}$, the above
isomorphism has degree $(-1)^{n}$.
\end{proof}

\nnpar{Products in $K$-theory with support.} 
We will now construct the product of $K$-theory with support. 
Let $X$ be a scheme in the site, let $U$, $V$ be open subsets 
of $X$, and set $Y=X\setminus U$, $Z=X\setminus V$. Then, 
there is a pairing 
\begin{align*}
K_{m}^{Y}(X)\otimes K_{n}^{Z}(X)=&
H^{-m}(C(U,X),\KK)\otimes H^{-n}(C(V,X),\KK)\longrightarrow \\ 
&H^{-m-n}(C(U,X)\wedge C(V,X),\KK).
\end{align*}
Using the fact that the smash product and the mapping cone commute, we obtain
that $C(U,X)\wedge C(V,X)$ is weakly equivalent to either of the
spaces  
\begin{align*}
&C(C(U\wedge V,U\wedge X),C(X\wedge V,X\wedge X)),\,\text{or} \\  
&C(C(U\wedge V,X\wedge V),C(U\wedge X,X\wedge X)).
\end{align*}
We will denote any of them as $C(U\wedge V;X\wedge V,U\wedge X;
X\wedge X)$. The diagonal map induces a morphism of diagrams
\begin{displaymath}
\begin{CD}
U\cap V@>>>U \\
@VVV @VVV \\
V@>>>X
\end{CD}
\longrightarrow
\begin{CD}
U\wedge V@>>>U\wedge X \\
@VVV @VVV \\
X\wedge V@>>>X\wedge X.
\end{CD}
\end{displaymath}
Since $\{U,V\}$ is an open covering of $U\cup V$, we find 
that $C(U\cap V;U,V;X)$ is weakly equivalent to $C(U\cup V,X)$. 
Therefore, we obtain a pairing 
\begin{displaymath}
K_{m}^{Y}(X)\otimes K_{n}^{Z}(X)\overset{\cup}
{\longrightarrow}H^{-m-n}(C(U\cup V,X),\KK)=
K^{Y\cap Z}_{m+n}(X).  
\end{displaymath}
A consequence of proposition \ref{prop:2} is the following

\begin{proposition}
Let $X$ be a scheme in the site, let $U$, $V$ be open subsets 
of $X$, and set $Y=X\setminus U$, $Z=X\setminus V$. Then, the 
following diagram commutes: 
\begin{displaymath}
\begin{CD}
K_{m}(U)\otimes K_{n}^{Z}(X)@>\delta\otimes\Id>>  
K_{m-1}^{Y}(X)\otimes K_{n}^{Z}(X) \\
@V\cup VV @VV\cup V \\
K_{m+n}^{U\cap Z}(X\setminus(Y\cap Z))@>\delta>>
K_{m+n-1}^{Y\cap Z}(X);    
\end{CD}
\end{displaymath}
here we have identified $K_{m+n}^{U\cap Z}(U)$ with $K_{m+n}^
{U\cap Z}(X\setminus(Y\cap Z))$ using excision, and $\delta $ 
is the connection morphism. 
\hfill $\square$
\end{proposition}

\nnpar{Direct images and the projection formula.} 
Let $f:X\longrightarrow X'$ be a proper morphism of regular 
schemes in the site. Let $Z\subset X$ be a closed subscheme. 
The isomorphisms $K'_{m}(Z)\longrightarrow K^{Z}_{m}(X)$ and
$K'_{m}(f(Z))\longrightarrow K^{f(Z)}_{m}(X')$ together with 
the direct image of the $K'$-groups induce direct image 
morphisms 
\begin{displaymath}
f_{\ast}:K^{Z}_{m}(X)\longrightarrow K^{f(Z)}_{m}(X').
\end{displaymath}

\noindent
These morphisms satisfy the projection formula.
  
\begin{proposition}
Let $f:X\longrightarrow X'$ be a proper morphism of regular 
schemes. Let $Z\subset X$, and $Z'\subset X'$ be closed subschemes. 
Then, for $\alpha\in K_{m}^{Z'}(X')$ and $\beta\in K_{n}^{Z}(X)$, 
we have 
\begin{displaymath}
f_{\ast}(f^{\ast}(\alpha)\cup\beta)=\alpha\cup f_{\ast}(\beta)
\in K_{m+n}^{f(Z)\cap Z'}(X').
\end{displaymath}
\hfill $\square$
\end{proposition}

\nnpar{Gersten-Quillen spectral sequence.} 
Let $X$ be a regular scheme in the site. Let us denote by 
$X^{(p)}$ the set of $p$-codimensional points of $X$. If 
$x\in X^{(p)}$, we will denote
\begin{displaymath}
K^{x}_{m}(X)=\lim_{\longrightarrow}K^{\overline{\{x\}}
\cap U}_{m}(U),
\end{displaymath}
where the limit is taken over all open subsets $U$ of $X$ 
containing the point $x$. 

\begin{theorem}
\label{thm:1}
Let $X$ be a regular scheme in $\ZAR$. Then, there exists a 
spectral sequence 
\begin{displaymath}
E^{p,q}_{r}(X)\Longrightarrow K_{-p-q}(X)
\end{displaymath}
with 
$$
E^{p,q}_{1}(X)=\bigoplus_{x\in X^{(p)}}K_{-p-q}^{x}(X).
$$
Moreover, this spectral sequence is convergent (since $X$ has 
finite Krull dimension).
\hfill $\square$
\end{theorem}

Observe moreover that, using purity, we can identify $K_{-p-q}^{x}(X)$ 
with $K_{-p-q}(k(x))$, where $k(x)$ is the residue field at the point
$x$. See \cite{Gillet:RRhK} for a proof of this theorem for generalized 
cohomology.

\subsection{$\lambda $-structure in $K$-theory and absolute 
cohomology}\label{sec:absc}

\nnpar{$\lambda$-rings and $\lambda$-algebras.} 
Recall that a \emph{$\lambda$-ring with involution} is a ring $R$ 
equipped with a family of operators $\{\lambda^{k}\}_{k\ge 0}$ 
and an involution $\psi^{-1}$ satisfying  certain relations (see
\cite{SGA6}). In particular, any $\lambda$-ring has a unit. Let 
$R$ be a $\lambda $-ring with involution, and $A$ an $R$-algebra. 
We call $A$ a \emph{$\lambda$-$R$-algebra with involution}, if it 
is equipped with a family of operators $\{\lambda^{k}\}_{k\ge 1}$
and an involution $\psi^{-1}$ such that $R\oplus A$ is a $\lambda$-ring 
(see \cite{Kratzer:lsKta} and \cite{HuberWildeshaus:cmp}).  

\nnpar{$\lambda$-structure in $K$-theory.} 
For any $N\ge 1$, one can define a family of operators $\{\lambda^{k}_
{N}\}_{k\ge 1}$ (see \cite{GilletSoule:fKt}, \S 4) 
\begin{displaymath}
\lambda^{k}_{N}:\mathbb{Z}\times\mathbb{Z}_{\infty}\BGL_{N}
\longrightarrow\mathbb{Z}\times\mathbb{Z}_{\infty}\BGL\,,
\end{displaymath}
which are compatible with the inclusions $\BGL_{N-1}\longrightarrow
\BGL_{N}$. Therefore, for any space $X$, there are induced operators  
\begin{displaymath}
\lambda^{k}_{N}:H^{-m}(X,\KK^{N})\longrightarrow H^{-m}(X,\KK).
\end{displaymath}
If the space $X$ is $K$-coherent, we have induced operators 
\begin{displaymath}
\lambda^{k}:H^{-m}(X,\KK)\longrightarrow H^{-m}(X,\KK).
\end{displaymath}
In the same way one can define an involution $\psi^{-1}$.

Let $S^{0}$ be the simplicial pointed $0$-sphere. Following
\cite{HuberWildeshaus:cmp}, we will write $K_{0}(sT)=H^{0}
(S^{0},\KK)$. This is a $\lambda$-ring. If the site $C$ equals 
$\ZAR(S)$, we have $K_{0}(sT)=K_{0}(S)$. A proof of the following 
theorem can be found in \cite{GilletSoule:fKt} (see also
\cite{HuberWildeshaus:cmp}).   

\begin{theorem}
Let $X$ be a $K$-coherent space. The family of operators $\{\lambda^
{k}\}_{k\ge 1}$ together with the involution $\psi^{-1}$ turn 
$H^{-m}(X,\KK)$ into a $\lambda$-$K_{0}(sT)$-algebra with involution.
Moreover, $H^{0}(X,\KK)$ is provided with a natural augmentation
\begin{displaymath}
\epsilon:H^{0}(X,\KK)\longrightarrow H^{0}(X,\mathbb{Z}).
\end{displaymath}
\end{theorem}

\nnpar{Adams operations and $\gamma$-filtration.} 
Let $\Phi $ be the monoid of natural operations on $\lambda$-rings
with an involution. Among the operations in $\Phi$, there are the 
$\gamma$-operations $\gamma^{k}$, $k\ge 1$, and the Adams operations 
$\psi^{k}$, $k\in\mathbb{Z}\setminus\{0\}$. These operations are 
defined by the usual formulae
\begin{gather*}
\gamma^{k}(x)=\lambda^{k}(x+k-1),\qquad k\ge 1, \\
\psi^{k}-\lambda^{1}\psi^{k-1}+\dots+(-1)^{k-1}\lambda^{k-1}\psi^{1}+
(-1)^{k}\lambda^{k}=0,\qquad k\ge 1, 
\end{gather*}
and $\psi^{-k}=\psi^{-1}\circ\psi^{k}$. Since these operations do 
not involve the unit, they are also operations of $\lambda$-algebras. 
By means of the $\gamma$-operations we can define a filtration in 
the $K$-groups. 

\begin{definition}
Let $X$ be a $K$-coherent space. Then $F^{i}_{\gamma}H^{-m}(X,\KK)$ 
is defined as the subgroup of $H^{-m}(X,\KK)$ generated by the 
products
\begin{displaymath}
\gamma^{i(1)}(x_{1})\cdots\gamma^{i(r)}(x_{r}),
\end{displaymath}
where $i(1)+\dots+i(r)\ge i$, $x_{j}\in H^{-m(j)}(X,\KK)$ with
$m(1)+\dots+m(r)=m$ and $\epsilon(x_{j})=0$.
\end{definition}

\noindent
The following properties of the Adams operations are well known.

\begin{proposition}
The Adams operations respect the additive structure of $K$-theory, 
the Loday product and pull-backs. Moreover, we have $\psi^{k}\circ
\psi^{k'}=\psi^{k+k'}$.
\hfill $\square$ 
\end{proposition}

Observe that the $\lambda$-structure in $K$-theory involves only 
the product between $K_{0}$ and higher $K$-theory and not the 
Loday product. Nevertheless, in order to see that the Adams 
operations are compatible with the Loday product one can argue 
as in \cite{Kratzer:lsKta} (see also \cite{Jeu:Zc}).  

\nnpar{Eigenspaces for the Adams operations.} 
Let us write 
\begin{displaymath}
H^{-m}(X,\KK)^{(i)}=\{x\in H^{-m}(X,\KK)\otimes\QQ\mid\psi^{k}(x)=
k^{i}x,\,\forall k\in\mathbb{Z}\setminus\{0\}\}.
\end{displaymath}
We will also write $K_{m}(X)^{(i)}=H^{-m}(X,\KK)^{(i)}$. The 
following proposition is proven in \cite{GilletSoule:fKt}.

\begin{proposition}
\label{prop:14}
Let $X$ be a $K$-coherent space of dimension at most $d$. Then, 
for $m\ge 0$, there is a decomposition 
\begin{displaymath}
H^{-m}(X,\KK)\otimes\QQ=\sum_{i=\alpha}^{m+d}H^{-m}(X,\KK)^{(i)}, 
\end{displaymath}
where $\alpha=\min(m,2)$. Moreover, we have
\begin{displaymath}
H^{-m}(X,\KK)^{(i)}=\Gr_{\gamma }^{i}H^{-m}(X,\KK)\otimes\QQ.
\end{displaymath}
\hfill $\square$
\end{proposition}

\nnpar{Riemann-Roch without denominators.} 
The Riemann-Roch theorem without denominators studies the effect of 
closed immersions on the operations in $K$-theory. Let $X$ and $Y$ 
be regular schemes in $\ZAR(S)$, let $j:Y\longrightarrow X$ be a 
closed immersion, and let us denote by $N=N_{X/Y}$ the class of the 
normal bundle of $Y$ in $X$ in $K_{0}(Y)$. For any operation $\tau 
\in\Phi$, let us write $\tau(N,x)$ as in \cite{Soule:oKa}, \S 4. 
Let $Z$ be a closed subset of $Y$. Using the identifications
$K_{m}'(Z)=K^{Z}_{m}(Y)$ and $K_{m}'(Z)=K^{Z}_{m}(X)$, we obtain 
an isomorphism 
\begin{displaymath}
j_{\ast}:K^{Z}_{m}(Y)\longrightarrow K^{Z}_{m}(X).
\end{displaymath}

\begin{theorem}[Soul\'e \cite{Soule:oKa}]
\label{thm:rrwd}
For $x\in K^{Z}_{m}(Y)$ and $\tau\in\Phi$, we have 
\begin{displaymath}
\tau(j_{\ast}(x))=j_{\ast}(\tau(N,x)).
\end{displaymath}
\hfill $\square$
\end{theorem}

\begin{remark}
This result is generalized in \cite{Jeu:Zc} and \cite{HuberWildeshaus:cmp}
to a certain class of closed immersions of spaces. Nevertheless in 
these papers the schemes are assumed to be smooth over the base $S$.
\end{remark}

\nnpar{Direct images and the gamma filtration.} 
\begin{proposition}
Let $f:X\longrightarrow X'$ be a proper morphism between regular
quasi-projective schemes of $\ZAR(S)$ of dimension $d$ and $d'$,
respectively. Let $Z$ be a closed subscheme of $X$. Then, there is
an inclusion  
\begin{displaymath}
f_{\ast}(F_{\gamma }^{j}K^{Z}_{m}(X)\otimes\QQ)\subset 
F_{\gamma}^{j+d'-d}K^{f(Z)}_{m}(X')\otimes\QQ.
\end{displaymath}
\end{proposition}
\begin{proof}
The Riemann-Roch theorem (\cite{Soule:oKa}, 7.2) states that there 
are increasing filtrations $F_{\cdot}$ of $K'_{m}(Z)\otimes\QQ$ and
$K'_{m}(f(Z))\otimes\QQ$ satisfying $f_{\ast}(F_{j})\subset F_{j}$. 
These filtrations are defined as follows. One chooses a closed 
immersion $i:Z\longrightarrow M$, where $M$ is an equidimensional 
scheme of dimension $n$, smooth and surjective over the base $S$. 
Then, there is an isomorphism $i_{\ast}:K_{m}(Z)\longrightarrow 
K^{Z}_{m}(M)$, and we have $F_{j}K'_{m}(Z)\otimes\QQ=i_{\ast}^{-1}
(F_{\gamma}^{n-j}K^{Z}_{m}(M))$. We can choose $M$ such that the 
inclusion $i$ factors as $Z\overset{i'}{\hookrightarrow}X\overset{j}
{\hookrightarrow}M$. Using the Riemann-Roch theorem without  
denominators (as in the proof of the independence of $M$ in 
\cite{Soule:oKa}, 7.2, 3) one obtains that $j_{\ast}(F^{j}_{\gamma}
K^{Z}_{m}(X)\otimes\QQ)=F^{j-d+n}_{\gamma}K^{Z}_{m}(M)$. This leads
to $F^{j}_{\gamma}K^{Z}_{m}(X)\otimes\QQ=F_{d-j}K'_{m}(Z)\otimes\QQ$ 
which implies the claim.  
\end{proof}   

\nnpar{Absolute cohomology.} 
The graded pieces for the $\gamma$-filtration tensored by $\QQ$ or,
put another way, the eigenspaces for the Adams operations, form a
cohomology theory. By the Riemann-Roch theorem for higher $K$-theory
given in \cite{Gillet:RRhK}, it is a universal cohomology. One can 
define it for coherent spaces satisfying certain technical conditions. 
Since we will not need it in the sequel, we will define it only for 
a pair of schemes.

\begin{definition}
Let $X$ be a regular scheme in the site and let $Y$ be a closed subset
of $X$. Then, the \emph{absolute cohomology groups of $X$ with support
on $Y$} are defined by
\begin{displaymath}
H^{2i-m}_{\mathcal{A},Y}(X,\mathbb{Q}(i))=K_{m}^{Y}(X)^{(i)}.
\end{displaymath}
\end{definition}

Let us summarize the basic properties of absolute cohomology. The
next theorem follows easily from the properties of $K$-theory and 
the Adams operations. 

\begin{theorem}
\label{thm:11}
\begin{enumerate}
\item[ (i)] 
Let $X$ be a regular scheme in $\ZAR$, and let $Y$, $Z$ be closed
subsets of $X$. Then, the product in $K$-theory induces a product
\begin{displaymath}
H^{m}_{\mathcal{A},Y}(X,\mathbb{Q}(i))\otimes
H^{n}_{\mathcal{A},Z}(X,\mathbb{Q}(j))\overset
{\cup}{\longrightarrow}H^{m+n}_{\mathcal{A},Y\cap Z}
(X,\mathbb{Q}(i+j)).
\end{displaymath}
\item[(ii)] 
Let $f:X'\longrightarrow X$ be a morphism of regular schemes
in $\ZAR$. Let $Y$ be a closed subset of $X$, and $Y'$ a closed
subset of $X'$ satisfying $f^{-1}(Y)\subset Y'$. Then, there is
pull-back morphism 
\begin{displaymath}
f^{\ast}:H^{n}_{\mathcal{A},Y}(X,\mathbb{Q}(i))\longrightarrow 
H^{n}_{\mathcal{A},Y'}(X',\mathbb{Q}(i))
\end{displaymath}
respecting the multiplicative structure and turning absolute
cohomology into a contravariant functor.
\item[(iii)] 
Let $X$ be a regular scheme in $\ZAR$, and $Y$ closed a subset 
of $X$. Then, there is a long exact sequence      
\begin{displaymath}
\dots\longrightarrow H^{n}_{\mathcal{A}}(X,\mathbb{Q}(i))
\longrightarrow H^{n}_{\mathcal{A}}(X\setminus Y,\mathbb{Q}(i))
\overset{\delta}{\longrightarrow} H^{n+1}_{\mathcal{A},Y}
(X,\mathbb{Q}(i))\longrightarrow \dots
\end{displaymath}
\item[(iv)]\label{part:abscomm} 
Let $X$ be a regular scheme in $\ZAR$, let $Y$, $Z$ be closed 
subsets of $X$, and set $U=X\setminus Y$. Then, the following 
diagram commutes:
\begin{displaymath}
\xymatrix
{H^{m}_{\mathcal{A}}(U,\mathbb{Q}(p))\otimes H^{n}_{\mathcal{A},Z}
(X,\mathbb{Q}(q))\ar[r]^-{\delta\otimes\Id}\ar[d]_{\cup}\!&
\!H^{m+1}_{\mathcal{A},Y}(X,\mathbb{Q}(p))\otimes H^{n}_{\mathcal{A},Z}
(X,\mathbb{Q}(q))\ar[d]^{\cup} \\
H^{m+n}_{\mathcal{A},U\cap Z}(X\setminus(Y\cap Z),\mathbb{Q}(p+q)) 
\ar[r]^-{\delta}\!&\!H^{m+n+1}_{\mathcal{A},Y\cap Z}(X,\mathbb{Q}(p+q)).}
\end{displaymath}
\item[(v)]
Let $f:X\longrightarrow X'$ be a proper morphism of regular,
equidimensional, quasi-projective schemes in $\ZAR(S)$ of dimension 
$d$ and $d'$, respectively. Let $Z$ be a closed subscheme of $X$. 
Then, there are morphisms        
\begin{displaymath}
f_{\ast}:H^{n}_{\mathcal{A},Z}(X,\mathbb{Q}(p))\longrightarrow 
H^{n-2d+2d'}_{\mathcal{A},f(Z)}(X',\mathbb{Q}(p-d+d')).
\end{displaymath}
If $g:X'\longrightarrow X''$ is another proper morphism, then
$(g\circ f)_{\ast}=g_{\ast}\circ f_{\ast}$. Moreover, if $Z'$ is 
a closed subscheme of $X'$, $\alpha\in H^{m}_{\mathcal{A},Z'}
(X',\mathbb{Q}(p))$, and $\beta\in H^{n}_{\mathcal{A},Z}(X,
\mathbb{Q}(q))$, we have
\begin{displaymath}
f_{\ast}(f^{\ast}(\alpha)\cup \beta)=\alpha\cup f_{\ast}(\beta)
\in H^{m+n-2d+2d'}_{\mathcal{A},f(Z)\cap Z'}(X',\mathbb{Q}(p+q-d+d')).
\end{displaymath}
\end{enumerate}
\hfill $\square$
\end{theorem}

\nnpar{The Adams operations and the Gersten-Quillen spectral 
sequence.} 
Since the Adams operations are defined on the level of sheaves of
simplicial spaces, they induce morphisms of Gersten-Quillen spectral
sequences. Therefore, after tensoring with $\mathbb{Q}$, the Quillen 
spectral sequence splits as a direct sum of spectral sequences
\begin{displaymath}
E^{p,q}_{r}(X)^{(i)}\Longrightarrow K_{-p-q}(X)^{(i)},
\end{displaymath}
where $E^{p,q}_{r}(X)^{(i)}\subset E^{p,q}_{r}(X)_{\mathbb{Q}}$ 
is the eigenspace of $\psi^{k}$ of eigenvalue $k^{i}$. By the 
Riemann-Roch theorem without denominators (see theorem \ref{thm:rrwd}, 
or \cite{Soule:oKa}, theorem 4), we have for $x\in X^{(p)}$  
\begin{displaymath}   
K^{x}_{m}(X)^{(i)}=K_{m}(k(x))^{(i-p)}.
\end{displaymath}
Therefore, we have the equality
\begin{displaymath} 
E^{p,q}_{1}(X)^{(i)}=\bigoplus_{x\in X^{(p)}}K_{-p-q}(k(x))^{(i-p)}.
\end{displaymath}
Since we have for a field $k$ (see \cite{Soule:oKa})
\begin{align}
K_{0}(k)_{\mathbb{Q}}&=K_{0}(k)^{(0)},\notag \\
K_{1}(k)_{\mathbb{Q}}&=K_{1}(k)^{(1)},\label{eq:3} \\
K_{2}(k)_{\mathbb{Q}}&=K_{2}(k)^{(2)},\notag
\end{align}
we obtain

\begin{proposition}
\label{prop:15} 
Let $X$ be a regular scheme in the site. Then, the lines $p=-q$ 
and $p=-q-1$ of the spectral sequence $E^{p,q}_{r}(X)_{\mathbb{Q}}$ 
degenerate at the term $E_{2}$. Thus there are natural isomorphisms
\begin{align*}
E^{p,-p}_{2}(X)_{\mathbb{Q}}&\cong K_{0}(X)^{(p)}, \\
E^{p-1,-p}_{2}(X)_{\mathbb{Q}}&\cong K_{1}(X)^{(p)}.
\end{align*}  
\end{proposition}
\begin{proof}
By \eqref{eq:3} we have $E^{p,-p}_{1}(X)_{\QQ}=E^{p,-p}_{1}(X)^{(p)}$. 
Again by \eqref{eq:3} we have for $r\ge 2$
\begin{displaymath}
E^{p-r,-p+r-1}_{1}(X)^{(p)}=\bigoplus_{x\in X^{(p-r)}}
K_{1}(k(x))^{(r)}=0. 
\end{displaymath}
This leads to the isomorphism $E^{p,-p}_{2}(X)_{\mathbb{Q}}\cong
K_{0}(X)^{(p)}$. The same argument applies to the line $p=-q-1$.
\end{proof}

Note that in \cite{Soule:oKa} it is also proven that the line
$p=-q-2$ degenerates at the term $E_{2}$. A consequence of 
proposition \ref{prop:15} is the following

\begin{proposition}
\label{prop:absspec} 
Let $X$ be a regular scheme in the site. Then, there are natural
isomorphisms 
\begin{align*}
E^{p,-p}_{2}(X)_{\mathbb{Q}}&\cong 
H^{2p}_{\mathcal{A}}(X,\mathbb{Q}(p)), \\
E^{p-1,-p}_{2}(X)_{\mathbb{Q}}&\cong
H^{2p-1}_{\mathcal{A}}(X,\mathbb{Q}(p)).
\end{align*}
\hfill $\square$
\end{proposition}

\subsection{Chow Groups of regular schemes}
\label{sec:chow-groups}
\nopagebreak
\nnpar{$K$-chains and Chow groups.} 
Let ${X}$ be a regular scheme in $\ZAR$. Let us denote
$$
R^{q}_p(X)=E^{q,-p}_{1}(X)=\bigoplus_{x\in X^{(q)}}K_{p-q}(k(x)),
$$
where $E^{q,-p}_{1}(X)$ are the terms of the Gersten-Quillen
spectral sequence.

Recall that $K_{0}(k(x))=\mathbb{Z}$, and $K_1(k(x))=k^{*}(x)$. 
Therefore, the group $R^{p}_{p}(X)$ equals the group ${\rm Z}^{p}(X)$ 
of $p$-codimensional cycles on $X$. Any element $f\in R^{p-1}_{p}(X)$ 
can be written as $f=\sum_{x\in X^{(p-1)}}f_{x} $ with $f_{x}\in 
k^{*}(x)$. The elements of $R^{p-1}_{p}(X)$ are called  $K_{1}$-chains. 
As usual we denote by $\CH^{p}(X)$ the $p$-th Chow group of $X$. In
\cite{Quillen:haKt}, section 5, it is proven: 

\begin{proposition}
\label{prop:9}
The differential in the Gersten-Quillen spectral sequence
$\dd:R^{p-1}_{p}(X)\longrightarrow R^{p}_{p}(X)$ is given by 
$\dd f=\dv(f)=\sum_{x\in X^{(p-1)}}\dv(f_{x})$.  
\hfill $\square$ 
\end{proposition}

\noindent
Therefore, we obtain an identification
\begin{displaymath}
\CH^{p}(X)=R^{p}_{p}(X)/\dd R^{p-1}_{p}(X)=E^{p,-p}_{2}(X).
\end{displaymath}
Hence, by proposition \ref{prop:absspec}, we obtain a natural
isomorphism
\begin{displaymath}
\CH^{p}(X)_{\QQ}\cong H^{2p}_{\mathcal{A}}(X,\mathbb{Q}(p)).
\end{displaymath}
The discrepancy in sign between proposition \ref{prop:9} and
\cite{Burgos:CDB}, p. 365, stems from the fact that we are using 
a different convention for the connecting morphism (see remark
\ref{rem:signs}).  

\noindent
Let us write
\begin{displaymath}
\CH^{p-1,p}(X)={\rm Z}R^{p-1}_{p}(X)/\dd R^{p-2}_{p}(X)=
E^{p-1,-p}_{2}(X),
\end{displaymath}
where ${\rm Z}R^{p-1}_{p}(X)=\Ker(\dd:R^{p-1}_{p}(X)\longrightarrow
R^{p}_{p}(X))$ is the subgroup of cycles. Then, again by proposition
\ref{prop:absspec}, we have an isomorphism 
\begin{displaymath}
\CH^{p-1,p}(X)_{\mathbb{Q}}\cong H^{2p-1}_{\mathcal{A}}(X,\mathbb{Q}(p)).
\end{displaymath}

\nnpar{Chow groups with support.} 
For a closed subset $Y \subset X$, we define 
\begin{displaymath}
R^{q}_{p,Y}(X)=\bigoplus_{x\in X^{(q)}\cap Y}K_{p-q}(k(x)).
\end{displaymath}
In particular, $R^{p}_{p,Y}(X)$ equals the group ${\rm Z}^{p}_{Y}(X)$ 
of $p$-codimensional cycles on $X$ with support on $Y$. We write 
\begin{align*}
\CH^{p}_{Y}(X)&=R^{p}_{p,Y}(X)/\dd R^{p-1}_{p,Y}(X), \\
\CH^{p-1,p}_{Y}(X)&={\rm Z}R^{p-1}_{p,Y}(X)/\dd R^{p-2}_{p,Y}(X).
\end{align*}
We call $\CH^{p}_{Y}(X)$ the Chow group of $p$-codimensional cycles 
on $X$ with support on $Y$. As before we have isomorphisms 
\begin{align*}
\CH^{p}_{Y}(X)_{\QQ}&\cong E^{p,-p}_{2,Y}(X)_{\QQ}\cong
H^{2p}_{\mathcal{A},Y}(X,\QQ(p)), \\
\CH^{p-1,p}_{Y}(X)_{\QQ}&\cong E^{p-1,-p}_{2,Y}(X)_{\QQ}\cong
H^{2p-1}_{\mathcal{A},Y}(X,\QQ(p)).
\end{align*}
Writing $U=X\setminus Y$, we have a long exact sequence 
\begin{displaymath}
\CH^{p-1,p}(U)\overset{\delta}{\longrightarrow}\CH^{p}_{Y}(X)
\longrightarrow\CH^{p}(X)\longrightarrow\CH^{p}(U)\longrightarrow 0.
\end{displaymath}
A family $\varphi$ of supports on $X$ is a family of closed
subsets of $X$ such that, if $Y,Z\in\varphi$, then $Y\cup Z\in
\varphi$. For any family $\varphi$ of supports on $X$, we define
\begin{align*}
R^{q}_{p,\varphi}(X)&=\lim_{\substack{\longrightarrow \\ 
Y\in\varphi}}R^{q}_{p,Y}(X), \\
\CH^{p}_{\varphi}(X)&=\lim_{\substack{\longrightarrow \\ 
Y\in\varphi}}\CH^{p}_{Y}(X), \\
\CH^{p-1,p}_{\varphi}(X)&=\lim_{\substack{\longrightarrow \\ 
Y\in\varphi}}\CH^{p-1,p}_{Y}(X), \\
H^{n}_{\mathcal{A},\varphi}(X,\QQ(p))&=\lim_{\substack{\longrightarrow \\
Y\in\varphi}}H^{n}_{\mathcal{A},Y}(X,\QQ(p)).
\end{align*}
We have the isomorphisms
\begin{align*}
\CH^{p}_{\varphi}(X)_{\QQ}&\cong 
H^{2p}_{\mathcal{A},\varphi }(X,\QQ(p)), \\
\CH^{p-1,p}_{\varphi}(X)_{\QQ}&\cong 
H^{2p-1}_{\mathcal{A},\varphi}(X,\QQ(p)).
\end{align*}

\nnpar{Intersection of cycles on regular schemes.} 
A method of defining products for the Chow groups of regular 
schemes, where no moving lemma is available and it is not possible 
to use the deformation to the normal cone technique, is given 
by means of the isomorphism between Chow groups and absolute 
cohomology groups. In this way we can transfer the multiplicative 
properties of $K$-theory to the Chow groups. Namely, we define 
a pairing 
\begin{displaymath}
\CH^{p}_{Y}(X)_{\mathbb{Q}}\otimes\CH_{Z}^{q}(X)_{\mathbb{Q}}
\overset{\cdot}{\longrightarrow}\CH^{p+q}_{Y\cap Z}(X)_{\mathbb{Q}}
\end{displaymath}
by means of the commutative diagram
\begin{displaymath}
\begin{CD}
\CH^{p}_{Y}(X)_{\mathbb{Q}}\otimes\CH^{q}_{Z}(X)_{\mathbb{Q}}
@>\cdot>>\CH^{p+q}_{Y\cap Z}(X)_{\mathbb{Q}} \\
@V\cong VV@VV\cong V \\
H^{2p}_{\mathcal{A},Y}(X,\mathbb{Q}(p))\otimes H^{2q}_
{\mathcal{A},Z}(X,\mathbb{Q}(q))@>\cup>>H^{2p+2q}_{\mathcal
{A},Y\cap Z}(X,\mathbb{Q}(p+q)).
\end{CD}
\end{displaymath}

\begin{theorem} [Gillet-Soul\'e \cite{GilletSoule:ItAo}]
\label{thm:intreg} 
If $X$ is a regular scheme of finite Krull dimension, and $Y,Z$ are
closed subsets, then the above pairing
\begin{displaymath}
\CH^{p}_{Y}(X)_{\mathbb{Q}}\otimes\CH^{q}_{Z}(X)_{\mathbb{Q}}
\overset{\cdot}{\longrightarrow}\CH^{p+q}_{Y\cap Z}(X)_{\mathbb{Q}}  
\end{displaymath}
satisfies the following properties: 
\begin{enumerate}
\item[(i)]
It turns $\bigoplus_{Y,p}\CH^{p}_{Y}(X)_{\mathbb{Q}}$ 
into a commutative ring with unit $[X]\in\CH^{0}(X)$.
\item[(ii)]
It is compatible with the change of support maps
\begin{displaymath}
\CH^{p}_{Y}(X)_{\mathbb{Q}}\longrightarrow\CH^{p}_{Y'}
(X)_{\mathbb{Q}},
\end{displaymath}
whenever $Y\subset Y'\subset X$.
\item[(iii)]
Let $y\in\CH^{p}_{Y}(X)_{\mathbb{Q}}$, $z\in\CH^{q}_{Z}(X)_
{\mathbb{Q}}$ with $Y=\supp y$, $Z=\supp Z$, and $Y,Z$ having 
proper intersection. Then, their intersection product $y\cdot 
z\in\CH^{p+q}_{Y\cap Z}(X)_{\mathbb{Q}}$ is given by Serre's
$Tor$-formula
\begin{displaymath}
y\cdot z=\left[\sum_{x\in Y\cap Z}\left(\sum_{i\ge 0}(-1)^{i}
\ell_{\mathcal{O}_{X,x}}(Tor_{i}^{\mathcal{O}_{X,x}}(\mathcal
{O}_{Y,x},\mathcal{O}_{Z,x}))\right)\overline{\{x\}}\right],
\end{displaymath}
where $\ell_{\mathcal{O}_{X,x}}$ denotes the length of a
${\mathcal{O}_{X,x}}$-module.
\end{enumerate}
\hfill $\square$ 
\end{theorem}

\noindent
It is clear from this result that, if $\varphi$ and $\psi$ are two
families of supports, then there is an induced pairing
\begin{displaymath}
\CH^{p}_{\varphi}(X)_{\mathbb{Q}}\otimes\CH^{q}_{\psi}(X)_
{\mathbb{Q}}\overset{\cdot}{\longrightarrow}\CH^{p+q}_
{\varphi\cap\psi}(X)_{\mathbb{Q}},
\end{displaymath}
where $\varphi\cap\psi$ is the family of supports given by
\begin{displaymath}
\varphi\cap\psi=\{Y\cap Z\mid Y\in \varphi,Z\in\psi\}.
\end{displaymath}
We will also need products between $\CH^{p-1,p}(X)$ and $\CH^{q}(X)$. 

\begin{theorem}
\label{thm:3}
Let $X$ be a regular scheme, let $Y$, $Z$ be closed subsets 
of $X$, and set $U=X\setminus Y$. Then, there exists a well-defined
pairing  
\begin{displaymath}
\CH^{p-1,p}(U)_{\mathbb{Q}}\otimes\CH^{q}_{Z}(X)_{\mathbb{Q}}
\overset{\cdot}{\longrightarrow}\CH^{p+q-1,p+q}_{U\cap Z}
(X\setminus(Y\cap Z))_{\mathbb{Q}}
\end{displaymath}
such that the diagram 
\begin{displaymath}
\begin{CD}
\CH^{p-1,p}(U)_{\mathbb{Q}}\otimes\CH^{q}_{Z}(X)_{\mathbb{Q}}
@>\cdot>>\CH^{p+q-1,p+q}_{U\cap Z}(X\setminus(Y\cap Z))_
{\mathbb{Q}} \\
@V\delta\otimes\Id VV@VV\delta V \\
\CH^{p}_{Y}(X)_{\mathbb{Q}}\otimes\CH^{q}_{Z}(X)_{\mathbb{Q}}
@>\cdot>>\CH^{p+q}_{Y\cap Z}(X)_{\mathbb{Q}}
\end{CD}
\end{displaymath}
commutes.
\end{theorem}
\begin{proof}
The pairing is defined by imposing the commutativity of the 
diagram
\begin{displaymath}
\begin{CD}
\CH^{p-1,p}(U)_{\mathbb{Q}}\otimes\CH_{Z}^{q}(X)_{\mathbb{Q}}
@>\cong>>H^{2p-1}_{\mathcal{A}}(U,\mathbb{Q}(p))\otimes H^{2q}_
{\mathcal{A},Z}(X,\mathbb{Q}(q)) \\
@V\cdot VV@VV\cup V \\
\CH^{p+q-1,p+q}_{U\cap Z}(X\setminus(Y\cap Z))_{\mathbb{Q}}
@>\cong>>H^{2p+2q-1}_{\mathcal{A},U\cap Z}(X\setminus(Y\cap Z),
\mathbb{Q}(p+q)).
\end{CD}
\end{displaymath}
The compatibility of the pairing with the connection morphism now
follows from theorem \ref{thm:11} (iv).
\end{proof}

For a $K_{1}$-chain $f\in R^{p-1}_{p}(X)$, put $y=\dv(f)$, $Y=
\supp y$, and $U=X\setminus Y$; observe that $f$ determines an 
element $[f]\in\CH^{p-1,p}(U)$. Let $z$ be a $q$-codimensional 
cycle with $Z=\supp z$. By the above theorem, there is a well-defined
element  
\begin{displaymath}
[f]\cdot[z]\in\CH^{p+q-1,p+q}_{U\cap Z}(X\setminus(Y\cap Z))_
{\mathbb{Q}}.
\end{displaymath}
Using that $R^{p+q-1}_{p+q}(X\setminus(Y\cap Z))$ is a direct 
summand of $R^{p+q-1}_{p+q}(X)$, the above theorem implies  

\begin{corollary}
\label{cor:2}
Assume that $U\cap Z\neq\emptyset$. Then, there exists a $K_{1}
$-chain $g\in R^{p+q-1}_{p+q}(X)_{\mathbb{Q}}$, determined up 
to $\dd R^{p+q-2}_{p+q}(X)_{\mathbb{Q}}+R^{p+q-1}_{p+q,Y\cap Z}
(X)_{\mathbb{Q}}$, such that its class in $\CH^{p+q-1,p+q}_{U
\cap Z}(X\setminus(Y\cap Z))_{\mathbb{Q}}$ is equal to $[f]
\cdot[z]$. Moreover, we have 
\begin{equation}    
\label{remarkcycle} 
[\dv(g)]=[\dv(f)]\cdot[z]
\end{equation}
in the group $\CH^{p+q}_{Y\cap Z}(X)_{\mathbb{Q}}$.
\hfill $\square$ 
\end{corollary}

\begin{remark}
If we have $\codim(Y\cap Z)=p+q$ in corollary \ref{cor:2}, the
equality $\CH^{p+q}_{Y\cap Z}(X)={\rm Z}^{p+q}_{Y\cap Z}(X)$
shows that identity (\ref{remarkcycle}) holds true on the level 
of cycles, i.e., we have 
\begin{displaymath}
\dv(g)=\dv(f)\cdot z.
\end{displaymath}
\end{remark}

\nnpar{Pull-back.} 
Here we recall the properties of the pull-back of algebraic cycles and
$K_{1}$-chains by a quasi-projective morphism of regular schemes. Let 
us start with the case of cycles. Since any such morphism can be factored 
as the composition of a closed immersion and a flat morphism, we will 
discuss both cases separately. For flat morphisms the pull-back is 
defined as in \cite{Fulton:IT}, 1.7, and for closed immersions one uses 
the deformation to the normal cone technique (see \cite{GilletSoule:ait},
4.4.1). Hence one can define the inverse image for a quasi-projective 
morphism using any factorization into a closed immersion followed by a 
flat morphism, and then proves that the result is independent of the 
chosen factorization. The main properties we need are:

\begin{theorem}
\label{thm:5} 
Let $f: X\longrightarrow Y$ be a quasi-projective morphism of regular 
schemes. Let $\varphi$ be a family of supports on $Y$, and let $\psi=
f^{-1}(\varphi)$ be the corresponding family of supports on $X$.
\begin{enumerate} 
\item[(i)]
Assuming that $f=g\circ i=h\circ j$, where $i$, $j$ are closed
immersions and $g$, $h$ are flat morphisms, we have
\begin{displaymath}
i^{\ast}\circ g^{\ast}=j^{\ast}\circ h^{\ast}.
\end{displaymath}
Therefore, for all $p\ge 0$, there is a well-defined homomorphism of 
Chow groups   
$$
f^{\ast}:\CH^{p}_{\varphi}(Y)\longrightarrow\CH^{p}_{\psi}(X).
$$
\item[(ii)]  
The map  $f^{\ast}$ induces a ring homomorphism
$$
f^{\ast}:\bigoplus_{p\ge 0}\CH^{p}(Y)_{\QQ}\longrightarrow
\bigoplus_{p\ge 0}\CH^{p}(X)_{\QQ}.
$$
\item[(iii)]
The pull-back map $f^{\ast}:\CH^{p}(Y)_{\QQ}\longrightarrow\CH^{p}
(X)_{\QQ}$ corresponds via the isomorphisms $\CH^{p}(Y)_{\QQ}\cong
H^{2p}_{\mathcal{A}}(Y,\mathbb{Q}(p))$ and $\CH^{p}(X)_{\QQ}\cong 
H^{2p}_{\mathcal{A}}(X,\mathbb{Q}(p))$ to the pull-back defined in 
absolute cohomology.
\item[(iv)]
If $g:Y\longrightarrow Z$ is another morphism, then we have
$(g\circ f)^{\ast}=f^{\ast}\circ g^{\ast}$.
\end{enumerate}
\hfill $\square$
\end{theorem}

To define the pull-back for $K_{1}$-chains, one observes first that
the Gersten-Quillen spectral sequence is contravariant for flat
morphisms. Moreover, in \cite{Gillet:RRhK}, there is a definition
for a pull-back of $K_{1}$-chains for closed immersions using the 
deformation to the normal cone technique. The next result follows 
from \cite{Gillet:RRhK} and \cite{GilletSoule:ait}, 4.4.2.

\begin{theorem}
\label{thm:12}
Let $f:X\longrightarrow Y$ be a quasi-projective morphism of regular
schemes. Let $f=g\circ i$ be a factorization of $f$ into a closed
immersion $i$ followed by a flat morphism $g$. Let $U\subset Y$ be 
an open subset and $Z\subset U$ a closed subset. Let us write $U'=
f^{-1}(U)$ and $Z'=f^{-1}(Z)$.   
\begin{enumerate}
\item[(i)] 
For any $p\ge 0$, there is a homomorphism
\begin{displaymath}
f^{\ast}=i^{\ast}\circ g^{\ast}:\CH^{p-1,p}_{Z}(U)\longrightarrow 
\CH^{p-1,p}_{Z'}(U')
\end{displaymath}
which does not depend on the factorization chosen.
\item[(ii)] 
Let us denote by $\delta$ the connection morphisms
\begin{align*}
\CH^{p-1,p}_{Z}(U)&\overset{\delta}{\longrightarrow}
\CH^{p}_{\overline{Z}\setminus U}(Y), \\
\CH^{p-1,p}_{Z'}(U')&\overset{\delta}{\longrightarrow}
\CH^{p}_{\overline{Z'}\setminus U'}(X).
\end{align*}
Then, we have $f^{\ast}\circ\delta=\delta\circ f^{\ast}$.
\item[(iii)] 
The induced map $f^{\ast}:\CH^{p-1,p}_{Z}(U)_{\QQ}\longrightarrow
\CH^{p-1,p}_{Z'}(U')_{\QQ}$ corresponds via the isomorphisms
$\CH^{p-1,p}_{Z}(U)_{\QQ}\cong H^{2p-1}_{\mathcal{A},Z}(U,\mathbb{Q}(p))$
and $\CH^{p-1,p}_{Z'}(U')_{\QQ}\cong H_{\mathcal{A},Z'}^{2p-1}(U',
\mathbb{Q}(p))$ to the pull-back defined in absolute cohomology.
\end{enumerate}
\end{theorem}
\begin{proof}
The definition of $i^{\ast}$ is given in \cite{GilletSoule:ait},
4.4.2. By definition, to show that $f^{\ast}$ is compatible with the
connection morphisms and with the pull-back in absolute cohomology,
it is enough to treat the flat morphisms and the closed immersions
separately. For flat morphisms the result is clear. The compatibility 
with the connection morphisms for closed immersions is proven in
\cite{GilletSoule:ait}, 4.4.2. To prove the compatibility with the 
pull-back in absolute cohomology, one observes that all the steps in 
the definition of $i^{\ast}$ are compatible with the isomorphism with
absolute cohomology.  
\end{proof}

\nnpar{Push-forward and the projection formula.} 
Let $f:X\longrightarrow X'$ be a proper morphism of regular,
equidimensional, quasi-projective schemes in $\ZAR(S)$ of relative
dimension $d$. Then, the Gersten-Quillen spectral sequence is 
covariant for proper morphisms (see \cite{Gillet:RRhK}). In particular, 
for any closed subset $Y\subset X$, we have morphisms
\begin{align*}
f_{\ast}:\CH^{p}_{Y}(X)&\longrightarrow\CH^{p-d}_{f(Y)}(X'), \\
f_{\ast}:\CH^{p-1,p}_{Y}(X)&\longrightarrow\CH^{p-d-1,p-d}_{f(Y)}(X').
\end{align*}
These morphisms satisfy the following properties.

\begin{theorem}
\label{thm:13}
Let $f:X\longrightarrow X'$ be a proper morphism of regular,
equidimensional, quasi-projective schemes in $\ZAR(S)$, of 
relative dimension $d$. 
\begin{enumerate}
\item[(i)]
Let $Y$ be a closed subset of $X'$. Let us write $U'=X'\setminus Y$ 
and $U=X\setminus f^{-1}(Y)$. Then, the following diagram commutes
\begin{displaymath}
\xymatrix{\CH^{p-1,p}(U)\ar[r]^{\delta}\ar[d]_{f_{\ast}}&
\CH^{p}_{f^{-1}(Y)}(X)\ar[d]^{f_{\ast}} \\
\CH^{p-1,p}(U')\ar[r]^{\delta}&\CH^{p}_{Y}(X').}
\end{displaymath}
\item[(ii)] 
If $g:X'\longrightarrow X''$ is another proper morphism of regular,
equidimensional, quasi-projective schemes in $\ZAR(S)$, we have
\begin{displaymath}
(g\circ f)_{\ast}=g_{\ast}\circ f_{\ast}.
\end{displaymath}
\item[(iii)]
The direct image of cycles and $K_{1}$-chains is compatible with
the direct image in absolute cohomology.
\item[(iv)] 
For $\alpha\in\CH^{p}(X')$ and $\beta\in\CH^{q}(X)$ we have
\begin{displaymath}
f_{\ast}(f^{\ast}(\alpha)\cdot\beta)=\alpha\cdot f_{\ast}(\beta)
\in\CH^{p+q-d}(X')_{\mathbb{Q}}.   
\end{displaymath}
\end{enumerate}
\end{theorem}
\begin{proof}
The first and second statement follow from the covariance of the
Gersten-Quillen spectral sequence. The third statement follows from
the Riemann-Roch theorem, and the fourth statement follows from the
third and the projection formula for absolute cohomology.
\end{proof}

\subsection{Gillet cohomology and characteristic classes}
\label{sec:ctcc} 

In this section we will recall some facts about characteristic
classes from higher $K$-theory to a suitable cohomology theory. 
We will follow the axiomatic approach of \cite{Gillet:RRhK}; hence,
we will ask the cohomology theory to satisfy all the properties 
stated in \cite{Gillet:RRhK}. Any such cohomology theory will be 
called a Gillet cohomology. Note however that the definition of
Gillet cohomology given here is slightly different from the definition 
in \cite{Gillet:RRhK}.

\nnpar{Gillet cohomologies.}
We fix a field $k$. Let $\mathscr{V}$ be the category of smooth 
schemes over $k$ (recall that by scheme we mean a noetherian, 
separated scheme of finite Krull dimension). Let $\mathscr{V}_
{\ast}$ be a category of schemes and proper morphisms that contains 
all closed immersions $Y\longrightarrow X$ with $X$ in $\mathscr{V}$. 
Let $C$ be the big Zariski site of the category $\mathscr{V}$. 

\begin{definition}
A \emph{Gillet cohomology} is given by the following data: 
\begin{enumerate}
\item[(1)] 
A graded complex of sheaves of abelian groups $\Gi^{\ast}(\ast)$ 
in the site $C$, together with a pairing in the derived category 
of graded complexes of abelian sheaves on $C$
\begin{displaymath}
\Gi^{\ast}(\ast)\underset{\mathbb{Z}}{\overset{L}{\otimes}}
\Gi^{\ast}(\ast)\longrightarrow\Gi^{\ast}(\ast),
\end{displaymath}
which is associative, (graded-)commutative, and has a unit. Given
such a complex $\Gi=\Gi^{\ast}(\ast)$, for each pair $(Y,X)$, with
$X$ in $\mathscr{V}$ and $Y$ a closed subscheme of $X$, we define 
the cohomology groups of $X$ with coefficients in $\Gi$ and support 
in $Y$ by     
\begin{displaymath}
H^{i}_{Y}(X,\Gi(j))=\mathbb{H}^{i}_{Y}(X,\Gi^{\ast}(j)).
\end{displaymath}
Since these cohomology groups are defined as the hypercohomology 
groups of a multiplicative Zariski sheaf, they satisfy the usual
multiplicative and functorial  properties (see \cite{Gillet:RRhK}, 
definition 1.1). 
\item[(2)]
A homology theory, that is, a covariant functor from the category 
$\mathscr{V}_{\ast}$ to bigraded abelian groups 
\begin{displaymath}
X\longmapsto\bigoplus_{\substack{i\ge 0\\[1mm]j\in\mathbb{Z}}}
H_{i}(X,\Gi(j)).
\end{displaymath}
\end{enumerate}
This homology-cohomology theory has to satisfy the properties 
(i)-(xi) of \cite{Gillet:RRhK}, definition 1.2 with $d=2$. Since 
we are distinguishing between objects of $\mathscr{V}$ and of
$\mathscr{V}_{\ast}$, we have to specify to which category the
objects in question in properties (i)-(xi) belong: In (i), 
(ii), (viii) $X,Y\in\mathscr{V}_{\ast}$; in (iii), (v) $X\in
\mathscr{V}$, $Y\in\mathscr{V}_{\ast}$; in (iv) $X\in\mathscr{V}_
{\ast}$; in (vi), (vii) $X,Y\in\mathscr{V}$; and in (ix), (x) 
$X\in\mathscr{V}$. In addition to Gillet's axioms, following 
Beilinson \cite{Beilinson:hr}, we will also assume    
\begin{enumerate}
\item[(xii)] 
$H^{i}(\cdot,\Gi(j))=0$ for all $i<0$ and all $j\in\mathbb{Z}$.  
\end{enumerate}
\end{definition}

\begin{definition}
A \emph{Gillet complex over $k$} is the graded complex of sheaves 
of abelian groups $\Gi=\Gi^{\ast}(\ast)$ of a Gillet cohomology.
\end{definition}

\begin{remark}
The main discrepancy between the definition here and that of
\cite{Gillet:RRhK} is that the objects of $\mathscr{V}$ and of
$\mathscr{V}_{\ast}$ are not the same. The reason for this is 
that we want a Gillet complex as a sheaf only over smooth schemes. 
On the other hand, we want to have cycle classes for arbitrary 
subschemes. Therefore, we need the homology also to be defined 
for singular schemes.    
\end{remark}

\nnpar{Purity.}
Property (vi) of the Gillet axioms includes the following purity 
condition: If $(X,Y)$ is a pair of schemes in $C$ with $Y$ a closed 
subscheme of $X$ of pure codimension $p$, then the natural map
\begin{displaymath}
H^{i}(Y,\Gi(j))\longrightarrow H_{Y}^{i+2p}(X,\Gi(j+p))
\end{displaymath}
is an isomorphism.

This purity condition, together with the vanishing assumption (xii),
implies:

\begin{proposition}
Let $\Gi$ be a Gillet complex over $k$, and $X$ a regular scheme 
in $\ZAR(\Spec(k))$. For any closed subset $Y$ of $X$ of codimension 
greater or equal to $p$, we then have
\begin{displaymath}
H^{i}_{Y}(X,\Gi(p))=0 
\end{displaymath}
for all $i<2p$.
\hfill $\square$
\end{proposition}

\nnpar{Sheaf cohomology as generalized cohomology.} 
Usual sheaf cohomology of abelian sheaves can also be realized as
generalized cohomology. The reader is referred to \cite{BrownGersten:haKt},
\cite{Gillet:RRhK}, \cite{Jeu:Zc}, and \cite{HuberWildeshaus:cmp},
appendix B, for details. 

Given a cohomological complex of abelian groups $G^{\ast}$ indexed by
non-positive integers, the Dold-Puppe functor associates to it a
simplicial abelian group $K(G)$ pointed by $0$ such that 
\begin{displaymath}
h^{-i}(G^{\ast})=\pi_{i}(K(G),0).
\end{displaymath}
Given a general cohomological complex $G$, we will denote by $\sigma_
{p}$ the b\^ete filtration (see section \ref{sec:tcdb}). Then, we can
define an infinite loop spectrum
\begin{displaymath}
K(G)_{N}=K(\sigma_{0}(G[N])).
\end{displaymath}
We can recover all the cohomology groups of $G$ from the stable
homotopy groups of this spectrum.

We can sheafify the above construction and then associate to each complex 
of abelian sheaves $\mathcal{G}$ an infinite loop spectrum of spaces 
$K(\mathcal{G})$. Moreover, we have (see \cite{HuberWildeshaus:cmp}): 

\begin{proposition}
Let $\mathcal{G}$ be a bounded below abelian sheaf over the site $C$. 
Then, we have for every $X$ in $C$,
\begin{displaymath}
H^{i}_{sT}(X,K(\mathcal{G}))=H^{i}_{\ZAR}(X,\mathcal{G}).
\end{displaymath}
\hfill $\square$
\end{proposition}

\nnpar{Characteristic classes.} 
To a Gillet cohomology theory, one can associate a theory of 
characteristic classes from higher $K$-theory (see \cite{Gillet:RRhK},
\cite{Beilinson:hr}). These characteristic classes can be constructed 
on the level of spaces (see \cite{Gillet:RRhK}, \cite{Jeu:Zc}). More
specifically, we have the following result:  

\begin{theorem}
Let $\bigoplus_{j\in\mathbb{Z}}\Gi^{\ast}(j)$ be a Gillet complex
over $k$. Moreover, assume that the sheaves $\Gi^{i}(j)$ are injective
for all $i,j$. Then, for every $j\ge 0$, there is a Chern class map 
of spaces
\begin{displaymath}    
c_{j}:\KK\longrightarrow K(\Gi(j)[2j]).
\end{displaymath}
These maps induce morphisms
\begin{displaymath}
c_{j}:H^{i}_{sT}(X,\KK)\longrightarrow H_{sT}^{2j+i}(X,K(\Gi(j)))=
H_{\ZAR}^{2j+i}(X,\Gi(j))
\end{displaymath}
for all spaces $X$.
\end{theorem}
\begin{proof}
The injectivity of $\Gi^{i}(j)$ implies that $\Gi^{\ast}(j)$ are 
pseudo-flasque complexes (see section \ref{sec:ktgc}). This result
is a particular case of \cite{HuberWildeshaus:cmp}, theorem B.3.7.
\end{proof}

This theorem implies in particular that there is a definition of
Chern classes with support. That is, if $Y$ is a closed subscheme
of $X$, then there are classes
\begin{displaymath}
c^{Y}_{j}:K_{m}^{Y}(X)\longrightarrow H^{2j-m}_{Y}(X,\Gi(j)).
\end{displaymath}
Moreover, these Chern classes are natural and satisfy

\begin{proposition}
\label{prop:8}
Let $Y$ be a closed subset of $X$ and $U=X\setminus Y$. Then, 
the Chern classes form a morphism of exact sequences  
\begin{displaymath}
\begin{CD}
K_{m}(X)@>>>K_{m}(U)@>\delta>>K_{m-1}^{Y}(X)@>>>\dots \\
@V c_{j} VV@V c_{j} VV@V c^{Y}_{j} VV \\
H^{2j-m}(X,\Gi(j))@>>>H^{2j-m}(U,\Gi(j))@>\delta>> 
H^{2j-m+1}_{Y}(X,\Gi(j))@>>>\dots 
\end{CD}
\end{displaymath}
\hfill $\square$
\end{proposition}

\nnpar{Bloch-Ogus spectral sequence.} 
Since the set of axioms in \cite{Gillet:RRhK} include in particular 
all the axioms of \cite{BlochOgus:Gchs}, we can form the Bloch-Ogus 
spectral sequence associated to the cohomology $\Gi$. Moreover, by 
the fact that characteristic classes are defined on the level of 
sheaves, characteristic classes induce morphisms from the Gersten-Quillen
spectral sequence of $K$-theory to the Bloch-Ogus spectral sequence 
of $\Gi$-cohomology (see \cite{Gillet:RRhK}, \S 3). 

Let $x\in X^{(p)}$, and let $\iota:\overline{\{x\}}\longrightarrow X$ 
be the inclusion. We will write
\begin{align*}
H^{i}(x,\Gi(j))=\lim_{\substack{\longrightarrow\\U}}
H^{i}(\overline{\{x\}}\cap U,\Gi(j)), \\
H^{i}_{x}(X,\Gi(j))=\lim_{\substack{\longrightarrow\\U}}
H^{i}_{\overline{\{x\}}\cap U}(X,\Gi(j)),  
\end{align*}
where the limit is taken over all open sets $U$ of $X$ containing
$x$. By the purity property of a Gillet cohomology, $\iota$ induces
isomorphisms 
\begin{displaymath}
\iota_{!}:H^{i}(x,\Gi(j))\longrightarrow H^{i+2p}_{x}(X,\Gi(j+p)).
\end{displaymath}
The Bloch-Ogus spectral sequence is a spectral sequence
\begin{displaymath}
E^{p,q}_{r}(j)\Longrightarrow H^{p+q}(X,\Gi(j))
\end{displaymath}
with 
\begin{align*}
E^{p,q}_{1}(j)&=\bigoplus_{x\in X^{(p)}}H^{p+q}_{x}(X,\Gi(j)) \\
&\cong\bigoplus_{x\in X^{(p)}}H^{q-p}(x,\Gi(j-p)). 
\end{align*}
Characteristic classes induce morphisms between the 
Gersten-Quillen spectral sequence and the Bloch-Ogus 
spectral sequence (see \cite{Gillet:RRhK}, \S3). 

\begin{proposition}
For each Chern class $c_{j}$ there is a morphism of spectral sequences
\begin{displaymath}
\bigoplus_{x\in X^{(p)}}K^{x}_{-p-q}(X)\overset{\bigoplus c_{j}^{x}}
{\longrightarrow}\bigoplus_{x\in X^{(p)}}H^{2j+p+q}_{x}(X,\Gi(j)).
\end{displaymath}
Moreover, this morphism is contravariant for flat morphisms.  
\hfill $\square$
\end{proposition}

We can form a commutative diagram
\begin{displaymath}
\begin{CD}
\bigoplus_{x\in X^{(p)}}K^{x}_{-p-q}(X)@>\bigoplus c_{j}^{x}>>
\bigoplus_{x\in X^{(p)}}H^{2j+p+q}_{x}(X,\Gi(j)) \\
@V\cong VV @VV\cong V \\
\bigoplus_{x\in X^{(p)}}K_{-p-q}(k(x))@>\bigoplus\gamma_{j}^{x}>>
\bigoplus_{x\in X^{(p)}}H^{2j-p+q}(x,\Gi(j-p)),
\end{CD}
\end{displaymath}
where the bottom arrow is determined as a consequence of the Riemann-Roch
theorem without denominators (see \cite{Gillet:RRhK}, theorem 3.9).

\begin{lemma}
\label{lemm:1}
Let $X$ be a smooth, noetherian $k$-scheme of finite Krull dimension. 
Let $x\in X^{(p)}$, and let $\eta_{x}\in H^{0}(x,\Gi(0))$ be the 
fundamental class. Then, there is a commutative diagram 
\begin{displaymath}
\begin{CD}
K^{x}_{-p-q}(X)@>c_{j}^{x}>>H^{2j+p+q}_{x}(X,\Gi (j)) \\
@A\iota_{*}A\cong A@A\cong A\iota_{!}A \\
K_{-p-q}(k(x))@>\gamma_{j}^{x}>>H^{2j-p+q}(x,\Gi(j-p))
\end{CD}
\end{displaymath}
such that $\alpha\in K_{-p-q}(k(x))$ satisfies
\begin{displaymath}
\gamma_{j}^{x}(\alpha )=
\begin{cases}
\frac{(-1)^{p}(j-1)!}{(j-p-1)!}c_{j-p}^{x}(\alpha),&\ \text{if }j>p, \\
\rk(\alpha)\eta_{x},&\ \text{if }j=p. 
\end{cases}
\end{displaymath}
\hfill $\square$
\end{lemma}

\nnpar{Classes for cycles and $K_{1}$-chains.} 
For the construction of arithmetic Chow rings we are mainly interested 
in the induced classes of algebraic cycles and $K_{1}$-chains.

\begin{definition}
\label{def:1}
Let $x\in X^{(p)}$, and let $\iota:\overline{\{x\}}\longrightarrow X$ 
be the inclusion. For $\alpha\in K_{0}(k(x))$, we will denote by 
$\cl_{\Gi}(\alpha)$ the class   
\begin{displaymath}
\cl_{\Gi}(\alpha)=\iota_{!}(\gamma_{p}^{x}(\alpha))=\iota_{!}
(\rk(\alpha)\eta_{x})\in H^{2p}_{x}(X,\Gi(p)).
\end{displaymath}
Furthermore, for $f\in K_{1}(k(x))$, we will denote by $\cl_{\Gi}(f)$ 
the class 
\begin{displaymath}
\cl_{\Gi}(f)=\iota_{!}((-1)^{p}\gamma_{p+1}^{x}(f)/p!)=\iota_{!}
(c_{1}^{x}(f))\in H^{2p+1}_{x}(X,\Gi(p+1)).
\end{displaymath}
By linearity, we obtain well-defined morphisms
\begin{align*}
\cl_{\Gi}:R^{p}_{p}(X)&\longrightarrow\bigoplus_{x\in X^{(p)}}
H^{2p}_{x}(X,\Gi(p))\,, \\ 
\cl_{\Gi}:R^{p-1}_{p}(X)&\longrightarrow\bigoplus_{x\in X^{(p-1)}}
H^{2p-1}_{x}(X,\Gi(p))\,.  
\end{align*}
\end{definition}

These classes should be interpreted as the Chern character. Namely,
using the isomorphism between Chow groups and absolute cohomology, 
the classes of definition \ref{def:1} induce classes between absolute
cohomology groups and the Gillet cohomology groups which agree with
the Chern character.
 
\begin{proposition}
There are commutative diagrams
\begin{displaymath}
\begin{CD}
\CH^{p}(X)@>>>H^{2p}_{\mathcal{A}}(X,\mathbb{Q}(p)) \\
@V\cl_{\Gi} VV @VV\ch V \\
H^{2p}(X,\Gi(p))@>>>H^{2p}(X,\Gi(p))_{\mathbb{Q}},
\end{CD}
\end{displaymath}
and
\begin{displaymath}
\begin{CD}
\CH^{p-1,p}(X)@>>>H^{2p-1}_{\mathcal{A}}(X,\mathbb{Q}(p)) \\
@V\cl_{\Gi} VV @VV\ch V \\
H^{2p-1}(X,\Gi(p))@>>>H^{2p-1}(X,\Gi(p))_{\mathbb{Q}}.
\end{CD}
\end{displaymath}
\end{proposition}
\begin{proof}
This follows from lemma \ref{lemm:1} and the power series expansion
of the Chern character.
\end{proof}

\nnpar{Basic properties of the morphism $\cl_{\Gi}$.} 
The next result follows from the purity of a Gillet cohomology,
proposition \ref{prop:8} and proposition \ref{prop:9}.
\begin{proposition}
\begin{enumerate}
\item[(i)] 
Let $y\in\bigoplus_{x\in X^{(p)}}K_{0}(k(x)) $ be a 
$p$-codimensional cycle with $Y=\supp y$. Then, $\cl_{\Gi}(y)$
determines a well-defined class $\cl_{\Gi}(y)\in H^{2p}_{Y}
(X,\Gi(p))$.
\item[(ii)] 
For a $K_{1}$-chain $f\in R^{p-1}_{p}(X)$, put $y=\dv(f)$, $Y=\supp y$, 
and $U=X\setminus Y$. Then, $\cl_{\Gi}(f)$ determines a well-defined 
class $\cl_{\Gi}(f)\in H^{2p-1}(U,\Gi(p)).$  
\item[(iii)] 
With $f$ as in (ii), the equality 
\begin{displaymath}
\cl_{\Gi}(\dv(f))=\delta\cl_{\Gi}(f)\in H^{2p}_{Y}(X,\Gi(p))
\end{displaymath}
holds, where $\delta$ is the connection morphism
\begin{displaymath}
\delta:H^{2p-1}(U,\Gi(p))\longrightarrow H^{2p}_{Y}(X,\Gi(p)).
\end{displaymath}
\item[(iv)] 
If $h\in R^{p-2}_{p}(X)$ is a $K_{2}$-chain, then the class
$\cl_{\Gi}(\dd h)$ vanishes in the group $H^{2p-1}(X,\Gi(p))$.
\end{enumerate}
\hfill $\square$
\end{proposition}

\begin{corollary}
Let $\varphi$ be any family of closed supports. Then, the morphisms
$\cl_{\Gi}$ induce morphisms
\begin{align*}
\cl_{\Gi}:\CH^{p}_{\varphi}(X)&\longrightarrow H^{2p}_{\varphi}
(X,\Gi(p)), \\
\cl_{\Gi}:\CH^{p-1,p}_{\varphi}(X)&\longrightarrow H^{2p-1}_
{\varphi}(X,\Gi(p)).  
\end{align*}
Moreover, if $Y\subset X$ is a closed subset with $U=X\setminus Y$, 
then there is a morphism of exact sequences
\begin{displaymath}
\begin{matrix}
\CH^{p-1,p}(U)&\rightarrow&\CH^{p}_{Y}(X)&\rightarrow &
\CH^{p}(X)&\rightarrow&\CH^{p}(U)&\!\rightarrow 0 \\
\downarrow&&\downarrow&&\downarrow&&\downarrow& \\
H^{2p-1}(U,\Gi(p))\!&\rightarrow&\!H^{2p}_{Y}(X,\Gi(p))
\!&\rightarrow&\!H^{2p}(X,\Gi(p))\!&\rightarrow&\!H^{2p}
(U,\Gi (p)).\!&
\end{matrix}
\end{displaymath}
\hfill $\square$
\end{corollary}

\nnpar{Product, pull-back and  push-forward.} 
Let us summarize the compatibility properties of characteristic 
classes from Chow groups with products, inverse images and direct 
images. Recall that in this section all schemes are defined over a 
field. 
 
\begin{proposition}
\label{prop:16}
The morphism
\begin{displaymath}
\cl_{\Gi}:\bigoplus_{p\ge 0}\CH^{p}(X)\longrightarrow 
\bigoplus_{p\ge 0}H^{2p}(X,\Gi(p))
\end{displaymath}
is multiplicative.
\hfill $\square$
\end{proposition}

\begin{proposition}
\label{prop:10}
For a $K_{1}$-chain $f\in R^{p-1}_{p}(X)$, put $y=\dv(f)$, 
$Y=\supp y$, and $U=X\setminus Y$; for $z\in R^{q}_{q}(X)$, 
put $Z=\supp z$. Then, for any $K_{1}$-chain $g\in R^{p+q-1}_
{p+q}(X)_{\mathbb{Q}}$ representing  $[f\cdot z]$, the equality 
\begin{displaymath}
\cl_{\Gi}(g)=\cl_{\Gi}(f)\cdot\cl_{\Gi}(z)\in H^{2p+2q-1}_
{U\cap Z}(X\setminus(Y\cap Z),\Gi(p+q))_{\mathbb{Q}}  
\end{displaymath}
holds. Moreover, we have $\cl_{\Gi}(\dv(g))=\cl_{\Gi}(\dv(f))
\cdot\cl_{\Gi}(z)$.
\hfill $\square$
\end{proposition}

\begin{proposition}
\label{prop:iicccg}
Let $f: X\longrightarrow Y$ be a quasi-projective morphism of regular 
schemes. Let $\varphi$ be a family of supports on $Y$, and let $\psi=
f^{-1}(\varphi)$ be the corresponding family of supports on $X$. Then,
we have a commutative diagram
\begin{displaymath}
\begin{CD}
\CH^{p}_{\varphi}(Y)@>\cl_{\Gi}>>H^{2p}_{\varphi}(Y,\Gi(p)) \\
@V f^{\ast}VV @VV f^{\ast}V \\
\CH^{p}_{\psi}(X)@>\cl_{\Gi}>>H^{2p}_{\psi}(X,\Gi(p))\,.
\end{CD}
\end{displaymath}
Furthermore, let $U\subset Y$ be an open subset and $Z\subset U$ a 
closed subset. Let us write $U'=f^{-1}(U)$ and $Z'=f^{-1}(Z)$. Then,
we have a commutative diagram   
\begin{displaymath}
\begin{CD}
\CH^{p-1,p}_{Z}(U)@>\cl_{\Gi}>>H^{2p-1}_{Z}(U,\Gi(p)) \\
@V f^{\ast}VV @VV f^{\ast}V \\
\CH^{p-1,p}_{Z'}(U')@>\cl_{\Gi}>>H^{2p-1}_{Z'}(U',\Gi(p))\,.
\end{CD}    
\end{displaymath}
\hfill $\square$
\end{proposition}

\begin{proposition}
\label{prop:dicccg}
Let $f:X\longrightarrow X'$ be a projective morphism of regular,
equidimensional schemes of relative dimension $d$. Let $Y$ be a
closed subset of $X$. Then, the diagrams   
\begin{displaymath}
\begin{CD}
\CH^{p}_{Y}(X)@>\cl_{\Gi}>>H^{2p}_{Y}(X,\Gi(p)) \\
@V f_{\ast}VV @VV f_{\ast}V \\
\CH^{p-d}_{f(Y)}(X')@>\cl_{\Gi}>>H^{2p-2d}_{f(Y)}(X',\Gi(p-d))
\end{CD}
\end{displaymath}
and
\begin{displaymath}
\begin{CD}
\CH^{p-1,p}_{Y}(X)@>\cl_{\Gi}>>H^{2p-1}_{Y}(X,\Gi(p)) \\
@V f_{\ast}VV @VV f_{\ast}V \\
\CH^{p-d-1,p-d}_{f(Y)}(X')@>\cl_{\Gi}>>H^{2p-2d-1}_{f(Y)}
(X',\Gi(p-d))
\end{CD}    
\end{displaymath}
are commutative.
\hfill $\square$
\end{proposition}

\newpage
\section{Some topics from homological algebra}
\label{sec:trcg}

\subsection{$k$-iterated complexes}
\label{sec:kic} 

In this section we will introduce the notion of $k$-iterated 
complexes and discuss its relationship with $k$-complexes. 
Let $R$ be a commutative ring and $\mathcal{A}$ the category 
of $R$-modules.

\nnpar{$k$-complexes and $k$-iterated complexes.} 
\begin{definition}
A \emph{$k$-complex} $A=(A^{\ast,\dots,\ast},\dd_{1},\dots,\dd_{k})$ 
is a $k$-graded $R$-module $A$ together with $k$ endomorphisms
$\dd_{1},\dots,\dd_{k}$ of multi-degrees 
$$
(1,0,\dots,0),\dots,(0,\dots,0,1),
$$ 
respectively, such that for all $i,j$  
\begin{displaymath}
\dd_{i}^{2}=0,\qquad\dd_{i}\dd_{j}=-\dd_{j}\dd_{i}.
\end{displaymath}
We will denote by $k\text{\,-\,}\mathcal{C}^{+}(\mathcal{A})$ the 
category of $k$-complexes which are bounded from below, i.e., there 
is an integer $m$ such that $A^{n_{1},\dots,n_{k}}=0$, whenever 
$n_{i}\le m$ ($i=1,\dots,k$).
\end{definition}

\begin{definition}
A \emph{$k$-iterated complex} $A=(A^{\ast,\dots,\ast},\dd_{1},\dots,
\dd_{k})$ is a $k$-graded $R$-module $A$ together with $k$ endomorphisms
$\dd_{1},\dots,\dd_{k}$ of multi-degrees 
$$
(1,0,\dots,0),\dots,(0,\dots,0,1),
$$ 
respectively, such that for all $i,j$  
\begin{displaymath}
\dd_{i}^{2}=0,\qquad\dd_{i}\dd_{j}=\dd_{j}\dd_{i}.
\end{displaymath}
We will denote by $(\mathcal{C}^{+})^{k}(\mathcal{A})$ the category
of $k$-iterated complexes which are bounded from below.
\end{definition}

By convention, in the sequel all $k$-complexes and $k$-iterated 
complexes will be bounded from below, even if it is not stated
explicitly.

\begin{definition} 
Let $\mathcal{C}_{k}:(\mathcal{C}^{+})^{k}(\mathcal{A})\longrightarrow
k\text{\,-\,}\mathcal{C}^{+}(\mathcal{A})$ be the functor given by
associating to a $k$-iterated complex $(A^{n_{1},\dots,n_{k}},\dd_{1},
\dd_{2},\dots,\dd_{k})$ the $k$-complex $(A^{n_{1},\dots,n_{k}},
\dd_{1},(-1)^{n_{1}}\dd_{2},\dots,(-1)^{n_{1}+\dots+n_{k-1}}\dd_{k})$.
\end{definition}

\noindent
This functor is an equivalence of categories.

\nnpar{The simple complex.} 
\begin{definition}
(i) Let $A=(A^{n_{1},\dots,n_{k}},\dd_{1},\dd_{2},\dots,\dd_{k})$ 
be a $k$-complex. Then, \emph{the associated simple complex $s(A)$} 
is defined by   
\begin{displaymath}
s(A)^{n}=\bigoplus_{n_{1}+\dots+n_{k}=n}A^{n_{1},\dots,n_{k}},
\qquad\dd=\sum_{i=1}^{k}\dd_{i}. 
\end{displaymath}
(ii) Let $A=(A^{n_{1},\dots,n_{k}},\dd_{1},\dd_{2},\dots,\dd_{k})$ 
be a $k$-iterated complex. Then, \emph{the associated simple complex} 
is defined by 
\begin{displaymath}
s(A)=s(\mathcal{C}_{k}(A)).
\end{displaymath}
\end{definition}

We have defined the simple complex of a $k$-iterated complex by 
contracting all the degrees into one degree. Sometimes it will be 
also useful to make only a partial simple. 

\begin{definition}
Let $A=(A^{\ast,\dots,\ast},\dd_{A,1},\dots,d_{A,k})$ be a $k$-iterated 
complex with $k\ge 2$. We denote by $s_{j,j+1}(A)$ the $(k-1)$-iterated
complex given by   
\begin{displaymath}
s_{j,j+1}(A)^{n_{1},\dots,n_{k-1}}=\bigoplus_{p+q=n_{j}}A^{n_{1},
\dots,n_{j-1},p,q,n_{j+1},\dots,n_{k-1}}
\end{displaymath}
with differentials
\begin{displaymath}
\dd_{i} x=  
\begin{cases}
\dd_{A,i}x,&\text{ if }i<j, \\
\dd_{A,i}x+(-1)^{n_{i}}\dd_{A,i+1}x,&\text{ if }i=j, \\
\dd_{A,i+1}x,&\text{ if }i>j,
\end{cases}
\end{displaymath}
where $x\in A^{n_{1},\dots,n_{j},\dots,n_{k}}$.
\end{definition}

\nnpar{Signs.} 
The simple complex of a $k$-iterated complex depends on the choice 
of an ordering of the degrees. If we choose a different ordering, 
we obtain a naturally isomorphic complex. In order to describe
explicitly these isomorphisms it suffices to treat the case of 
the transposition of two adjacent degrees.

\begin{definition}
Let $T_{i,i+1}$ be the functor on $(\mathcal{C}^{+})^{k}(\mathcal{A})$ 
to itself given, for $A\in(\mathcal{C}^{+})^{k}(\mathcal{A})$, by
\begin{displaymath}
T_{i,i+1}(A)^{n_{1},\dots,n_{i},n_{i+1},\dots,n_{k}}=
A^{n_{1},\dots,n_{i+1},n_{i},\dots,n_{k}}.
\end{displaymath}
\end{definition}

\begin{proposition} 
\label{prop:7}
Let $A$ be a $k$-iterated complex. The map 
$$
s(A)\longrightarrow s(T_{i,i+1}(A))
$$ 
given by
\begin{displaymath}
x\longmapsto(-1)^{n_{i}n_{i+1}}x\qquad(x\in A^{n_{1},\dots,n_{k}})
\end{displaymath}  
is an isomorphism of simple complexes. 
\hfill $\square$
\end{proposition}

\nnpar{External product.} 
\begin{definition}
Let $A$ be a $k$-iterated complex, and let $B$ be an $l$-iterated
complex. Then, \emph{the external product} $A\boxtimes B$ is the
$(k+l)$-iterated complex with groups
\begin{displaymath}
(A\boxtimes B)^{n_{1},\dots,n_{k+l}}=A^{n_{1},\dots,n_{k}}\otimes 
B^{n_{k+1},\dots,n_{k+l}},
\end{displaymath}
and differentials
\begin{displaymath}
\dd_{j}=
\begin{cases}
\dd_{i}\otimes\Id,&\text{ if }i\le k, \\
\Id\otimes\dd_{i-k},&\text{ if }i> k. 
\end{cases}
\end{displaymath}
\end{definition}

The definitions of the functors $\mathcal{C}_{k}$ and of the simple 
of a $k$-iterated complex are chosen in order to have compatibility 
with the usual sign convention in the tensor product of complexes.

\begin{lemma} 
\label{lem:simple-iso}
Let $A$ be a $k$-iterated and $B$ be an $l$-iterated complex. The map 
\begin{displaymath}
s(A)\otimes s(B)\longrightarrow s(A\boxtimes B)
\end{displaymath}
given by $a\otimes b\mapsto a\otimes b$ is an isomorphism of 
complexes.
\hfill $\square$
\end{lemma}

\nnpar{Tensor product.} 
Let us show how to define the tensor product in the category of
$2$-iterated complexes.

\begin{definition}
\label{def:2-tensor}
Let $A,B$ be a pair of $2$-iterated complexes. Then, the
\emph{tensor product} $A\otimes B$ of $A$ and $B$ is the
$2$-iterated complex
$$
A\otimes B=s_{1,2}s_{3,4}(T_{2,3}(A\boxtimes B)).
$$  
\end{definition}

\begin{remark} 
This definition is justified by the fact that $A\otimes B$ satisfies 
the universal properties of a tensor product of $A$ and $B$ in the 
category of $2$-iterated complexes. In particular, there are canonical
isomorphisms $A\otimes B\cong B\otimes A$ and $A\otimes\left( B\otimes
C\right)\cong\left(A\otimes B\right)\otimes C$; moreover, if $A$, $B$, 
$C$, $D$ are complexes of $R$-modules, then the identity map on the
level of elements induces an isomorphism
\begin{displaymath}
(A\boxtimes B)\otimes(C\boxtimes D)\cong
(A\otimes C)\boxtimes(B\otimes D).
\end{displaymath}
\end{remark}

\nnpar{Examples.} 
Let us  discuss some fundamental examples of $2$-iterated complexes.

\begin{example} 
\label{ex:fundamental}
Let $f:(A^{*},\dd_{A})\longrightarrow (B^{*},\dd_{B})$ be a 
morphism of complexes. Then, we can consider the $2$-iterated 
complex $(f^{*,*},\dd',\dd'')$ defined by
\begin{displaymath}
f^{0,q}=A^{q},\qquad f^{1,q}=B^{q},\qquad\dd'=f,\qquad
\dd''|_{f^{0,*}}=\dd_{A},\qquad\dd''|_{f^{1,*}}=\dd_{B}.
\end{displaymath}
By abuse of notation we denote this $2$-iterated complex by $f$. 
By definition, the simple complex $s(f)=s(A,B)$ of $f$ is given 
by $s(f)=s(\mathcal{C}_{2}(f^{*,*}))$. Explicitly, we have
\begin{displaymath}
s(f)^{n}=s(A,B)^n=A^{n}\oplus B^{n-1},\quad\dd(a,b)=
(\dd_{A}a,f(a)-\dd_{B}b). 
\end{displaymath}
\end{example}

\begin{definition}
Given a complex $(A^{*},\dd)$, and an integer $k$, the shifted
complex $A[k]$ is given by $A[k]^{n}=A^{n+k}$ with differential
$(-1)^{k}\dd$.
\end{definition}

\begin{remark} 
Recall that the mapping cone of a morphism $f:A\longrightarrow 
B$ of complexes is defined as 
\begin{displaymath}
\cone(f)^{n}=A^{n+1}\oplus B^{n},\quad\dd(a,b)=
(-\dd_{A}a,f(a)+\dd_{B}b). 
\end{displaymath}
Therefore, we have $s(f)=\cone(-f)[-1]$.
\end{remark}

\begin{example} 
\label{ex:comp1}
Let $f:A\longrightarrow B$ and $g:B\longrightarrow C $ be morphisms
of complexes with $g\circ f=0$. We may consider the diagram
\begin{displaymath}
\begin{CD}
\xi =(A@>f>>B@>g>>C)     
\end{CD}
\end{displaymath}
as a $2$-iterated complex $(\xi^{*,*},\dd',\dd'')$, where the groups
are 
\begin{displaymath}
\xi^{0,q}=A^{q},\qquad\xi^{1,q}=B^{q},\qquad\xi^{2,q}=C^{q},
\end{displaymath}
the differential $\dd'$ is either $f$ or $g$ and the differential
$\dd''$ is the differential of the complexes $A$, $B$ or $C$. By 
abuse of notation, we denote this $2$-iterated complex by $\xi$.
\end{example}

\begin{proposition}
\label{prop:comp1}
With the notations of example \ref{ex:comp1}, let $\iota:A
\longrightarrow s(-g)$ be the morphism given by $\iota(a)=
(f(a),0)$. Then, there is a natural isomorphism of complexes 
\begin{displaymath}
s(\xi)\longrightarrow s(A\overset{\iota}{\longrightarrow}s(-g))
\end{displaymath}
given by $(a,b,c)\longmapsto(a,(b,c))$. 
\hfill $\square$
\end{proposition}

\nnpar{The simple complex of a tensor product.} 
For the rest of this section $f_{1}:A_{1}\longrightarrow B_{1}$ and
$f_{2}:A_{2}\longrightarrow B_{2}$ will be morphisms of complexes.
Considering these morphisms as $2$-iterated complexes as in example
\ref{ex:fundamental}, we have by definition \ref{def:2-tensor} 
\begin{align} 
\label{eq:2-tensor}
f_{1}\otimes f_{2}=s_{1,2}s_{3,4}(T_{2,3}(f_{1}\boxtimes f_{2})).
\end{align} 
This $2$-iterated complex is naturally identified with the
$2$-iterated complex associated to the diagram
\begin{align} 
\label{eq:2-diagram}
\begin{CD}
\xi= \big(A_{1}\otimes A_{2}@>(f_{1}\otimes\Id,\Id\otimes f_{2})>> 
B_{1}\otimes A_{2}\oplus A_{1}\otimes B_{2}@>-\Id\otimes f_{2}+
f_{1}\otimes \Id>>B_{1}\otimes B_{2}\big)
\end{CD}
\end{align}
as in example \ref{ex:comp1}.

\begin{remark}
Note that, if we consider $f_{1}:A_{1}\longrightarrow B_{1}$ and 
$f_{2}:A_{2}\longrightarrow B_{2}$ as morphisms and not as
$2$-iterated complexes, then $f_{1}\otimes f_{2}$ is the morphism 
$f_{1}\otimes f_{2}:A_{1}\otimes A_{2}\longrightarrow B_{1}\otimes 
B_{2}$. It will be clear from the context which point of view is 
adopted in each case. If there is a danger of confusion, we will 
point out whether $f_{1}\otimes f_{2}$ is a $2$-iterated complex 
or a morphism. 
\end{remark}

\begin{proposition}
\label{prop:funiso}
There is an isomorphism of complexes
\begin{displaymath}
s(f_{1})\otimes s(f_{2})\longrightarrow s(f_{1}\otimes f_{2})
\end{displaymath}
given by
\begin{equation}
\label{eq:16}
(a_{1},b_{1})\otimes(a_{2},b_{2})\longmapsto(a_{1}\otimes a_{2},
(b_{1}\otimes a_{2},(-1)^{n}a_{1}\otimes b_{2}),(-1)^{n-1}b_{1}
\otimes b_{2}),
\end{equation}  
where $(a_{1},b_{1})\in s(f_{1})^{n}$ and $(a_{2},b_{2})\in 
s(f_{2})^{m}$.
\end{proposition}
\begin{proof} 
By lemma \ref{lem:simple-iso}, we have an isomorphism $s(f_{1})
\otimes s(f_{2})\cong s(f_{1}\boxtimes f_{2})$. Since $s(f_{1}
\boxtimes f_{2})\cong s(T_{2,3}(f_{1}\boxtimes f_{2}))$ by
proposition \ref{prop:7}, and since $s(T_{2,3}(f_{1}\boxtimes
f_{2}))=s(s_{1,2}s_{3,4}(T_{2,3}(f_{1}\boxtimes f_{2})))$, the
existence of the claimed isomorphism follows by definition
\ref{def:2-tensor} of $f_{1}\otimes f_{2}$.
  
Let us compute this isomorphism explicitly: For this we take
$(a_{1},b_{1})\in s(f_{1})^{n}$ and $(a_{2},b_{2})\in s(f_{2})^{m}$.
Viewing $a_{1}$, $b_{1}$, resp. $a_{2}$, $b_{2}$, as elements of 
the $2$-iterated complex $f_{1}$, resp. $f_{2}$, the elements
$a_{1}$, $b_{1}$, $a_{2}$, $b_{2}$ have bidegree equal to $(0,n)$,
$(1,n-1)$, $(0,m)$, $(1,m-1)$, respectively.  Then, the element
$(a_{1},b_{1})\otimes(a_{2},b_{2})\in s(f_{1})\otimes s(f_{2})$ 
is mapped to the element
$$
a_{1}\otimes a_{2}+b_{1}\otimes a_{2}+a_{1}\otimes b_{2}+ 
b_{1}\otimes b_{2}
$$
of $s(f_{1}\boxtimes f_{2})$, where the summands have the following 
$4$-degree:
\begin{center}
\begin{tabular}{c|c}
element & 4-degree \\
\hline
$a_{1}\otimes a_{2}$ & $(0,n,0,m)$ \\
$b_{1}\otimes a_{2}$ & $(1,n-1,0,m)$ \\
$a_{1}\otimes b_{2}$ & $(0,n,1,m-1)$ \\
$b_{1}\otimes b_{2}$ & $(1,n-1,1,m-1)$
\end{tabular}
\end{center}
Applying $T_{2,3}$ to $f_{1}\boxtimes f_{2}$, the latter element 
is mapped to
$$
a_{1}\otimes a_{2}+b_{1}\otimes a_{2}+(-1)^{n}a_{1}\otimes b_{2}+
(-1)^{n-1}b_{1}\otimes b_{2},
$$
which by the identification \ref{eq:2-diagram} is mapped to the element
claimed in formula \ref{eq:16}.
\end{proof}

\nnpar{Commutativity of the external product.} 
The following observations will be needed in the discussion of
commutativity and associativity of the product in relative 
cohomology defined in the next section.

\begin{proposition} 
\label{prop:commten}
There is a commutative diagram of complexes
\begin{displaymath}
\begin{CD}
s(f_{1})\otimes s(f_{2})@>\alpha_{1}>>s(f_{1}\otimes f_{2}) \\
@V\alpha_{2}VV@VV\beta_{2}V \\
s(f_{2})\otimes s(f_{1})@>\beta_{1}>>s(f_{2}\otimes f_{1}), 
\end{CD}
\end{displaymath}
where $\alpha_{1}$ and $\beta_{1}$ are the isomorphisms 
determined  in \eqref{eq:16}, $\alpha_2$ is given by
\begin{align}
\label{eq:com-iso1}
\alpha_{2}((a_{1},b_{1})\otimes(a_{2},b_{2}))&=
(-1)^{nm}(a_{2},b_{2})\otimes(a_{1},b_{1}), 
\end{align}
and $\beta_{2}$ is given by 
\begin{align} 
\label{eq:com-iso2}
\!&\beta_{2}( a_{1}\otimes a_{2},(b_{1}\otimes a_{2},  
a_{1}\otimes b_{2}),b_{1}\otimes b_{2})=\notag \\
&((-1)^{nm}a_{2}\otimes a_{1},((-1)^{n(m-1)}b_{2}\otimes a_{1},
(-1)^{(n-1)m}a_{2}\otimes b_{1}),(-1)^{nm-n-m}b_{2}\otimes b_{1})
\end{align}
with $\deg(a_{1})=n$, $\deg(b_{1})=n-1$, $\deg(a_{2})=m$ and 
$\deg(b_{2})=m-1$.
\end{proposition}
\begin{proof} First we note that $\beta_{2}$ is well defined by 
proposition \ref{prop:funiso}. The commutativity can now be checked 
directly. The signs in the definition of the morphisms $\alpha_{2}$ 
and $\beta_{2}$ are obtained as in the proof of proposition
\ref{prop:funiso}.
\end{proof}

\nnpar{Associativity of the external product.} 
Furthermore, let $f_{3}:A_{3}\longrightarrow B_{3}$ be a third
morphism of complexes, and put $\phi_{1}=f_{1}\otimes\Id\otimes\Id$, 
$\phi_{2}=\Id\otimes f_{2}\otimes\Id$, $\phi_{3}=\Id\otimes\Id 
\otimes f_{3}$. Then, the $2$-iterated complexes $(f_{1}\otimes
f_{2})\otimes f_{3}$ and $f_{1}\otimes(f_{2}\otimes f_{3})$ are 
both naturally identified with the $2$-iterated complex associated 
to the diagram     
\begin{align}
A_{1}\otimes A_{2}\otimes A_{3}\overset{\delta_{1}}{\longrightarrow}&
B_{1}\otimes A_{2}\otimes A_{3}\oplus A_{1}\otimes B_{2}\otimes A_{3}
\oplus A_{1}\otimes A_{2}\otimes B_{3}\notag \\
\overset{\delta_{2}}{\longrightarrow}&B_{1}\otimes B_{2}\otimes A_{3}
\oplus B_{1}\otimes A_{2}\otimes B_{3}\oplus A_{1}\otimes B_{2}\otimes 
B_{3}\notag \\
\overset{\delta_{3}}{\longrightarrow}&B_{1}\otimes B_{2}\otimes B_{3},
\label{eq:secass} 
\end{align}
where
\begin{align*}
\delta_{1}(a)&=(\phi_{1}a,\phi_{2}a,\phi_{3}a), \\
\delta_{2}(a,b,c)&=(\phi_{1}b-\phi_{2}a,\phi_{1}c-\phi_{3}a,
\phi_{2}c-\phi_{3}b), \\ 
\delta_{3}(a,b,c)&=\phi_{3}a-\phi_{2}b+\phi_{1}c.
\end{align*}
In the sequel we will denote any of the two $2$-iterated complexes
$(f_{1}\otimes f_{2})\otimes f_{3}$ or $f_{1}\otimes(f_{2}\otimes
f_{3})$ by $f_{1}\otimes f_{2} \otimes f_{3}$.

\begin{proposition}
\label{prop:assten}
There is a commutative diagram of complexes
\begin{displaymath}
\begin{CD}
s(f_{1})\otimes s(f_{2})\otimes s(f_{3})@>\alpha_{1}>>
s(f_{1}\otimes f_{2})\otimes s(f_{3}) \\
@V\alpha_{2}VV@VV\beta_{2}V \\
s(f_{1})\otimes s(f_{2}\otimes f_{3})@>\beta_{1}>>
s(f_{1}\otimes f_{2}\otimes f_{3}), 
\end{CD}
\end{displaymath}
where $\alpha_{1}$ and $\alpha_{2}$ are isomorphisms induced by
\eqref{eq:16}, $\beta_{1}$ is given by
\begin{align}
\beta_{1}((a,b)\otimes(c,(d,e),f))&=(a\otimes c,(b\otimes c, 
(-1)^{n}a\otimes d,(-1)^{n}a\otimes e),\notag \\
&((-1)^{n-1}b\otimes d,(-1)^{n-1}b\otimes e,a\otimes f),b\otimes f),
\end{align}
with $n=\deg(a,b)$, and $\beta_{2}$ is given by
\begin{align}
\beta_{2}((a,(b,c),d)\otimes(e,f))&=(a\otimes e,(b\otimes e, 
c\otimes e,(-1)^{n}a\otimes f),\notag \\
&(d\otimes e,(-1)^{n-1}b\otimes f,(-1)^{n-1}c\otimes f),
(-1)^{n-2}d\otimes f),
\end{align}
with $n=\deg(a,(b,c),d)$.
\end{proposition}
\begin{proof}
The commutativity can be checked directly. The signs in the definition
of the morphisms $\beta_{1}$ and $\beta_{2}$ are obtained as in the 
proof of proposition \ref{prop:funiso}.
\end{proof}

\subsection{Relative cohomology groups}
\label{sec:rcg}

\nnpar{Relative cohomology.}
\begin{definition}
Let $f:A\longrightarrow B$ be a morphism of complexes. Then, the
\emph{relative cohomology groups (of the pair $A,B$)} are defined 
by
\begin{displaymath}
H^{*}(A,B)=H^{*}(s(f)).
\end{displaymath}
For a cycle $(a,b)\in s(f)^{n}$ we will denote by $[a,b]$ its class
in $H^{n}(A,B)$. Observe that a cocycle of $s(f)$ is a pair $(a,b)$
with $\dd_{A}a=0$ and $\dd_{B}b=f(a)$. Moreover, the subspace of
coboundaries is generated by elements of the form $(\dd_{A}a,f(a))$ 
and $(0,\dd_{B}b)$.
\end{definition}

\nnpar{Natural maps.} 
There are two natural maps relating the simple $s(f)$ of a morphism 
$f:A\longrightarrow B$ with the underlying complexes $A,B$: 
\begin{displaymath}
\alpha:s(f)\longrightarrow A,\,\,\,\,{\rm given\,\,by}\,\,\,\,
(a,b)\longmapsto a,
\end{displaymath}
and
\begin{displaymath}
\beta:B\longrightarrow s(f)[1],\,\,\,\,{\rm given\,\,by}\,\,\,\,
b\longmapsto (0,b).
\end{displaymath}
We will denote also by $\alpha$ and $\beta $ the corresponding 
morphisms for the cone given by the same formulas.

\nnpar{Distinguished exact triangles.} 
Given a short exact sequence of complexes there are two possible
choices for the connecting morphism; one is the opposite of the 
other. To choose the sign of the connecting morphism is equivalent 
to choose a set of distinguished exact triangles. 

Given a morphism of complexes $f:A\longrightarrow B$, one has 
chosen in \cite{SGA4_1/2}, p. 269, the exact triangle 
\begin{equation}
\label{eq:31}
A\overset{f}{\longrightarrow}B\overset{\beta}{\longrightarrow}
\cone{(f)}\overset{-\alpha}{\longrightarrow}A[1]
\end{equation}
as a distinguished exact triangle. 

There are two basic principles that can be applied to obtain a
distinguished exact triangle from another. The first one is that 
any exact triangle isomorphic (in the derived category) to a 
distinguished exact triangle is also distinguished. For instance, 
by means of the isomorphism $\cone(f)\cong s(f)[1]$ given by
mapping $(a,b)$ to $(-a,b)$, we find that the exact triangle 
\begin{equation}
\label{eq:32}
A\overset{f}{\longrightarrow}B\overset{\beta}{\longrightarrow}
s(f)[1]\overset{\alpha[1]}{\longrightarrow}A[1]
\end{equation}
is also distinguished. As another example we can change the sign 
of two out of the three morphisms of a distinguished exact triangle, 
and the resulting exact triangle will again be distinguished; this 
results from the following isomorphism of exact triangles
\begin{displaymath}
\xymatrix{
A\ar[r]^{u}\ar[d]_{\Id}&B\ar[r]^{v}\ar[d]_{-\Id}&
C\ar[r]^{w}\ar[d]_{\Id}&... \\
A\ar[r]^{-u}&B\ar[r]^{-v}&C\ar[r]^{w}&...
}
\end{displaymath}
The second principle is that a distinguished exact triangle can 
be shifted by one and the resulting exact triangle is again exact
provided that the sign of one of the morphisms is changed. 

Applying these two principles to the distinguished exact triangle 
\eqref{eq:31}, we see that the exact triangles  
\begin{gather}
\label{eq:37}
s(f)\overset{\alpha}{\longrightarrow}A\overset{f}{\longrightarrow}
B\overset{-\beta}{\longrightarrow}...\,\,,\\
\label{eq:38}
s(-f)\overset{\alpha}{\longrightarrow}A\overset{f}{\longrightarrow}
B\overset{\beta}{\longrightarrow}...
\end{gather}
are also distinguished. We note that these exact triangles are 
isomorphic by means of the isomorphism
\begin{equation}
\label{eq:42}
s(f)\longrightarrow s(-f),\,\,\,\,{\rm given\,\,by}\,\,\,\,
(a,b)\longmapsto(a,-b).
\end{equation}

\nnpar{Kernel, cokernel, and the simple.} 
Let 
\begin{displaymath}
0\longrightarrow A\overset{f}{\longrightarrow}B\overset{g}
{\longrightarrow}C\longrightarrow 0
\end{displaymath}
be a short exact sequence of complexes. The sign of the connecting
homomorphism $\delta$ is determined by the isomorphisms (in the 
derived category) of exact triangles 
\begin{equation}
\label{eq:39}
\xymatrix{
A\ar[r]^{f}\ar[d]_{\Id}&B\ar[r]^{\beta}\ar[d]_{\Id}&
s(f)[1]\ar[r]^{\alpha[1]}\ar[d]_{\pi}&... \\
A\ar[r]^{f}&B\ar[r]^{g}&C\ar[r]^{\delta}&...
}
\end{equation}
with $\pi(a,b)=g(b)$, and
\begin{equation}
\label{eq:40}
\xymatrix{
s(-g)\ar[r]^{\alpha}&B\ar[r]^{g}&C\ar[r]^{\beta}&... \\
A\ar[r]^{f}\ar[u]^{\iota}&B\ar[r]^{g}\ar[u]^{\Id}&
C\ar[r]^{\delta}\ar[u]^{\Id}&...
}
\end{equation}
with $\iota(a)=(f(a),0)$.

\begin{definition} 
\label{def:MVP}
The morphism $\pi$ will be called the \emph{simple-cokernel
quasi-isomorphism}, and the morphism $\iota$ will be called 
the \emph{kernel-simple quasi-isomorphism}.
\end{definition}

The sign of the connecting morphism determined by the above
isomorphisms is given as follows. If $c$ is a cycle in $C^{n}$
satisfying $c=g(b)$ with $b\in B^{n}$, then we have
\begin{equation}
\label{eq:41}
\delta(c)=f^{-1}(\dd_{B}b).
\end{equation}

\begin{remark}
\label{rem:6}
Observe that the above connecting morphism is minus the connecting 
morphism considered in \cite{Burgos:CDB}. This forces us to adjust
many of the signs. 
\end{remark}

\nnpar{The connecting morphism in relative cohomology.} 
Since the exact triangle \eqref{eq:37} is distinguished, we are 
led to the following definition 

\begin{definition}
\label{def:lo-ex-seq}
Let $f:A\longrightarrow B$ be a morphism of complexes. The
\emph{connecting homomorphism in relative cohomology} is the 
morphism 
\begin{displaymath}
\delta:H^{n-1}(B)\longrightarrow H^{n}(A,B)
\end{displaymath}
induced by the natural map $-\beta$.
\end{definition}

\nnpar{Split exact sequences.} 
Let 
\begin{displaymath}
0\longrightarrow A\overset{f}{\longrightarrow}B\overset{g}
{\longrightarrow}C\longrightarrow 0
\end{displaymath}
be a short exact sequence of complexes. Assume that the short 
exact sequence  
\begin{displaymath}
0\longrightarrow A^{n}\overset{f}{\longrightarrow}B^{n}
\overset{g}{\longrightarrow}C^{n}\longrightarrow 0
\end{displaymath}
is split for all $n$. Thus there are sections $\sigma:C^{n}
\longrightarrow B^{n}$ such that $g\circ\sigma=\Id_{C^{n}}$. 
These sections allow us to give a quasi-inverse of the 
simple-cokernel quasi-isomorphism $\pi$ and the kernel-simple 
quasi-isomorphism $\iota$ 

\begin{proposition}
\label{prop:quisinv}
Let  
\begin{displaymath}
0\longrightarrow A\overset{f}{\longrightarrow}B\overset{g}
{\longrightarrow}C\longrightarrow 0
\end{displaymath}
be a short exact sequence as above such that for each $n$ there
is a section $\sigma:C^{n}\longrightarrow B^{n}$ of $g$. Then,
we have the following statements:
\begin{enumerate}
\item[(i)]
The map $\iota':C\longrightarrow s(f)[1]$ given by
\begin{displaymath}
\iota'(c)=(f^{-1}(\dd_{B}\sigma(c)-\sigma(\dd_{C}c)),\sigma(c)) 
\end{displaymath}
is a morphism of complexes, and satisfies $\pi\circ\iota'=
\Id_{C}$ and $\iota'\circ\pi\sim\Id_{s(f)[1]}$, where $\sim$ 
means homotopically equivalent.
\item[(ii)] 
The map $\pi':s(-g)\longrightarrow A$ given by
\begin{displaymath}
\pi'(b,c)=f^{-1}\left(b-\sigma(g(b))-\sigma(\dd_{C}c)+
\dd_{B}\sigma(c)\right) 
\end{displaymath}
is a morphism of complexes satisfying $\pi'\circ\iota=\Id_{A}$,
and $\iota\circ\pi'\sim\Id_{s(-g)}$, where $\sim$ means again
homotopically equivalent.   
\end{enumerate}
\hfill $\square$
\end{proposition}

\nnpar{Change of complex.} 
The following result is obvious. We quote it for future reference.

\begin{proposition}
\label{prop:quis2}
Let $f:A\longrightarrow B $ and $g:B\longrightarrow C $ be morphisms 
of complexes such that $g$ is a quasi-isomorphism. Then, there is a 
natural quasi-isomor\-phism $s(f)\longrightarrow s(g\circ f)$. 
\hfill $\square$
\end{proposition}

\nnpar{De Rham cohomology with support.}
\begin{example}
\label{ex:coh-supp}
Let $X$ be a differentiable manifold and $Y$ a closed subset of 
$X$. Let us denote by $E(X)$ the complex of $\mathbb{C}$-valued 
differential forms and by $\sigma:E(X)\longrightarrow E(X\setminus 
Y)$ the restriction morphism. Then, the relative de Rham cohomology 
groups of the pair $(X,X\setminus Y)$ are the cohomology groups of 
the simple $s(\sigma)$. By standard arguments in sheaf theory these 
groups can be naturally identified with the cohomology groups of 
the constant sheaf $\mathbb{C}$ with support in $Y$, denoted by 
$H^{\ast}_{Y}(X,\mathbb{C})$. Now let $Z$ be another closed subset 
of $X$ and $j$ the morphism  
\begin{align*}
j:&\,E(X\setminus Y)\oplus E(X\setminus Z)\longrightarrow 
E(X\setminus(Y\cup Z))
\end{align*}
given by $j(\omega,\eta)=-\omega+\eta$. The restriction of 
differential forms induces the Mayer-Vietoris sequence 
\begin{displaymath}
0\rightarrow  E(X\setminus(Y\cap Z))\rightarrow E(X\setminus Y)
\oplus E(X\setminus Z)\overset{j}{\rightarrow}E(X\setminus(Y\cup Z))
\rightarrow 0,
\end{displaymath}
where the first map assigns to $\omega$ the pair $(\omega,\omega)$.
In this case the kernel-simple quasi-isomorphism of definition
\ref{def:MVP} is the quasi-isomorphism
\begin{align}
\label{cor:1}
\iota:E(X\setminus(Y\cap Z))\longrightarrow s(-j).
\end{align}
Using the quasi-isomorphism \eqref{cor:1} together with proposition
\ref{prop:quis2}, we obtain  an isomorphism  
\begin{align}
\label{eq:coh-supp-sim}
H^{n+m}_{Y\cap Z}(X,\mathbb{C})\cong H^{n+m}(E(X),s(-j)).
\end{align}
Using partitions of unity and proposition \ref{prop:quisinv}, one 
can construct a quasi-inverse of the kernel-simple quasi-isomorphism 
$\iota$, hence an inverse of the isomorphism \eqref{eq:coh-supp-sim}.
\end{example}

\subsection{Products in relative cohomology}
\label{sec:prc}

\nnpar{External product.} 
If $A$ and $B$ are complexes of $R$-modules, there is a well-defined 
external product  
\begin{displaymath}
H^{\ast}(A)\otimes H^{\ast}(B)\longrightarrow H^{\ast}(A\otimes B). 
\end{displaymath}
In particular, there is an external product in relative cohomology
\begin{displaymath}
H^{n}(s(f_{1}))\otimes H^{m}(s(f_{2}))\longrightarrow 
H^{n+m}(s(f_{1})\otimes s(f_{2})).
\end{displaymath}

\begin{proposition} 
\label{thm:pairing-relative}
The external product in relative cohomology can be identified with 
the pairing
\begin{displaymath}
H^{n}(s(f_{1}))\otimes H^{m}(s(f_{2}))\longrightarrow 
H^{n+m}\left(A_{1}\otimes A_{2},s(\Id\otimes f_{2}- f_{1}\otimes 
\Id)\right)
\end{displaymath}
given by formula \eqref{eq:16}.
\end{proposition}
\begin{proof} 
The result is a direct consequence of proposition \ref{prop:funiso}, 
the identification \eqref{eq:2-diagram}, and proposition \ref{prop:comp1}.  
\end{proof} 

\nnpar{Cup product.} 
We have defined an external product in relative cohomology. Moreover, 
if there is a product defined on the level of complexes, then we 
obtain a cup product in relative cohomology.

\begin{theorem} 
\label{prop:prc} 
Let
\begin{equation}
\label{eq:bullet}
\begin{CD}
{\begin{CD}
A_{1}\otimes B_{2}@>f_{1}\otimes\Id>>B_{1}\otimes B_{2} \\
@A\Id\otimes f_{2}AA@AA\Id\otimes f_{2}A \\
A_{1}\otimes A_{2}@>f_{1}\otimes\Id>>B_{1}\otimes A_{2}
\end{CD}}
@>\bullet>>
{\begin{CD}
E_{0,1}@>\delta>>E_{1,1} \\
@A\alpha AA@AA\gamma A \\
E_{0,0}@>\beta>>E_{1,0}
\end{CD}}
\end{CD}
\end{equation}
be a morphism of commutative diagrams in the category of complexes. 
Let $i:E_{0,0}\longrightarrow E_{1,0}\oplus E_{0,1}$, resp. $j:E_{1,0}
\oplus E_{0,1}\longrightarrow E_{1,1}$, be the morphism given by 
$$
i(x)=(\beta(x),\alpha(x)),\quad\text{resp.}\quad j(x,y)=-\gamma(x)+
\delta(y).
$$
Then, there is a pairing
\begin{displaymath}
H^{n}(A_{1},B_{1})\otimes H^{m}(A_{2},B_{2})\overset{\bullet}
{\longrightarrow}H^{n+m}(E_{0,0},s(-j))
\end{displaymath}
given by
\begin{displaymath}
[a_{1},b_{1}]\bullet[a_{2},b_{2}]=[a_{1}\bullet a_{2},((b_{1}\bullet 
a_{2},(-1)^{n}a_{1}\bullet b_{2}),(-1)^{n-1}b_{1}\bullet b_{2})];
\end{displaymath}
note that the map $E_{0,0}\longrightarrow s(-j)$ is given by $x\mapsto 
(i(x),0)$.
\end{theorem}
\begin{proof} 
By proposition \ref{prop:funiso} there is a pairing
\begin{displaymath}
H^{n}(A_{1},B_{1})\otimes H^{m}(A_{2},B_{2})\longrightarrow
H^{n+m}(s(\xi)),  
\end{displaymath}
where $\xi$ is given by \eqref{eq:2-diagram}. Let $\eta$ be the 
$2$-iterated complex 
\begin{displaymath}
\eta=(E_{0,0}\overset{i}{\longrightarrow}E_{1,0}\oplus E_{0,1}
\overset{j}{\longrightarrow}E_{1,1}).
\end{displaymath}
From the data given in the assumptions, we obtain a morphism $\bullet$
of $2$-iterated complexes 
$$
\bullet:\xi\longrightarrow\eta,
$$ 
hence we obtain a morphism $H^{n+m}(s(\xi))\longrightarrow H^{n+m}
(s(\eta))$. By proposition \ref{prop:comp1} the latter cohomology 
group is isomorphic to $H^{n+m}(E_{0,0},s(-j))$, which completes
the proof. 
\end{proof}

\begin{remark} 
\label{prevrmk}
Letting $A_{1}=0$ in the above notations and observing that $H^{n}
(0,B_{1})\cong H^{n-1}(B_{1})$, we obtain from theorem \ref{prop:prc} 
with $E_{0,0}=0$ a pairing, also denoted by $\bullet$,
\begin{equation}
\label{eq:pairingb}   
H^{n-1}(B_{1})\otimes H^{m}(A_{2},B_{2})\overset{\bullet}
{\longrightarrow}H^{n+m-1}(s(-j))
\end{equation}
given by
\begin{displaymath}
[b_{1}]\bullet[a_{2},b_{2}]=[(b_{1}\bullet a_{2},0),(-1)^{n-1}b_{1}
\bullet b_{2}].
\end{displaymath}
\end{remark}

\begin{corollary} 
With the assumptions of theorem \ref{prop:prc} and remark 
\ref{prevrmk} we have a commutative diagram 
\begin{displaymath}
\begin{CD}
H^{n-1}(B_{1})\otimes H^{m}(A_{2},B_{2})@>\bullet>>H^{n+m-1}
(s(-j)) \\
@V\delta\otimes\Id VV@VV\delta V \\
H^{n}(A_{1},B_{1})\otimes H^{m}(A_{2},B_{2})@>\bullet>>H^{n+m}
(E_{0,0},s(-j)), 
\end{CD}
\end{displaymath}
where $\delta$ is the connecting morphism in the long exact sequence 
of relative cohomology described in \ref{def:lo-ex-seq}. 
\hfill $\square$
\end{corollary}

\nnpar{The cup product in de Rham cohomology with support.} 
\begin{example}
\label{ex:cup-product-de}
Recall the notations of example \ref{ex:coh-supp}. Let us illustrate 
how the formalism developed in theorem \ref{prop:prc} allows us to 
compute the product   
\begin{displaymath} 
H^{n}_{Y}(X,\mathbb{C})\otimes H^{m}_{Z}(X,\mathbb{C})\longrightarrow 
H^{n+m}_{Y\cap Z}(X,\mathbb{C})
\end{displaymath}
by means of differential forms. Let us write
\begin{displaymath}
\begin{array}{ll}
A_{1}=A_{2}=E_{0,0}=E(X),&B_{1}=E_{1,0}=E(X\setminus Y), \\
B_{2}=E_{0,1}=E(X\setminus Z),&E_{1,1}=E(X\setminus(Y\cup Z)),
\end{array}
\end{displaymath}
and let $\bullet$ be the wedge product $\land$. Then, the assumptions 
of theorem \ref{prop:prc} are satisfied. Hence the wedge product induces 
a pairing  
\begin{equation}
\label{eq:17}
H^{n}_{Y}(X,\mathbb{C})\otimes H^{m}_{Z}(X,\mathbb{C})\longrightarrow 
H^{n+m}(E(X),s(-j))
\end{equation}
given, for appropriate cocycles $(\omega,\eta)$ and $(\omega',\eta')$, 
by
\begin{equation}
\label{eq:cycle}
(\omega\land\omega',((\eta\land\omega',(-1)^{n}\omega\land\eta'), 
(-1)^{n-1}\eta\land\eta'));
\end{equation}
this is a cocycle in $s(E(X),s(-j))$. Then, the desired pairing 
is defined by composing the pairing \eqref{eq:17} with the inverse 
of the isomorphism \eqref{eq:coh-supp-sim}. In particular, the 
cycle \eqref{eq:cycle} represents a class in $H^{n+m}_{Y\cap Z}
(X,\mathbb{C})$.

Using standard arguments from sheaf theory one can see that the 
product constructed above agrees with the sheaf theoretic product 
in cohomology with support (see, e.g., \cite{Burgos:Gftp}, proposition 
2.5). One can also compare this construction with the definition 
of the Loday product in $K$-theory with support given in section 
\ref{sec:ktgc}.
\end{example}

\subsection{Truncated relative cohomology groups}
\label{sec:tcdb}

\nnpar{Definition of truncated relative cohomology groups.} 
Let $A=(A^{\ast},\dd_{A})$ be a complex of $R$-modules. We will 
denote by ${\rm Z}^{n}(A)$ the submodule of cocycles of $A^{n}$ 
and by $\widetilde{A}^{n}= A^{n}/\Img\dd_A$. For $a\in A^{n}$ 
we write $\widetilde{a}$ for its class in $\widetilde{A}^{n}$. 
Recall that the b\^ete filtration $\sigma_{p}A$ of $A$ is given 
by
\begin{displaymath}
(\sigma_{p}A)^{n}=
\begin{cases}
A^{n},&\text{ if }n\ge p, \\
0,&\text{ if }n<p. 
\end{cases}
\end{displaymath}

\begin{definition}
Let $f: (A^{\ast},\dd_{A})\longrightarrow(B^{\ast},\dd_{B}) $ 
be a morphism of complexes of $R$-modules. The \emph{truncated 
relative cohomology groups} associated with $f$ are defined by
\begin{align*}
\widehat{H}^{n}(A,B)=H^{n}(\sigma_{n}A,B).
\end{align*} 
In other words,
\begin{align*}
\widehat{H}^{n}(A,B)=\left\{(a,\widetilde{b})\in{\rm Z}^{n}(A)
\oplus\widetilde{B}^{n-1}\big|f(a)=\dd_{B}b\right\}.
\end{align*} 
\end{definition}

\nnpar{Basic properties.} 
If $(a,b)\in{\rm Z}^{n}\left(s(A,B)\right)$, we denote by $(a, 
\widetilde{b})$ its class in $\widehat{H}^{n}(A,B)$. The relative
and truncated relative cohomology groups are $R$-modules in a 
natural way, and there are various natural maps relating them. 
The most important ones are the \emph{class map} 
\begin{alignat*}{2}
\cl&:\widehat{H}^{n}(A,B)\longrightarrow  H^{n}(A,B),&&\qquad 
\cl(a,\widetilde{b})=[a,b], \\
\intertext{and the \emph{cycle map}}
\omega&:\widehat{H}^{n}(A,B)\longrightarrow{\rm Z}^{n}(A),
&&\qquad\omega(a,\widetilde{b})=a. \\
\intertext{We also recall the maps $\amap$ and $\bmap$, namely}
\amap&:\widetilde{A}^{n-1}\longrightarrow\widehat{H}^{n}(A,B), 
&&\qquad\amap(\widetilde{a})=(-\dd_{A}a,-\widetilde{f(a)}), \\
\intertext{and}
\bmap&:H^{n-1}(B)\longrightarrow\widehat{H}^{n}(A,B),&&\qquad 
\bmap([b])=(0,-\widetilde{b}).  
\end{alignat*}
Denoting the induced morphism $\amap|_{H^{n-1}(A)}:H^{n-1}(A) 
\longrightarrow\widehat{H}^{n}(A,B)$ also by $\amap$, we observe 
for $[a]\in H^{n-1}(A)$ the relation   
\begin{equation} 
\label{rem:abf}
\amap([a])=\bmap(f(a)).
\end{equation}
Finally, we note the following commutative diagram:
$$
\xymatrix
{{\rm Z}^{n}\left(s(A,B)\right)\ar[dr]\ar[r]&\widehat{H}^{n}(A,B) 
\ar[d]^{\cl} \\
&H^{n}( A,B).}
$$

\begin{remark}
The signs of the maps $\amap$ and $\bmap$ are the reverse of the 
signs in \cite{Burgos:CDB}. This is due to the fact that the map 
$\bmap$ represents the connecting morphism in relative cohomology 
and we have changed the convention on the sign of the connecting 
morphism in order to be compatible with \cite{SGA4_1/2}, p. 269. 
The signs in the definition of the map $\amap$ change accordingly.
\end{remark}

For more details about truncated relative cohomology groups and 
the proof of the following proposition we refer to \cite{Burgos:CDB}, 
p. 352.

\begin{proposition}
\label{prop:exacttrunc}
The following sequences are exact:
\begin{align}
&H^{n-1}(A,B)\longrightarrow\widetilde{A}^{n-1}\overset{\amap}
{\longrightarrow}\widehat{H}^{n}(A,B)\overset{\cl}{\longrightarrow} 
H^{n}(A,B)\longrightarrow 0,  
\label{eq:seq1} \\[3mm]
&0\longrightarrow H^{n-1}(B)\overset{\bmap}{\longrightarrow}  
\widehat{H}^{n}(A,B)\overset{\omega}{\longrightarrow}{\rm Z}^{n}(A)
\longrightarrow H^{n}(B),
\label{eq:seq2} \\[3mm]
&H^{n-1}(A,B)\longrightarrow H^{n-1}(A)\overset{\amap}
{\longrightarrow}\widehat{H}^{n}(A,B)\overset{\cl\oplus\omega}
{\longrightarrow}\notag \\
&\phantom{H^{n-1}(A,B)\longrightarrow H^{n-1}(A)}H^{n}(A,B)\oplus  
{\rm Z}^{n}(A)\longrightarrow H^{n}(A)\longrightarrow 0.
\end{align}
\hfill $\square$
\end{proposition}

\nnpar{Functorial properties.} 
Let $f:A\longrightarrow B$ and $f':A'\longrightarrow B'$ be morphisms
of complexes and let $g:f\longrightarrow f'$ be a morphism of the
$2$-iterated complexes $f$, $f'$, i.e., $g$ is a pair of morphisms
$(g_{A},g_{B})$ such that the diagram
\begin{displaymath}
\begin{CD}
A@>f>>B \\
@V g_{A}VV@VV g_{B}V \\
A'@>f'>>B'
\end{CD}
\end{displaymath}
commutes.

\begin{definition}
\label{def:functru}
The \emph{morphism induced (in truncated relative cohomology) by $g$} 
is the morphism   
\begin{displaymath}
\widehat{g}:\widehat{H}^{n}(A,B)\longrightarrow\widehat{H}^{n}(A',B')
\end{displaymath}
defined by $\widehat{g}(a,\widetilde{b})=(g_{A}(a),\widetilde
{g_{B}(b)})$. 
\end{definition}

\noindent
The exact sequence \eqref{eq:seq2} implies: 

\begin{corollary}
\label{cor:tci} 
Let $g$ be a morphism of $2$-iterated complexes as in definition
\ref{def:functru} such that $g_{A}$ is an isomorphism and $g_{B}$ 
is a quasi-isomorphism. Then, the morphism $\widehat{g}$ induced 
by $g$ is an isomorphism.  
\hfill$\square$ 
\end{corollary}

This fundamental property of truncated relative
cohomology groups will be of great use later. 
We also note that the requirement of $g_{A}$ being an isomorphism is 
essential. If we replace $A$ by a quasi-isomorphic complex, the 
truncated relative cohomology groups are no longer isomorphic.

\subsection{Products in truncated relative cohomology}
\label{sec:ptrc}

\nnpar{The $*$-product.} 
We now specialize the discussion about products in relative
cohomology groups to the case of truncated relative cohomology.
Let $f_{1}$ and $f_{2}$ be as in section \ref{sec:prc}. Moreover, 
suppose that we have the morphism of commutative diagrams  
\eqref{eq:bullet}. Then, there is an induced morphism 
\begin{displaymath}
\begin{CD}
{\begin{CD}
\sigma_{n}A_{1}\otimes B_{2}@>>>B_{1}\otimes B_{2} \\
@AAA@AAA \\
\sigma_{n}A_{1}\otimes\sigma_{m}A_{2}@>>>B_{1}\otimes\sigma_{m}A_{2}
\end{CD}}
@>\bullet >>
{\begin{CD}
E_{0,1}@>>>E_{1,1} \\
@AAA@AAA \\
\sigma_{n+m}E_{0,0}@>>>E_{1,0}.
\end{CD}}
\end{CD}
\end{displaymath}

\begin{definition}
\label{def:starprod}
The $*$-\emph{product in truncated relative cohomology (induced
by $\bullet$)} is the pairing  
$$
\begin{CD}
\widehat{H}^{n}(A_{1},B_{1})\otimes\widehat{H}^{m}(A_{2},B_{2})
@>*>>\widehat{H}^{n+m}(E_{0,0},s(-j))
\end{CD}
$$
provided by theorem \ref{prop:prc}. In particular, for $(a_{1},
\widetilde{b}_{1})\in\widehat{H}^{n}(A_{1},B_{1})$ and $(a_{2},
\widetilde{b}_{2})\in\widehat{H}^{m}(A_{2},B_{2})$, it is given 
by
\begin{equation}
\label{eq:stf}
(a_{1},\widetilde{b}_{1})*(a_{2},\widetilde{b}_{2})=(a_{1}
\bullet a_{2},\left((b_{1}\bullet a_{2},(-1)^{n}a_{1}\bullet b_{2}), 
(-1)^{n-1}b_{1}\bullet b_{2}\right)^{\sim}).  
\end{equation}

Notice that, since the $*$-product is induced by a pairing 
of complexes, the right-hand side of the above equation is
independent of the choice of representatives $b_{1}$ and $b_{2}$ 
of $\widetilde{b}_{1}$ and $\widetilde{b}_{2}$. 
\end{definition}

\nnpar{Properties of the $*$-product.} We now summarize the basic
properties of the $\ast$-product.
\begin{proposition}
\label{prop:comp}
With the above notations, there are commutative diagrams
$$
\begin{CD}
\widehat{H}^{n}(A_{1},B_{1})\otimes\widehat{H}^{m}(A_{2},B_{2})
@>*>>\widehat{H}^{n+m}(E_{0,0},s(-j)) \\
@V\cl\otimes\cl VV@VV\cl V \\
H^{n}(A_{1},B_{1})\otimes H^{m}(A_{2},B_{2})@>\bullet>>
H^{n+m}(E_{0,0},s(-j)),
\end{CD}
$$
and
$$
\begin{CD}
\widehat{H}^{n}(A_{1},B_{1})\otimes\widehat{H}^{m}(A_{2},B_{2})
@>*>>\widehat{H}^{n+m}(E_{0,0},s(-j)) \\
@V\omega\otimes\omega VV@VV\omega V \\
{\rm Z}^{n}(A_{1})\otimes{\rm Z}^{m}(A_{2})@>\bullet>>
{\rm Z}^{n+m}(E_{0,0}).
\end{CD}
$$
\hfill $\square$
\end{proposition} 

\begin{proposition}
\label{prop:restrlem}
With the above notations we have the following statements:
\begin{enumerate}
\item[(i)] 
$\bmap(H(B_{1}))*\amap(\widetilde{A}_{2})=0$.
\item[(ii)] 
If $[b_{1}]\in H^{n-1}(B_{1})$, and $g=(a_{2},\widetilde{b}_{2})
\in\widehat{H}^{m}(A_{2},B_{2})$, then, in the group $\widehat{H}^
{n+m}(E_{0,0},s(-j))$, the following  equality holds:
\begin{displaymath}
\bmap([b_{1}])*g=\bmap([b_{1}]\bullet\cl(g))
\end{displaymath}
with $[b_{1}]\bullet\cl(g)\in H^{n+m-1}(s(-j))$ given by the 
pairing \eqref{eq:pairingb}.
\item[(iii)]
If $\widetilde{a}_{1}\in\widetilde{A}_{1}^{n-1}$, and $g=(a_{2},
\widetilde{b}_{2})\in\widehat{H}^{m}(A_{2},B_{2})$, then, 
in the group $\widehat{H}^{n+m}(E_{0,0},s(-j))$, we have the equality
\begin{displaymath}
\amap(\widetilde{a}_{1})*g=\amap(\widetilde{a_{1}\bullet\omega(g)})
\end{displaymath}
with $\widetilde{a_{1}\bullet\omega(g)}\in\widetilde{E}_{0,0}^
{n+m-1}$ given by the pairing $A_{1}\otimes A_{2}\overset{\bullet}
{\longrightarrow}E_{0,0}$.
\item[(iv)]
If $[a_{1}]\in H^{n-1}(A_{1})$, and $g=(a_{2},\widetilde{b}_{2})\in 
\widehat{H}^{m}(A_{2},B_{2})$, then, again in the group $\widehat{H}^{n+m}
(E_{0,0},s(-j))$, we have the equality
\begin{displaymath}
\amap([a_{1}])*g=\amap([a_{1}]\bullet[a_{2}]);
\end{displaymath}
here $[a_{2}]$ is the class of $a_{2}$ in $H^{m}(A_{2})$.
\end{enumerate}
\end{proposition}
\begin{proof}
(i) For $[b_{1}]\in H^{n-1}(B_{1})$ and $\widetilde{a}_{2}\in 
\widetilde{A}_{2}$, we have
\begin{align*}
\bmap([b_{1}])*\amap(\widetilde{a}_{2})&=(0,-\widetilde{b}_{1})*
(-\dd_{A_{2}}a_{2},-\widetilde{f_{2}(a_{2})}) \\
&=\left(0,\left((b_{1}\bullet\dd_{A_{2}}a_{2},0),(-1)^{n-1}b_{1}
\bullet f_{2}(a_{2})\right)^{\sim}\right).
\end{align*}
Since $b_{1}$ is a cycle, we have  
\begin{displaymath}
(-1)^{n-1}\dd((b_{1}\bullet a_{2},0),0)=\left((b_{1}\bullet 
\dd_{A_{2}}a_{2},0),(-1)^{n-1}b_{1}\bullet f_{2}(a_{2})\right),
\end{displaymath}
which immediately proves the first statement. \\
(ii) The formula we seek follows by a straightforward 
computation. \\
(iv) For a cycle $a_{1}\in A_{1}^{n-1}$, we have
\begin{align*}
\amap(\widetilde{a}_{1})*(a_{2},\widetilde{b}_{2})&=(0,-\widetilde
{f_{1}(a_{1})})*(a_{2},\widetilde{b}_{2}) \\
&=\left(0,-\left((f_{1}(a_{1})\bullet a_{2},0),(-1)^{n-1}f_{1}(a_{1})
\bullet b_{2}\right)^{\sim}\right).
\end{align*}
Now we note
\begin{displaymath}
\dd((0,a_{1}\bullet b_{2}),0)=\left((0,\dd_{A_{1}}a_{1}\bullet
b_{2}+(-1)^{n-1}a_{1}\bullet\dd_{B_{2}}b_{2}),-f_{1}(a_{1})\bullet
b_{2}\right),
\end{displaymath}
hence
\begin{displaymath}
(-1)^{n}\dd((0,a_{1}\bullet b_{2}),0)=\left((0,-a_{1}\bullet
f_{2}(a_{2})),(-1)^{n-1}f_{1}(a_{1})\bullet b_{2}\right).
\end{displaymath}
This immediately leads to
\begin{displaymath}
\amap(\widetilde{a}_{1})*(a_{2},\widetilde{b}_{2})=\left(0,
-((f_{1}(a_{1})\bullet a_{2},a_{1}\bullet f_{2}(a_{2})),0)^
{\sim}\right)=\amap(\widetilde{a_{1}\bullet a_{2}}),
\end{displaymath}
which proves the claim. \\
(iii) Since $\omega(g)=a_{2}$, the equality follows
immediately from the proof given for part (iv).
\end{proof}

\nnpar{Example.} 
The following gives an example of the formalism developed for 
truncated cohomology.  

Let $A$ be an associative, graded commutative, differential algebra 
over $R$, and $B$ a graded commutative, differential algebra over $R$, 
which is associative up to homotopy, i.e., there exist morphisms
\begin{displaymath}
h:B^{n}\otimes B^{m}\otimes B^{l}\longrightarrow B^{n+m+l-1}
\end{displaymath}
satisfying
\begin{displaymath}
(b_{1}b_{2})b_{3}-b_{1}(b_{2}b_{3})=\dd h(b_{1}\otimes b_{2}
\otimes b_{3})+h\dd(b_{1}\otimes b_{2}\otimes b_{3}).
\end{displaymath}
Furthermore, let $f:A\longrightarrow B$ be a morphism of graded 
differential algebras satisfying   
\begin{displaymath}
h(f(a_{1})\otimes f(a_{2})\otimes f(a_{3}))=0
\end{displaymath}
for all elements $a_{1},a_{2},a_{3}$ in $A$.

\begin{proposition}
\label{pro:prodtrun}
Let $A$, $B$, $f:A\longrightarrow B$ be as above. Then, we have
the following statements:
\begin{enumerate}
\item[(i)] 
The truncated cohomology groups $\widehat{H}^{\ast}(A,B)$ have a 
natural structure of an associative, graded commutative algebra.
\item[(ii)] 
In the particular case assumed we have $\bmap(H^{\ast}(B))*
\widehat{H}^{\ast}(A,B)=0$.
\item[(iii)] 
The class map $\cl:\widehat{H}^{\ast}(A,B)\longrightarrow
H^{\ast}(A,B)$ is a morphism of algebras.
\item[(iv)] 
The cycle map $\omega:\widehat{H}^{\ast}(A,B)\longrightarrow 
{\rm Z}^{\ast}(A)$ is a morphism of algebras.
\end{enumerate}
\end{proposition}
\begin{proof}
The key point of the proof is to show (i). Then (iii), (iv) will
follow from proposition \ref{prop:comp1}; part (ii) is a direct
consequence of the explicit description of the $*$-product under
the above assumptions and will be shown at the end of the proof.

By assumption, we have a morphism of commutative diagrams
\begin{displaymath}
\begin{CD}
{\begin{CD}
A\otimes B@>>>B\otimes B \\
@AAA@AAA \\
A\otimes A@>>>B\otimes A
\end{CD}}
@>\bullet>>
{\begin{CD}
B@>>>B \\
@AAA@AAA \\
A@>>>B,
\end{CD}}
\end{CD}
\end{displaymath}
where $\bullet$ is given by the underlying algebra structure.
Therefore, if $j:B\oplus B\longrightarrow B$ is the morphism 
given by $j(x,y)=-x+y$, we obtain a pairing
\begin{displaymath}
\widehat{H}^{n}(A,B)\otimes\widehat{H}^{m}(A,B)\overset{*}
{\longrightarrow}\widehat{H}^{n+m}(A,s(-j)).
\end{displaymath}
Now consider the short exact sequence   
\begin{equation}
\label{eq:tses}
\begin{CD}
0@>>>B@>i>>B\oplus B@>j>>B@>>>0,
\end{CD}
\end{equation}
where the morphism $i:B\longrightarrow B\oplus B$ is given 
by $i(x)=(x,x)$. By the kernel-simple quasi-isomorphism (see
definition \ref{def:MVP}), $i$ induces a quasi-isomorphism
between $B$ and $s(-j)$. By means of corollary \ref{cor:tci}
we then obtain a natural isomorphism $\widehat{H}^{n+m}(A,s(-j))
\cong\widehat{H}^{n+m}(A,B)$. Therefore, the $*$-product 
defines an inner product on $\widehat{H}^{\ast}(A,B)$.

We now give an explicit formula for the $*$-product. Since 
the exact sequence $\ref{eq:tses}$ is split, we easily find 
an explicit quasi-inverse for the quasi-isomorphism induced 
by $i$. Namely, by applying proposition \ref{prop:quisinv}, 
we can choose a quasi-inverse from $s(-j)$ to $B$ by the 
morphism given by
\begin{displaymath}
((x,y),z)\longmapsto x.
\end{displaymath}
With this choice we derive from (\ref{eq:stf}) the formula
\begin{equation}
\label{eq:exa-star}
(a_{1},\widetilde{b}_{1})*(a_{2},\widetilde{b}_{2})=
(a_{1}a_{2},\widetilde{b_{1}f(a_{2})}).
\end{equation}
By corollary \ref{cor:tci} this product does not depend on the 
choice of the quasi-inverse. 

We now show that the $*$-product is graded commutative. To do so,
let $(a_{1},\widetilde{b}_{1})\in\widehat{H}^{n}(A,B)$,
$(a_{2},\widetilde{b}_{2})\in\widehat{H}^{m}(A,B)$, and let
$(a_{1},b_{1})$,$(a_{2},b_{2})$ be respective representatives. 
Then, the difference 
\begin{displaymath}
(a_{1},\widetilde{b}_{1})*(a_{2},\widetilde{b}_{2})-(-1)^{nm}
(a_{2},\widetilde{b}_{2})*(a_{1},\widetilde{b}_{1})
\end{displaymath}
is represented by
\begin{displaymath}  
(a_{1}a_{2}-(-1)^{nm}a_{2}a_{1},b_{1}f(a_{2})-(-1)^{nm}b_{2}f(a_{1}))= 
(0,b_{1}f(a_{2})-(-1)^{n}f(a_{1})b_{2}).
\end{displaymath}
Since $\dd b_{1}=f(a_{1})$ and $\dd b_{2}=f(a_{2})$, we obtain
\begin{displaymath}
(-1)^{n-1}\dd(b_{1}b_{2})=b_{1}f(a_{2})-(-1)^{n}f(a_{1})b_{2}.
\end{displaymath}
This shows
\begin{displaymath}
(a_{1},\widetilde{b}_{1})*(a_{2},\widetilde{b}_{2})=(-1)^{nm}
(a_{2},\widetilde{b}_{2})*(a_{1},\widetilde{b}_{1}). 
\end{displaymath}

Finally, we show that the $*$-product is associative. For this
purpose, let $(a_{3},\widetilde{b}_{3})\in\widehat{H}^{l}(A,B)$ 
be a third element. Consider the difference 
\begin{displaymath}
\left((a_{1},\widetilde{b}_{1})*(a_{2},\widetilde{b}_{2})\right)*
(a_{3},\widetilde{b}_{3})-(a_{1},\widetilde{b}_{1})*\left((a_{2},
\widetilde{b}_{2})*(a_{3},\widetilde{b}_{3})\right),
\end{displaymath}
which is represented by
\begin{displaymath}
((a_{1}a_{2})a_{3},(b_{1}f(a_{2}))f(a_{3}))-
(a_{1}(a_{2}a_{3}),b_{1}(f(a_{2})f(a_{3}))).
\end{displaymath}
We have to show that the class of the latter difference 
in $\widehat{H}^{n+m+l}(A,B)$ is zero.
By the associativity of $A$, we have $(a_{1}a_{2})a_{3}=a_{1}
(a_{2}a_{3})$. Moreover, by the associativity of $B$
up to homotopy we have 
\begin{align*}
(b_{1}f(a_{2}))f(a_{3})-b_{1}(f(a_{2})f(a_{3}))&= \\
\dd h(b_{1}\otimes f(a_{2})\otimes f(a_{3}))&+h\dd(b_{1}
\otimes f(a_{2})\otimes f(a_{3})). 
\end{align*}
Since $f(a_{2})$ and $f(a_{3})$ are cycles, and since
$\dd b_{1}=f(a_{1})$, we obtain by assumption
\begin{displaymath}
h\dd(b_{1}\otimes f(a_{2})\otimes f(a_{3}))=h(f(a_{1})
\otimes f(a_{2})\otimes f(a_{3}))=0.
\end{displaymath}
From this the associativity of the $*$-product follows.
  
Statement (ii) is obvious from the explicit formula \eqref{eq:exa-star} 
given for the $*$-product. 
\end{proof}

\begin{remark} 
\label{rem:signs}
There are many discrepancies in the signs between this paper and the 
previous papers \cite{Burgos:Gftp}, \cite{Burgos:CDB}. The main 
source for these discrepancies is that the convention for the sign 
of the connecting morphism in this paper is minus the connecting 
morphism in \cite{Burgos:CDB}. In particular, compare formula 
\eqref{eq:41} with \cite{Burgos:CDB}, \S2. But there are also some 
differences in the sign of the product.
\end{remark}

\newpage
\section{Green objects}
\label{sec:green-objects}
The aim of this section is to develop an abstract theory of Green
objects. It is a generalization of \cite{Burgos:CDB} and some arguments,
like the proof of the associativity and the commutativity 
of \cite{Burgos:CDB}, carry over to this abstract setting directly.

\subsection{$\Gi$-complexes}
\label{sec:ac}

\nnpar{The Mayer-Vietoris principle.}
\begin{definition}
Let $X$ be a topological space. We say that a sheaf $\mathcal{F}$ is 
a \emph{totally acyclic sheaf}, if the restriction $\mathcal{F}|_{U}$ 
of $\mathcal{F}$ to $U$ is acyclic for all open subsets $U$ of $X$.
\end{definition}

\noindent
We note that for instance every flasque or fine sheaf on $X$ is 
totally acyclic. 

\begin{definition}
We say that a presheaf $\mathcal{F}$ on $X$ \emph{satisfies the
Mayer-Vietoris principle}, if the sequence   
\begin{displaymath}
0\longrightarrow\mathcal{F}(U\cup V)\overset{i}{\longrightarrow} 
\mathcal{F}(U)\oplus\mathcal{F}(V)\overset{j}{\longrightarrow} 
\mathcal{F}(U\cap V)\longrightarrow 0 
\end{displaymath}
with $i(\eta)=(\eta,\eta)$ and $j(\omega,\eta)=-\omega+\eta$ is 
exact for any pair of open subsets $U$ and $V$ of $X$.
\end{definition}

\begin{proposition} 
\label{prop:17}
Let $X$ be a noetherian space. Then, a presheaf $\mathcal{F}$ on 
$X$ satisfies the Mayer-Vietoris principle, if and only if it is 
a totally acyclic sheaf.
\end{proposition}
\begin{proof}
Since $X$ is noetherian, the exactness of 
\begin{displaymath}
0\longrightarrow\mathcal{F}(U\cup V)\overset{i}{\longrightarrow} 
\mathcal{F}(U)\oplus\mathcal{F}(V)\overset{j}{\longrightarrow} 
\mathcal{F}(U\cap V)
\end{displaymath}
is equivalent to the fact that $\mathcal{F}$ is a sheaf. If $\mathcal
{F}$ is totally acyclic, then the map $j$ is surjective. Hence
$\mathcal{F}$ 
satisfies the Mayer-Vietoris principle. Conversely, if $\mathcal{F}$
satisfies the Mayer-Vietoris principle, one can use induction on
the number of open sets to show that, for any finite open covering  
$\mathfrak{U}$ of $X$, the \v{C}ech cohomology group $H^{1}(\mathfrak
{U},\mathcal{F})$ vanishes. This implies that, if $\mathcal{F}'$ 
satisfies the Mayer-Vietoris principle, then for any short exact 
sequence   
\begin{displaymath}     
0\longrightarrow\mathcal{F}'\longrightarrow\mathcal{F}     
\longrightarrow\mathcal{F}''\longrightarrow 0  
\end{displaymath}   
and any open set $U$, the sequence    
\begin{displaymath}
0\longrightarrow\mathcal{F}'(U)\longrightarrow\mathcal{F}(U)
\longrightarrow \mathcal{F}''(U)\longrightarrow 0
\end{displaymath}
is exact. Then, the nine lemma implies that, if the sheaves $\mathcal
{F}'$ and $\mathcal{F}$ in an exact sequence as above satisfy the
Mayer-Vietoris principle, then also $\mathcal{F}''$ does. One finally
concludes the proof as in \cite{Hartshorne:ag}, III 2.5. 
\end{proof}

\nnpar{Definition of $\Gi$-complexes.} 
In order to define arithmetic Chow groups with Green objects in
a certain cohomology theory, the main property we have to require of 
the cohomology theory is that it receives characteristic classes 
from $K$-theory (at least from $K_{0}$ and $K_{1}$) with natural 
properties. Since the arithmetic Chow groups will depend on the 
choice of the complex used to compute the cohomology, it might be 
more important to have a particular property of this complex than 
to have all the properties of a Gillet cohomology. Thus, in order
to have characteristic classes, but also to retain flexibility in 
the choice of the complexes, we will ask only that the cohomology 
theory factors through a Gillet cohomology. Thus, we fix a field 
$k$ and an auxiliary Gillet complex $\Gi=\Gi^{\ast}(\ast)$ for 
regular schemes of finite type over $k$. Furthermore, until the 
discussion of the functoriality, we also fix a regular separated 
scheme $X$ of finite type over $k$. 

\begin{definition} 
A \emph{$\Gi$-complex over $X$} is a graded complex $\cc=(\cc^{\ast}
(\ast),\dd)$ of sheaves of abelian groups over $X$ together with a 
morphism 
\begin{displaymath}
\mathfrak{c}_{\cc}:\Gi\longrightarrow\cc
\end{displaymath}
in the derived category of graded complexes of sheaves of abelian 
groups over $X$ such that the sheaves $\cc^{n}(p)$ satisfy the 
Mayer-Vietoris principle for all $n,p$. For a $\Gi$-complex 
$(\cc,\mathfrak{c}_{\cc})$ over $X$ we will simply write $\cc$ 
as a shorthand.
\end{definition}

\begin{definition}
A \emph{morphism of $\Gi$-complexes over $X$} is a morphism $f:
\cc\longrightarrow\cc'$ of graded complexes of sheaves of abelian 
groups over $X$ such that the diagram     
\begin{displaymath}
\begin{CD}
\Gi@>\mathfrak{c}_{\cc}>>\cc \\
@V\Id VV@VV f V \\
\Gi@>\mathfrak{c}_{\cc'}>>\cc'
\end{CD}    
\end{displaymath}
commutes in the derived category of graded complexes of sheaves of 
abelian groups over $X$.
\end{definition}
 
\nnpar{Cohomology groups of a $\Gi$-complex.}
\begin{notation}
For an open subset $U$ of $X$, the sections of $\cc^{n}(p)$ over $U$ 
will be denoted by $\cc^{n}(U,p)$. If $\omega\in\cc^{n}(X,p)$, we
will write $\omega|_{U}$ for the restriction of $\omega$ to $U$;
moreover, if the open set $U$ is clear from the context, we will
simply write $\omega$ instead of $\omega|_{U}$. Furthermore, if $Y$ 
is a closed subset of $X$ and $U=X\setminus Y$, we introduce the 
following notation:
\begin{align*}
H^{\ast}_{\cc}(X,p)&=H^{\ast}({\cc}(X,p)), \\
H^{\ast}_{\cc}(U,p)&=H^{\ast}({\cc}(U,p)), \\
H^{\ast}_{\cc,Y}(X,p)&=H^{\ast}\left({\cc}(X,p),{\cc}(U,p)\right).
\end{align*}
Analogously, for any family $\varphi$ of supports on $X$, we define
the cohomology groups with support in $\varphi$ by $H^{\ast}_{\cc,
\varphi}(X,p)$.  
\end{notation}

Since the sheaves $\cc^{\ast}(p)$ are totally acyclic, the cohomology
of global sections agrees with the hypercohomology of the complex.
Hence, there are induced morphisms
\begin{align*}
H^{\ast}(X,\Gi(p))&\longrightarrow H^{\ast}_{\cc}(X,p), \\
H^{\ast}(U,\Gi(p))&\longrightarrow H^{\ast}_{\cc}(U,p), \\
H^{\ast}_{Y}(X,\Gi(p))&\longrightarrow H^{\ast}_{\cc,Y}(X,p),
\end{align*}
which, by abuse of notation, will be denoted again by $\mathfrak{c}_
{\cc}$. 

\begin{lemma}
\label{lemma:cohsup}
Let $Y$ be a closed subset of $X$, and $U=X\setminus Y$. Then, there 
is a morphism of exact sequences 
\begin{equation*}
\begin{CD}
H^{n}(X,\Gi(p))@>>>H^{n}(U,\Gi(p))@>\delta>>H^{n+1}_{Y}(X,\Gi(p))
@>>>\dots \\
@V\mathfrak{c}_{\cc}VV@V\mathfrak{c}_{\cc}VV@V\mathfrak{c}_{\cc}VV \\
H^{n}_{\cc}(X,p)@>>>H^{n}_{\cc}(U,p)@>\delta>>H^{n+1}_{\cc,Y}(X,p)
@>>>\dots 
\end{CD}
\end{equation*}
\end{lemma}
\begin{proof}
The proof follows immediately from the fact that $\mathfrak{c}_{\cc}:
\Gi\longrightarrow\cc$ is a morphism of sheaves.
\end{proof}

\nnpar{Characteristic classes in $\cc$-cohomology.} 
Using the fact that there are well-defined characteristic classes in 
$\Gi$-cohomology, we can also define classes for cycles and 
$K_{1}$-chains in $\cc$-cohomology.

\begin{definition}
\label{def:3.8}
Let $y$ be a $p$-codimensional cycle of $X$, and $Y=\supp y$. 
We define $\cl_{\cc}(y)\in H^{2p}_{\cc,Y}(X,p)$ to be the class
\begin{displaymath}
\cl_{\cc}(y)=\mathfrak{c}_{\cc}(\cl_{\Gi}(y)).
\end{displaymath}  
Let $f\in R^{p-1}_{p}(X)$ be a $K_{1}$-chain, $y=\dv(f)$, $Y=
\supp y$, and $U=X\setminus Y$. We define $\cl_{\cc}(f)\in 
H^{2p-1}_{\cc}(U,p)$ to be the class     
\begin{displaymath}
\cl_{\cc}(f)=\mathfrak{c}_{\cc}(\cl_{\Gi}(f)).
\end{displaymath}
Here $\cl_{\Gi}$ is given in definition \ref{def:1}.  
\end{definition}


If $Y=\supp y$ and $Y\subset Z$ for a closed subset $Z$ of $X$, by
abuse of notation, we will also denote by $\cl_{\cc}(y)$ the class 
in $H^{2p}_{\cc,Z}(X,p)$. If there is a danger of confusion, we will 
write explicitly $\cl_{\cc}(y)\in H^{2p}_{\cc,Z}(X,p)$ to indicate 
in which cohomology group the class is considered. We will use the 
same convention for $\cl_{\cc}(f)$. In particular, if $y=\dv(f)$, 
$Y=\supp y$, $U=X\setminus Y$, and $W=\supp f$, we will consider 
the classes $\cl_{\cc}(f)\in H^{2p-1}_{\cc}(U,p)$, and $\cl_{\cc}(f)
\in H^{2p-1}_{\cc,W\setminus Y}(U,p)$.

\medskip

Observe that the cohomology $\Gi$ is only an auxiliary device 
to ensure the existence of classes for cycles and $K_{1}$-chains 
having good properties. The properties which will be needed to 
define arithmetic Chow groups are summarized in the following 
lemma. This lemma can be taken as an alternative starting point 
for the definition of arithmetic Chow groups. It is an immediate 
consequence of the  properties of characteristic classes in 
$\Gi$-cohomology.

\begin{lemma}
\label{lemma:asc}
Let $\cc$ be a $\Gi$-complex over $X$. Then, we have the following
statements:
\begin{enumerate}
\item[(i)]
\label{lemma:asc1} 
For any family $\varphi$ of supports on $X$, the map $\cl_{\cc}$ 
induces a morphism of groups   
\begin{displaymath}
{\rm Z}^{p}_{\varphi}(X)\longrightarrow H^{2p}_{\cc,\varphi}(X,p)
\end{displaymath}
which is compatible with change of support; here ${\rm Z}^{p}_
{\varphi}(X)=R^{p}_{p,\varphi}(X)$ is the group of $p$-codimensional 
cycles on $X$ with support on $\varphi$.
\item[(ii)]
\label{lemma:asc2} 
For a $K_{1}$-chain $f\in R^{p-1}_{p}(X)$, put $y=\dv(f)$, 
$Y=\supp y$, and $U=X\setminus Y$. Then, the equality
\begin{displaymath}
\cl_{\cc}(\dv(f))=\delta\cl_{\cc}(f)\in H^{2p}_{\cc,Y}(X,p)
\end{displaymath}
holds, where $\delta$ is the connection morphism
\begin{displaymath}
\delta:H^{2p-1}_{\cc}(U,p)\longrightarrow H^{2p}_{\cc,Y}(X,p).
\end{displaymath}
Moreover, if $g\in R^{p-1}_{p}(X)$ is another $K_{1}$-chain, 
$z=\dv(g)$, and $Z=\supp z$, we have
\begin{displaymath}
\cl_{\cc}(f)+\cl_{\cc}(g)=\cl_{\cc}(f+g)\in H^{2p-1}_{\cc}
(X\setminus(Y\cup Z),p).
\end{displaymath}
\item[(iii)] 
\label{lemma:asc3} 
If $h\in R^{p-2}_{p}(X)$ is a $K_{2}$-chain, then the class  
$\cl_{\cc}(\dd h)$ vanishes in the group $H^{2p-1}_{\cc}(X,p)$.
\end{enumerate}
\hfill $\square$
\end{lemma}

\noindent
Lemma \ref{lemma:asc} has the following direct consequences.

\begin{corollary}
For any family $\varphi$ of supports on $X$, there are well-defined 
morphisms
\begin{align*}
\cl_{\cc}&:\CH^{p}_{\varphi}(X)\longrightarrow 
H^{2p}_{\cc,\varphi}(X,p), \\
\cl_{\cc}&:\CH^{p-1,p}_{\varphi}(X)\longrightarrow
H^{2p-1}_{\cc,\varphi}(X,p). 
\end{align*}
\hfill $\square$
\end{corollary}

Following \cite{BlochOgus:Gchs}, let us denote by $\mathcal{Z}^{p}=
\mathcal{Z}^{p}(X)$ the set of all closed subsets of $X$ of 
codimension greater or equal to $p$ ordered by inclusion. We then 
write
\begin{align*}
H^{n}_{\cc}(X\setminus\mathcal{Z}^{p},p)&=\lim_{\substack{\longrightarrow\\  
Y\in\mathcal{Z}^{p}}}H^{n}_{\cc}(X\setminus Y,p), \\
H^{n}_{\cc,\mathcal{Z}^{p}}(X,p)&=\lim_{\substack{\longrightarrow\\ 
Y\in\mathcal{Z}^{p}}}H^{n}_{\cc,Y}(X,p).
\end{align*}

\begin{corollary}
The following diagram is commutative and has exact rows:
$$
\xymatrix{\CH^{p-1,p}(X)\ar[r]\ar[d]_{\cl_{\cc}}&R^{p-1}_{p}(X)/
\dd R^{p-2}_{p}(X)\ar[r]^-{\dv}\ar[d]_{\cl_{\cc}}&R^p_{p}(X)   
\ar[r]\ar[d]_{\cl_{\cc}}&\CH^{p}(X)\ar[d]_{\cl_{\cc}} \\ 
H^{2p-1}_{\cc}(X,p)\ar[r]&H^{2p-1}_{{\cc}}(X\setminus
\mathcal{Z}^{p},p)\ar[r]^-{\delta}&H^{2p}_{{\cc},\mathcal{Z}^p}
(X,p)\ar[r]&H^{2p}_{\cc}(X,p).} 
$$  
\hfill $\square$
\end{corollary}

\nnpar{Purity.} 
In many cases the cohomology theory associated with a $\Gi$-complex
$\cc$ will satisfy a purity axiom, e.g., if $\cc$-cohomology and
$\Gi$-cohomology agree. Later on we will see that the arithmetic 
Chow groups obtained in this case have better properties.

\begin{definition}
\label{def.311}
We say that the $\Gi$-complex $\cc$ satisfies the \emph{weak 
purity condition}, if for any closed subset $Y$ of $X$ of 
codimension greater or equal to $p$, we have
\begin{displaymath}
H^{2p-1}_{\cc,Y}(X,p)=0.
\end{displaymath}
\end{definition}

Recall that for any complex $A$, we let $\widetilde A=A/\Img\dd$.
In order to be able to deal with complexes which may not satisfy 
the weak purity condition, we introduce the following notation
\begin{align*}
\widetilde{\cc}^{2p-1}(X,p)^{\pure}&=\widetilde{\cc}^{2p-1}(X,p)\big/ 
\Img\left(H^{2p-1}_{\cc,\mathcal{Z}^{p}}(X,p)\longrightarrow\widetilde
{\cc}^{2p-1}(X,p)\right), \\ 
H^{2p-1}_{\cc}(X,p)^{\pure}&=H^{2p-1}_{\cc}(X,p)\big/
\Img\left(H^{2p-1}_{\cc,\mathcal{Z}^{p}}(X,p)\longrightarrow
H^{2p-1}_{\cc}(X,p)\right). 
\end{align*}  

Note that, if $\cc$ satisfies the weak purity condition, then
\begin{align*}
\widetilde{\cc}^{2p-1}(X,p)^{\pure}&=\widetilde{\cc}^{2p-1}(X,p), \\
H^{2p-1}_{\cc}(X,p)^{\pure}&=H^{2p-1}_{\cc}(X,p).
\end{align*}

\subsection{Definition of Green objects}
\label{sec:go}

\nnpar{Preliminaries.} 
Given a $\Gi$-complex $\cc=(\cc^{\ast}(\ast),\dd)$ over $X$, 
we will define Green objects with values in this complex. For 
concrete examples we refer to \cite{Burgos:CDB}, or chapters
6 and 7 below. 

\begin{notation}
Let $Y$ be a closed subset of $X$, and $U=X\setminus Y$. We will
then write
\begin{displaymath}
\widehat{H}^{n}_{\cc,Y}(X,p)=\widehat{H}^{n}(\cc(X,p),\cc(U,p)), 
\end{displaymath}
where $\widehat{H}^{n}(\cc(X,p),\cc(U,p))$ are the truncated
relative cohomology groups defined in section \ref{sec:tcdb}. 
We will also write 
\begin{displaymath}
\widehat{H}^{n}_{\cc,\mathcal{Z}^{p}}(X,p)=\lim_{\substack{\longrightarrow 
\\ Y\in\mathcal{Z}^{p}}}\widehat{H}^{n}_{\cc,Y}(X,p).
\end{displaymath}
By writing
\begin{displaymath}
\cc(X\setminus\mathcal{Z}^{p},p)=\lim_{\substack{\longrightarrow 
\\ Y\in\mathcal{Z}^{p}}}\cc(X\setminus Y,p),
\end{displaymath}
we have
\begin{displaymath}
\widehat{H}^{n}_{\cc,\mathcal{Z}^{p}}(X,p)=\widehat{H}^{n}
(\cc(X,p),\cc(X\setminus\mathcal{Z}^{p},p)).
\end{displaymath}
\end{notation}

From the definition of the truncated relative cohomology groups 
we recall the morphisms
\begin{align*}
\cl:\widehat{H}^{n}_{\cc,\mathcal{Z}^{p}}(X,p)&\longrightarrow 
H^{n}_{\cc,\mathcal{Z}^{p}}(X,p), \\
\omega:\widehat{H}^{n}_{\cc,\mathcal{Z}^{p}}(X,p)&\longrightarrow
{\rm Z}(\cc^{n}(X,p)), \\
\amap:\widetilde{\cc}^{n-1}(X,p)&\longrightarrow\widehat{H}^{n}_
{\cc,\mathcal{Z}^{p}}(X,p), \\
\amap:H^{n-1}_{\cc}(X,p)&\longrightarrow\widehat{H}^{n}_
{\cc,\mathcal{Z}^{p}}(X,p), \\ 
\bmap:H^{n-1}_{\cc}(X\setminus\mathcal{Z}^{p},p)&\longrightarrow 
\widehat{H}^{n}_{\cc,\mathcal{Z}^{p}}(X,p);
\end{align*}
analogous morphisms exist for $\widehat{H}^{n}_{\cc,Y}(X,p)$ with
$Y$ a closed subset of $X$. We leave it to the reader to 
write down the exact sequences of proposition \ref{prop:exacttrunc} 
in these cases. In particular, for $Y\in\mathcal{Z}^{p}$ we note 
the commutative diagram with exact rows    
\begin{equation}
\label{eq:18}
\xymatrix{H^{n-1}_{\cc,Y}(X,p)\ar[r]\ar[d]&\widetilde{\cc}^{n-1}(X,p)
\ar[r]^{\amap}\ar@{=}[d]&\widehat{H}^{n}_{\cc,Y}(X,p)\ar[r]^{\cl}
\ar[d]&H^{n}_{\cc,Y}(X,p)\ar[r]\ar[d]&0 \\
H^{n-1}_{\cc,\mathcal{Z}^{p}}(X,p)\ar[r]&\widetilde{\cc}^{n-1}(X,p)
\ar[r]^{\amap}&\widehat{H}^{n}_{\cc,\mathcal{Z}^{p}}(X,p)\ar[r]^
{\cl}&H^{n}_{\cc,\mathcal{Z}^{p}}(X,p)\ar[r]&0.}
\end{equation}

While dealing with $p$-codimensional cycles, we will be mainly 
interested in the groups with $n=2p$. In this case we obtain as
a direct consequence of proposition \ref{prop:exacttrunc} 

\begin{proposition}
\label{prop:exactgreen}
There are exact sequences
\begin{align}
&0\longrightarrow\widetilde{\cc}^{2p-1}(X,p)^{\rm pure}\overset
{\amap}{\longrightarrow}\widehat{H}^{2p}_{\cc,\mathcal{Z}^{p}}(X,p)
\overset{\cl}{\longrightarrow}H^{2p}_{\cc,\mathcal{Z}^{p}}(X,p)  
\longrightarrow 0,  
\label{exseq:go4} \\[5mm]
&0\longrightarrow H^{2p-1}_{\cc}(X,p)^{\rm pure}\overset{\amap}
{\longrightarrow}\widehat{H}^{2p}_{\cc,\mathcal{Z}^{p}}(X,p)
\overset{\cl\oplus\omega}{\longrightarrow}\notag \\
&\phantom{H^{2p-1}_{\cc,\mathcal{Z}^{p}}(X,p)\longrightarrow}
H^{2p}_{\cc,\mathcal{Z}^{p}}(X,p)\oplus{\rm Z}\cc^{2p}(X,p) 
\longrightarrow H^{2p}_{\cc}(X,p)\longrightarrow 0.
\label{exseq:go5}
\end{align}
\hfill $\square$    
\end{proposition}

\nnpar{Green objects.}
\begin{definition} Let $y$ be a $p$-codimensional cycle on $X$. 
A \emph{Green object for the class of $y$ (with values in $\cc$)} 
is an element $\mathfrak{g}_{y}\in\widehat{H}^{2p}_{\cc,\mathcal{Z}^
{p}}(X,p)$ such that
\begin{displaymath}
\cl(\mathfrak{g}_{y})=\cl_{\cc}(y)\in H^{2p}_{\cc,\mathcal{Z}^{p}}
(X,p).
\end{displaymath}
In other words, a Green object for the class of $y$ is a pair
$\mathfrak{g}_{y}=(\omega_{y},\widetilde{g}_{y})$ with $\omega_{y} 
\in{\rm Z}\cc^{2p}(X,p)$ and $\widetilde{g}_{y}\in\widetilde{\cc}^
{2p-1}(X\setminus\mathcal{Z}^{p},p)$ such that $\omega_{y}=\dd
g_{y}$ and such that this pair represents the class of $y$ in
$H^{2p}_{\cc,\mathcal{Z}^{p}}(X,p)$; note that $g_{y}$ here
denotes a representative of $\widetilde{g}_{y}$ in $\cc^{2p-1}
(X\setminus\mathcal{Z}^{p},p)$.
\end{definition}

Observe that the Green objects for the class of any $p$-codimensional
cycle belong to the same space $\widehat{H}^{2p}_{\cc,\mathcal{Z}^{p}}
(X,p)$. This is the reason for taking the limit over all 
$p$-codimensional cycles of $X$. Nevertheless, sometimes it will be 
necessary to have a preciser control on the group in which a Green 
object is defined.

\begin{definition} 
Let $y$ be a $p$-codimensional cycle on $X$ with $Y=\supp y$. 
A \emph{Green object for the cycle $y$} is an element $\mathfrak
{g}_{y}\in\widehat{H}^{2p}_{\cc,Y}(X,p)$ such that 
\begin{displaymath}
\cl(\mathfrak{g}_{y})=\cl_{\cc}(y)\in H^{2p}_{\cc,Y}(X,p).
\end{displaymath}
\end{definition}
  
\begin{definition} 
Let $Z$ be a closed subset of $X$, and $y\in\CH^{p}_{Z}(X)$.
A \emph{weak Green object for $y$ with support in $Z$} is an 
element $\mathfrak{g}_{y}\in\widehat{H}^{2p}_{\cc,Z}(X,p)$ 
such that
\begin{displaymath} 
\cl(\mathfrak{g}_{y})=\cl_{\cc}(y)\in H^{2p}_{\cc,Z}(X,p).
\end{displaymath} 
\end{definition}

\begin{remark} 
A Green object for the cycle $y$ represents a cohomology class with
support exactly equal to the support of $y$. A Green object for the
class of $y$ represents a cohomology class whose support has the
same codimension as $y$, but may be bigger than the support of $y$.
From proposition \ref{prop:exgo} (iii) below, it will be clear that,
if the $\Gi$-complex $\cc$ satisfies the weak purity condition,
there is no difference between a Green object for the class of $y$
and a Green object for the cycle $y$. As we will see, even in the
case of a complex which does not satisfy the weak purity condition,
the distinction between these two kinds of Green objects is a minor
technical point. In contrast to this, a weak Green object for $y$
represents a cohomology class whose support may be bigger than the
support of $y$, and may even have codimension smaller than the
codimension of $y$. In
general, a weak Green object $\mathfrak{g}_{y}$ carries less
information than a Green object for $y$ and is not a Green object
for the class of $y$. Weak Green objects appear naturally when
handling non proper intersections and are therefore useful in
intermediate steps.
\end{remark}

\nnpar{Additivity.}
\begin{remark}
If $\mathfrak{g}_{y}$ is a Green object for $y$, and $\mathfrak
{g}_{z}$ a Green object for $z$ such that $\supp(y+z)=\supp(y)
\cup\supp(z)$, lemma \ref{lemma:asc} implies that $\mathfrak
{g}_{y}+\mathfrak{g}_{z}$ is a Green object for $y+z$. On the 
other hand, if $\supp(y+z)\subsetneq\supp(y)\cup\supp(z)$, then  
$\mathfrak{g}_{y}+\mathfrak{g}_{z}$ determines only a Green object 
for the class of $y+z$.    
\end{remark}

\nnpar{Existence of Green objects.}
\begin{proposition} 
\label{prop:exgo}
Let $\cc$ be a $\Gi$-complex on $X$, and let $y$ be any 
$p$-codimen\-sional cycle on $X$ with $Y=\supp y$. Then,
we have the following statements:
\begin{enumerate}
\item[(i)]
There exists a Green object $\mathfrak{g}_{y}\in\widehat{H}^{2p}_
{\cc,\mathcal{Z}^{p}}(X,p)$ for the class of $y$.
\item[(ii)]
Any Green object $\mathfrak{g}_{y}$ for the class of $y$ can be
lifted to a Green object $\mathfrak{g}'_{y}\in\widehat{H}^{2p}_
{\cc,Y}(X,p)$ for $y$.
\item[(iii)]
\label{prop:exgo3} 
If the $\Gi$-complex $\cc$ satisfies the weak purity condition, 
the morphism
\begin{displaymath} 
\widehat{H}^{2p}_{\cc,Y}(X,p)\longrightarrow\widehat{H}^{2p}_
{\cc,\mathcal{Z}^{p}}(X,p)
\end{displaymath}  
is injective. Therefore, the Green object $\mathfrak{g}'_{y}$ 
for $y$ is uniquely determined by the Green object $\mathfrak{g}_
{y}$ for the class of $y$ in this case.
\end{enumerate}
\end{proposition}
\begin{proof}
(i) The first statement is an immediate consequence of the 
surjectivity of the class map $\cl$ in the exact sequence 
\eqref{exseq:go4}. \\
(ii) Since the class of the cycle $y$ lies in $H^{2p}_{\cc,Y}(X,p)$,
the second statement follows from a standard diagram chase in the
commutative diagram \eqref{eq:18}. \\ 
(iii) If $\cc$ satisfies the weak purity condition, the 
commutative diagram \eqref{eq:18} gives rise to the diagram
\begin{displaymath}
\xymatrix{0\ar[r]&\widetilde{\cc}^{2p-1}(X,p)\ar[r]^{\amap}\ar@{=}[d] 
&\widehat{H}^{2p}_{\cc,Y}(X,p)\ar[r]^{\cl}\ar[d]&H^{2p}_{\cc,Y}(X,p)  
\ar[r]\ar[d]&0 \\
0\ar[r]& \widetilde{\cc}^{2p-1}(X,p)\ar[r]^{\amap}&\widehat{H}^{2p}_
{\cc,\mathcal{Z}^{p}}(X,p)\ar[r]^{\cl}&H^{2p}_{\cc,\mathcal{Z}^{p}}
(X,p)\ar[r]& 0,}
\end{displaymath}
where the vertical map on the right is injective. Therefore, by the 
five lemma, the vertical map in the middle is also injective.
\end{proof}

\nnpar{The Green object associated to a $K_{1}$-chain.}
\begin{definition} 
\label{def:8} 
For a $K_{1}$-chain $f\in R^{p-1}_{p}(X)$ with $y=\dv(f)$, and
$Y=\supp y$, we denote by $\mathfrak{g}(f)$ the distinguished 
Green object  
\begin{displaymath}
\mathfrak{g}(f)=\bmap(\cl_{\cc}(f))\in\widehat{H}^{2p}_{\cc,Y}
(X,p).
\end{displaymath}
If we need to specify the $\Gi$-complex $\cc$, we will write
$\mathfrak{g}_{\cc}(f)$ instead of $\mathfrak{g}(f)$.
\end{definition}

\begin{remark}
\label{rem:1} 
By lemma \ref{lemma:asc} (ii), $\mathfrak{g}(f)$ is a well-defined
Green object for $y=\dv f$. By lemma \ref{lemma:asc}  
(iii), we have $\mathfrak{g}(\dd h)=0$ for any $K_{2}$-chain 
$h$. This means in particular that the distinguished Green 
object $\mathfrak{g}(f)$ depends only on the class of $f$ in  
$\CH^{p-1,p}(X\setminus Y)$. In other words, for any closed
subset $Z$ of $X$ and $U=X\setminus Z$, there is a well-defined 
morphism
\begin{displaymath}
\CH^{p-1,p}(U)\longrightarrow\widehat{H}^{2p}_{\cc,Z}(X,p)
\end{displaymath}
which we also denote by $\mathfrak{g}$.
\end{remark}

\subsection{The $*$-product of Green objects}
\label{sec:sp}

In this section we will define a pairing of Green objects, the 
$*$-product, following the strategy of proposition \ref{pro:prodtrun}. 
We will define the $*$-product for any pairing of $\Gi$-complexes 
which is compatible with the product given in $\Gi$-cohomology. This 
will ensure the compatibility with the intersection product of cycles.  

\nnpar{Pairing of $\Gi$-complexes.}
\begin{definition}
\label{def:arpasc}
Let $\cc$, $\cc'$ and $\cc''$ be $\Gi$-complexes over $X$ with 
the same auxiliary cohomology $\Gi$. A \emph{$\Gi$-pairing} is 
a pairing of graded complexes of sheaves of abelian groups
\begin{displaymath}
\cc\otimes\cc'\overset{\bullet}{\longrightarrow }\cc''
\end{displaymath}
such that, in the derived category of graded complexes of sheaves 
of abelian groups over $X$, there is a commutative diagram
\begin{displaymath}
\begin{CD}
\Gi\overset{L}{\otimes}\Gi@>\mathfrak{c}_{\cc}\otimes
\mathfrak{c}_{\cc'}>>\cc\overset{L}{\otimes}\cc' \\
@VVV@VV\bullet V \\
\Gi@>\mathfrak{c}_{\cc''}>>\cc''\,,
\end{CD}
\end{displaymath}
where the vertical arrow on the left is the product in 
$\Gi$-cohomology (see \cite{Gillet:RRhK}, 1.1).
\end{definition}

\nnpar{The $*$-product of Green objects.}
\begin{definition}
\label{def:2}
Let $\cc$ be a $\Gi$-complex over $X$. Let $Y$, $Z$ be closed 
subsets of $X$, $U=X\setminus Y$, $V=X\setminus Z$, respectively, 
and let
\begin{displaymath}
j:\cc^{n}(U,p)\oplus\cc^{n}(V,p)\longrightarrow  
\cc^{n}(U\cap V,p)
\end{displaymath}
be the map given by $j(\omega,\eta)=-\omega+\eta$; we then put
\begin{displaymath} 
\cc(X;Y,Z,p)=s(-j).
\end{displaymath}
Observe that this complex is, up to the sign of the morphism, 
the \v{C}ech complex associated to the sheaf $\cc$ and the open
covering $\{U,V\}$. 

The kernel-simple quasi-isomorphism (see definition \ref{def:MVP})
associated to the short exact sequence    
\begin{equation}
\label{eq:20}
0\longrightarrow\cc^{n}(U\cup V,p)\overset{i}{\longrightarrow} 
\cc^{n}(U,p)\oplus\cc^{n}(V,p)\overset{j}{\longrightarrow } 
\cc^{n}(U\cap V,p)\longrightarrow 0 
\end{equation}
with $i(\eta)=(\eta,\eta)$, is the induced morphism 
\begin{equation}
\label{eq:2}
\iota:\cc(U\cup V,p)\longrightarrow\cc(X;Y,Z,p).
\end{equation}
\end{definition}

Throughout this section we will assume that there are given three
$\Gi$-complexes $\cc$, $\cc'$, $\cc''$ over $X$, and that $\bullet$ 
is a $\Gi$-pairing between $\cc$ and $\cc'$ to $\cc''$. Furthermore,
let $Y$, $Z$ be closed subsets of $X$, and $U=X\setminus Y$, $V=X
\setminus Z$, respectively. Since $\bullet $ is a pairing of sheaves,
and therefore compatible with restrictions, it induces for all $p$, 
$q$ a morphism of commutative diagrams
\begin{equation}
\label{eq:mcdaa}
\begin{matrix}
\cc(X,p)\otimes\cc'(V,q)\hspace{-2mm}&\!\rightarrow\!&\hspace{-2mm}
\cc(U,p)\otimes\cc'(V,q) \\
\uparrow&&\uparrow \\
\cc(X,p)\otimes\cc'(X,q)\hspace{-2mm}&\!\rightarrow\!&\hspace{-2mm}
\cc(U,p)\otimes\cc'(X,q)
\end{matrix}
\overset{\bullet}{\rightarrow}
\begin{matrix}
\cc''(V,p+q)\hspace{-2mm}&\!\rightarrow\!&\hspace{-2mm}
\cc''(U\cap V,p+q) \\
\uparrow&&\uparrow \\
\cc''(X,p+q)\hspace{-2mm}&\!\rightarrow\!&\hspace{-2mm}\cc''(U,p+q)
\end{matrix}
\end{equation}

\begin{theorem}
\label{thm:ap-com} 
There exist well-defined pairings and a commutative diagram    
\begin{displaymath}
\begin{CD}
\widehat{H}^{n}_{\cc,Y}(X,p)\otimes\widehat{H}^{m}_{\cc',Z}(X,q)
@>*>>\widehat{H}^{n+m}_{\cc'',Y\cap Z}(X,p+q) \\
@V\cl\otimes\cl VV@VV\cl V \\
H^{n}_{\cc,Y}(X,p)\otimes H^{m}_{\cc',Z}(X,q)@>\bullet>> 
H^{n+m}_{\cc'',Y\cap Z}(X,p+q).
\end{CD}
\end{displaymath}
\end{theorem}
\begin{proof}
By proposition \ref{prop:prc} and proposition \ref{prop:comp}, 
we have with $r=p+q$, well-defined pairings and a commutative 
diagram
\begin{equation}
\label{eq:ap-nat}
\begin{CD}
\widehat{H}^{n}_{\cc,Y}(X,p)\otimes\widehat{H}^{m}_{\cc',Z}(X,q)
@>*>>\widehat{H}^{n+m}(\cc''(X,r),\cc''(X;Y,Z,r)) \\
@V\cl\otimes\cl VV@VV\cl V \\   
H^{n}_{\cc,Y}(X,p)\otimes H^{m}_{\cc',Z}(X,q)@>\bullet>> 
H^{n+m}(\cc''(X,r),\cc''(X;Y,Z,r)).
\end{CD}
\end{equation}
Using the kernel-simple quasi-isomorphism \eqref{eq:2} and 
corollary \ref{cor:tci}, we obtain the isomorphisms 
\begin{align}
\label{eq:ap-is1}
H^{n+m}_{\cc'',Y\cap Z}(X,r)&\overset{\cong}{\longrightarrow}
H^{n+m}(\cc''(X,r),\cc''(X;Y,Z,r)), \\
\label{eq:ap-is2}
\widehat{H}^{n+m}_{\cc'',Y\cap Z}(X,r)&\overset{\cong}
{\longrightarrow}\widehat{H}^{n+m}(\cc''(X,r),\cc''(X;Y,Z,r)).
\end{align}
This completes the proof of the claim.
\end{proof}

\begin{remark} 
Let $(\omega,\widetilde{g})\in\widehat{H}^{n}_{\cc,Y}(X,p)$, 
and $(\omega',\widetilde{g}')\in\widehat{H}^{m}_{\cc',Z}(X,q)$. 
A representative of $(\omega,\widetilde{g})*(\omega',
\widetilde{g}')$ in the group $\widehat{H}^{n+m}(\cc''(X,r),
\cc''(X;Y,Z,r))$ is then given by formula \eqref{eq:stf}. 
Observe that, in order to obtain an explicit representative for 
this class in the group 
\begin{displaymath}
\widehat{H}^{n+m}_{\cc'',Y\cap Z}(X,p+q)=\widehat{H}^{n+m}
(\cc''(X,p+q),\cc''(X\setminus(Y\cap Z),p+q)), 
\end{displaymath}
we have to use an explicit quasi-inverse of the kernel-simple
quasi-isomor\-phism. This will be done in the examples (see, e.g.,
section \ref{sec:arithm-chow-groups}).    
\end{remark}

\nnpar{The $*$-product and the intersection product.} 
The above pairing induced by $\bullet$ in cohomology with support 
is compatible with the product defined in $\Gi$-cohomology.

\begin{theorem}
\label{thm:7}
There is a commutative diagram 
\begin{displaymath}
\begin{CD}
H^{n}_{Y}(X,\Gi(p))\otimes H^{m}_{Z}(X,\Gi(q))@>\cdot>>
H^{n+m}_{Y\cap Z}(X,\Gi(p+q)) \\
@V\mathfrak{c}_{\cc}\otimes\mathfrak{c}_{\cc'}VV
@VV\mathfrak{c}_{\cc''}V \\
H^{n}_{\cc,Y}(X,p)\otimes H^{m}_{\cc',Z}(X,q)@>\bullet>> 
H^{n+m}_{\cc'',Y\cap Z}(X,p+q).
\end{CD}
\end{displaymath}
\end{theorem}
\begin{proof}
The proof follows from standard arguments in sheaf theory 
(see \cite{Godement:Tf}, 6.2; \cite{Iversen:Cs}, II.10; 
\cite{Burgos:Gftp}, 2.5; and example \ref{ex:cup-product-de}).
\end{proof}

\begin{theorem} 
\label{thm:wstp} 
Let $y\in\CH^{p}_{Y}(X)$, $z\in\CH^{q}_{Z}(X)$ with closed
subsets $Y$, $Z$ of $X$, respectively. Let $\mathfrak{g}_{y}
\in\widehat{H}^{2p}_{\cc,Y}(X,p)$, $\mathfrak{g}_{z}\in  
\widehat{H}^{2q}_{\cc',Z}(X,q)$ be weak Green objects for 
$y$, $z$, respectively, and let $w=y\cdot z\in\CH^{p+q}_
{Y\cap Z}(X)$. Then, the $*$-product $\mathfrak{g}_{y}*
\mathfrak{g}_{z}$ is a weak Green object for $w$. Moreover, 
if $\mathfrak{g}_{y}=(\omega_{y},\widetilde{g}_{y})$, and 
$\mathfrak{g}_{z}=(\omega_{z},\widetilde{g}_{z})$, then 
$\mathfrak{g}_{y}*\mathfrak{g}_{z}$ is represented in the 
group $\widehat{H}^{2p+2q}(\cc''(X,p+q),\cc''(X;Y,Z,p+q))$ 
by the element
\begin{equation}
\label{eq:stfg}   
\left(\omega_{y}\bullet\omega_{z},\left((g_{y}\bullet\omega_{z},
\omega_{y}\bullet g_{z}),-g_{y}\bullet g_{z}\right)\widetilde
{\phantom{a}}\right).
\end{equation}  
\end{theorem}
\begin{proof}
We have to show that $\cl(\mathfrak{g}_{y}*\mathfrak{g}_{z})=
\cl_{\cc''}(w)\in H^{2p+2q}_{\cc'',Y\cap Z}(X,p+q)$. But by 
theorem \ref{thm:ap-com}, we have
\begin{displaymath}
\cl(\mathfrak{g}_{y}*\mathfrak{g}_{z})=\cl(\mathfrak{g}_{y})
\bullet\cl(\mathfrak{g}_{z})=\cl_{\cc}(y)\bullet\cl_{\cc'}(z)=
\mathfrak{c}_{\cc}(\cl_{\Gi}(y))\bullet\mathfrak{c}_{\cc'}
(\cl_{\Gi}(z)).
\end{displaymath}
Now, by theorem \ref{thm:7}, $\mathfrak{c}_{\cc}(\cl_{\Gi}(y))
\bullet\mathfrak{c}_{\cc'}(\cl_{\Gi}(z))=\mathfrak{c}_{\cc''}
(\cl_{\Gi}(y)\cdot\cl_{\Gi}(z))$. Finally, by proposition 
\ref{prop:16}, $\cl_{\Gi}(y)\cdot\cl_{\Gi}(z)=\cl_{\Gi}(w)$. 
The claimed formula \eqref{eq:stfg} follows from \eqref{eq:stf} 
taking into account the degrees of $\omega_{y}$ and $g_{y}$. 
\end{proof}

\begin{theorem} 
\label{thm:stcl} 
Let $y$ be a $p$-codimensional, resp. $z$ a $q$-codimensional 
cycle of $X$, and $Y=\supp y$, resp. $Z=\supp z$. Let $\mathfrak
{g}_{y}\in\widehat{H}^{2p}_{\cc,\mathcal{Z}^{p}}(X,p)$, resp. 
$\mathfrak{g}_{z}\in\widehat{H}^{2q}_{\cc',\mathcal{Z}^{q}}(X,q)$ 
be Green objects for the class of $y$, resp. $z$. Let $\mathfrak
{g}_{y}'\in\widehat{H}^{2p}_{\cc,Y}(X,p)$, resp. $\mathfrak{g}_
{z}'\in\widehat{H}^{2q}_{\cc',Z}(X,q)$ be representatives of
$\mathfrak{g}_{y}$, resp. $\mathfrak{g}_{z}$. If $y$ and $z$ 
intersect properly, we have the following statements:   
\begin{enumerate}
\item[(i)] 
The product $\mathfrak{g}'_{y}*\mathfrak{g}'_{z}$ is a Green
object for the cycle $y\cdot z$.
\item[(ii)]
The image of $\mathfrak{g}'_{y}*\mathfrak{g}'_{z}$ in 
$\widehat{H}^{2p+2q}_{\cc'',\mathcal{Z}^{p+q}}(X,p+q)$ does 
not depend on the choice of the representatives $\mathfrak{g}'_
{y}$ and $\mathfrak{g}'_{z}$.
\end{enumerate}
\end{theorem}
\begin{proof}
The first statement is a direct consequence of theorem
\ref{thm:wstp}. Observe that, if the complexes $\cc$ and 
$\cc'$ satisfy the weak purity condition, the second statement 
is automatically fulfilled. \\
To prove the second statement in the case in which the weak 
purity condition is not satisfied, we let $\mathfrak{g}''_{y}$ 
be another choice for a lift of $\mathfrak{g}_{y}$. We then
consider the commutative diagram \eqref{eq:18}, i.e.,       
\begin{displaymath}
\xymatrix{H^{2p-1}_{\cc,Y}(X,p)\ar[r]\ar[d]&
\widetilde{\cc}^{2p-1}(X,p)\ar[r]^{\amap}\ar@{=}[d]& 
\widehat{H}^{2p}_{\cc,Y}(X,p)\ar[r]^{\cl}\ar[d]&
H^{2p}_{\cc,Y}(X,p)\ar[r]\ar[d]&0 \\
H^{2p-1}_{\cc,\mathcal{Z}^{p}}(X,p)\ar[r]&
\widetilde{\cc}^{2p-1}(X,p)\ar[r]^{\amap}& 
\widehat{H}^{2p}_{\cc,\mathcal{Z}^{p}}(X,p)\ar[r]^{\cl}&
H^{2p}_{\cc,\mathcal{Z}^{p}}(X,p)\ar[r]&0.} 
\end{displaymath}
Since we are not assuming weak purity, the vertical map on 
the right of this diagram need not be injective. But, by the 
definition of a Green object for a cycle, we know that the 
images of $\mathfrak{g}'_{y}$ and $\mathfrak{g}''_{y}$ by
the class map agree in $H^{2p}_{\cc,Y}(X,p)$. Therefore, 
we have 
\begin{displaymath}
\mathfrak{g}'_{y}-\mathfrak{g}''_{y}=\amap(\widetilde{x})
\end{displaymath}
for some $x\in\cc^{2p-1}(X,p)$. Moreover, since $\mathfrak{g}'_
{y}$ and $\mathfrak{g}''_{y}$ have the same image in $\widehat{H}^
{2p}_{\cc,\mathcal{Z}^{p}}(X,p)$, we have $\amap(\widetilde{x})=0$
in the group $\widehat{H}^{2p}_{\cc,\mathcal{Z}^{p}}(X,p)$ which
implies $\omega(\amap(\widetilde{x}))=0$, and hence shows $\dd x=0$. 
Therefore, $\widetilde{x}\in H^{2p-1}_{\cc,\mathcal{Z}^{p}}(X,p)
\subset\widetilde{\cc}^{2p-1}(X,p)$ which implies that $\widetilde
{x}$ arises as the image of a class in $H^{2p-1}_{\cc,W}(X,p)$
for a suitable closed subset $W$ of $X$ of codimension greater or
equal to $p$. Writing $\mathfrak{g}'_{z}=(\omega_{z},\widetilde
{g}'_{z})$, we obtain by proposition \ref{prop:restrlem} (iii), 
the equality 
\begin{displaymath}  
\amap(\widetilde{x})*\mathfrak{g}'_{z}=\amap(\widetilde{x\bullet 
\omega_{z}}) 
\end{displaymath}
in the group $\widehat{H}^{2p+2q}_{\cc'',W\cap Z}(X,p+q)$. By
the moving lemma for cycles, there exists a cycle $z'$ linearly 
equivalent to $z$ which intersects $W$ properly. Since the class 
of $z'$ in $H^{2q}_{\cc'}(X,q)$ can also be represented by 
$\omega_{z}$, the class
\begin{displaymath} 
\widetilde{x\bullet\omega_{z}}\in\widetilde{\cc}''^{2p+2q-1}
(X,p+q)
\end{displaymath}
arises as the image of a class with support in the closed subset 
$W\cap\supp z'$ which has codimension greater or equal to $p+q$; 
this shows that this class also arises as the image of a class 
in $H^{2p+2q-1}_{\cc'',\mathcal{Z}^{p+q}}(X,p+q)$. By exactness, 
the element $\amap(\widetilde{x\bullet\omega_{z}})$ therefore 
vanishes in $\widehat{H}^{2p+2q}_{\cc'',\mathcal{Z}^{p+q}}(X,p+q)$. 
This proves the claim. 
\end{proof}

\begin{definition} 
\label{def:5} 
Let $y$ be a $p$-codimensional cycle on $X$, $\mathfrak{g}_
{y}\in\widehat{H}^{2p}_{\cc,\mathcal{Z}^{p}}(X,p)$ a Green 
object for the class of $y$; let $z$ be a $q$-codimensional 
cycle on $X$, $\mathfrak{g}_{z}\in\widehat{H}^{2q}_{\cc',
\mathcal{Z}^{q}}(X,q)$ a Green object for the class of $z$. 
Then, the \emph{$*$-product of $\mathfrak{g}_{y}$ with 
$\mathfrak{g}_{z}$}, denoted by $\mathfrak{g}_{y}*\mathfrak 
{g}_{z}$, is defined as the image of $\mathfrak{g}'_{y}*
\mathfrak{g}'_{z}$ in the group $\widehat{H}^{2p+2q}_{\cc'',
\mathcal{Z}^{p+q}}(X,p+q)$, where $\mathfrak{g}'_{y}$, resp. 
$\mathfrak{g}'_{z}$ is any lift of $\mathfrak{g}_{y}$, resp. 
$\mathfrak{g}_{z}$. 
\end{definition}

Note that, in the above definition, the image of $\mathfrak{g}_
{y}*\mathfrak{g}_{z}\in\widehat{H}^{2p+2q}_{\cc'',\mathcal{Z}^
{p+q}}(X,p+q)$ is a well-defined Green object for the class of 
$y\cdot z$ by theorem \ref{thm:stcl}.

\subsection{Associativity and commutativity}
\label{sec:ca}

\nnpar{Associativity.} 
Let $\cc_{1}$, $\cc_{2}$, $\cc_{3}$, $\cc_{12}$, $\cc_{23}$,
$\cc_{123}$ be $\Gi$-complexes and assume that there are
$\Gi$-pairings   
\begin{align*}
\cc_{1}\otimes\cc_{2}&\longrightarrow\cc_{12} \\
\cc_{2}\otimes\cc_{3}&\longrightarrow\cc_{23} \\
\cc_{12}\otimes\cc_{3}&\longrightarrow\cc_{123} \\
\cc_{1}\otimes\cc_{23}&\longrightarrow\cc_{123}
\end{align*}
which we will denote collectively by $\bullet$. Recall that 
the product $\bullet$ is called associative up to homotopy,
if there exist sheaf morphisms 
\begin{displaymath}
h_{a}:\cc_{1}^{n}(\cdot,p)\otimes\cc_{2}^{m}(\cdot,q)\otimes 
\cc_{3}^{l}(\cdot,r)\longrightarrow\cc_{123}^{n+m+l-1}(\cdot,p+q+r)
\end{displaymath}
satisfying
\begin{displaymath}
(\alpha\bullet\beta)\bullet\gamma -\alpha\bullet(\beta
\bullet\gamma)=\dd h_{a}(\alpha\otimes\beta\otimes\gamma)+
h_{a}\dd(\alpha\otimes\beta\otimes\gamma).
\end{displaymath}
In order to have an associative product in the arithmetic 
Chow groups, we will need a slightly stronger condition 
than associativity up to homotopy. 

\begin{definition}
\label{def:9}
The product $\bullet$ is called \emph{pseudo-associative}, 
if it is associative up to homotopy and the equality
\begin{displaymath}
h_{a}(\alpha\otimes\beta\otimes\gamma)=0
\end{displaymath}
holds for all $\alpha\in{\rm Z}\cc_{1}^{2p}(X,p)$, $\beta 
\in{\rm Z}\cc_{2}^{2q}(X,q)$, $\gamma\in{\rm Z}\cc_{3}^{2r}
(X,r)$.
\end{definition}

\begin{theorem} 
\label{thm:assst} 
Assume that the product $\bullet$ is associative up to homotopy.
Furthermore, let $Y$, $Z$, $W$ be closed subsets of $X$, and 
$\mathfrak{g}_{1}\in\widehat{H}^{2p}_{\cc_{1},Y}(X,p)$, 
$\mathfrak{g}_{2}\in\widehat{H}^{2q}_{\cc_{2},Z}(X,q)$, 
$\mathfrak{g}_{3}\in\widehat{H}^{2r}_{\cc_{3},W}(X,r)$ with
$\mathfrak{g}_{j}=(\omega_{j},\widetilde{g}_{j})$ for $j=1,2,3$. 
Then, we have the following statements:
\begin{enumerate}
\item[(i)] 
In the group $\widehat{H}^{2p+2q+2r}_{\cc_{123},Y\cap Z\cap W}
(X,p+q+r)$, there is an equality 
\begin{equation}\label{eq:45}
(\mathfrak{g}_{1}*\mathfrak{g}_{2})*\mathfrak{g}_{3}-
\mathfrak{g}_{1}*(\mathfrak{g}_{2}*\mathfrak{g}_{3})=
\amap(h_{a}(\omega_{1}\otimes\omega_{2}\otimes\omega_{3})).
\end{equation}
\item[(ii)] 
For $\mathfrak{g}_{1}$, or $\mathfrak{g}_{2}$, or $\mathfrak{g}_{3}
\in\Ker(\omega)$, the equality
\begin{displaymath}
(\mathfrak{g}_{1}*\mathfrak{g}_{2})*\mathfrak{g}_{3}=
\mathfrak{g}_{1}*(\mathfrak{g}_{2}*\mathfrak{g}_{3})
\end{displaymath}
holds.
\item[(iii)]
Moreover, if the product $\bullet$ is pseudo-associative, 
again the equality
\begin{displaymath}
(\mathfrak{g}_{1}*\mathfrak{g}_{2})*\mathfrak{g}_{3}=
\mathfrak{g}_{1}*(\mathfrak{g}_{2}*\mathfrak{g}_{3})
\end{displaymath}
holds.
\end{enumerate}
\end{theorem}
\begin{proof} 
The statements (ii) and (iii) follow immediately from (i). 
To prove (i), let us write $t=p+q+r$, and let
\begin{align*}
f_{1}&:\cc_{1}(X,p)\longrightarrow\cc_{1}(X\setminus Y,p), \\
f_{2}&:\cc_{2}(X,p)\longrightarrow\cc_{2}(X\setminus Z,p), \\
f_{3}&:\cc_{3}(X,p)\longrightarrow\cc_{3}(X\setminus W,p)
\end{align*}
be the restriction morphisms. For simplicity, let us denote 
$\cc_{123}$ by $\cc$. By means of the identification 
\eqref{eq:secass} we now observe that the product $\bullet$ 
induces a morphism of $2$-iterated complexes between $f_{1}
\otimes f_{2}\otimes f_{3}$ and the $2$-iterated complex
\begin{align}
\xi=(\cc(X,t)\overset{\delta_{1}}{\longrightarrow}&
\cc(X\setminus Y,t)\oplus\cc(X\setminus Z,t)\oplus 
\cc(X\setminus W,t)\notag \\
\overset{\delta_{2}}{\longrightarrow}&\cc(X\setminus
(Y\cup Z),t)\oplus\cc(X\setminus(Y\cup W),t)\oplus
\cc(X\setminus(Z\cup W),t)\notag \\
\overset{\delta_{3}}{\longrightarrow}&\cc(X\setminus 
(Y\cup Z\cup W),t));
\end{align}
here 
\begin{align*}
\delta_{1}(a)&=(a,a,a), \\
\delta_{2}(a,b,c)&=(b-a,c-a,c-b), \\
\delta_{3}(a,b,c)&=a-b+c.
\end{align*}
Let us denote by $\cc(X;Y,Z,W,t)$ the simple of the $2$-iterated
complex 
\begin{align}
&\cc(X\setminus Y,t)\oplus\cc(X\setminus Z,t)\oplus
\cc(X\setminus W,t)\overset{-\delta_{2}}{\longrightarrow}
\notag \\
&\cc(X\setminus(Y\cup Z),t)\oplus\cc(X\setminus(Y\cup W),t)
\oplus\cc(X\setminus(Z\cup W),t)\overset{-\delta_{3}}
{\longrightarrow}\notag \\
&\cc(X\setminus(Y\cup Z\cup W),t).
\label{eq:7}
\end{align}
Except for the signs of the morphisms this is the \v{C}ech 
complex associated to the open covering $\{X\setminus Y,
X\setminus Z,X\setminus W\}$. As in proposition \ref{prop:comp1}, 
there is a natural isomorphism    
\begin{displaymath}
s(\xi)\longrightarrow s(\cc(X,t)\overset{\delta_{1}}
{\longrightarrow}\cc(X;Y,Z,W,t))   
\end{displaymath}
which is given by the identity on the level of elements.

Let us consider the commutative diagram
\begin{equation}
\label{eq:comass}
\begin{CD}
\cc(X\setminus(Y\cap Z\cap W),t)@>i_{1}>>\cc(X;Y\cap Z,W,t) \\
@V i_{2}VV@VV k_{1}V \\
\cc(X;Y,Z\cap W,t)@>k_{2}>>\cc(X;Y,Z,W,t),
\end{CD}
\end{equation}
where
\begin{align*}
i_{1}(a)&=((a,a),0), \\
i_{2}(a)&=((a,a),0), \\
k_{1}((a,b),c)&=((a,a,b),(0,c,c),0), \\
k_{2}((a,b),c)&=((a,b,b),(c,c,0),0).
\end{align*}
The Mayer-Vietoris principle of the sheaf $\cc$ implies that 
all the arrows in the above diagram are quasi-isomorphisms.

Recalling $\cc(X;Y,Z,t)$ from definition \ref{def:2}, we
introduce the notation
\begin{align*}
\widehat{H}^{2t}_{\cc,Y,Z}(X,t)&=\widehat{H}^{2t}(\cc(X,t),
\cc(X;Y,Z,t)), \\
\widehat{H}^{2t}_{\cc,Y,Z,W}(X,t)&=\widehat{H}^{2t}(\cc(X,t),
\cc(X;Y,Z,W,t)),
\end{align*}
and similar notations for the other closed subsets under
consideration. We then obtain from the commutative diagram 
\eqref{eq:comass} and corollary \ref{cor:tci} the commutative 
diagram 
\begin{displaymath}
\begin{CD}
\widehat{H}^{2t}_{\cc,Y\cap Z\cap W}(X,t)@>i_{1}>>
\widehat{H}^{2t}_{\cc,Y\cap Z,W}(X,t) \\
@V i_{2}VV@VV k_{1}V \\ 
\widehat{H}^{2t}_{\cc,Y,Z\cap W}(X,t)@>k_{2}>>
\widehat{H}^{2t}_{\cc,Y,Z,W}(X,t),
\end{CD}
\end{displaymath}
where all the arrows are isomorphisms. Therefore, in order 
to prove equation \eqref{eq:45}, it is enough to compare 
explicit representatives of both elements in the group  
$\widehat{H}^{2t}_{\cc,Y,Z,W}(X,t)$. 

To obtain these representatives, we consider the commutative
diagram 
\begin{displaymath}
\begin{CD}
\widehat{H}^{2p}_{\cc_{1},Y}(X,p)\otimes\widehat{H}^{2q+2r}_
{\cc_{23},Z\cap W}(X,q+r)@>*>>\widehat{H}^{2t}_{\cc,Y,Z\cap W}
(X,t) \\
@V\cong VV@VV k_{2}V \\
\widehat{H}^{2p}_{\cc_{1},Y}(X,p)\otimes\widehat{H}^{2q+2r}_
{\cc_{23},Z,W}(X,q+r)@>*>>\widehat{H}^{2t}_{\cc,Y,Z,W}(X,t),
\end{CD}
\end{displaymath}
where the bottom arrow is induced from the morphism $\beta_{1}$ 
given in proposition \ref{prop:assten}. We thus obtain
\begin{align*}
(a,b)*(c,(d,e),f)=&(a\bullet c,(b\bullet c,a\bullet d,
a\bullet e), \\
&(-b\bullet d,-b\bullet e,a\bullet f),b\bullet f).
\end{align*}
Analogously, there is another commutative diagram 
\begin{displaymath}
\begin{CD}
\widehat{H}^{2p+2q}_{\cc_{12},Y\cap Z}(X,p+q)\otimes\widehat
{H}^{2r}_{\cc_{3},W}(X,r)@>*>>\widehat{H}^{2t}_{\cc,Y\cap Z,W}
(X,t) \\
@V\cong VV@VV k_{1}V \\
\widehat{H}^{2p+2q}_{\cc_{12},Y,Z}(X,p+q)\otimes\widehat{H}^{2r}_
{\cc_{3},W}(X,r)@>*>>\widehat{H}^{2t}_{\cc,Y,Z,W}(X,t),
\end{CD}
\end{displaymath}
where the bottom arrow is induced from the morphism $\beta_{2}$ 
given in proposition \ref{prop:assten}. We obtain
\begin{align*}
(a,(b,c),d)*(e,f)=&(a\bullet e,(b\bullet e,c\bullet e, 
a\bullet f), \\
&(d\bullet e,-b\bullet f,-c\bullet f),d\bullet f).
\end{align*}
Summarizing, we find the following identities in the group
$\widehat{H}^{2t}_{\cc,Y,Z,W}(X,t)$ 
\begin{align*}
(\mathfrak{g}_{1}*\mathfrak{g}_{2})*&\mathfrak{g}_{3}=
((\omega_{1}\bullet\omega_{2})\bullet\omega_{3},
(((g_{1}\bullet\omega_{2})\bullet\omega_{3},(\omega_{1}
\bullet g_{2})\bullet\omega_{3},(\omega_{1}\bullet\omega_{2})
\bullet g_{3}), \\
&(-(g_{1}\bullet g_{2})\bullet\omega_{3},-(g_{1}\bullet 
\omega_{2})\bullet g _{3},-(\omega_{1}\bullet g_{2})\bullet 
g_{3}),-(g_{1}\bullet g_{2})\bullet g_{3})\widetilde{\phantom{a}}),
\end{align*}
and
\begin{align*}
\mathfrak{g}_{1}*(\mathfrak{g}_{2}*&\mathfrak{g}_{3})=
(\omega_{1}\bullet(\omega_{2}\bullet\omega_{3}),((g_{1}
\bullet(\omega_{2}\bullet\omega_{3}),\omega_{1}\bullet 
(g_{2}\bullet\omega_{3}),\omega_{1}\bullet(\omega_{2}
\bullet g_{3})), \\
&(-g_{1}\bullet(g_{2}\bullet\omega_{3}),-g_{1}\bullet 
(\omega_{2}\bullet g _{3}),-\omega_{1}\bullet(g_{2}\bullet 
g_{3})),-g_{1}\bullet(g_{2}\bullet g_{3}))\widetilde{\phantom{a}}).
\end{align*}
With the element $y\in\cc^{2t-2}(X;Y,Z,W,t)$ given by
\begin{align*}
y=&((h_{a}(g_{1}\otimes\omega_{2}\otimes\omega_{3}),h_{a}
(\omega_{1}\otimes g_{2}\otimes\omega_{3}),h_{a}(\omega_{1}
\otimes\omega_{2}\otimes g_{3})), \\
&(h_{a}(g_{1}\otimes g_{2}\otimes\omega_{3}),h_{a}(g_{1}
\otimes\omega_{2}\otimes g_{3}),h_{a}(\omega_{1}\otimes 
g_{2}\otimes g_{3})), \\
&-h_{a}(g_{1}\otimes g_{2}\otimes g_{3})),
\end{align*}
we derive the identity 
\begin{equation}
\label{eq:19}
(\mathfrak{g}_{1}*\mathfrak{g}_{2})*\mathfrak{g}_{3}-
\mathfrak{g}_{1}*(\mathfrak{g}_{2}*\mathfrak{g}_{3})-
(0,\dd y)=\amap(h_{a}(\omega_{1}\otimes\omega_{2}\otimes 
\omega_{3})),
\end{equation}
which proves (i).
\end{proof}
 

\nnpar{Commutativity.} 
Let $\cc_{1}$ and $\cc$ be $\Gi$-complexes and let  
\begin{displaymath}
\bullet:\cc_{1}\otimes\cc_{1}\longrightarrow\cc
\end{displaymath}
be a $\Gi$-pairing. Recall that the pairing $\bullet$ 
is commutative up to homotopy, if there exist sheaf 
morphisms
\begin{displaymath}
h_{c}:\cc_{1}^{n}(\cdot,p)\otimes\cc_{1}^{m}(\cdot,q)
\longrightarrow \cc^{n+m-1}(\cdot,p+q)
\end{displaymath}
satisfying
\begin{displaymath}
\alpha\bullet\beta -(-1)^{n m}\beta\bullet\alpha=
\dd h_{c}(\alpha\otimes\beta)+ h_{c}\dd(\beta\otimes\alpha ).
\end{displaymath}

As in the case of the associativity, in order to have a 
commutative product on the level of arithmetic groups, 
we need a slightly stronger condition than commutativity
up to homotopy.

\begin{definition}
The product $\bullet$ is called \emph{pseudo-commutative},
if it is commutative up to homotopy and if the equality
\begin{displaymath} 
h_{c}(\alpha\otimes\beta)=0
\end{displaymath}
holds for all $\alpha\in{\rm Z}\cc_{1}^{2p}(X,p)$ and
$\beta\in{\rm Z}\cc_{1}^{2q}(X,q)$. 
\end{definition}

\begin{theorem}
\label{thm:commst}
Let $Y$ and $Z$ be closed subsets of $X$, and $\mathfrak{g}_{1}
\in\widehat{H}^{2p}_{\cc_{1},Y}(X,p)$, $\mathfrak{g}_{2}\in
\widehat{H}^{2q}_{\cc_{1},Z}(X,q)$ with $\mathfrak{g}_{j}=(\omega_
{j},\widetilde{g}_{j})$ for $j=1,2$. Then, we have the following
statements:
\begin{enumerate}
\item[(i)]
In the group $\widehat{H}^{2p+2q}_{\cc,Y\cap Z}(X,p+q)$, there 
is an equality 
\begin{displaymath}
\mathfrak{g}_{1}*\mathfrak{g}_{2}-\mathfrak{g}_{2}*\mathfrak{g}_{1}=
\amap(h_{c}(\omega_{1}\otimes\omega_{2})). 
\end{displaymath}
\item[(ii)] 
For $\mathfrak{g}_{1}$ or $\mathfrak{g}_{2}\in\Ker(\omega)$, the
equality
\begin{displaymath}
\mathfrak{g}_{1}*\mathfrak{g}_{2}=\mathfrak{g}_{2}*\mathfrak{g}_{1}
\end{displaymath} 
holds.
\item[(iii)]
Moreover, if the product $\bullet$ is pseudo-commutative, again
the equality 
\begin{displaymath}
\mathfrak{g}_{1}*\mathfrak{g}_{2}=\mathfrak{g}_{2}*\mathfrak{g}_{1}
\end{displaymath}
holds.
\end{enumerate}
\end{theorem}
\begin{proof} 
The statements (ii) and (iii) follow immediately from (i). To prove
(i), let us write $s=p+q$. We recall that $\mathfrak{g}_{1}*\mathfrak
{g}_{2}$ is defined as an element of 
\begin{displaymath}
\widehat{H}^{2s}_{\cc,Y,Z}(X,s)=\widehat{H}^{2s}(\cc(X,s),\cc(X;Y,Z,s))
\end{displaymath}
with $\cc(X;Y,Z,s)$ given in definition \ref{def:2}. By means of 
the inverse of the isomorphism
\begin{displaymath}
\widehat{H}^{2s}_{\cc,Y\cap Z}(X,s)\overset{\cong}{\longrightarrow} 
\widehat{H}^{2s}_{\cc,Y,Z}(X,s),
\end{displaymath}
which is induced by the kernel-simple quasi-isomorphism, the element
$\mathfrak{g}_{1}*\mathfrak{g}_{2}$ is sent to $\widehat{H}^{2s}_{\cc,
Y\cap Z}(X,s)$. We now consider the commutative diagram of complexes
\begin{displaymath}
\xymatrix
{&\cc(X;Y,Z,s)\ar[dd]^{T} \\
\cc(X\setminus(Y\cap Z),s)\ar[ur]^{i}\ar[dr]^{i}& \\
&\cc(X;Z,Y,s),}
\end{displaymath}
where $T((\alpha,\beta),\gamma)=((\beta,\alpha),-\gamma)$. We note 
that the isomorphism $T$ induces an isomorphism, also denoted by $T$,
\begin{displaymath}
T:\widehat{H}^{2s}_{\cc,Y,Z}(X,s)\overset{\cong}{\longrightarrow} 
\widehat{H}^{2s}_{\cc,Z,Y}(X,s).
\end{displaymath}
We are left to compare the element
\begin{align*}
T((\omega_{1},\widetilde{g}_{1})*(\omega_{2},\widetilde{g}_{2}))&=
T\left(\omega_{1}\bullet\omega_{2},\left((g_{1}\bullet\omega_{2},
\omega_{1}\bullet g_{2}),-g_{1}\bullet g_{2}\right)\widetilde
{\phantom {a}}\right) \\
&=\left(\omega_{1}\bullet\omega_{2},\left((\omega _{1}\bullet 
g_{2},g_{1}\bullet\omega_{2}),g_{1}\bullet g_{2}\right)\widetilde
{\phantom {a}}\right)
\end{align*}
with the element
\begin{align*}
(\omega_{2},\widetilde{g}_{2})*(\omega_{1},\widetilde{g}_{1})&=
\left(\omega_{2}\bullet\omega_{1},\left((g_{2}\bullet\omega_{1},
\omega_{2}\bullet g_{1}),-g_{2}\bullet g_{1}\right)\widetilde
{\phantom {a}}\right)
\end{align*}
in the group $\widehat{H}^{2s}_{\cc,Z,Y}(X,s)$. To do this,
let us consider the element $y\in\cc^{2s-2}(X;Z,Y,s)$ given by
\begin{align*}
y=((h_{c}(\omega_{1}\otimes g_{2}),h_{c}(g_{1}\otimes\omega_{2})),
h_{c}(g_{1}\otimes g_{2})).
\end{align*}
With this element we derive
\begin{equation}
T(\mathfrak{g}_{1}*\mathfrak{g}_{2})-\mathfrak{g}_{2}*\mathfrak
{g}_{1}-(0,\dd y)=\amap(h_{c}(\omega_{1}\otimes\omega_{2}));
\end{equation}
here we have used the fact that the elements $\omega_{j}$ have even degree,
and the elements $g_{j}$ have odd degree. This completes the proof
of (i).
\end{proof}

\nnpar{$\Gi$-algebras and $\Gi$-modules.} 
The above results can be applied to the case of $\Gi$-algebras and
$\Gi$-modules. 

\begin{definition}
\label{def:4}
A \emph{$\Gi$-algebra over $X$} is a $\Gi$-complex $\cc$ together 
with a $\Gi$-pairing
\begin{displaymath}
\bullet: \cc\otimes\cc\longrightarrow\cc
\end{displaymath}
making $(\cc,\bullet)$ into a graded differential algebra with
unit element which is associative up to homotopy and commutative 
up to homotopy. We call $\cc$ \emph{pseudo-associative} or
\emph{pseudo-commutative}, if the product $\bullet$ is.
\end{definition}

\noindent
A direct consequence of theorems \ref{thm:assst} and \ref{thm:commst}
is the following theorem

\begin{theorem}
\label{thm:asscomm}
Let $\cc$ be a pseudo-associative and pseudo-commutative algebra 
over $X$. Then, the direct sum  
\begin{displaymath}
\bigoplus_{\substack{Y\text{\rm closed}\\p\ge 0}}\widehat{H}^
{2p}_{\cc,Y}(X,p)
\end{displaymath}
is an associative and commutative algebra.
\hfill$\square$
\end{theorem}

Once we have the notion of $\Gi$-algebra, we can define the concept 
of a $\Gi$-module.

\begin{definition} 
Let $\cc$ be a $\Gi$-algebra over $X$. A \emph{$\Gi$-module 
over $\cc$} is a $\Gi$-complex $\cc'$ over $X$ together with 
a $\Gi$-pairing 
\begin{displaymath}
\cc\otimes\cc'\longrightarrow\cc',
\end{displaymath}
which is associative up to homotopy. We call such a pairing 
a $\Gi$-action. We call the $\Gi$-action \emph{pseudo-associative}, 
if the pairing is pseudo-associative.
\end{definition}

\noindent
As a consequence of theorem \ref{thm:assst}, we have the theorem.

\begin{theorem}
Let $\cc$ be a pseudo-associative $\Gi$-algebra over $X$, and 
$\cc'$ a $\Gi$-complex provided with a pseudo-associative
$\Gi$-action. Then, the direct sum
\begin{displaymath}
\bigoplus_{\substack{Y\text{\rm closed}\\p\ge 0}}\widehat{H}^
{2p}_{\cc',Y}(X,p) 
\end{displaymath}
is an associative module over the associative algebra $\bigoplus_{Y,p}
\widehat{H}^{2p}_{\cc,Y}(X,p)$.
\hfill$\square$
\end{theorem}

\nnpar{Multiple products.}
\begin{proposition}
Let $y_{j}$ be $p_{j}$-codimensional cycles of $X$, $Y_{j}=
\supp y_{j}$, and $\mathfrak{g}_{j}$ Green objects for the 
class of $y_{j}$ ($j=1,...,r$) such that $\codim(Y_{1}\cap...
\cap Y_{r})=p_{1}+...+p_{r}$. Then, the $r$-fold $*$-product 
$\mathfrak{g}_{1}*...*\mathfrak{g}_{r}$ is a well-defined 
Green object for the class of $y_{1}\cdot...\cdot y_{r}$ 
even though partial intersections of the $Y_{j}$'s need not 
be proper.
\end{proposition}
\begin {proof}
The proof of the proposition is an immediate consequence of
theorems \ref{thm:wstp}, \ref{thm:stcl}, and \ref{thm:asscomm}.
\end{proof}

\subsection{Functorial properties of Green objects}
\label{sec:funcgreen}

Throughout this section, \emph{$k$-scheme} will mean regular separated
scheme of finite type over $k$.

\nnpar{Direct image of a $\Gi$-complex.} 
Let $f:X\longrightarrow Y$ be a morphism of $k$-schemes, and 
$\cc_{X}$ a $\Gi$-complex over $X$. Since $\cc_{X}$ satisfies the 
Mayer-Vietoris principle, it is clear that $f_{*} \cc_{X}$ also 
satisfies the Mayer-Vietoris principle. Since $\Gi$ is a sheaf
in the big Zariski site over $k$, there is a morphism of sheaves 
$\Gi_{Y}\longrightarrow f_{*}\Gi_{X}$ on $Y$. Therefore, the 
composition 
\begin{displaymath}
\Gi_{Y}\longrightarrow f_{*}\Gi_{X}\longrightarrow f_{*}\cc_{X} 
\end{displaymath}
determines a structure of $\Gi$-complex on $f_{*}\cc_{X}$.

\nnpar{Pull-back of Green objects.} 
\begin{definition} 
\label{def:353}
Let $f:X\longrightarrow Y$ be a morphism of $k$-schemes, and let
$\cc_{X}$, $\cc_{Y}$ be $\Gi$-complexes over $X$, $Y$, respectively.
A \emph{contravariant $f$-morphism (of $\Gi$-complexes)} is a
morphism of $\Gi$-complexes over $Y$ 
\begin{displaymath}
f^{\#}:\cc_{Y}\longrightarrow f_{*}\cc_{X}.
\end{displaymath}
\end{definition}

If $g:Y\longrightarrow Z$ is another morphism of $k$-schemes,
and $g^{\#}:\cc_{Z}\longrightarrow g_{*}\cc_{Y}$ a contravariant
$g$-morphism, then $g_{*}(f^{\#})\circ g^{\#}$ is a contravariant
$(g\circ f)$-morphism.

Let $f^{\#}$ be as in definition \ref{def:353}. Then, for any
open subset $U\subset Y$, any closed subset $Z\subset U$, and 
any integer $p$, we have a morphism     
\begin{align*}
s(\cc_{Y}(U,p),\cc_{Y}(U\setminus Z,p))&\longrightarrow
s(\cc_{X}(f^{-1}(U),p),\cc_{X}(f^{-1}(U\setminus Z),p) )
\end{align*}
given by the assignment
\begin{displaymath}
(\omega,g)\longmapsto(f^{\#}\omega,f^{\#}g).
\end{displaymath}
Therefore, we have induced morphisms, again denoted by $f^{\#}$,    
\begin{align*}
f^{\#}:H^{n}_{\cc_{Y},Z}(U,p)&\longrightarrow H^{n}_{\cc_{X},
f^{-1}(Z)}(f^{-1}(U),p), \\
f^{\#}:\widehat{H}^{n}_{\cc_{Y},Z}(U,p)&\longrightarrow\widehat 
{H}^{n}_{\cc_{X},f^{-1}(Z)}(f^{-1}(U),p).
\end{align*}

\begin{theorem}
\label{thm:8}
Let $f:X\longrightarrow Y$ be a morphism of $k$-schemes, 
let $\cc_{X}$, $\cc_{Y}$ be $\Gi$-complexes over $X$, $Y$, 
respectively, and let $f^{\#}$ be a contravariant $f$-morphism. 
Then, we have the following statements:
\begin{enumerate}
\item[(i)]
There is a commutative diagram
\begin{displaymath}
\begin{CD}
H^{n}_{Z}(U,\Gi(p))@>f^{\ast}>>H^{n}_{f^{-1}(Z)}(f^{-1}(U),
\Gi(p)) \\
@V\mathfrak{c}_{\cc_Y}VV@VV\mathfrak{c}_{\cc_X}V \\
H^{n}_{\cc_{Y},Z}(U,p)@>f^{\#}>>H^{n}_{\cc_{X},f^{-1}(Z)}
(f^{-1}(U),p).
\end{CD}
\end{displaymath}
\item[(ii)] 
Setting $U'=f^{-1}(U)$ and $Z'=f^{-1}(Z)$, we have the 
commutative diagram
\begin{align*}
\makebox[-1cm]{}& 
\begin{CD}
H^{n-1}_{\cc_{Y},Z}(U,p)@>>>\widetilde{\cc}_{Y}^{n-1}(U,p) 
@>\amap>>\widehat{H}^{n}_{\cc_{Y},Z}(U,p)@>\cl>>H^{n}_
{\cc_{Y},Z}(U,p) \\
@Vf^{\#}VV@Vf^{\#}VV@Vf^{\#}VV@Vf^{\#}VV \\    
H^{n-1}_{\cc_{X},Z'}(U',p)@>>>\widetilde{\cc}_{X}^{n-1}
(U',p)@>\amap>>\widehat{H}^{n}_{\cc_{X},Z'}(U',p)@>\cl>>
H^{n}_{\cc_{X},Z'}(U',p).
\end{CD}
\end{align*}
\hfill$\square$
\end{enumerate}
\end{theorem}

\begin{corollary}
Let $z\in{\rm Z}^{p}(Y)$ be a $p$-codimensional cycle on $Y$
with $Z=\supp z$, and $\mathfrak{g}_{z}=(\omega_{z},\widetilde
{g}_{z})$ a Green object for the cycle $z$. Then, we have:
\begin{enumerate}
\item[(i)] 
The element $f^{\#}(\mathfrak{g}_{z})$ is a weak Green object for 
the class of $f^{*}(z)$ in the group $\CH^{p}_{f^{-1}(Z)}(X)$.   
\item[(ii)]
If the cycle $f^{-1}(Z)$ has codimension $p$ in $X$, the element
$f^{\#}(\mathfrak{g}_{z})$ is a Green object for the cycle $f^{*}(z)$.
\end{enumerate}
\hfill$\square$  
\end{corollary}

By a similar argument to the one used in the course of the proof 
of theorem \ref{thm:stcl}, we derive:

\begin{proposition}
\label{prop:356}
Let $z\in{\rm Z}^{p}(Y)$ be a $p$-codimensional cycle on $Y$
with $Z=\supp z$ and $\codim f^{-1}(Z)=p$. Let $\mathfrak{g}_{z}$ 
be a Green object for the class of $z$, and $\mathfrak{g}'_{z}
\in\widehat{H}^{2p}_{\cc_{Y},Z}(Y,p)$ be any representative of
$\mathfrak{g}_{z}$. Then, the image of $f^{\#}(\mathfrak{g}'_{z})$ 
in the group $\widehat{H}^{2p}_{\cc_{X},\mathcal{Z}^{p}}(X,p)$ is
independent of the choice of $\mathfrak{g}_{z}'$.
\hfill$\square$
\end{proposition}

\begin{definition}
\label{def:6} 
With the hypothesis of proposition \ref{prop:356}, we will 
denote the image of $f^{\#}(\mathfrak{g}'_{z})$ in the group 
$\widehat H^{2p}_{\cc_{X},\mathcal{Z}^{p}}(X,p)$ by $f^{\#}
(\mathfrak{g}_{z})$.
\end{definition}

\begin{definition}
\label{def:358}
Let $\cc_{X}\otimes\cc'_{X}\overset{\bullet_{X}}{\longrightarrow}
\cc''_{X}$, and $\cc_{Y}\otimes\cc'_{Y}\overset{\bullet_{Y}}
{\longrightarrow}\cc''_{Y}$ be $\Gi$-pairings, and let 
\begin{align*}
f^{\#}_{\cc}:\cc_{Y}&\longrightarrow f_{*}\cc_{X}, \\
f^{\#}_{\cc'}:\cc'_{Y}&\longrightarrow f_{*}\cc'_{X}, \\
f^{\#}_{\cc''}:\cc''_{Y}&\longrightarrow f_{*}\cc''_{X}
\end{align*}
be contravariant $f$-morphisms. We will call the $f$-morphisms  
$f^{\#}_{\cc}$, $f^{\#}_{\cc'}$, $f^{\#}_{\cc''}$ compatible with 
the pairings $\bullet_{X}$ and $\bullet_{Y}$, if the equality
\begin{displaymath}
f^{\#}_{\cc}(y)\bullet_{X} f^{\#}_{\cc'}(y')=
f^{\#}_{\cc''}(y\bullet_{Y} y')
\end{displaymath}
holds for all sections $y$, $y'$ of $\cc_{Y}$, $\cc'_{Y}$,
respectively.
\end{definition}

\begin{proposition} 
\label{prop:11}
Let $f^{\#}_{\cc}$, $f^{\#}_{\cc'}$, and $f^{\#}_{\cc''}$ 
be contravariant $f$-morphisms which are compatible with 
the pairings $\bullet_{X}$ and $\bullet_{Y}$ as in definition 
\ref{def:358}. Then, for any pair of Green objects $\mathfrak
{g}_{1}$, $\mathfrak{g}_{2}$ the equality
\begin{displaymath}
f^{\#}(\mathfrak{g}_{1}*\mathfrak{g}_{2})=f^{\#}
(\mathfrak{g}_{1})*f^{\#}(\mathfrak{g}_{2})
\end{displaymath}
holds.
\hfill$\square$
\end{proposition}

\nnpar{Change of $\Gi$-complex.}
\begin{remark}
We point out that all the results concerning inverse images 
are compatible with the change of $\Gi$-complex, i.e., they
apply to the case in which we have two $\Gi$-complexes $\cc$ and 
$\cc'$ over $X$ together with a morphism $f^{\#}:\cc\longrightarrow 
\cc'$ of $\Gi$-complexes.   
\end{remark}

By corollary \ref{cor:tci}, there is a case, when the change of
$\Gi$-complex does not change the space of Green objects.

\begin{proposition}
\label{prop:21}
Let $X$ be a $k$-scheme, and $\cc$, $\cc'$ two $\Gi$-complexes
over $X$. Let $f^{\#}:\cc\longrightarrow\cc'$ be a morphism of
$\Gi$-complexes such that the morphism  
\begin{displaymath}
f^{\#}:\cc(X,p)\longrightarrow\cc'(X,p)
\end{displaymath}
is an isomorphism for any integer $p$, and such that the morphism
\begin{displaymath}
f^{\#}:\cc(U,p)\longrightarrow\cc'(U,p)
\end{displaymath}
is a quasi-isomorphism for any open subset $U\subset X$ and any   
integer $p$. Then, the induced morphism 
\begin{displaymath}
f^{\#}:\widehat{H}^{n}_{\cc,Z}(X,p)\longrightarrow\widehat
{H}^{n}_{\cc',Z}(X,p)
\end{displaymath}
is an isomorphism for any closed subset $Z\subset X$ and any pair 
of integers $n,p$. 
\end{proposition}
\begin{proof}
This follows immediately from corollary \ref{cor:tci}.
\end{proof}

\nnpar{Push-forward of Green objects.}
\begin{definition}
Given a graded complex $A=(A^{*}(*),\dd)$, the \emph{twisted 
complex $A(d)$} is given by $A(d)^{n}(p)=A^{n}(p+d)$ with the 
same differential.   
\end{definition}

\begin{definition} 
\label{def:3}
Let $f:X\longrightarrow Y$ be a proper morphism of equidimensional 
$k$-schemes of relative dimension $d$. Let $\cc_{X}$, $\cc_{Y}$ be
$\Gi$-complexes over $X$, $Y$, respectively. A \emph{covariant 
$f$-morphism (of $\Gi$-complexes)} is a morphism of graded complexes 
of sheaves over $Y$ 
\begin{displaymath}
f_{\#}:f_{*}\cc_{X}\longrightarrow\cc_{Y}(-d)[-2d]
\end{displaymath}
such that for any open subset $U\subset Y$, and any closed subset 
$Z\subset f^{-1}(U)$, the induced diagram      
\begin{displaymath}
\begin{CD}
H^{n}_{Z}(f^{-1}(U),\Gi(p))@>f_{!}>>H^{n-2d}_{f(Z)}(U,\Gi(p-d)) \\
@V\mathfrak{c}VV@VV\mathfrak{c}V \\
H^{n}_{\cc_{X},Z}(f^{-1}(U),p)@>f_{\#}>>H^{n-2d}_{\cc_{Y},f(Z)}
(U,p-d)
\end{CD}
\end{displaymath}
is commutative. In the above diagram, the arrow $f_{\#}$ is given 
by the composition of morphisms
\begin{displaymath}
H^{n}_{\cc_{X},Z}(f^{-1}(U),p)\longrightarrow H^{n}_{\cc_{X},
f^{-1}(f(Z))}(f^{-1}(U),p)\longrightarrow H^{n-2d}_{\cc_{Y},f(Z)}
(U,p-d);
\end{displaymath}
here the first map is the restriction morphism, and the second map 
is induced by the morphism of complexes
\begin{align}
\label{eq:21}
s(\cc_{X}(f^{-1}(U),p),\cc_{X}&( f^{-1}(U\setminus f(Z)),p))
\longrightarrow\notag \\
&s(\cc_{Y}(U,p-d),\cc_{Y}(U\setminus f(Z),p-d))[-2d]
\end{align}
which is given by the assignment
\begin{displaymath}
(\omega,g)\longmapsto(f_{\#}\omega,f_{\#}g).
\end{displaymath}
\end{definition}

With the notations and assumptions of definition \ref{def:3}, 
a covariant $f$-morphism $f_{\#}$ induces a morphism by means of the  
restriction morphism followed by the morphism of complexes 
\eqref{eq:21}; this morphism is again denoted $f_{\#}$, 
\begin{displaymath}
f_{\#}:\widehat{H}^{n}_{\cc_{X},Z}(f^{-1}(U),p)\longrightarrow 
\widehat{H}^{n-2d}_{\cc_{Y},f(Z)}(U,p-d).
\end{displaymath}

If $g:Y\longrightarrow Z$ is another proper morphism of
equidimensional $k$-schemes of relative dimension $e$,
and $g_{\#}:g_{*}\cc_{Y}\longrightarrow\cc_{Z}(-e)[-2e]$ 
is a covariant $g$-morphism, then $g_{\#}\circ g_{*}(f_{\#})$ 
is a covariant $(g\circ f)$-morphism.

\begin{theorem}
\label{thm:9} 
Let $f:X\longrightarrow Y$ be a proper morphism of equidimensional 
$k$-schemes of relative dimension $d$, and let $f_{\#}$ be a
covariant $f$-morphism. Given $U$, $Z$ as in definition \ref{def:3}, 
set $U'=f^{-1}(U)$, $Z'=f(Z)$, respectively. Then, for every $n,p$, we
have the 
commutative diagram
\begin{displaymath}
\begin{CD}
H^{n-1}_{\cc_{X},Z}(U',p)@>>>\widetilde{\cc}_{X}^{n-1}(U',p) 
@>\amap>>\widehat{H}^{n}_{\cc_{X},Z}(U',p)@>\cl>>H^{n}_
{\cc_{X},Z}(U',p) \\
@Vf_{\#}VV@Vf_{\#}VV@Vf_{\#}VV@Vf_{\#}VV \\    
H^{m-1}_{\cc_{Y},Z'}(U,q)@>>>\widetilde{\cc}_{Y}^{m-1}(U,q) 
@>\amap>>\widehat{H}^{m}_{\cc_{Y},Z'}(U,q)@>\cl>>H^{m}_
{\cc_{Y},Z'}(U,q),
\end{CD}
\end{displaymath}
where $q=p-d$ and $m=n-2d$.
\hfill $\square$
\end{theorem}

Given a proper morphism $f:X\longrightarrow Y$ of equidimensional 
$k$-schemes of relative dimension $d$, and a closed subset $Z\subset 
X$, we have $\codim f(Z)\ge\max(\codim Z-d,0)$. Therefore, the
following result is simpler to prove than the corresponding statement 
for the inverse image of Green objects. 

\begin{proposition} 
\label{prop:366}
Let $f:X\longrightarrow Y$ be a proper morphism of equidimensional 
$k$-schemes of relative dimension $d$, and let $f_{\#}$ be a covariant
$f$-morphism. Let $z\in{\rm Z}^{p}(X)$ be a $p$-codimensional cycle 
on $X$ with $Z=\supp z$. Let $\mathfrak{g}_{z}$ be a Green object 
for the class of $z$, and $\mathfrak{g}'_{z}\in\widehat{H}^{2p}_
{\cc_{X},Z}(X,p)$ be any representative of $\mathfrak{g}_{z}$. Then,
the image of $f_{\#}(\mathfrak{g}'_{z})$ in the group $\widehat
{H}^{2p-2d}_{\cc_{Y},\mathcal{Z}^{p-d}}(Y,p-d)$ is independent of 
the choice of $\mathfrak{g}_{z}'$; it is a Green object for the 
class of the cycle $f_{*}(z)$.
\hfill$\square$
\end{proposition}

\begin{definition} 
With the hypothesis of proposition \ref{prop:366}, we will 
denote the image of $f_{\#}(\mathfrak{g}'_{z})$ in the group 
$\widehat{H}^{2p-2d}_{\cc_{Y},\mathcal{Z}^{p-d}}(Y,p-d)$ by
$f_{\#}(\mathfrak{g}_{z})$.
\end{definition}

\nnpar{Covariant pseudo-morphisms.} In some cases, the notion of
covariant $f$-morphism is too restrictive. Moreover, it is not 
strictly necessary in order to define direct images of Green objects. 
Hence, we introduce the notion of covariant $f$-pseudo-morphisms.

\begin{definition} 
\label{def:21}
With the hypothesis of definition \ref{def:3}, a \emph{covariant
$f$-pseudo-morphism (of $\Gi$-complexes)} is a diagram of morphisms 
$f_{\#}$ of graded complexes of sheaves over $Y$ 
\begin{displaymath}
f_{\ast}\cc_{X}\overset{u}{\longleftarrow}\mathcal{F}\longrightarrow 
\cc_{Y}(-d)[-2d],
\end{displaymath}
where $\mathcal{F}$ is an auxiliary $\Gi$-complex satisfying the
following conditions:
\begin{enumerate}
\item[(1)]
The morphism $u$ is a quasi-isomorphism such that the morphism 
of global sections 
\begin{displaymath}
u:\Gamma(Y,\mathcal{F})\longrightarrow\Gamma(Y,f_{\ast}\cc_{X})=
\Gamma(X,\cc_{X})
\end{displaymath}
is an isomorphism. Then, as in definition \ref{def:3}, the 
pseudo-morphism $f_{\#}$ induces a morphism 
\begin{displaymath}
H^{n}_{\cc_{X},Z}(f^{-1}(U),p)\overset{f_{\#}}{\longrightarrow}
H^{n-2d}_{\cc_{Y},f(Z)} (U,p-d)
\end{displaymath}
for any subset $U\subset Y$, and any closed subset $Z\subset 
f^{-1}(U)$.  
\item[(2)]
The diagram 
\begin{displaymath}
\begin{CD}
H^{n}_{Z}(f^{-1}(U),\Gi(p))@>f_{!}>>H^{n-2d}_{f(Z)}(U,\Gi(p-d)) \\
@V\mathfrak{c}VV@VV\mathfrak{c}V \\
H^{n}_{\cc_{X},Z}(f^{-1}(U),p)@>f_{\#}>>H^{n-2d}_{\cc_{Y},f(Z)}
(U,p-d)
\end{CD}
\end{displaymath}
is commutative.
\end{enumerate}
\end{definition}

Observe that the definition of a covariant $f$-pseudo-morphism is
stronger than the notion of a morphism in the derived category,
because we need a well-defined morphism of groups on the level 
of global sections. Note also that a covariant $f$-morphism 
determines a covariant $f$-pseudo-morphism by taking for $u$ the 
identity map.

Since $\mathcal{F}$ and $f_{\ast}\cc_{X}$ satisfy the Mayer-Vietoris
principle, they are totally acyclic. Therefore, the fact that $u$ 
is a quasi-isomorphism of sheaves implies that the induced morphism 
\begin{displaymath}
\mathcal{F}(U,p)\longrightarrow\cc_{X}(f^{-1}(U),p)
\end{displaymath}
is a quasi-isomorphism for every open subset $U\subset Y$, and 
every $p\in\mathbb{Z}$. Taking into account that the induced 
morphism is an isomorphism on the level of global sections by 
corollary \ref{cor:tci}, we obtain

\begin{lemma}
Let $f:X\longrightarrow Y$ be a proper morphism of equidimensional 
$k$-schemes of relative dimension $d$, and let $f_{\#}$ be a
covariant $f$-pseudo-morphism. If $Z$ is a closed subset of $X$,
the induced morphism
\begin{displaymath}
\widehat{H}^{n}(\mathcal{F}(Y,p),\mathcal{F}(Y\setminus f(Z),p))
\longrightarrow\widehat{H}^{n}(\cc_{X}(X,p),\cc_{X}(X\setminus 
f^{-1}(f(Z)),p))
\end{displaymath}
is an isomorphism.
\hfill $\square$
\end{lemma}

\noindent
This result enables us to define direct images.

\begin{definition} 
\label{def:22}
Let $f:X\longrightarrow Y$ be a proper morphism of equidimensional 
$k$-schemes of relative dimension $d$, and let $f_{\#}$ be a 
covariant $f$-pseudo-morphism. If $Z$ is a closed subset of $X$, 
the induced morphism 
\begin{displaymath}
f_{\#}:\widehat{H}^{n}_{\cc_{X},Z}(X,p)\longrightarrow 
\widehat{H}^{n-2d}_{\cc_{Y},f(Z)}(Y,p-d)
\end{displaymath}
is defined by the composition
\begin{align*}
\widehat{H}^{n}_{\cc_{X},Z}(X,p)&=\widehat{H}^{n}(\cc_{X}(X,p),
\cc_{X}(X\setminus Z,p)) \\
&\longrightarrow\widehat{H}^{n}(\cc_{X}(X,p),\cc_{X}(X\setminus 
f^{-1}(f(Z)),p)) \\
&\overset{\cong}{\longleftarrow}\widehat{H}^{n}(\mathcal{F}(Y,p),
\mathcal{F}(Y\setminus f(Z),p)) \\
&\longrightarrow\widehat{H}^{n-2d}(\cc_{Y}(Y,p-d),\cc_{Y}(Y\setminus 
f(Z),p-d)) \\
&=\widehat{H}^{n-2d}_{\cc_{Y},f(Z)}(Y,p-d).
\end{align*}
\end{definition}

\begin{remark} 
\label{rem:4}
The analogues of theorem \ref{thm:9} with $U=Y$, and of proposition
\ref{prop:366} then also hold for covariant $f$-pseudo-morphisms. 
\end{remark}

\nnpar{Composition of pseudo-morphisms.} 
It is not clear from the definition whether we can always define 
the composition of covariant pseudo-morphisms. Nevertheless, we 
can determine when a pseudo-morphism should be considered as the 
composition of two pseudo-morphisms. 

For this, let $g:Y\longrightarrow Z$ be another proper morphism 
of equidimensional $k$-schemes of relative dimension $e$, and 
let $g_{\#}:g_{\ast}\cc_{Y}\overset{u'}{\longleftarrow}\mathcal{F}'
\longrightarrow\cc_{Z}(-e)[-2e]$ be a covariant $g$-pseudo-morphism.
Furthermore, let
\begin{displaymath}
h_{\#}:(g\circ f)_{\ast}\cc_{X}\overset{u''}{\longleftarrow}
\mathcal{F}''\longrightarrow\cc_{Z}(-d-e)[-2d-2e].
\end{displaymath}

\begin{definition}
We say that $h_{\#}$ is the composition of $f_{\#}$ and $g_{\#}$,
if there is a commutative diagram
\begin{displaymath}
\xymatrix{ 
(g\circ f)_{\ast}\cc_{X}&& \\
g_{\ast}\mathcal{F}\ar[u]_{g_{\ast}u}\ar[r]&g_{\ast}\cc_{Y}(-d)
[-2d]& \\
\mathcal{F}''\ar@(dl,ul)[uu]^{u''}\ar[u]\ar[r]\ar@(dl,dr)[rr]
&\mathcal{F}'(-d)[-2d]\ar[u]_{u'}\ar[r]&\cc_{Z}(-d-e)[-2d-2e]\,. \\
}
\end{displaymath}
\end{definition}

\begin{proposition} 
\label{prop:20}
Let $h_{\#}$ be a composition of the pseudo-morphisms $f_{\#}$ 
and $g_{\#}$. If $W$ is a closed subset of $X$, the induced 
morphisms on the level of Green objects   
\begin{align*}
f_{\#}&:\widehat{H}^{n}_{\cc_{X},W}(X,p)\longrightarrow
\widehat{H}^{n-2d}_{\cc_{Y},f(W)}(Y,p-d)\,, \\
g_{\#}&:\widehat{H}^{n-2d}_{\cc_{Y},f(W)}(Y,p-d)\longrightarrow
\widehat{H}^{n-2d-2e}_{\cc_{Z},g(f(W))}(Z,p-d-e)\,, \\
h_{\#}&:\widehat{H}^{n}_{\cc_{X},W}(X,p)\longrightarrow
\widehat{H}^{n-2d-2e}_{\cc_{Z},h(W)}(Z,p-d-e)
\end{align*}
satisfy $h_{\#}=g_{\#}\circ f_{\#}$.
\end{proposition}
\begin{proof}
We may assume that $W$ satisfies $W=h^{-1}(h(W))$. Writing
$U=X\setminus W$, $U'=Y\setminus f(W)$, $U''=Z\setminus 
h(W)$, and $n'=n-2d$, $n''=n-2d-2e$, we obtain the following 
commutative diagram (where we have omitted the corresponding 
twists for ease of notation)
\begin{displaymath}
\xymatrix{ 
\widehat{H}^{n}(\cc_{X}(X),\cc_{X}(U))&& \\
\widehat{H}^{n}(\mathcal{F}(Y),\mathcal{F}(U'))\ar[u]_{\cong}
\ar[r]&\widehat{H}^{n'}(\cc_{Y}(Y),\cc_{Y}(U'))& \\
\widehat{H}^{n}(\mathcal{F}''(Z),\mathcal{F}''(U''))
\ar@<1.2cm>@(l,l)[uu]_{\cong}\ar[u]\ar[r]\ar@(dl,dr)[rr]& 
\widehat{H}^{n'}(\mathcal{F}'(Z),\mathcal{F}'(U''))\ar[u]_
{\cong}\ar[r]&\widehat{H}^{n''}(\cc_{Z}(Z),\cc_{Z}(U''))\,.
}
\end{displaymath}
This proves the proposition.
\end{proof}

\nnpar{Projection formula.} 
In order to have a projection formula for Green objects, we need 
that the complexes under consideration also satisfy a projection
formula.

\begin{definition}
\label{def:3.76}
Let $f:X\longrightarrow Y$ be a proper morphism of equidimensional 
$k$-schemes of relative dimension $d$. Let $\cc_{X}$, $\cc_{X}'$,
$\cc''_{X}$ be $\Gi$-complexes over $X$, and let $\cc_{Y}$, $\cc'_
{Y}$, $\cc''_{Y}$ be $\Gi$-complexes over $Y$. Let     
\begin{displaymath}
f^{\#}:\cc_{Y}\longrightarrow f_{\ast}\cc_{X}
\end{displaymath}
be a contravariant $f$-morphism, and let
\begin{displaymath}
f'_{\#}:f_{\ast}\cc'_{X}\longrightarrow\cc'_{Y}(-d)[-2d]\,,\qquad
f''_{\#}:f_{\ast}\cc''_{X}\longrightarrow\cc''_{Y}(-d)[-2d]
\end{displaymath}
be covariant $f$-morphisms. Finally, let
\begin{displaymath}
\cc_{X}\otimes\cc'_{X}\overset{\bullet_{X}}{\longrightarrow}
\cc''_{X}\,,\qquad
\cc_{Y}\otimes\cc'_{Y}\overset{\bullet_{Y}}{\longrightarrow}
\cc''_{Y}
\end{displaymath} 
be $\Gi$-pairings. We call $(f^{\#},f'_{\#},f''_{\#},\bullet_{X},
\bullet_{Y})$ a \emph{projection five-tuple}, if the equality    
\begin{displaymath}      
f''_{\#}(f^{\#}(y)\bullet_{X}x)=y\bullet_{Y}f'_{\#}(x) 
\end{displaymath}
holds for all sections $x$, $y$ of $f_{*}\cc'_{X}$, $\cc_{Y}$,
respectively.
\end{definition}

\begin{proposition}
\label{prop:12}
Let $f:X\longrightarrow Y$ be a proper morphism of equidimensional 
$k$-schemes of relative dimension $d$, and let $(f^{\#},f'_{\#},
f''_{\#},\bullet_{X},\bullet_{Y})$ be a projection  five-tuple. 
Then, for any pair of Green objects $\mathfrak{g}_{1}$, $\mathfrak
{g}_{2}$ the equality  
\begin{displaymath}     
f''_{\#}(f^{\#}(\mathfrak{g}_{1})*_{X}\mathfrak{g}_{2})=     
\mathfrak{g}_{1}*_{Y}f'_{\#}(\mathfrak{g}_{2}) 
\end{displaymath}
holds.
\hfill$\square$
\end{proposition}

\nnpar{Pseudo-morphism and the projection formula.} 

\begin{remark} 
\label{rem:5}
Assume that $f'_{\#}$ and $f''_{\#}$ in definition \ref{def:3.76}  
are only covariant $f$-pseudo-morphisms given by
\begin{align*}
&f_{\ast}\cc'_{X}\overset{u'}{\longleftarrow}\mathcal{F}'\overset{v'}
{\longrightarrow}\cc'_{Y}(-d)[-2d]\,, \\
&f_{\ast}\cc''_{X}\overset{u''}{\longleftarrow}\mathcal{F}''\overset{v''}
{\longrightarrow}\cc''_{Y}(-d)[-2d]\,,
\end{align*}
respectively.
In order to have a projection formula on the level of Green
objects in this case, we need the following three $\Gi$-pairings
\begin{displaymath}
\cc_{X}\otimes\cc'_{X}\overset{\bullet_{X}}{\longrightarrow}
\cc''_{X}\,,\qquad\cc_{Y}\otimes\mathcal{F}'\overset{\bullet_
{\mathcal{F}}}{\longrightarrow}\mathcal{F}''\,,\qquad\cc_{Y}
\otimes\cc'_{Y}\overset{\bullet_{Y}}{\longrightarrow}\cc''_{Y}
\end{displaymath} 
satisfying
\begin{displaymath}
u''(y\bullet_{\mathcal{F}}x)=f^{\#}(y)\bullet_{X}u'(x)\,,
\qquad v''(y\bullet_{\mathcal{F}}x)=y\bullet_{Y}v'(x)
\end{displaymath}
for all sections $x$, $y$ of $\mathcal{F}'$, $\cc_{Y}$,
respectively.
\end{remark}

\newpage
\section{Abstract arithmetic Chow groups}
\label{sec:AC}

The idea behind the definition of arithmetic Chow groups for a 
variety $X$ over the ring of integers of a number field is that 
the variety in question can be ``compactified'' by adding the 
complex variety $X_{\infty}=X\times\mathbb{C}$, or more precisely, 
by adding a certain cohomology theory on the complex variety 
$X_{\infty}$. For instance, the cohomology theory involved in 
the definition of arithmetic Chow groups by H. Gillet and C. 
Soul\'e is the real Deligne-Beilinson cohomology (see 
\cite{Burgos:CDB}). However, the use of secondary characteristic
invariants (Green's forms, Bott-Chern forms) implies that the
arithmetic Chow groups depend not only on the cohomology theory, 
but also on the particular complex used to compute the cohomology. 
Therefore, the properties of this complex are reflected in the 
properties of the arithmetic Chow groups. The objective of this
section is to develop an abstract version of the construction of 
arithmetic Chow rings given in \cite{Burgos:CDB} emphasizing how 
the properties of the complexes under consideration are reflected 
by the properties of the arithmetic Chow groups.

\subsection{Arithmetic varieties}
\label{sec:av}

In this and the next section we will fix frequently used
notations. Unless stated otherwise, we will adhere to the 
notations of \cite{GilletSoule:ait}. For more details 
we refer the reader to loc. cit., and also to \cite{Soule:lag},
\cite{Burgos:CDB}. 

\nnpar{Arithmetic rings.}
\begin{definition}
An \emph{arithmetic ring} is a triple $(A,\Sigma,F_{\infty})$
consisting of an excellent regular noetherian integral domain 
$A$, a finite non-empty set $\Sigma$ of monomorphisms $\sigma:
A\longrightarrow\mathbb{C}$, and an antilinear involution $F_
{\infty}:\mathbb{C}^{\Sigma}\longrightarrow\mathbb{C}^{\Sigma}$ 
of $\mathbb{C}$-algebras such that the diagram
\begin{displaymath}
\xymatrix
{&&\mathbb{C}^{\Sigma}\ar[dd]^{F_{\infty}} \\
A\ar[rru]^{\delta}\ar[rrd]_{\delta}&& \\
&&\mathbb{C}^{\Sigma}}
\end{displaymath}
commutes; here $\delta$ is induced by the set of monomorphisms 
in $\Sigma$.
\end{definition}

The main examples of arithmetic rings are subrings of a number 
field $F$ which have $F$ as fraction field, subrings of $\mathbb{R}$,
and $\mathbb{C}$ itself (see \cite{GilletSoule:ait}, 3.2.1, for 
details).

\nnpar{Arithmetic varieties.}
\begin{notation}
Let $(A,\Sigma,F_{\infty})$ be an arithmetic ring with fraction 
field $F$, and let $X$ be a scheme over $\Spec(A)$. We will write 
$X_{F}$ for the generic fiber of $X$. If $\sigma\in\Sigma $, we 
will write $X_{\sigma}=X\underset{\sigma}{\otimes}\mathbb{C}$, 
and $X_{\Sigma}=X\underset{A}{\otimes}\mathbb{C}^{\Sigma }$. We
denote by $X_{\infty}$ the analytic space 
\begin{displaymath}
X_{\infty}=X_{\Sigma}(\mathbb{C})=\coprod_{\sigma\in\Sigma}
X_{\sigma}(\mathbb{C}).
\end{displaymath}
We observe that the antilinear involution $F_{\infty}$ of
$\mathbb{C}^{\Sigma}$ induces an antilinear involution
$F_{\infty}$ of $X_{\Sigma}$. We denote the real variety 
$(X_{\Sigma},F_{\infty})$  by $X_{\mathbb{R}}$.   
\end{notation}

Let us fix an arithmetic ring $A$, and a Gillet complex $\Gi$
for schemes over the field of real numbers $\mathbb{R}$.

\begin{definition}
An \emph{arithmetic variety $X$ over $A$} is a regular scheme 
$X$, which is flat and quasi-projective over $\Spec(A)$. If $X$
is equidimensional, we mean by the \emph{dimension of $X$ over
$A$}, the relative dimension of the scheme $X$ over $\Spec(A)$.
If $A$ is fixed and clear from the context, we call $X$ simply 
an arithmetic variety.
\end{definition}

\begin{definition}
\label{def:20}
A \emph{$\Gi$-arithmetic variety over $A$} is a pair 
\begin{displaymath}
\widehat{X}=(X,\cc)
\end{displaymath}
consisting of an arithmetic variety $X$ over $A$ and a 
$\Gi$-complex $\cc$ on the real variety $X_{\mathbb{R}}$.   
\end{definition}

\begin{definition}
Let $\widehat{X}=(X,\cc_{X})$ and $\widehat{Y}=(Y,\cc_{Y})$ be
$\Gi$-arithmetic varieties over $A$. A \emph{morphism}   
\begin{displaymath}
f:\widehat{X}\longrightarrow\widehat{Y}
\end{displaymath} 
\emph{of $\Gi$-arithmetic varieties over $A$} is a pair $f=
(f,f^{\#})$, where $f:X\longrightarrow Y$ is a morphism of 
$A$-schemes, and $f^{\#}:\cc_{Y}\longrightarrow f_{*}\cc_{X}$ 
is a contravariant $f_{\mathbb{R}}$-morphism of $\Gi$-complexes.
\end{definition}

The class of arithmetic varieties over $A$ together with 
their morphisms forms a category. We denote the category of 
$\Gi$-arithmetic varieties over $A$ by $\mathfrak{A}_{A,\Gi}$. 
If the cohomology theory $\Gi$ is fixed, we will write simply
$\mathfrak{A}_{A}$ instead of $\mathfrak{A}_{A,\Gi}$, and we 
will call a $\Gi$-arithmetic variety simply an arithmetic 
variety.   

\begin{definition}
A \emph{covariant morphism of $\Gi$-arithmetic varieties}
\begin{displaymath}  
f:\widehat{X}=(X,\cc_{X})\longrightarrow\widehat{Y}=(Y,\cc_{Y})
\end{displaymath}
\emph{over $A$} is a pair $f=(f,f_{\#})$, where $f:X\longrightarrow 
Y$ is a proper morphism of equidimensional $A$-schemes of relative
dimension $d$, and $f_{\#}:f_{*}\cc_{X}\longrightarrow\cc_{Y}(-d)
[-2d]$ is a covariant $f_{\mathbb{R}}$-morphism of $\Gi$-complexes. 
We call $f$ a \emph{covariant pseudo-morphism of $\Gi$-arithmetic 
varieties over $A$}, if $f_{\#}$ is a covariant $f_{\mathbb{R}}
$-pseudo-morphism.  
\end{definition}

\begin{remark} 
Let $X$ be an arithmetic variety over $A$. If $\cc$, $\cc'$ are two 
$\Gi$-complexes on $X_{\mathbb{R}}$, and $f:\cc\longrightarrow\cc'$ 
is a morphism of $\Gi$-complexes, then $f$ can be seen as a morphism 
of $\Gi$-arithmetic varieties $(X,\cc')\longrightarrow(X,\cc)$, as 
well as a covariant morphism of $\Gi$-arithmetic varieties $(X,\cc)
\longrightarrow(X,\cc')$. Therefore, $f$ will enjoy the properties 
of both kinds of morphisms.   
\end{remark}

\subsection{Arithmetic Chow groups}
Throughout this section $\widehat{X}=(X,\cc)$ will be a 
$\Gi$-arithmetic variety over an arithmetic ring $A$. 

\nnpar{Arithmetic cycles.}
\begin{notation} 
A cycle $y\in{\rm Z}^{p}(X)$ on $X$ determines a cycle on 
$X_{\Sigma}$. Clearly, this cycle is invariant under the 
action of $F_{\infty}$. Therefore, it is a cycle on $X_
{\mathbb{R}}$, which will be denoted by $y_{\mathbb{R}}$.  
Analogously, any $K_{1}$-chain $f$ on $X$ determines a 
$K_{1}$-chain $f_{\mathbb{R}}$ on $X_{\mathbb{R}}$.  

Since the $\Gi$-complex $\cc$ depends only on $X_{\mathbb{R}}$, 
we will write as a shorthand
\begin{align*}
\cc^{\ast}(X,p)&=\cc^{\ast}(X_{\mathbb{R}},p), \\
H^{*}_{\cc,\mathcal{Z}^{p}}(X,p)&=H^{*}_{\cc,\mathcal{Z}^{p}}
(X_{\mathbb{R}},p), \\
\widehat{H}^{*}_{\cc,\mathcal{Z}^{p}}(X,p)&=\widehat{H}^{*}_
{\cc,\mathcal{Z}^{p}}(X_{\mathbb{R}},p),
\end{align*}
and similar notations for cohomology with supports.  

For any cycle $y\in{\rm Z}^{p}(X)$, we write $\cl(y)$ for the 
class $\cl_{\cc}(y_{\mathbb{R}})$. For any $K_{1}$-chain $f\in
R^{p-1}_{p}(X)$, we write $\cl(f)$ for the class $\cl_{\cc}
(f_{\mathbb{R}})$, and $\mathfrak{g}(f)=\mathfrak{g}(f_{\mathbb
{R}})$ (see definition \ref{def:8}). 
\end{notation}

\begin{definition}
We define the \emph{group of $p$-codimensional arithmetic cycles 
of $\widehat{X}=(X,\cc)$} by 
\begin{displaymath}
\widehat{{\rm Z}}^{p}({X},\cc)=\left\{(y,\mathfrak{g}_{y})\in
{\rm Z}^{p}(X)\oplus\widehat{H}^{2p}_{\cc,\mathcal{Z}^{p}}(X,p)
\,\Big|\,\cl(y)=\cl(\mathfrak{g}_{y})\right\}
\end{displaymath}
with the obvious group structure.
\end{definition}

\nnpar{Rational equivalence.} 
We now define rational equivalence in this set-up.  

\begin{definition} 
Let $f\in R^{p-1}_{p}(X)$ be a $K_{1}$-chain. We set 
\begin{displaymath}
\widehat{\dv}(f)=(\dv(f),\mathfrak{g}(f)). 
\end{displaymath}
By remark \ref{rem:1}, we have $\widehat{\dv}(f)\in\widehat
{{\rm Z}}^{p}(X,\cc)$. We define 
\begin{displaymath} 
\rata^{p}(X,\cc)=\left\{\widehat{\dv}(f)\,\Big|\, f\in
R^{p-1}_{p}(X)\right\}\subseteq\widehat{{\rm Z}}^{p}(X,\cc).
\end{displaymath} 
\end{definition}

\noindent
It is easy to see that $\rata^{p}(X,\cc)$ is a subgroup of
$\widehat{{\rm Z}}^{p}(X,\cc)$. 

\nnpar{Definition and basic properties.}
\begin{definition}
\label{def:gacg} 
The \emph{$p$-th arithmetic Chow group of $\widehat{X}=(X,\cc)$} 
is defined by
\begin{displaymath}
\cha^{p}(X,\cc)=\widehat{{\rm Z}}^{p}(X,\cc)\big/\rata^{p}
(X,\cc).
\end{displaymath}
The groups $\cha^{p}(X,\cc)$ will also be called \emph{$\cc
$-arithmetic Chow groups of $X$}. The class of the arithmetic
cycle $(y,\mathfrak{g}_{y})$ in $\cha^{p}(X,\cc)$ will be
denoted by $[y,\mathfrak{g}_{y}]$.
\end{definition}

\begin{notation}
There are well-defined maps
\begin{alignat*}{2}
\zeta&:\cha^{p}(X,\cc)\longrightarrow\CH^{p}(X),&&\qquad 
\zeta[y,\mathfrak{g}_{y}]=[y], \\  
\rho&:\CH^{p,p-1}(X)\longrightarrow H^{2p-1}_{\cc}(X,p)
\subseteq\widetilde{\cc}^{2p-1}(X,p),&&\qquad\rho[f]= 
\cl(f), \\
\amap&:\widetilde{\cc}^{2p-1}(X,p)\longrightarrow\cha^{p}
(X,\cc),&&\qquad\amap(\widetilde{a})=[0,\amap(\widetilde{a})], \\
\omega&:\cha^{p}(X,\cc)\longrightarrow{\rm Z}\cc^{2p}(X,p), 
&&\qquad\omega[y,\mathfrak{g}_{y}]=\omega(\mathfrak{g}_{y}), \\
h&:{\rm Z}\cc^{2p}(X,p)\longrightarrow H^{2p}_{\cc}(X,p), 
&&\qquad h(\alpha)=[\alpha]. 
\end{alignat*}  
We will also write 
\begin{align*}
\cha^{p}(X,\cc)_{0}=&\Ker(\omega:\cha^{p}(X,\cc)\longrightarrow 
{\rm Z}\cc^{2p}(X,p)), \\ 
\CH^{p}(X)_{0}=&\Ker(\cl:\CH^{p}(X)\longrightarrow H^{2p}_{\cc}
(X,p)). 
\end{align*}
\end{notation}

\begin{theorem}
\label{thm:15}
Using the notations following definition \ref{def.311}, we have
the following exact sequences:
\begin{align}
&\CH^{p-1,p}(X)\overset{\rho}{\longrightarrow}\widetilde{\cc}^
{2p-1}(X,p)^{\pure}\overset{\amap}{\longrightarrow}\cha^{p}(X,\cc) 
\overset{\zeta}{\longrightarrow}\CH^{p}(X)\longrightarrow 0,
\label{eq:8} \\[3mm]
&\CH^{p-1,p}(X)\overset{\rho}{\longrightarrow}H^{2p-1}_{\cc}
(X,p)^{\pure}\overset{\amap}{\longrightarrow}\cha^{p}(X,\cc)
\overset{(\zeta,-\omega)}{\longrightarrow}\notag \\ 
&\phantom{CH^{p-1,p}(X)\overset{\rho}{\longrightarrow}}\CH^{p}(X) 
\oplus{\rm Z}\cc^{2p}(X,p)\overset{\cl+h}{\longrightarrow} H^{2p}_
{\cc}(X,p)\longrightarrow 0,
\label{eq:9} \\[3mm]
&\CH^{p-1,p}(X)\overset{\rho}{\longrightarrow}H^{2p-1}_{\cc}(X,p)^
{\pure}\overset{\amap}{\longrightarrow}\cha^{p}(X,\cc)_{0}\overset
{\zeta}{\longrightarrow}\CH^{p}(X)_{0}\longrightarrow 0.
\label{eq:10}
\end{align}
\end{theorem} 
\begin{proof} 
Let us prove the exactness of the first sequence. By proposition
\ref{prop:exgo}, for any $p$-codimensional cycle $y$ on $X$ there
exists a Green object $\mathfrak{g}_{y}$ for the class of $y$.
Therefore, $\zeta$ is surjective.

It is clear that $\zeta\circ\amap=0$. On the other hand, assume 
that $\zeta(\alpha)=0$, and let $(y,\mathfrak{g}_{y})$ be any 
representative of $\alpha$. Then, we have $y=\dv(f)$ for some 
$K_{1}$-chain $f$, and we obtain $(y,\mathfrak{g}_{y})\sim(0,
\mathfrak{g}_{y}-\mathfrak{g}(f))$. Since $\mathfrak{g}_{y}-
\mathfrak{g}(f)$ is a Green object for the trivial cycle, we have
$\cl(\mathfrak{g}_{y}-\mathfrak{g}(f))=0$. By the exact sequence
\eqref{exseq:go4}, we find $\mathfrak{g}_{y}-\mathfrak{g}(f)\in
\amap(\widetilde{\cc}^{2p-1}(X,p)^{\pure})$, i.e., $\alpha\in
\Img(\amap)$.  

Recalling $\mathfrak{g}(f)=\bmap(\cl(f))$ together with equation 
\eqref{rem:abf}, we find $\mathfrak{g}(f)=\amap(\cl(f))$, which 
shows $\amap\circ\rho=0$. On the other hand, if $\amap(\widetilde
{a})=0$, we have $(0,\amap(\widetilde{a}))=\diva(f)$ for some 
$K_{1}$-chain $f$. Therefore, we have $\dv(f)=0$, and $\mathfrak
{g}(f)=\amap(\widetilde{a})$. Since $\dv(f)=0$, we have $\cl(f)
\in H^{2p-1}_{\cc}(X,p)$, from which we derive as before $\mathfrak
{g}(f)=\bmap(\cl(f))=\amap(\cl(f))$. Therefore, $\widetilde{a}-
\cl(f)$ lies in the kernel of $\amap$. By the exact sequence 
\eqref{exseq:go4}, this implies $\widetilde{a}=\cl(f)=\rho[f]$.      

The proof of the other two exact sequences follows  the same 
lines, but uses the exact sequence \eqref{exseq:go5} instead.
\end{proof}

Note that at this abstract level, the above theorem is a
formal consequence of the definitions.

\subsection{Arithmetic intersection pairing}
\label{sec:aip}

The aim of this section is to define an arithmetic intersection 
pairing. We will follow the strategy of \cite{GilletSoule:ait}. 
Therefore, we divide out by finite rational equivalence before 
taking Green objects into account. 

\nnpar{Finite support.}
\begin{definition}
Let $\fin$, resp. $\tau(p)$ be the following families of supports
\begin{align*}
\fin&=\{Y\subset X\mid Y\text{closed},\,Y\cap X_{F}=\emptyset\},\,
\text{resp.} \\
\tau(p)&=\{Y\subset X\mid Y\text{closed},\codim(Y\cap X_{F})\ge p\}. 
\end{align*}
\end{definition}

\noindent
With these notations, we have the exact sequences 
\begin{align*}
&\bigoplus_{\substack{x\in X^{(p-1)} \\ \overline {\{x\}}\cap 
X_{F}=\emptyset}}k(x)^{*}\overset{\dv}{\longrightarrow}{\rm Z}^
{p}(X)\longrightarrow{\rm Z}^{p}(X_{F})\oplus\CH^{p}_{\fin}(X)
\longrightarrow 0, \\
&\bigoplus_{x\in X_{F}^{(p-1)}}k(x)^{*}\overset{\dv}{\longrightarrow} 
{\rm Z}^{p}(X_{F})\oplus\CH^{p}_{\fin}(X)\longrightarrow\CH^{p}(X)
\longrightarrow 0,
\end{align*}
and a canonical isomorphism
\begin{displaymath}
\CH^{p}_{\tau(p)}(X)\cong{\rm Z}^{p}(X_{F})\oplus\CH^{p}_{\fin}(X).
\end{displaymath}
Moreover, if $x\in X^{(p-1)}$ with $\overline{\{x\}}\cap X_{F}=
\emptyset$, and $f\in k(x)^{*}$, we have $\dv(f)\cap X_{F}=\emptyset$,
and $\mathfrak{g}(f)=0$. Therefore, we find
\begin{displaymath}
\widehat{\dv}(f)=(\dv(f),0).
\end{displaymath}
Thus, for any $\Gi$-arithmetic variety, there is a short exact sequence
\begin{displaymath}
\bigoplus_{x\in X_{F}^{(p-1)}}k(x)^{*}\overset{\widehat{\dv}}
{\longrightarrow}\widehat{{\rm Z}}^{p}(X_{F},\cc)\oplus\CH^{p}_
{\fin}(X)\longrightarrow\cha^{p}(X,\cc)\longrightarrow 0.    
\end{displaymath}
Observe that the notation $\widehat{{\rm Z}}^{p}(X_{F},\cc)$ makes
sense, since $(X_{F},\cc)$ is a $\Gi$-arithmetic variety over the 
arithmetic ring $F$.

\nnpar{The intersection pairing.} 
Let $y$ be a $p$-codimensional, and $z$ a $q$-codimensional cycle 
of $X$. Put $Y=\supp y$, and $Z=\supp z$. Assume that $y_{F}$ and 
$z_{F}$ intersect properly; this shows $Y\cap Z\in\tau(p+q)$. By 
theorem \ref{thm:intreg}, there is a well-defined class  
\begin{displaymath}
[y]\cdot[z]\in\CH^{p+q}_{Y\cap Z}(X)_{\mathbb{Q}},
\end{displaymath}
and hence a well-defined class
\begin{displaymath}
[y]\cdot[z]\in\CH^{p+q}_{\tau(p+q)}(X)_{\mathbb{Q}}\cong{\rm Z}^
{p+q}(X_{F})_{\mathbb{Q}}\oplus\CH^{p+q}_{\fin}(X)_{\mathbb{Q}}. 
\end{displaymath}
Let $(X,\cc)$, $(X,\cc')$, $(X,\cc'')$ be three $\Gi$-arithmetic 
varieties with the same underlying scheme $X$, and let $\bullet:
\cc\otimes\cc'\longrightarrow\cc''$ be a $\Gi$-pairing. Let 
$\mathfrak{g}_{y}$ be a Green object for the class of $y$ with 
values in $\cc$, and $\mathfrak{g}_{z}$ a Green object for the 
class of $z$ with values in $\cc'$. Then, the $*$-product 
$\mathfrak{g}_{y}*\mathfrak{g}_{z}$ is a Green object for 
the class of $y_{F}\cdot z_{F}$ with values in $\cc''$ (see 
definition \ref{def:5}), and we write 
\begin{equation}
\label{eq:6}
(y,\mathfrak{g}_{y})\cdot(z,\mathfrak{g}_{z})=((y_{F}\cdot z_{F},
\mathfrak{g}_{y}*\mathfrak{g}_{z}),[y\cdot z]_{\fin})\in\widehat
{{\rm Z}}^{p+q}(X_{F},\cc'')_{\mathbb{Q}}\oplus\CH^{p+q}_{\fin}(X)_
{\mathbb{Q}}.
\end{equation}
We will denote the image of the product $(y,\mathfrak{g}_{y})\cdot
(z,\mathfrak{g}_{z})$ in $\cha^{p+q}(X,\cc'')_{\mathbb{Q}}$ by
$[(y,\mathfrak{g}_{y})\cdot(z,\mathfrak{g}_{z})]$.

\begin{theorem}
\label{thm:2}
Let $(X,\cc)$, $(X,\cc')$, $(X,\cc'')$ be three $\Gi$-arithmetic 
varieties with the same underlying scheme $X$, and let $\bullet:
\cc\otimes\cc'\longrightarrow\cc''$ be a $\Gi$-pairing. Then, we 
have the following statements:
\begin{enumerate}
\item[(i)]
There exists a pairing
\begin{displaymath}
\cha^{p}(X,\cc)\otimes\cha^{q}(X,\cc')\overset{\cdot}{\longrightarrow}
\cha^{p+q}(X,\cc'')_{\mathbb{Q}}\,,
\end{displaymath}
which is determined by formula \eqref{eq:6} for cycles intersecting
properly in the generic fiber.
\item[(ii)] 
The following diagrams commute:
\begin{displaymath}
\begin{CD}
\cha^{p}(X,\cc)\otimes\cha^{q}(X,\cc')@>\cdot>>\cha^{p+q}(X,\cc'')_
{\mathbb{Q}} \\
@V\zeta\otimes\zeta VV@VV\zeta V \\
\CH^{p}(X)\otimes\CH^{q}(X)@>\cdot>>\CH^{p+q}(X)_{\mathbb{Q}}\,,
\end{CD}
\end{displaymath}
and 
\begin{displaymath}
\begin{CD}
\cha^{p}(X,\cc)\otimes\cha^{q}(X,\cc')@>\cdot>>\cha^{p+q}(X,\cc'')_
{\mathbb{Q}} \\
@V\omega\otimes\omega VV@VV\omega V \\
{\rm Z}\cc^{2p}(X,p)\otimes{\rm Z}\cc'{}^{2q}(X,q)@>\bullet>> 
{\rm Z}\cc''{}^{2p+2q}(X,p+q)_{\mathbb{Q}}.
\end{CD}
\end{displaymath} 
\end{enumerate}
\end{theorem}
\begin{proof}
Let $\alpha\in\cha^{p}(X,\cc)$, and $\beta\in\cha^{q}(X,\cc')$. 
By the moving lemma for cycles on a regular variety over a field, 
there are representatives $(y,\mathfrak{g}_{y})$ of $\alpha$, and 
$(z,\mathfrak{g}_{z})$ of $\beta $ such that $y_{F}$ and $z_{F}$ 
intersect properly. By means of formula \eqref{eq:6}, we define 
\begin{displaymath}
\alpha\cdot\beta=[(y,\mathfrak{g}_{y})\cdot(z,\mathfrak{g}_{z})].
\end{displaymath}
We have to show that this definition does not depend on the 
choice of representatives. Let $(y',\mathfrak{g}'_{y'})$ be 
another representative of $\alpha $ such that $y'_{F}$ and 
$z_{F}$ also intersect properly. Then, there exists a $K_{1}
$-chain $f$ such that 
\begin{displaymath}
\diva(f)=(y',\mathfrak{g}'_{y'})-(y,\mathfrak{g}_{y}).
\end{displaymath}
Therefore, it suffices to show that whenever $f$ is a $K_{1}$-chain 
whose divisor $\dv(f)_{F}$ intersects $z_{F}$ properly, then there 
exists a $K_{1}$-chain $g$ such that
\begin{displaymath}
\diva(g)=\diva(f)\cdot(z,\mathfrak{g}_{z})\in\za^{p+q}(X_{F},\cc'')_
{\mathbb{Q}}\oplus\CH^{p+q}_{\fin}(X)_{\mathbb{Q}}.
\end{displaymath}
Let $S=\supp\dv(f)$, $U=X\setminus S$, and $Z=\supp z$. By 
corollary \ref{cor:2} and proposition \ref{prop:10}, there 
exists a $K_{1}$-chain $g\in R^{p+q-1}_{p+q}(X)_{\mathbb{Q}}$ 
satisfying
\begin{align}
[\dv(g)]&=[\dv(f)]\cdot[z]\in\CH^{p+q}_{S\cap Z}(X)_{\mathbb{Q}},
\label{eq:4} \\
\cl(g)&=\cl(f)\bullet\cl(z)\in H^{2p+2q-1}_{\cc'',U\cap Z}
(X\setminus(S\cap Z),p+q)_{\mathbb{Q}}.
\label{eq:5}
\end{align}
Since $S\cap Z$ lies in $\tau(p+q)$, equation \eqref{eq:4} is also
valid in $\CH_{\tau(p+q)}^{p+q}(X)_{\mathbb{Q}}$ which shows
\begin{displaymath} 
[\dv(g)]=[\dv(f)]\cdot[z]=(\dv(f)_{F}\cdot z_{F},[\dv(f)\cdot
z]_{\fin}).
\end{displaymath}
Furthermore, we derive from equation \eqref{eq:5} and proposition 
\ref{prop:restrlem}
\begin{displaymath}
\mathfrak{g}(g)=\bmap(\cl(g))=\bmap(\cl(f)\bullet\cl(z))=\bmap
(\cl(f))*\mathfrak{g}_{z}=\mathfrak{g}(f)*\mathfrak{g}_{z}.
\end{displaymath}
Summing up, we find the equality
\begin{displaymath}
\diva(g)=\left((\dv(f)_{F}\cdot z_{F},\mathfrak{g}(f)*\mathfrak
{g}_{z}),[\dv(f)\cdot z]_{\fin}\right)=\diva(f)\cdot(z,\mathfrak
{g}_{z}) 
\end{displaymath}
in the group $\za^{p+q}(X_{F},\cc'')_{\mathbb{Q}}\oplus\CH^{p+q}_
{\fin}(X)_{\mathbb{Q}}$, as desired. 

The compatibility of the product with the morphisms $\zeta$ and
$\omega $ follows directly from the definitions. 
\end{proof}

\begin{remark}\label{rem:7}
  It may be interesting to explain more clearly why we do not need 
  the $K_1$-chain moving lemma. A close look at the proof of the
  well-definedness of the arithmetic intersection product in
  \cite{GilletSoule:ait}  shows that
  the $K_1$-chain moving lemma is used, first, to define a
  $K_{1}$-chain up to boundaries
  $f.[z]$ (which is the $g$ we are using in our proof) and second, to
  prove two 
  compatibility properties  
  \begin{align*}
    \dv(f.[z])&=\dv(f).[z],\\
    \log|f.[z]|^2&=\log|f|^2\land \delta _{z}.
  \end{align*}
  We observe that in the proof, the $K_{1}$-chain $f.[z]$ and the
  first compatibility condition is used on the whole arithmetic
  variety $X$ where no $K_{1}$-chain moving lemma is
  available. Therefore a result similar to corollary \ref{cor:2} is
  used implicitly. This means that the essential use of the moving lemma
  for $K_{1}$-chains is to prove the compatibility with Green
  currents: it is the concrete nature of the Green currents that
  forces the use of the $K_{1}$-chain moving lemma. 

  In our case the compatibility that we need is equation \eqref{eq:5},
  which follows directly from the properties of a Gillet
  cohomology. Roughly speaking, since the map that sends a
  $K_{1}$-chain to the associated Green object factors through a
  cohomology group, we are able to avoid the $K_{1}$-chain moving lemma.
\end{remark}

The fact that the product is defined only after tensoring with 
$\mathbb{Q}$ is due to the lack of a satisfactory intersection 
theory on general regular schemes, and the use of $K$-theory. 
In the cases in which an intersection product of cycles can be 
defined without using $K$-theory, it is possible to work over 
the integers. For instance, this is the case when $y$ and 
$z$ intersect properly on the whole of $X$, and not just on
$X_{F}$; then, the product $y\cdot z$ is a well-defined class 
in $\CH^{p+q}_{Y\cap Z}(X)$. Therefore, formula \eqref{eq:6} 
gives rise to an element of $\widehat{Z}^{p+q}(X,\cc'')\oplus
\CH^{p+q}(X)$. In particular, by the same argument as in 
\cite{GilletSoule:ait}, 4.2.3.iii, one can prove the following 
refinement of theorem \ref{thm:2}. 

\begin{theorem}
If $p=1$, or $q=1$, there is a unique pairing
\begin{displaymath}
\cha^{p}(X,\cc)\otimes\cha^{q}(X,\cc')\overset{\cdot}
{\longrightarrow}\cha^{p+q}(X,\cc'')
\end{displaymath}
given by formula \eqref{eq:6} when the cycles intersect 
properly. Moreover, this pairing induces the pairing of 
theorem \ref{thm:2}.
\hfill $\square$
\end{theorem}

\nnpar{$\Gi$-algebras.} 
Let us consider the case when $(X,\cc)$ is a $\Gi$-arithmetic variety
with $\cc$ being a $\Gi$-algebra.

\begin{theorem}
\label{thm:10} 
Let $(X,\cc)$ be a $\Gi$-arithmetic variety, and assume that 
$\cc$ is a pseudo-associative and pseudo-commutative 
$\Gi$-algebra. Then, we have the following statements:
\begin{enumerate}
\item[(i)] 
$\bigoplus_{p}\cha^{p}(X,\cc)_{\QQ}$ is a commutative and 
associative $\mathbb{Q}$-algebra with unit.
\item[(ii)]
For $\alpha\in\cha^{1}(X,\cc)$, $\beta\in\cha^{1}(X,\cc)$,
and $\gamma\in\cha^{p}(X,\cc)$, we have 
\begin{align*}
\alpha\cdot\gamma&=\gamma\cdot\alpha\in\cha^{p+1}(X,\cc), \\ 
(\alpha\cdot\beta )\cdot\gamma&=\alpha\cdot(\beta\cdot\gamma)
\in\cha^{p+2}(X,\cc).
\end{align*} 
\item[(iii)]
$\bigoplus_{p}\cha^{p}(X,\cc)_{0,\mathbb{Q}}=(\Ker\omega)_
{\mathbb{Q}}$ is an ideal of $\bigoplus_{p}\cha^{p}(X,\cc)_
{\mathbb{Q}}$.
\item[(iv)] 
$(\Ker\zeta)_{\mathbb{Q}}$ is an ideal of $\bigoplus_{p}
\cha^{p}(X,\cc)_{0,\mathbb{Q}}$.
\end{enumerate}
\end{theorem}
\begin{proof}
Since the product of algebraic cycles is commutative and 
associative, the first two statements of the theorem follow 
from the commutativity and associativity of the $*$-product 
of Green objects (see theorems \ref{thm:commst} and
\ref{thm:assst}). The last two statements are immediate.
\end{proof}

\begin{theorem}
\label{thm:17}
Let $(X,\cc)$ be a $\Gi$-arithmetic variety, and assume that 
$\cc$ is a pseudo-associative and pseudo-commutative 
$\Gi$-algebra. The $\bigoplus_{p}\cha^{p}(X,\cc)_
{\mathbb{Q}}$-module structure of $\bigoplus_{p}\cha^{p}
(X,\cc)_{0,\mathbb{Q}}$ is then induced by a $\bigoplus_
{p}\CH^{p}(X)_{\mathbb{Q}}$-module structure, i.e., 
there is a commutative diagram
\begin{displaymath}
\xymatrix{\cha^{p}(X,\cc)_{0}\otimes\cha^{q}(X,\cc)
\ar[r]\ar[d]_{\Id\otimes\zeta}&\cha^{p+q}(X,\cc)_
{0,\mathbb{Q}} \\
\cha^{p}(X,\cc)_{0}\otimes\CH^{q}(X).\ar[ur].}  
\end{displaymath}
\end{theorem}
\begin{proof}
We have to show that $\Ker(\zeta)\cdot\Img(\amap)=0$. This 
follows from proposition \ref{prop:restrlem} (i).
\end{proof}

\begin{corollary}
\label{cor:4}
With the assumptions of theorem \ref{thm:17} there is a 
well-defined product 
\begin{displaymath}
\CH^{p}(X)_{0}\otimes\CH^{q}(X)_{0}\overset{\cdot}
{\longrightarrow}\cha^{p+q}(X,\cc)_{0,\mathbb{Q}}.
\end{displaymath}
\end{corollary}
\begin{proof}
Let $x\in\CH^{p}(X)_{0}$ and $y\in\CH^{q}(X)_{0}$. We choose
arithmetic cycles $\widehat{x}\in\cha^{p}(X,\cc)_{0}$, resp. 
$\widehat{y}\in\cha^{q}(X,\cc)_{0}$ such that $\zeta(\widehat
{x})=x$, resp. $\zeta(\widehat{y})=y$. By theorem \ref{thm:17}, 
the product $\widehat{x}\cdot\widehat{y}$ does not depend on 
the choice of $\widehat{x}$ and $\widehat{y}$.
\end{proof}

\begin{remark}
If $\cc$ is only a $\Gi$-algebra, then the product in $\cha^{*}
(X,\cc)_{\mathbb{Q}}$ does not need to be associative or 
commutative. Nevertheless, theorems \ref{thm:assst} and 
\ref{thm:commst} imply that the product induced in $\bigoplus_
{p}\cha^{\ast}(X,\cc)_{0,\mathbb{Q}}$ is associative and 
commutative.
\end{remark}

\subsection{Inverse images}
\label{sec:ii}
Let $f:(X,\cc_{X})\longrightarrow(Y,\cc_{Y})$ be a morphism of 
$\Gi$-arithmetic varieties. In this section we will construct 
a pull-back morphism for arithmetic Chow groups. 

\nnpar{Definition and basic properties of inverse images.} 
Let $(z,\mathfrak{g}_{z})\in\za^{p}(Y,\cc_{Y})$ be a 
$p$-codimensional arithmetic cycle such that $f^{-1}(\supp
z_{\mathbb{R}})$ has codimension $p$ in $X_{\mathbb{R}}$. 
By theorem \ref{thm:5}, there is a well-defined cycle $f^
{*}(z)\in\CH^{p}_{\tau(p)}(X)$. Moreover, by definition 
\ref{def:6}, there is a well-defined Green object $f^{\#}
(\mathfrak{g}_{z})$ for the class of $f^{*}(z_{\mathbb{R}})$. 
In this case we will write   
\begin{equation}
\label{eq:11}
f^{*}(z,\mathfrak{g}_{z})=(f^{*}(z),f^{\#}(\mathfrak{g}_{z})).
\end{equation}

\begin{theorem}
\label{thm:6}
Let $f:(X,\cc_{X})\longrightarrow(Y,\cc_{Y})$ be a morphism 
of $\Gi$-arithmetic varieties. Then, the following statements 
hold:
\begin{enumerate}
\item[(i)]   
There is a well-defined morphism of graded groups  
\begin{displaymath}
f^{*}:\cha^{p}(Y,\cc_{Y})\longrightarrow\cha^{p}(X,\cc_{X})
\end{displaymath}
induced by equation \eqref{eq:11} for arithmetic cycles
$(z,\mathfrak{g}_{z})\in\za^{p}(Y,\cc_{Y})$ such that 
$f^{-1}(\supp z_{\mathbb{R}})$ has codimension $p$ in 
$X_{\mathbb{R}}$.
\item[(ii)] 
If $g:(Y,\cc_{Y})\longrightarrow(Z,\cc_{Z})$ is another
morphism of $\Gi$-arithmetic varieties, then the equality 
$(g\circ f)^{*}=f^{*}\circ g^{*}$ holds.  
\item[(iii)] 
The pull-back of arithmetic cycles satisfies the following   
relations: 
\begin{align*}
\omega\circ f^{*}&=f^{\#}\circ\omega, \\
\zeta\circ f^{*}&=f^{*}\circ\zeta, \\
f^{*}\circ\amap&=\amap\circ f^{\#}.
\end{align*}
\end{enumerate}
\end{theorem}
\begin{proof}
Let $\alpha\in\cha^{p}(Y,\cc_{Y})$ be an arithmetic cycle. 
By the moving lemma over $\mathbb{R}$, there is a representative 
$(z,\mathfrak{g}_{z})$ of $\alpha$ such that $f^{-1}(\supp z_
{\mathbb{R}})$ has codimension $p$ in $X_{\mathbb{R}}$. We
define $f^{*}(\alpha)$ as the class of $f^{*}(z,\mathfrak{g}_
{z})$. It remains to show that $f^{*}(\alpha)$ does not depend 
on the choice of the representative $(z,\mathfrak{g}_{z})$. 
In order to do this, let $(z',\mathfrak{g}_{z'})$ be another 
representative, i.e.,
\begin{equation}
\label{eq:12}
(z,\mathfrak{g}_{z})-(z',\mathfrak{g}_{z'})=\diva(h)
\end{equation}
for some $K_{1}$-chain $h$. Writing $U=Y\setminus(\supp z\cup
\supp z')$, the $K_{1}$-chain $h$ determines an element $[h]
\in\CH^{p,p-1}(U)$. By theorem \ref{thm:12} and proposition
\ref{prop:iicccg}, there exists an element $f^{*}[h]\in\CH^
{p,p-1}(f^{-1}(U))$ satisfying  
\begin{align*}
\dv(f^{*}[h])&=f^{*}(\dv(h))\in\CH^{p}_{\tau(p)}(X), \\
\cl(f^{*}[h])&=f^{\#}(\cl(h))\in H^{2p}_{\cc_{X}}(f^{-1}(U),p).
\end{align*}
These two equations then imply 
\begin{displaymath}
f^{*}(z,\mathfrak{g}_{z})-f^{*}(z',\mathfrak{g}_{z'})=
\diva(f^{*}(h)),
\end{displaymath}
which proves the first statement. The statements (ii), (iii)
are now shown easily.
\end{proof}

\nnpar{Multiplicativity of inverse images.} 
Since the pull-back of cycles is compatible with the 
intersection product of cycles by theorem \ref{thm:5}, 
we obtain the following result using proposition 
\ref{prop:11}.

\begin{theorem}
\label{thm:4.30}
Let $(f,f^{\#}_{\cc}):(X,\cc_{X})\longrightarrow(Y,\cc_{Y})$,
$(f,f^{\#}_{\cc'}):(X,\cc'_{X})\longrightarrow(Y,\cc'_{Y})$,
$(f,f^{\#}_{\cc''}):(X,\cc''_{X})\longrightarrow(Y,\cc''_{Y})$ 
be morphisms of $\Gi$-arithmetic varieties, and let $\cc_{X}\otimes
\cc'_{X}\overset{\bullet_{X}}{\longrightarrow}\cc''_{X}$, 
$\cc_{Y}\otimes\cc'_{Y}\overset{\bullet_{Y}}{\longrightarrow}
\cc''_{Y}$ be $\Gi$-pairings which are compatible with the 
$f$-morphisms $f^{\#}_{\cc}$, $f^{\#}_{\cc'}$ and $f^{\#}_
{\cc''}$. Then, the diagram    
\begin{displaymath}
\begin{CD}
\cha^{p}(Y,\cc_{Y})\otimes\cha^{q}(Y,\cc'_{Y})@>\cdot>>
\cha^{p+q}(Y,\cc''_{Y})_{\mathbb{Q}} \\
@Vf^{*}VV@VV f^{*}V \\
\cha^{p}(X,\cc_{X})\otimes\cha^{q}(X,\cc'_{X})@>\cdot>>
\cha^{p+q}(X,\cc''_{X})_{\mathbb{Q}}
\end{CD}
\end{displaymath}
is commutative. If $p=1$, or $q=1$ the same result holds
true without tensoring with $\mathbb{Q}$. 
\hfill$\square$
\end{theorem}

\nnpar{Change of $\Gi$-complex}. 
The results for inverse images apply in particular to the 
case of a change of complexes on $X$, i.e., when there is a 
morphism of $\Gi$-arithmetic varieties 
\begin{displaymath}
f=(\Id,f^{\#}):(X,\cc')\longrightarrow(X,\cc).
\end{displaymath}
Moreover, a direct consequence of proposition \ref{prop:21} 
is the following result.

\begin{proposition}
\label{prop:22}
Let $f=(\Id,f^{\#}):(X,\cc')\longrightarrow(X,\cc)$ be a
morphism of $\Gi$-arithmetic varieties such that the morphism   
\begin{displaymath}
f^{\#}:\cc(X,p)\longrightarrow\cc'(X,p)
\end{displaymath}
is an isomorphism for any integer $p$, and such that the 
morphism
\begin{displaymath}
f^{\#}:\cc(U,p)\longrightarrow\cc'(U,p)
\end{displaymath}
is a quasi-isomorphism for any open subset $U\subset X$
and any integer $p$. Then, the induced morphism
\begin{displaymath}
f^{*}:\cha^{p}(X,\cc)\longrightarrow\cha^{p}(X,\cc')
\end{displaymath}
is an isomorphism for any integer $p$.
\hfill$\square$
\end{proposition}

\subsection{Proper push-forward} 
\label{sec:pp}
Let $f:(X,\cc_{X})\longrightarrow(Y,\cc_{Y})$ be a covariant 
morphism or, more generally, a covariant pseudo-morphism of 
$\Gi$-arithmetic varieties of relative dimension $d$. In this 
section we will construct a push-forward morphism for arithmetic 
Chow groups. 

\nnpar{Push-forward.} 
Let $(z,\mathfrak{g}_{z})\in\za^{p}(X,\cc_{X})$ be a 
$p$-codimensional arithmetic cycle. By theorem \ref{thm:13},
there is a well-defined cycle $f_{*}(z)\in{\rm Z}^{p-d}(Y)$. 
Moreover, there is a well-defined Green object $f_{\#}
(\mathfrak{g}_{z})$ for the class of $f_{*}(z_{\mathbb{R}})$ 
(see definition \ref{def:3}). We will write  
\begin{equation}
\label{eq:22}
f_{*}(z,\mathfrak{g}_{z})=(f_{*}(z),f_{\#}(\mathfrak{g}_{z})).
\end{equation}

\begin{theorem} 
\label{thm:18}
Let $f:(X,\cc_{X})\longrightarrow(Y,\cc_{Y})$ be a covariant 
morphism or a covariant pseudo-morphism of $\Gi$-arithmetic 
varieties of relative dimension $d$. Then, the following 
statements hold:
\begin{enumerate}
\item[(i)]   
There is a well-defined morphism of graded groups  
\begin{displaymath}
f_{*}:\cha^{p}(X,\cc_{X})\longrightarrow\cha^{p-d}(Y,\cc_{Y}) 
\end{displaymath}
induced by equation \eqref{eq:22} for arithmetic cycles
$(z,\mathfrak{g}_{z})\in\za^{p}(X,\cc_{X})$. 
\item[(ii)] If $g:(Y,\cc_{Y})\longrightarrow(Z,\cc_{Z})$ is 
another covariant morphism or covariant pseudo-morphism of 
$\Gi$-arithmetic varieties, then the equality $(g\circ f)_
{*}=g_{*}\circ f_{*}$ holds.
\item[(iii)]
The push-forward of arithmetic cycles satisfies the following
relations: 
\begin{align*}
\omega\circ f_{*}&=f_{\#}\circ\omega, \\
\zeta\circ f_{*}&=f_{*}\circ\zeta, \\
f_{*}\circ\amap&=\amap\circ f_{\#}.
\end{align*}
\end{enumerate}
\end{theorem}
\begin{proof}
We will treat here only the case of a covariant morphism 
of $\Gi$-arith\-metic varieties; the case of a covariant 
pseudo-morphism can be treated analogously.

The only thing that remains to be shown is the compatibility 
with rational equivalence. For this purpose, let $h\in R^{p-1}_
{p}(X)$ be a $K_{1}$-chain. By theorem \ref{thm:13} and 
proposition \ref{prop:dicccg}, there is a well-defined $K_{1}
$-chain $f_{*}(h)\in R^{p-d-1}_{p-d}(Y)$ satisfying 
\begin{align}
\dv(f_{*}(h))&=f_{*}(\dv(h)), \\
\cl_{\Gi}(f_{*}(h))&=f_{\#}(\cl_{\Gi}(h)).
\label{eq:23}
\end{align}
By the definition of a covariant morphism of $\Gi$-arithmetic 
varieties, equation \eqref{eq:23} implies that $ \cl_{\cc_{Y}}
(f_{*}(h))=f_{\#}(\cl_{\cc_{X}}(h))$. Therefore, we have 
$\mathfrak{g}(f_{*}(h))=f_{\#}(\mathfrak{g}(h))$, and consequently   
\begin{displaymath}     
f_{*}(\diva(h))=\diva(f_{*}(h)).
\end{displaymath}
\end{proof}

\nnpar{Projection formula.} 
The next result is a direct consequence of the projection 
formula for algebraic cycles and proposition \ref{prop:12}.

\begin{proposition}
Let $(f,f^{\#}):(X,\cc_{X})\longrightarrow(Y,\cc_{Y})$ be 
a morphism of $\Gi$-arithmetic varieties, and let 
\begin{align*}
(f,f'_{\#}):(X,\cc'_{X})&\longrightarrow(Y,\cc'_{Y}), \\
(f,f''_{\#}):(X,\cc''_{X})&\longrightarrow(Y,\cc''_{Y})
\end{align*}
be covariant morphisms of $\Gi$-arithmetic varieties. Let
\begin{displaymath}
\cc_{X}\otimes\cc'_{X}\overset{\bullet_{X}}{\longrightarrow}
\cc''_{X}\qquad\cc_{Y}\otimes\cc'_{Y}\overset{\bullet_{Y}}
{\longrightarrow}\cc''_{Y}
\end{displaymath}
be $\Gi$-pairings such that $(f^{\#},f'_{\#},f''_{\#},
\bullet_{X},\bullet_{Y})$ is a projection five-tuple. Then, 
the projection formula
\begin{displaymath}
f''_{*}(f^{*}(\alpha)\cdot\beta)=\alpha\cdot f'_{*}(\beta)
\end{displaymath}
holds for all $\alpha\in\cha^{p}(Y,\cc_{Y})$ and $\beta\in
\cha^{q}(X,\cc'_{X})$.
\hfill$\square$
\end{proposition}

\begin{remark}
An analogous result also holds when $f'_{\#}$ and 
$f''_{\#}$ are covariant pseudo-morphisms. We leave it 
to the reader to make this result explicit.
\end{remark}

\newpage
\section{Deligne-Beilinson cohomology as a Gillet cohomology}
\label{sec:deligne-beil-cohom}

A particular example for a Gillet cohomology is Deligne-Beilinson 
cohomology for real and complex algebraic varieties. In this 
section we recall the definition and some facts on Deligne-Beilinson 
cohomology and homology which we will use in the sequel. 
The main references for this section are \cite{Beilinson:hr}, 
\cite{EsnaultViehweg:DBc}, and \cite{Jannsen:DcHD}. We start 
by recalling some properties of general Deligne-Beilinson 
cohomology; later we will shift to real Deligne-Beilinson 
cohomology. We will show that it satisfies most of the properties 
of a Gillet cohomology. Moreover, we will construct an explicit
Gillet complex for real Deligne-Beilinson cohomology. Note that the
results of this chapter are well known and we include them for the
convenience of the reader.

\begin{notation}
By a complex algebraic manifold we mean the analytic variety
associated to a smooth separated scheme of finite type over 
$\mathbb{C}$.   
\end{notation}

\subsection{Review of Deligne-Beilinson cohomology}

\nnpar{The definition of Deligne-Beilinson cohomology.}
Let $X$ be a complex algebraic manifold. Let $\Lambda$ be a subring 
of $\mathbb{R}$, and set $\Lambda(p)=(2\pi i)^{p}\Lambda\subseteq 
\mathbb{C}$. We will denote the corresponding constant sheaves on 
$X$ also by $\Lambda$, resp. $\Lambda(p)$.

Let us choose a smooth compactification $\overline{X}$ of $X$ 
with $D=\overline{X}\setminus X$ a normal crossing divisor, and 
denote by $j:X\longrightarrow\overline{X}$ the natural inclusion.
Let $\Omega^{\ast}_{X}$ be the sheaf of holomorphic differential 
forms on $X$, and let $\Omega^{\ast}_{\overline{X}}(\log D)$ be 
the sheaf of holomorphic differential forms on $\overline{X}$ with
logarithmic poles along $D$ (see \cite{Deligne:THII}). Let $F$ 
be the Hodge filtration of $\Omega^{\ast}_{\overline{X}}(\log D)$, 
i.e., 
\begin{displaymath}
F^{p}\Omega^{\ast}_{\overline{X}}(\log D)=\bigoplus_{p'\geq p}
\Omega^{p'}_{\overline{X}}(\log D). 
\end{displaymath}
Then, the Deligne-Beilinson complex of the pair $(X,\overline{X})$ 
is given by the simple complex
\begin{displaymath}
\Lambda(p)_{\mathcal{D}}=s\left(Rj_{\ast}\Lambda(p)\oplus F^{p}
\Omega^{\ast}_{\overline{X}}(\log D)\overset{u}{\longrightarrow}
j_{\ast}\Omega_{X}^{\ast}\right),
\end{displaymath}
where the morphism $u$ is defined by $u(a,f)=-a+f$.

\begin{definition}
The \emph{Deligne-Beilinson cohomology groups} are the hypercohomology 
groups of the sheaf $\Lambda(p)_{\mathcal{D}}$, i.e.,
\begin{displaymath}
H_{\mathcal{D}}^{\ast}(X,\Lambda(p))=\mathbb{H}^{\ast}(\overline{X}, 
\Lambda(p)_{\mathcal{D}}).
\end{displaymath}
\end{definition}

It can be shown that these groups are independent of the 
compactification $\overline{X}$ of $X$. The Deligne-Beilinson
cohomology groups can also be constructed as the hypercohomology
groups of a complex of sheaves of graded abelian groups in the 
Zariski topology; see \cite{EsnaultViehweg:DBc}, \S5, for details. 
We will denote the corresponding sheaf in the Zariski topology
by $\Lambda(\ast)_{\mathcal{D},\Zar}$; we will not need the precise 
definition of this sheaf for general $\Lambda$. In contrast, in 
definition \ref{def:18}, we will give a concrete construction in 
the case $\Lambda =\mathbb{R}$.   

\nnpar{An exact sequence.} 
The definition of Deligne-Beilinson cohomology as the cohomology 
of a simple complex associated to a morphism of complexes implies 
the existence of some exact sequences. One of these exact sequences 
relates Deligne-Beilinson cohomology to the usual cohomology and 
shows that the integral Deligne-Beilinson cohomology is an extension 
of the group of Hodge cycles by an intermediate Jacobian. We recall 
the construction of this exact sequence. There is a natural map 
\begin{equation}
\Lambda(p)_{\mathcal{D}}\longrightarrow Rj_{\ast}\Lambda(p)
\end{equation}
given by $(a,f,\omega)\longmapsto a$. The kernel of this morphism is
the simple 
\begin{displaymath}
s\left(F^{p}\Omega^{\ast}_{\overline{X}}(\log D)\longrightarrow
j_{\ast}\Omega_{X}^{\ast}\right).
\end{displaymath}
Since the differential is strict with respect to the Hodge filtration
(see \cite{Deligne:THII}), the cohomology of this simple complex is
given by $H^{\ast}(X,\mathbb{C})/F^{p}H^{\ast}(X,\mathbb{C})$. Therefore, 
there is a long exact sequence
\begin{equation}
\label{eq:25}
...\to\left.H^{n-1}(X,\mathbb{C})\right/F^{p}H^{n-1}(X,\mathbb{C})     
\to H_{\mathcal{D}}^{n}(X,\Lambda(p))\to H^{n}(X,\Lambda(p))\to...  
\end{equation}

\nnpar{The product in Deligne-Beilinson cohomology.}
Beilinson has introduced a product in this cohomology. More specifically, 
he has introduced a family of pairings indexed by the interval $[0,1]$. 
For $0\leq\alpha\leq 1$, the pairing $\cup_{\alpha}$ is given by 
\begin{align*}
(a,f,\omega)\cup_{\alpha}(a',f',\omega')=&(aa',f\land f',\alpha(\omega 
\land a'+(-1)^{n}f\land\omega')+ \\
&(1-\alpha)(\omega\land f'+(-1)^{n}a\land\omega')),
\end{align*}
where $(a,f,\omega)\in\Lambda(p)_{\mathcal{D}}^{n}$ and $(a',f',\omega')
\in\Lambda(q)_{\mathcal{D}}^{m}$. All these pairings turn out to be
homotopically equivalent. If $\alpha=0,1$, the pairing is associative 
and, if $\alpha =1/2$, the pairing is graded commutative. Therefore,
we obtain a well-defined commutative and associative pairing in
Deligne-Beilinson cohomology.

\nnpar{Deligne-Beilinson cohomology and the Picard group.}
The axiom 1.2 (xi) of \cite{Gillet:RRhK} states for Deligne-Beilinson
cohomology that there is a natural transformation of functors
\begin{displaymath}
\Pic(\cdot)\longrightarrow H^{2}_{\mathcal{D}}(\cdot,\Lambda (1)). 
\end{displaymath}
This natural transformation is realized by a morphism in the derived
category of the category of graded complexes of Zariski sheaves of 
abelian groups
\begin{equation}
\label{eq:29}
{\rm c}_{1}:\mathcal{O}^{\times}_{\alg}[-1]\longrightarrow
\Lambda(1)_{\mathcal{D},\Zar}\,;
\end{equation}
we refer to \cite{EsnaultViehweg:DBc}, \S5, for the precise construction. 
Here we only sketch a slightly different version in the case $\Lambda=
\mathbb{Z}$; the other cases then follow by functoriality. Since, for 
$p\ge 1$ and $n\le 0$, we have $H^{n}_{\mathcal{D}}(X,\mathbb{Z}(p))=0$, 
we can assume that the complex of sheaves $\mathbb{Z}(1)_{\mathcal{D},
\Zar}$ starts in degree $1$. Therefore, ${\rm c}_{1}$ is determined by 
a functorial isomorphism 
\begin{equation}
\label{eq:15}
{\rm c}_{1}:H^{0}(X,\mathcal{O}^{\times}_{\alg})\longrightarrow H^{1}_
{\mathcal{D}}(X,\mathbb{Z}(1)),
\end{equation}
which can be described as follows: Using the quasi-isomorphism
\begin{displaymath}
\mathbb{Z}(1)\longrightarrow s(\mathcal{O}_{X}\overset{\exp}
{\longrightarrow}\mathcal{O}^{\times}_{X}[-1]),
\end{displaymath}
one sees that $\mathbb{Z}(1)_{\mathcal{D}}$ is quasi-isomorphic to
the complex 
\begin{displaymath}
s\left(Rj_{\ast}\mathcal{O}^{\times}_{X}\oplus F^{1}\Omega^{\ast}_
{\overline{X}}(\log D)\overset{u'}{\longrightarrow }Rj_{\ast}(\Omega^
{\ast}_{X}/\mathcal{O}_{X})\right),
\end{displaymath}
where $u'(a,f)=-\dd\log a+f$. Hence, the group $H^{1}_{\mathcal{D}}
(X,\mathbb{Z}(1))$ is isomorphic to 
\begin{displaymath}
\{(f,\omega)\in H^{0}(X,\mathcal{O}^{\times}_{X})\oplus H^{0}
(\overline{X},\Omega^{1}_{\overline{X}}(\log D))\mid\dd\log
f=\omega\}.
\end{displaymath}
By the GAGA principle the morphism $H^{0}(X,\mathcal{O}^{\times}_
{\alg})\longrightarrow H^{1}_{\mathcal{D}}(X,\mathbb{Z}(1))$ given 
by $\omega\mapsto(\omega,\dd\log\omega)$ is an isomorphism.

\subsection{Review of Deligne algebras}
\label{sec:RDBc}
In this section we will recall some definitions and results from
\cite{Burgos:CDB} about Dolbeault algebras and their associated
Deligne algebras. The interest in Deligne algebras is explained 
by the fact that they are very simple objects which compute real
Deligne-Beilinson cohomology.

\nnpar{Dolbeault complexes.}
\begin{definition}
\label{def:12}
A \emph{Dolbeault complex} $A=(A^{\ast}_{\mathbb{R}},\dd_{A})$ is 
a graded complex of real vector spaces, which is bounded from below 
and equipped with a bigrading on $A_{\mathbb{C}}=A_{\mathbb{R}}
\otimes_{\mathbb{R}}{\mathbb{C}}$, i.e.,
\begin{displaymath}
A^{n}_{\mathbb{C}}=\bigoplus_{p+q=n}A^{p,q},
\end{displaymath}  
satisfying the following properties:
\begin{enumerate}
\item[(i)]
The differential $\dd_{A}$ can be decomposed as the sum $\dd_{A}=
\partial+\bar{\partial}$ of operators $\partial$ of type $(1,0)$, 
resp. $\bar{\partial}$ of type $(0,1)$.
\item[(ii)] 
It satisfies the symmetry property $\overline{A^{p,q}}=A^{q,p}$,
where $\overline{\phantom{M}}$ denotes complex conjugation.
\end{enumerate}
\end{definition}

\begin{notation}
\label{def:13}
Given a Dolbeault complex $A=(A^{\ast}_{\mathbb{R}},\dd_{A})$, we 
will use the following notations. The Hodge filtration $F$ of $A$ 
is the decreasing filtration of $A_{\mathbb{C}}$ given by
\begin{displaymath}
F^{p}A^{n}=F^{p}A^{n}_{\mathbb{C}}=\bigoplus_{p'\geq p}A^{p',n-p'}.
\end{displaymath}
The filtration $\overline F$ of $A$ is the complex conjugate of $F$, 
i.e.,
\begin{displaymath}
\overline{F}^{p}A^{n}=\overline{F}^{p}A^{n}_{\mathbb{C}}=\overline
{F^{p}A^{n}_{\mathbb{C}}}.
\end{displaymath}
For an element $x\in A_{\mathbb{C}}$, we write $x^{i,j}$ for its 
component in $A^{i,j}$. For $k,k'\geq 0$, we define an operator 
$F^{k,k'}:A_{\mathbb{C}}\longrightarrow A_{\mathbb{C}}$ by the 
rule 
\begin{displaymath}
F^{k,k'}(x):=\sum_{l\geq k,l'\geq k'}x^{l,l'}.
\end{displaymath}
We note that the operator $F^{k,k'}$ is the projection of $A^{\ast}_
{\mathbb{C}}$ onto the subspace $F^{k}A^{\ast}\cap\overline{F}^{k'}
A^{\ast}$. We will write $F^{k}=F^{k,-\infty}$. 

We denote by $A^{n}_{\mathbb{R}}(p)$ the subgroup $(2\pi i)^{p}\cdot
A^{n}_{\mathbb{R}}\subseteq A^{n}_{\mathbb{C}}$, and we define the 
operator
\begin{displaymath}
\pi_{p}:A_{\mathbb{C}}\longrightarrow A_{\mathbb{R}}(p)
\end{displaymath}
by setting $\pi_{p}(x):=\frac{1}{2}(x+(-1)^{p}\bar{x})$.
\end{notation}

\nnpar{The Deligne complex.} 
To any Dolbeault complex we can associate a Deligne complex.

\begin{definition}
Let $A$ be a Dolbeault complex. We denote by $A^{\ast}(p)_
{\mathcal{D}}$ the complex $s(A_{\mathbb{R}}(p)\oplus F^{p}A
\overset{u}{\longrightarrow}A_{\mathbb{C}})$, where $u(a,f)=
-a+f$.
\end{definition}

\begin{definition}
Let $A$ be a Dolbeault complex. Then, the \emph{Deligne complex 
$(\mathcal{D}^{\ast}(A,\ast),\dd_{\mathcal{D}})$ associated to $A$} 
is the graded complex given by   
\begin{align*}
&\mathcal{D}^{n}(A,p)=
\begin{cases}
A_{\mathbb{R}}^{n-1}(p-1)\cap F^{n-p,n-p}A^{n-1}_{\mathbb{C}},
&\qquad\text{if}\quad n\leq 2p-1, \\
A_{\mathbb{R}}^n(p)\cap F^{p,p}A^{n}_{\mathbb{C}},
&\qquad\text{if}\quad n\geq 2p,  
\end{cases}  
\intertext{with differential given by ($x\in\mathcal{D}^{n}(A,p)$)}
&\dd_{\mathcal{D}}x=
\begin{cases}
-F^{n-p+1,n-p+1}\dd_{A}x,
&\qquad\text{if}\quad n<2p-1, \\
-2\partial\bar{\partial}x, 
&\qquad\text{if}\quad n=2p-1, \\
\dd_{A}x, 
&\qquad\text{if}\quad n\geq 2p.
\end{cases}  
\end{align*}  
\end{definition}

For instance, if $A$ is a Dolbeault complex satisfying $A^{p,q}=0$ 
for $p<0$, $q<0$, $p>n$, or $q>n$, then the complex $\mathcal{D}(A,0)$ 
agrees with the real complex $A^{\ast}_{\mathbb{R}}$; for $p>0$, we
have represented $\mathcal{D}(A,p)$ in figure \ref{fig:1} below, where 
the lower left square is shifted by one; this means in particular that
$A^{0,0}$ sits in degree 1 and $A^{p-1,p-1}$ sits in degree $2p-1$. 

\begin{figure}[htb]
\begin{center}
\begin{displaymath}
\xymatrix{&
\makebox{$\begin{pmatrix}
A^{p,n}&\to&\cdots&\to&A^{n,n} \\
\uparrow &&&& \uparrow \\
\vdots &&&& \vdots \\
\uparrow &&&& \uparrow \\
A^{p,p}&\to&\cdots&\to&A^{n,p} 
\end{pmatrix}$}_{\makebox[0pt]{$\mathbb{R}$}}
\makebox[0pt][l]{$(p)$} \\
\makebox{$\begin{pmatrix}
A^{0,p-1}&\to&\cdots&\to&A^{p-1,p-1} \\
\uparrow &&&& \uparrow \\
\vdots &&&& \vdots \\
\uparrow &&&& \uparrow \\
A^{0,0}&\to&\cdots&\to&A^{p-1,0}
\end{pmatrix}$}_{\makebox[0pt]{$\mathbb{R}$}}
\makebox[0pt][l]{$(p-1)$}
\ar[ur]^{-2\partial\overline{\partial}}
&}
\end{displaymath}
\caption{$\mathcal{D}(A,p)$}
\label{fig:1}
\end{center}
\end{figure}

\begin{remark} 
It is clear from the definition that the functor $\mathcal{D}
(\cdot,p)$ is exact for all $p$.
\end{remark}

The main property of the Deligne complex is expressed by the 
following proposition; for a proof see \cite{Burgos:CDB}.

\begin{proposition}
\label{prop:32}
The complexes $A^{\ast}(p)_{\mathcal{D}}$ and $\mathcal{D}^{\ast}
(A,p)$ are homotopically equivalent. The homotopy equivalences
$\psi:A^{n}(p)_{\mathcal{D}}\longrightarrow\mathcal{D}^{n}(A,p)$,
and $\varphi:\mathcal{D}^{n}(A,p)\longrightarrow A^{n}(p)_{\mathcal 
{D}}$ are given by
\begin{displaymath}
\psi(a,f,\omega)=
\begin{cases}
\pi(\omega),\qquad&\text{if }n\le 2p-1, \\
F^{p,p}a+2\pi_{p}(\partial\omega^{p-1,n-p+1}),\quad&\text{if }n\ge 2p,
\end{cases}
\end{displaymath}  
where $\pi(\omega)=\pi_{p-1}(F^{n-p,n-p}\omega)$, i.e., $\pi$ is
the projection of $A_{\mathbb{C}}$ over the co\-kernel of $u$, and
\begin{displaymath}
\varphi(x)=
\begin{cases}
(\partial x^{p-1,n-p}-\bar{\partial}x^{n-p,p-1},2\partial
x^{p-1,n-p},x),\quad&\text{if }n\le 2p-1, \\
(x,x,0),&\text{if }n\ge 2p.
\end{cases}
\end{displaymath}
Moreover, $\psi\circ\varphi=\Id$, and $\varphi\circ\psi-\Id=\dd h+
h\dd$, where $h:A^{n}(p)_{\mathcal{D}}\longrightarrow A^{n-1}(p)_
{\mathcal{D}}$ is given by
\begin{displaymath}
h(a,f,\omega)=
\begin{cases}
(\pi_{p}(\overline{F}^{p}\omega+\overline{F}^{n-p}\omega),-2F^{p}
(\pi_{p-1}\omega),0),\quad&\text{if }n\le 2p-1, \\
(2\pi_{p}(\overline{F}^{n-p}\omega),-F^{p,p}\omega-2F^{n-p}
(\pi_{p-1}\omega),0),\quad&\text{if }n\ge 2p.
\end{cases}
\end{displaymath}
\end{proposition}
\hfill $\square$

\nnpar{Example.} 
Let $X$ be a complex projective manifold, and $\mathscr{E}^{\ast}_
{X}$ the sheaf of smooth, complex differential forms on $X$. Then, 
the complex of global sections $E^{\ast}_{X}$ of $\mathscr{E}^
{\ast}_{X}$ has a natural structure of Dolbeault complex. The 
cohomology of the complex $\mathcal{D}^{\ast}(E_{X},p)$ is naturally
isomorphic to the real Deligne-Beilinson cohomology $H^{\ast}_
{\mathcal{D}}(X,\mathbb{R}(p))$ of $X$ (see \cite{Deligne:dc}).

\nnpar{The product in the Deligne complex.} 
The multiplicative structure of a Dolbeault algebra induces a 
product in the Deligne complex which is graded commutative and
associative up to homotopy.  

\begin{definition}
A \emph{Dolbeault algebra} $A=(A^{\ast}_{\mathbb{R}},\dd_{A},\land)$
is a Dolbeault complex equipped with an associative and graded 
commutative product 
\begin{displaymath}
\land:A^{\ast}_{\mathbb{R}}\times A^{\ast}_{\mathbb{R}}\longrightarrow
A^{\ast}_{\mathbb{R}} 
\end{displaymath}
such that the induced multiplication on $A^{\ast}_{\mathbb{C}}$ is 
compatible with the bigrading, i.e.,
\begin{displaymath}  
A^{p,q}\land A^{p',q'}\subseteq A^{p+p',q+q'}.
\end{displaymath}
\end{definition}

\begin{definition}
\label{def:17}
Let $A$ be a Dolbeault algebra. The \emph{Deligne algebra associated 
to $A$} is the Deligne complex $\mathcal{D}^{\ast}(A,\ast)$ together 
with the graded commutative product
$\bullet:\mathcal{D}^{n}(A,p)\times\mathcal{D}^{m}(A,q)
\longrightarrow\mathcal{D}^{n+m}(A,p+q)$ given by   
\begin{align*}
&x\bullet y= \\
&\begin{cases}
(-1)^{n}r_{p}(x)\land y+x\land r_{q}(y),
&\text{if }n<2p,\,m<2q, \\
F^{l-r,l-r}(x\land y),
&\text{if }n<2p,\,m\geq 2q,\,l<2r, \\
F^{r,r}(r_{p}(x)\land y)+2\pi_{r}(\partial(x\land y)^{r-1,l-r}),
&\text{if }n<2p,\,m\geq 2q,\,l\geq 2r, \\
x\land y,
&\text{if }n\geq 2p,\,m\geq 2q,
\end{cases}
\end{align*}
where we have written $l=n+m$, $r=p+q$, and $r_{p}(x)=2\pi_{p}
(F^{p}\dd_{A}x)$.  
\end{definition}

\noindent
For a proof of the next result we refer to \cite{Burgos:CDB}.

\begin{proposition}
\label{prop:19}
Let $\mathcal{D}^{\ast}(A,\ast)$ be the Deligne algebra associated 
to a Dolbeault algebra $A$. Then, we have the following statements:    
\begin{enumerate}
\item[(i)]
The product $\bullet$ is associative up to a natural homotopy $h_{a}$.
\item[(ii)] 
For $x\in\mathcal{D}^{2p}(A,p)$, $y \in\mathcal{D}^{2q}(A,q)$,
$z\in\mathcal{D}^{2r}(A,r)$, we have
\begin{displaymath}
h_{a}(x\otimes y\otimes z)=0,
\end{displaymath}
i.e., the product $\bullet$ is pseudo-associative in the sense of 
definition \ref{def:9}; in particular, the direct sum $\bigoplus_
{n\in\mathbb{Z}}\mathcal{D}^{2n}(A,n)$ is an associative subalgebra.
\item[(iii)] 
\label{item:rp} 
There is a natural morphism of graded complexes
$r_{p}:\mathcal{D}^{\ast}(A,p)\longrightarrow A^{\ast}(p)$ given by
\begin{displaymath}
r_{p}(x)= 
\begin{cases} 
\partial(x^{p-1,n-p})-\bar{\partial}(x^{n-p,p-1}),
&\quad\text{if  }n\leq 2p-1, \\
x, 
&\quad\text{if  }n\geq 2p;
\end{cases}
\end{displaymath}
here $x\in\mathcal{D}^{n}(A,p)$. This morphism is multiplicative 
up to homotopy and induces a morphism of graded algebras
\begin{displaymath}
H^{*}(\mathcal{D}(A,p))\longrightarrow H^{\ast}(A(p)).
\end{displaymath}
\end{enumerate}
\hfill $\square$
\end{proposition}

\begin{remark}
\label{remark:fulea}
In his thesis \cite{Fulea:Thesis}, D. Fulea has introduced a new
product structure for the Deligne algebra $\mathcal{D}^{\ast}(A,\ast)$, 
which is graded commutative and associative. Using this product 
structure one may simplify some of the constructions of the present
paper;  
in particular, one can avoid the use of pseudo-associativity. 
\end{remark}

\nnpar{Specific degrees.}
In the sequel we will be interested in some specific degrees
where we can give simpler formulas. Namely, we consider
\begin{align*}
&\mathcal{D}^{2p}(A,p)=A^{2p}_{\mathbb{R}}(p)\cap A^{p,p}, \\
&\mathcal{D}^{2p-1}(A,p)=A^{2p-2}_{\mathbb{R}}(p-1)\cap A^{p-1,p-1}, \\
&\mathcal{D}^{2p-2}(A,p)=A^{2p-3}_{\mathbb{R}}(p-1)\cap
(A^{p-2,p-1}\oplus A^{p-1,p-2}).
\end{align*}
The corresponding differentials are given by
\begin{alignat*}{2}
&\dd_{\mathcal{D}}x=\dd_{A}x, 
&\quad &\text{if }x\in\mathcal{D}^{2p}(A,p), \\
&\dd_{\mathcal{D}}x=-2\partial\bar{\partial}x, 
&\quad &\text{if }x\in\mathcal{D}^{2p-1}(A,p), \\
&\dd_{\mathcal{D}}(x,y)=-\partial x-\bar{\partial}y, 
&\quad &\text{if }(x,y)\in\mathcal{D}^{2p-2}(A,p). 
\end{alignat*}
Moreover, the product is given as follows: for $x\in\mathcal{D}^
{2p}(A,p)$, $y\in\mathcal{D}^{2q}(A,q)$ or $y\in\mathcal{D}^{2q-1}
(A,q)$, we have
\begin{align*}
x\bullet y&=x\land y, \\
\intertext{and for $x\in\mathcal{D}^{2p-1}(A,p)$, $y\in
\mathcal{D}^{2q-1}(A,q)$, we have}
x\bullet y&=-\partial x\land y+\bar{\partial}x\land y+
x\land\partial y-x\land\bar{\partial}y.
\end{align*}

\nnpar{Dolbeault modules.} Once we have defined Dolbeault algebras we
may introduce the concept of Dolbeault modules. The main example is
the space of currents over a variety as a module over the space of
differential forms. 

\begin{definition}
Let $A$ be a Dolbeault algebra and $M$ a Dolbeault complex. We 
say that $M$ is a \emph{Dolbeault module over $A$}, if $M$ is a 
differential graded module satisfying
\begin{displaymath}
A^{p,q} M^{p',q'}\subset M^{p+p',q+q'}.
\end{displaymath}
\end{definition}

\noindent
The following proposition is straightforward.

\begin{proposition} 
\label{prop:27}
Let $A$ be a Dolbeault algebra and $M$ a Dolbeault module over 
$A$. Then, $\mathcal{D}^{\ast}(M,\ast)$ is a differential graded 
module over $\mathcal{D}^{\ast}(A,\ast)$. Moreover, the action 
is pseudo-associative.  
\hfill $\square$
\end{proposition}

\nnpar{Deligne complexes and Deligne-Beilinson cohomology.}
The main interest in Deligne complexes is expressed by the 
following theorem which is proven in \cite{Burgos:CDB} in a 
particular case, although the proof is valid in general. It 
is a consequence of proposition \ref{prop:32}.

\begin{theorem} 
\label{thm:4} 
Let $X$ be a complex algebraic manifold, $\overline{X}$ a smooth 
compactification of $X$ with $D=\overline{X}\setminus X$ a normal 
crossing divisor, and denote by $j:X\longrightarrow\overline{X}$ 
the natural inclusion. Let $\mathscr{A}^{\ast}$ be a sheaf of 
Dolbeault algebras over $\overline{X}^{\an}$ such that, for every 
$n,p$, the sheaves $\mathscr{A}^{n}$ and $F^{p}\mathscr{A}^{n}$ 
are acyclic, $\mathscr{A}^{\ast}$ is a multiplicative resolution 
of $Rj_{\ast}\mathbb{R}$, and $(\mathscr{A}^{\ast}_{\mathbb{C}},
F)$ is a multiplicative filtered resolution of $(\Omega^{\ast}_
{\overline{X}}(\log D),F)$. Putting $A^{\ast}=\Gamma(\overline{X},
\mathscr{A}^{\ast})$, we have a natural isomorphism of graded 
algebras   
\begin{displaymath}
H^{\ast}_{\mathcal{D}}(X,\mathbb{R}(p))\cong H^{\ast}
(\mathcal{D}(A,p)).
\end{displaymath}
Moreover, the morphism $r_{p}$ of proposition \ref{prop:19} 
(iii) induces the natural morphism of graded algebras   
\begin{displaymath}
H^{\ast}_{\mathcal{D}}(X,\mathbb{R}(p))\longrightarrow 
H^{\ast}(X,\mathbb{R}(p)).
\end{displaymath}
\hfill $\square$
\end{theorem}

\subsection{A Gillet complex: The Deligne algebra $\mathcal{D}_{\log}$}
\label{sec:real-deligne-beil}

\nnpar{Smooth differential forms with logarithmic singularities.}
Let $W$ be a complex algebraic manifold and $D$ a normal crossing
divisor in $W$. We put $X=W\setminus D$, and denote by $j:X
\longrightarrow W$ the natural inclusion. We recall that $\mathscr
{E}^{\ast}_{W}$  denotes the sheaf of smooth, complex differential
forms on $W$. 

\begin{definition}
\emph{The complex of sheaves $\mathscr{E}^{\ast}_{W}(\log D)$ 
of differential forms with logarithmic singularities along $D$}  
is the $\mathscr{E}^{\ast}_{W}$-subalgebra of $j_{\ast}\mathscr
{E}^{\ast}_{X}$, which is locally generated by the sections 
\begin{displaymath}
\log z_{i}\bar{z}_{i},\,\frac{\dd z_{i}}{z_{i}},\,\frac
{\dd\bar{z}_{i}}{\bar{z}_{i}}\quad\text{for }i=1,\dots,m,
\end{displaymath}
where $z_{1}\cdots z_{m}=0$ is a local equation for $D$ (see 
\cite{Burgos:CDc}).
\end{definition}

\begin{notation}
\label{not:1} 
In the sequel we will adhere to the following convention. Sheaves in
the analytic topology 
will usually be denoted by script letters, whereas the group of 
global sections will be denoted by the corresponding roman letters.
For instance, we will write    
\begin{align*}
&E^{\ast}_{W}=\Gamma(W,\mathscr{E}^{\ast}_{W}), \\
&E^{\ast}_{W}(\log D)=\Gamma(W,\mathscr{E}^{\ast}_{W}(\log D)), \\
&E^{\ast}_{W,\mathbb{R}}(\log D)=\Gamma(W,\mathscr{E}^{\ast}_{W,
\mathbb{R}}(\log D)).
\end{align*}
With these notations, the complex $E_{W}(\log D)=(E^{\ast}_
{W,\mathbb{R}}(\log D),\dd)$ is a Dolbeault algebra.    
\end{notation}

Suppose now that $W$ is proper. Then, there are multiplicative 
isomorphisms in the derived category of abelian sheaves over $W$
\begin{align*}
&Rj_{\ast}\mathbb{R}(p)\longrightarrow\mathscr{E}^{\ast}_
{W,\mathbb{R}}(\log D)(p), \\
&j_{\ast}\Omega^{\ast}_{X}\longrightarrow\mathscr{E}^{\ast}_
{W}(\log D), \\
&F^{p}\Omega^{\ast}_{W}(\log D)\longrightarrow F^{p}\mathscr
{E}^{\ast}_{W}(\log D).
\end{align*}
We therefore obtain by theorem \ref{thm:4} 

\begin{theorem} 
\label{thm:14}
There is a natural multiplicative isomorphism 
\begin{displaymath}
H^{\ast}_{\mathcal{D}}(X,\mathbb{R}(p))\longrightarrow
H^{\ast}(\mathcal{D}(E_{W}(\log D),p)).
\end{displaymath}
\hfill $\square$
\end{theorem}

\nnpar{Logarithmic singularities at infinity.} 
We want to obtain a description of Deligne-Beilinson cohomology 
which is independent of a given compactification. Given a complex 
algebraic manifold $X$, let $I$ be the category of all smooth
compactifications of $X$ with a normal crossing divisor as its 
complement. This means that an element $(\overline{X}_{\alpha},
j_{\alpha})$ of $I$ consists of a proper complex algebraic 
manifold $\overline{X}_{\alpha}$ together with an immersion 
$j_{\alpha}:X\longrightarrow\overline{X}_{\alpha}$ such that 
$D_{\alpha}=\overline{X}_{\alpha}\setminus j_{\alpha}(X)$ is 
a normal crossing divisor. The morphisms of $I$ are the maps 
$f:\overline{X}_{\alpha}\longrightarrow\overline{X}_{\beta}$ 
satisfying $f\circ j_{\alpha}=j_{\beta}$. It can be shown that 
the opposite category $I^{o}$ is directed (see \cite{Deligne:THII}).

We put 
\begin{displaymath}
E^{\ast}_{\log}(X)^{\circ}=\lim_{\substack{\longrightarrow \\ \alpha
\in I^{o}}}E^{\ast}_{\overline{X}_{\alpha}}(\log D_{\alpha});
\end{displaymath}
it is clear that the vector spaces $E^{\ast}_{\log}(X)^{\circ}$ form a
complex of presheaves $E^{\ast}_{\log}{}^{\circ}$ in the Zariski topology. 

The corresponding real subcomplex
$E_{\log}(X)^{\circ}=(E^{\ast}_{\log,\mathbb{R}}(X)^{\circ},\dd)$
is  
a presheaf of Dolbeault algebras. Moreover,
by the results of \cite{Burgos:CDc}, 
if $f:\overline{X}_{\alpha}\longrightarrow\overline{X}_{\beta}$
is a morphism of $I$, the induced morphism 
\begin{equation} \label{eq:47}
f^{\ast}:(E^{\ast}_{\overline{X}_{\beta}}(\log D_{\beta}),F)
\longrightarrow(E^{\ast}_{\overline{X}_{\alpha}}(\log D_
{\alpha}),F)
\end{equation}
is a real filtered quasi-isomorphism. Since $I^{o}$ is directed, 
all the induced morphisms 
\begin{displaymath}
(E^{\ast}_{\overline{X}_{\alpha}}(\log D_{\alpha}),F)
\longrightarrow(E^{\ast}_{\log}(X)^{\circ},F)
\end{displaymath}
are also real filtered quasi-isomorphisms.

We will denote by
$E^{\ast}_{\log}$ the complex of sheaves in the Zariski
topology associated to the complex of presheaves $E^{\ast}_{\log}{}^{\circ}$.

\begin{definition}
\emph{The complex of differential forms with logarithmic 
singularities along infinity} is defined by
\begin{displaymath}
E^{\ast}_{\log}(X)=\Gamma (X,E^{\ast}_{\log});
\end{displaymath}
it is a subalgebra of $E^{\ast}_{X}=\Gamma(X,\mathscr{E}^{\ast}_
{X})$. We denote the corresponding real subcomplex by $E^{\ast}_
{\log,\mathbb{R}}(X)$.
\end{definition}

\begin{remark}
The natural map $E^{\ast}_{\log}(X)^{\circ}\longrightarrow
E^{\ast}_{\log}(X)$ is injective, but not surjective, in general.
\end{remark}

\nnpar{Pseudo-flasque complexes of presheaves.} In order to understand
the relationship between the cohomology of the complexes
$E^{\ast}_{\log}(X)$ and $E^{\ast}_{\log}(X)^{\circ}$, we introduce
pseudo-flasque complexes of presheaves. 

\begin{definition}
  Let $X$ be a scheme and $\mathcal{F}^{\ast}$ a complex of presheaves
  in the Zariski topology of $X$. Then, $\mathcal{F}^{\ast}$ is called
  \emph{pseudo-flasque}, if $\mathcal{F}^{\ast}(\emptyset)=0$ and if,
  for every pair of open subsets $U$, $V$, the natural map  
  \begin{displaymath}
    \mathcal{F}(U\cup V)\longrightarrow s(\mathcal{F}(U)\oplus
    \mathcal{F}(V) \rightarrow \mathcal{F}(U\cap V))
  \end{displaymath}
  is a quasi-isomorphism.
\end{definition}

For instance, any complex of totally acyclic sheaves is a
pseudo-flasque complex of presheaves. 

The basic property of pseudo-flasque complexes of presheaves is the
analogue of proposition \ref{prop:13}.

\begin{proposition}\label{prop:34}
  Let $\mathcal{F}$ be a pseudo-flasque complex of presheaves, and let
  $\mathcal{F}'$ be the  
  associated complex of sheaves. Then, for any scheme $X$, we have
\begin{displaymath}
\mathbb{H}^{n}(X,\mathcal{F}')=H^{n}(\Gamma (X,\mathcal{F})).
\end{displaymath}
\hfill $\square$
\end{proposition}

\begin{remark}
  It is possible to generalize the notion of $\mathcal{G}$-complexes
  using pseudo-flasque complexes of presheaves. All the theory
  developed in 
  section \ref{sec:green-objects} and section \ref{sec:AC} can be
  generalized with minor modifications.
\end{remark}

\begin{proposition}\label{prop:35}
  The natural morphism 
  \begin{displaymath}
    (E^{\ast}_{\log}(X)^{\circ},F)\longrightarrow (E^{\ast}_{\log}(X),F)
  \end{displaymath}
  is a filtered quasi-isomorphism.
\end{proposition}
\begin{proof}
  By the filtered quasi-isomorphism \eqref{eq:47}, the filtered
  complex $(E^{\ast}_{\log}(X)^{\circ},F)$ computes the cohomology of
  $X$ with complex coefficients with its Hodge filtration.  
  Since the Mayer-Vietoris sequence for the cohomology of $X$ with complex
  coefficients is an exact sequence of
  mixed  Hodge structures, it induces Mayer-Vietoris sequences for the
  graded pieces with respect to the Hodge filtration. Therefore,
  the complexes of presheaves
  $E^{p,\ast}_{\log}(X)^{\circ}$ are pseudo-flasque. By a partition of
  unity argument (see the discussion
  before lemma \ref{lem:quisinvdb}) it is easy to see that the sheaves 
  $E^{p,q}_{\log}(X)$ are totally acyclic. Hence, the result follows from
  proposition \ref{prop:34}.
\end{proof}

\nnpar{A Deligne complex with logarithmic singularities.} 

\begin{definition}
Let $X$ be a complex algebraic manifold. For any integer $p$, 
we put
\begin{displaymath}
\mathcal{D}^{\ast}_{\log}(X,p)=\mathcal{D}^{\ast}(E_{\log}(X),p).
\end{displaymath}
\end{definition}

\noindent
By the exactness of the functor $\mathcal{D}(\cdot,p)$, theorem
\ref{thm:14} and proposition \ref{prop:35}, we obtain 

\begin{corollary} 
\label{cor:3}
There is a natural multiplicative isomorphism 
\begin{displaymath}
H^{\ast}_{\mathcal{D}}(X,\mathbb{R}(p))\longrightarrow
H^{\ast}(\mathcal{D}_{\log}(X,p)).
\end{displaymath}
\hfill $\square$  
\end{corollary}

\nnpar{A Gillet complex for real Deligne-Beilinson cohomology.} 
We denote by $C$ the site of regular schemes in $\ZAR(\Spec(
\mathbb{C}))$. In particular, we recall that all schemes in $C$ 
are separated and of finite type over $\Spec(\mathbb{C})$. For 
a scheme $X$ in $C$, the set of complex points $X(\mathbb{C})$ 
is a complex algebraic manifold. 

\begin{definition} 
\label{def:18}
For any integers $n,p$, let $\mathcal{D}^{n}_{\log}(p)$ denote 
the presheaf (in fact, a sheaf) over $C$, which assigns to $X$ the group  
\begin{displaymath}
\mathcal{D}^{n}_{\log}(X,p)=\mathcal{D}^{n}_{\log}(X(\mathbb
{C}),p)=\mathcal{D}^{n}(E_{\log}(X(\mathbb{C})),p).
\end{displaymath}
For any scheme $X$ in $C$, we will denote the induced presheaf 
of graded complexes of real vector spaces on $X$ by $\mathcal{D}_
{\log,X}=\mathcal{D}^{\ast}_{\log,X}(\ast)$.
\end{definition}

\begin{proposition} 
\label{prop:24}
For any integers $n,p$, the presheaf $\mathcal{D}^{n}_{\log,X}(p)$ 
is a totally acyclic sheaf on $X$. 
\end{proposition}
\begin{proof}
  This follows from the facts that the functor $\mathcal{D}(\cdot,p)$ is
  exact and that the sheaves $E^{p,q}_{\log}$ are totally acyclic.
\end{proof}

\begin{theorem} 
\label{thm:26}
The graded complex of sheaves of abelian groups $\mathcal{D}_
{\log}$ is a Gillet complex for regular schemes over $\mathbb{C}$, 
which computes real Deligne-Beilinson cohomology. Moreover, the 
pair $(\mathcal{D}_{\log},\bullet)$ is a graded commutative and
pseudo-associative algebra for real Deligne-Beilinson cohomology.
\end{theorem}
\begin{proof}
Since $\mathcal{D}_{\log}$ is a totally acyclic sheaf in the 
Zariski topology, the hypercohomology agrees with the cohomology 
of the complex of global sections. By theorem \ref{thm:4}, the 
complex of global sections computes real Deligne-Beilinson 
cohomology. Since \cite{Jannsen:DcHD} shows that real Deligne-Beilinson
cohomology satisfies the Gillet axioms, the complex $\mathcal{D}_
{\log}$ is a Gillet complex. Finally, proposition \ref{prop:19} and 
theorem \ref{thm:4} imply the claimed multiplicative properties.
\end{proof}

\subsection{Deligne-Beilinson homology of proper smooth varieties} 
\label{sec:deligne-homology}

The homology theory associated to a Gillet cohomology determines 
direct images and classes for algebraic cycles. For this reason 
we have to discuss the construction and basic properties of 
Deligne-Beilinson homology. In general, Deligne-Beilinson 
homology is defined by means of currents and smooth singular 
chains (for details, see \cite{Jannsen:DcHD}). But since we are 
only interested  in real Deligne-Beilinson homology, we do not need  
to use singular chains. Apart from some minor changes, we 
will follow \cite{Jannsen:DcHD}. In particular, unless stated 
otherwise, we will follow the conventions therein. 

\nnpar{Currents.} 
Let $X$ be a complex algebraic manifold. The sheaf ${}'\mathscr
{E}^{n}_{X}$ of currents of degree $n$ on $X$ is defined as 
follows. For any open subset $U$ of $X$, the group ${}'\mathscr
{E}^{n}_{X}(U)$ is the topological dual of the group of sections 
with compact support $\Gamma_{c}(U,\mathscr{E}^{-n}_{X})$. The 
differential 
\begin{displaymath}
\dd:{}'\mathscr{E}^{n}_{X}\longrightarrow{}'\mathscr{E}^{n+1}_{X}
\end{displaymath}
is defined by
\begin{displaymath}
\dd T(\varphi)=(-1)^{n}T(\dd\varphi);
\end{displaymath}
here $T$ is a current and $\varphi$ a corresponding test form.

The real structure and the bigrading of $\mathscr{E}^{\ast}_{X}$ 
induce a real structure and a bigrading of ${}'\mathscr{E}^{\ast}_
{X}$. Furthermore, there is a pairing 
\begin{displaymath}
\mathscr{E}^{n}_{X}\otimes{}'\mathscr{E}^{m}_{X}\longrightarrow 
{}'\mathscr{E}^{n+m}_{X},\quad\omega\otimes T\longmapsto\omega 
\land T,
\end{displaymath}
where the current $\omega \land T$ is defined by
\begin{displaymath}
(\omega\land T)(\eta)=T(\eta\land\omega).
\end{displaymath}
This pairing, the real structure and the bigrading equip ${}'\mathscr
{E}^{\ast}_{X}$ with the structure of a Dolbeault module ${}'\mathscr
{E}_{X}=({}'\mathscr{E}^{\ast}_{X,\mathbb{R}},\dd)$ over the Dolbeault 
algebra $\mathscr{E}_{X}=(\mathscr{E}^{\ast}_{X,\mathbb{R}},\dd)$.

\nnpar{Deligne-Beilinson homology.}
\begin{definition}
Let $X$ be a proper complex algebraic manifold. Then, the 
\emph{real Deligne-Beilinson homology groups of $X$} are 
defined by  
\begin{displaymath}
{}'H^{\ast}_{\mathcal{D}}(X,\mathbb{R}(p))=H^{\ast}
(\mathcal{D}({}'E_{X},p)),
\end{displaymath}
where, as fixed in notation \ref{not:1}, we have written 
${}'E^{\ast}_{X}=\Gamma(X,{}'\mathscr{E}^{\ast}_{X})$. We 
will also write
\begin{displaymath}
H_{n}^{\mathcal{D}}(X,\mathbb{R}(p))={}'H^{-n}_{\mathcal{D}}
(X,\mathbb{R}(-p)).
\end{displaymath}
\end{definition}

\begin{remark}
The proof that the above definition of the real Deligne-Beilinson
homology groups agrees with Beilinson's definition as given in
\cite{Beilinson:hr} or \cite{Jannsen:DcHD}, is completely analogous 
to the proof of theorem \ref{thm:4} given in \cite{Burgos:CDB}. 
\end{remark}

\nnpar{Equidimensional varieties.}
\begin{definition} 
\label{def:11} 
Let $X$ be an equidimensional complex algebraic manifold of 
dimension $d$. Then, we put   
\begin{displaymath}
\mathscr{D}^{\ast}_{X}={}'\mathscr{E}^{\ast}_{X}[-2d](-d),
\end{displaymath}
and, according to notation \ref{not:1},
\begin{displaymath}
D^{\ast}_{X}=\Gamma(X,\mathscr{D}^{\ast}_{X}).
\end{displaymath}
In particular, we note
\begin{align*}
&\mathscr{D}^{p,q}_{X}={}'\mathscr{E}^{p-d,q-d}_{X}, \\
&\mathscr{D}^{\ast}_{X,\mathbb{R}}(p)=(2\pi i)^{p-d}\cdot
{}'\mathscr{E}^{\ast}_{X,\mathbb{R}};    
\end{align*}
thus, the subcomplex of real currents $\mathscr{D}^{n}_
{X,\mathbb{R}}(U,p)$ is the topological dual of $\Gamma_
{c}(U,\mathscr{E}^{2d-n}_{X,\mathbb{R}}(p-d))$.
\end{definition}

\nnpar{The current associated to a differential form.} 
Let $X$ be an equidimensional complex algebraic manifold of 
dimension $d$. Then, there is a natural morphism of sheaves  
\begin{displaymath}
\mathscr{E}_{X}^{\ast}\longrightarrow\mathscr{D}_{X}^{\ast}
\end{displaymath}
given by $\omega\longmapsto[\omega]$, where the current $[\omega]$ 
is defined by
\begin{equation}
\label{eq:26}
[\omega](\eta)=\frac{1}{(2\pi i)^{d}}\int_{X}\eta\land\omega, 
\end{equation}
for a corresponding test form $\eta$ with compact support. This 
morphism is compatible with the structure of the Dolbeault module
$\mathscr{D}_{X}=(\mathscr{D}^{\ast}_{X,\mathbb{R}},\dd)$ over
the Dolbeault algebra $\mathscr{E}_{X}=(\mathscr{E}^{\ast}_{X,
\mathbb{R}},\dd)$. 


We will use the same notation and normalization  for any locally
integrable differential form $\omega $.

\nnpar{Poincar\'e duality.} 
The local version of Poincar\'e duality is the following theorem.
For a proof we refer, e.g., to \cite{GriffithsHarris:pag}, p. 384.

\begin{theorem}[Local Poincar\'e Duality]
\label{thm:22}
Let $X$ be an equidimensional complex algebraic manifold. Then,
the morphism $(\mathscr{E}^{\ast}_{X},F)\longrightarrow
(\mathscr{D}^{\ast}_{X},F)$ is a filtered quasi-isomorphism,
which is compatible with the underlying real structures.
\hfill $\square$
\end{theorem}

\noindent
From the local version of Poincar\'e duality, one derives

\begin{corollary}[Poincar\'e Duality]
Let $X$ be a proper equidimensional complex algebraic manifold of 
dimension $d$. Then, there is a natural isomorphism
\begin{displaymath}
H^{n}_{\mathcal{D}}(X,\mathbb{R}(p))\longrightarrow{}'H^{n-2d}_
{\mathcal{D}}(X,\mathbb{R}(p-d))=H_{2d-n}^{\mathcal{D}}(X,\mathbb
{R}(d-p)).
\end{displaymath}
\hfill $\square$
\end{corollary}

\nnpar{Direct images.}
Let $f:X\longrightarrow Y$ be a proper morphism between complex
algebraic manifolds. Then, we define a morphism of sheaves
\begin{displaymath}
f_{!}:f_{\ast}{}'\mathscr{E}^{\ast}_{X}\longrightarrow 
{}'\mathscr{E}^{\ast}_{Y}
\end{displaymath}
by setting $(f_{!}T)(\eta)=T(f^{\ast}\eta)$ for a test form 
$\eta$. When $X$ and $Y$ are proper, this morphism induces a morphism 
\begin{displaymath}
f_{!}:{}'H^{\ast}_{\mathcal{D}}(X,p)\longrightarrow{}'H^{\ast}_
{\mathcal{D}}(Y,p).
\end{displaymath}
If furthermore $X$ and $Y$ are equidimensional, and $f$ has relative
dimension $e$, 
the morphism $f_{!}$ sends  
$\mathscr{D}^{n}_{X,\mathbb{R}}(p)$ to
$\mathscr{D}^{n-2e}_{Y,\mathbb{R}}(p-e)$. Therefore, by Poincar\'e
duality,  for each $n,p$, we obtain an induced morphism
\begin{displaymath}
f_{!}:H^{n}_{\mathcal{D}}(X,\mathbb{R}(p))
\longrightarrow H^{n-2e}_{\mathcal{D}}(Y,\mathbb{R}(p-e)).
\end{displaymath}
If $f$ is smooth and $\omega$ is a differential form on $X$, we 
will write 
\begin{displaymath}
f_{!}\omega =\frac{1}{(2\pi i)^{e}}\int_{f}\omega  
\end{displaymath}
for the integration of $\omega$ along the fiber. It turns out 
that the definition of a current associated to a differential 
form and of push-forwards are compatible, i.e.,
\begin{displaymath}
f_{!}[\omega]=[f_{!}\omega].
\end{displaymath}

\nnpar{The fundamental class.}
Let $X$ be a proper equidimensional complex algebraic manifold 
of dimension $d$. Then, the fundamental class of axiom 1.2 (iv) 
of \cite{Gillet:RRhK} in
\begin{displaymath}
H^{0}_{\mathcal{D}}(X,\mathbb{R}(0))\cong H_{2d}^{\mathcal{D}}
(X,\mathbb{R}(d))={}'H^{-2d}_{\mathcal{D}}(X,\mathbb{R}(-d))
\end{displaymath}
is determined by the constant function $1$. Therefore, the
fundamental class can be represented by the current
\begin{displaymath}
\delta_{X}=[1].
\end{displaymath}

\nnpar{The class of a cycle.}
\begin{definition}
\label{def:10}
Let $X$ be a complex algebraic manifold, and $Y$ an 
$e$-dimensional irreducible subvariety of $X$. Let $\widetilde
{Y}$ be a resolution of singularities of $Y$, and $\imath:
\widetilde{Y}\longrightarrow X$ the induced map. Then, \emph
{the current integration along $Y$}, denoted by $\delta_{Y}$,
is defined by   
\begin{displaymath}
\delta_{Y}=\imath_{!}\delta_{\widetilde{Y}}\in\mathcal{D}^{-2e}
({}'E_{X},-e)={}'E^{-2e}_{X,\mathbb{R}}(-e)\cap{}'E^{-e,-e}_{X}.
\end{displaymath}
Therefore, it satisfies
\begin{displaymath}
\delta_{Y}(\eta)=\frac{1}{(2\pi i)^{e}}\int_{\widetilde{Y}}
\imath^{\ast}\eta.
\end{displaymath}
If $X$ is equidimensional of dimension $d$ and $p=d-e$, we 
obtain by the convention on the real structure of $\mathscr{D}^
{\ast}_{X}$ that $\delta_{Y}$ is an element of $\mathcal{D}^
{2p}(D_{X},p)=D^{2p}_{X,\mathbb{R}}(p)\cap D^{p,p}_{X}$. By 
linearity, we define $\delta_{y}$ for any algebraic cycle $y$. 
\end{definition}

\noindent
Since the current $\delta_{Y}$ is closed, the following definition
makes sense.

\begin{definition}
Let $X$ be a proper complex algebraic manifold, and $y$ an 
$e$-dimensional algebraic cycle of $X$. The \emph{homology class 
of $y$}, denoted by $\cl_{\mathcal{D}}(y)$ or simply by $\cl(y)$, 
is the class in ${}'H^{-2e}_{\mathcal{D}}(X,\mathbb{R}(-e))$ 
represented by $\delta_{y}$. If $X$ is equidimensional of dimension 
$d$ and $p=d-e$, the \emph{(cohomology) class of $y$} in $H^{2p}_
{\mathcal{D}}(X,\mathbb{R}(p))$, also denoted by $\cl_{\mathcal{D}}
(y)$ or $\cl(y)$, is defined by Poincar\'e duality.
\end{definition}

\subsection{Deligne-Beilinson homology of arbitrary varieties}
\label{sec:deligne-homol-arbitr}

Up to now, we have discussed Deligne-Beilinson homology only for
proper complex algebraic manifolds. Therefore, we still have to treat 
the case of non proper complex algebraic manifolds. Moreover, 
even if a Gillet cohomology is only defined for smooth varieties, 
the homology should also be defined for singular varieties. In 
this section we will discuss Deligne-Beilinson homology for general 
varieties. By abuse of notation, we will denote in this section 
a complex variety and its associated analytic space by the same 
letter. This will not cause confusion, since we consider only sheaves in 
the analytic topology in this section.

Note that the complexes of currents which we will use are 
slightly different from the complexes used in \cite{Jannsen:DcHD}. 

\nnpar{Currents on a subvariety.} 
Let $X$ be a complex algebraic manifold, $Y$ a closed subvariety
of $X$, and $j:Y\longrightarrow X$ the natural inclusion. We put
\begin{displaymath}
\Sigma_{Y}\mathscr{E}^{\ast}_{X}=\{\omega\in\mathscr{E}^{\ast}_{X}
\mid j^{\ast}\omega=0\}.
\end{displaymath}

\begin{definition} 
\label{def:16}
\emph{The sheaf of currents on $Y$} is defined by 
\begin{displaymath}
{}'\mathscr{E}^{\ast}_{Y}=\{T\in{}'\mathscr{E}^{\ast}_{X}\mid 
T(\omega)=0\,\,\forall\omega\in\Sigma_{Y}\mathscr{E}^{\ast}_{X}\}.
\end{displaymath}
Furthermore, we put
\begin{displaymath}
{}'\mathscr{E}^{\ast}_{X/Y}={}'\mathscr{E}^{\ast}_{X}/
{}'\mathscr{E}^{\ast}_{Y}.
\end{displaymath}
\end{definition}

The sheaf $'\mathscr{E}_{Y}^{\ast}$ was introduced by 
T. Bloom and M. Herrera in \cite{BloomHerrera:drcas}. If $Y$ 
is smooth, it agrees with the usual definition. The sheaf 
${}'\mathscr{E}^{n}_{X/Y}$ is the sheaf of distributions for
$\Sigma_{Y}\mathscr{E}^{-n}_{X}$. We point out that when $Y$ is a
normal crossing divisor, the complex of sheaves
${}'\mathscr{E}^{\ast}_{X/Y}$ does not  
agree with the complex denoted by $'{}\Omega^{\ast}_{X^{\infty}}\left< 
Y\right>$ in \cite{Jannsen:DcHD} because this last complex is not
defined over $\mathbb{R}$ (see also \cite{King:lccaj}). Note however that,
as a consequence of theorem \ref{thm:21}, both complexes are filtered
quasi-isomorphic.   

If $X$ is equidimensional of dimension $d$, we put as before
\begin{align*}
&\mathscr{D}^{\ast}_{X/Y}={}'\mathscr{E}^{\ast}_{X/Y}[-2d](-d), \\
&\overline{\mathscr{D}}^{\ast}_{Y}={}'\mathscr{E}^{\ast}_{Y}[-2d](-d).
\end{align*}
We observe that the grading for the complex $\overline{\mathscr
{D}}^{\ast}_{Y}$ is relative to the dimension of $X$ and not to 
the dimension of $Y$. Therefore, if $Y$ is smooth, the complex
$\overline{\mathscr{D}}^{\ast}_{Y}$ does not agree with the 
complex $\mathscr{D}^{\ast}_{Y}$.   

The complex ${}'\mathscr{E}^{\ast}_{X/Y}$ defines a Dolbeault 
complex. If $X$ is equidimensional of dimension $d$, the complex 
$\mathscr{D}^{\ast}_{X/Y}$ also defines a Dolbeault complex. In 
particular, both complexes have a well-defined Hodge filtration.

\nnpar{Normal crossing divisors.} 

\begin{theorem}
\label{thm:21}
Let $X$ be an equidimensional complex algebraic manifold 
of dimension $d$, and $Y$ a normal crossing divisor in 
$X$. Then, there is a well-defined morphism of complexes
\begin{displaymath}
\mathscr{E}^{\ast}_{X}(\log Y)\longrightarrow\mathscr{D}^
{\ast}_{X/Y},\quad\omega\longmapsto[\omega],
\end{displaymath}
satisfying
\begin{displaymath}
[\omega](\eta)=\frac{1}{(2\pi i)^{d}}\int_{U}\eta\land\omega, 
\end{displaymath}
where $U$ is an open subset of $X$, $\omega\in\Gamma(U,\mathscr
{E}^{\ast}_{X}(\log Y))$, and $\eta\in\Gamma_{c}(U,\Sigma_{Y}
\mathscr{E}^{\ast}_{X})$. Moreover, this morphism is a 
quasi-isomorphism with respect to the Hodge filtration.
\end{theorem}
\begin{proof}
The fact that we have a well-defined morphism of complexes is
proven in \cite{Burgos:Gftp}, 3.3. 

All the ingredients for the proof that the morphism in question 
is a filtered quasi-isomorphism are contained in \cite{Fujiki:dmHs}. 
We put
\begin{align*}
\mathscr{E}^{\ast}_{X}(\hol\log Y)&=\Omega^{\ast}_{X}(\log Y)
\otimes_{\Omega^{\ast}_{X}}\mathscr{E}^{\ast}_{X}, \\
\mathscr{Q}^{\ast}_{X}(\hol\log Y)&=\mathscr{E}^{\ast}_{X}
(\hol\log Y)/\mathscr{E}^{\ast}_{X}.
\end{align*}
Since these complexes are not defined over $\mathbb{R}$, they do
not define Dolbeault complexes. Nevertheless, we can define the 
Hodge filtration in the usual way.

By \cite{Burgos:CDc}, the natural inclusion 
\begin{displaymath}
\mathscr{E}^{\ast}_{X}(\hol\log Y)\longrightarrow\mathscr{E}^
{\ast}_{X}(\log Y)
\end{displaymath}
is a filtered quasi-isomorphism with respect to the Hodge
filtration. Thus, it is enough to show that the composition
\begin{displaymath}
\mathscr{E}^{\ast}_{X}(\hol\log Y)\longrightarrow\mathscr{D}^
{\ast}_{X/Y}
\end{displaymath}
is a filtered quasi-isomorphism.

Following \cite{Fujiki:dmHs} and \cite{HerreraLieberman:rpvcs}, 
we define a map $\PV:\mathscr{E}^{\ast}_{X}(\hol\log Y)
\longrightarrow\mathscr{D}^{\ast}_{X}$ by setting
\begin{displaymath}
\PV(\omega)(\eta)=\frac{1}{(2\pi i)^{d}}\int_{U}\eta\land\omega,    
\end{displaymath}
where $U$ is an open subset of $X$, $\omega\in\Gamma(U,\mathscr{E}^
{\ast}_{X}(\hol\log Y))$, and $\eta\in\Gamma_{c}(U,\mathscr{E}^
{\ast}_{X})$. We note that the map $\PV$ is not a morphism of 
complexes.

Starting with the commutative diagram
\begin{displaymath}
\xymatrix{&0\ar[r]&\mathscr{E}^{\ast}_{X}\ar[r]\ar[d]&
\mathscr{E}^{\ast}_{X}(\hol\log Y)\ar[r]\ar[d]\ar@{.>}
[dl]^{\PV}&\mathscr{Q}^{\ast}_{X}(\hol\log Y)\ar[r]&0 \\
0\ar[r]&\overline{\mathscr{D}}^{\ast}_{Y}\ar[r]&\mathscr
{D}^{\ast}_{X}\ar[r]&\mathscr{D}^{\ast}_{X/Y}\ar[r]&0\,,}
\end{displaymath}
we define a morphism of complexes
\begin{displaymath}
\Res:\mathscr{Q}^{\ast}_{X}(\hol\log Y)[-1]\longrightarrow 
\overline{\mathscr{D}}^{\ast}_{Y}
\end{displaymath}
by
\begin{displaymath}
\Res(\omega)=\dd\PV(\overline{\omega})-\PV(\dd\overline{\omega}),
\end{displaymath}
where $\overline{\omega}$ is any representative of $\omega$ 
in $\mathscr{E}^{\ast}_{X}(\hol\log Y)[-1]$. It is easy to see 
that $\Res$ is well defined and its image lies in $\overline
{\mathscr{D}}^{\ast}_{Y}$. By the local description of $\Res$ 
given in \cite{Fujiki:dmHs}, 3.6, it is a homogeneous morphism 
of bidegree $(0,1)$. Thus, it is compatible with the Hodge
filtration, and we obtain a map of distinguished exact triangles    
\begin{displaymath}
\xymatrix{\Gr_{F}^{p}\mathscr{Q}^{\ast}_{X}(\hol\log Y)[-1]
\ar[r]\ar[d]_{\Res}&\Gr_{F}^{p}\mathscr{E}^{\ast}_{X}\ar[r]
\ar[d]&\Gr_{F}^{p}\mathscr{E}^{\ast}_{X}(\hol\log Y)\ar[r]
\ar[d]& \\
\Gr_{F}^{p}\overline{\mathscr{D}}^{\ast}_{Y}\ar[r]&\Gr_{F}^{p} 
\mathscr{D}^{\ast}_{X}\ar[r]&\Gr_{F}^{p}\mathscr{D}^{\ast}_
{X/Y}\ar[r]&}
\end{displaymath}
From this we conclude that, if two of the above vertical maps are 
quasi-isomorphisms, so is the third. By theorem \ref{thm:22}, the 
middle vertical arrow is a quasi-isomorphism. By means of an 
auxiliary complex $\mathscr{K}^{\ast}_{X}\left<Y\right>$ equipped 
with a Hodge filtration, the following commutative triangle, where 
all the morphisms are filtered, is established in \cite{Fujiki:dmHs}: 
\begin{displaymath}
\xymatrix{&\mathscr{K}^{\ast}_{X}\left<Y\right>\ar[dd]^
{\eta_{X}} \\
\mathscr{Q}^{\ast}_{X}(\hol\log Y)[-1]\ar[ur]^{H_{0}}
\ar[dr]_{\Res}& \\
&\overline{\mathscr{D}}^{\ast}_{Y}}
\end{displaymath}
Now it is proven in \cite{Fujiki:dmHs} that all the morphisms 
in this diagram are quasi-isomorphisms and that $H_{0}$ is a 
filtered quasi-isomorphism with respect to the Hodge filtration. 
It is also immediate from the proof in \cite{Fujiki:dmHs} that 
$\eta_{X}$ is a filtered quasi-isomorphism with respect to the 
Hodge filtration. This proves that $\Res$ is also a filtered
quasi-isomorphism with respect to the Hodge filtration from
which the claim follows.
\end{proof}

\nnpar{Deligne-Beilinson homology of non proper varieties.}

\begin{definition}
Let $X$ be a complex algebraic manifold, and $\overline{X}$ a 
smooth compactification of $X$ with $D=\overline{X}\setminus X$ 
a normal crossing divisor. Then, the \emph{real Deligne-Beilinson 
homology groups of $X$} are defined by 
\begin{displaymath}
{}'H^{\ast}_{\mathcal{D}}(X,\mathbb{R}(p))=H^{\ast}(\mathcal{D}
({'}E_{\overline{X}/D},p)).
\end{displaymath}
We will also write
\begin{displaymath}
H_{n}^{\mathcal{D}}(X,\mathbb{R}(p))={}'H^{-n}_{\mathcal{D}}
(X,\mathbb{R}(-p)).
\end{displaymath}
\end{definition}

\begin{theorem} \label{thm:20}
  Real Deligne-Beilinson homology is well defined. That is, it does not
  depend on the 
  choice of a compactification. Moreover it is covariant for proper
  morphisms between smooth complex varieties.   
\end{theorem}
\begin{proof}
  The main ingredients of the proof of the first statement are: first,
  that given 
  two compactifications $\overline {X}_{1}$ and $\overline {X}_{2}$ as
  above, there is always a third compactification $\overline {X}_{3}$
  that dominates both; and second, that given a diagram of
  compactifications 
  \begin{displaymath}
    \xymatrix{
      \overline{X}_{1}\ar[r]^{\varphi}&\overline{X}_{2} \\
      &X\ar[ul]\ar[u],}
  \end{displaymath}
  with normal crossing divisors $D_{1}$ and $D_{2}$, the morphism
  $\varphi $ induces a filtered quasi-isomorphism
  \begin{displaymath}
    \varphi_{!}:\varphi_{\ast} {}'\mathscr{E}^{\ast}_{\overline
      {X}_{1}/D_{1}}\longrightarrow 
{}'\mathscr{E}^{\ast}_{\overline
      {X}_{2}/D_{2}}.
  \end{displaymath}
  The main ingredient to show the covariance is that, given a proper
  morphism $f:X\longrightarrow Y$ between smooth complex varieties,
  it is possible to  construct a commutative 
  diagram 
  \begin{displaymath}
       \xymatrix{
      X\ar[r]\ar[d]^{f}&\overline{X}\ar[d]^{\varphi} \\
      Y\ar[r]&\overline {Y},}
  \end{displaymath}
  where $\overline X$ and $\overline Y$ are smooth compactifications
  with normal crossing divisor $D_{X}=\overline X\setminus X$
  and $D_{Y}=\overline Y\setminus Y$. Then there is an induced
  morphism 
  \begin{equation}\label{eq:46}
    \varphi_{!}:\varphi_{\ast} {}'\mathscr{E}^{\ast}_{\overline
      {X}/D_{X}}\longrightarrow 
{}'\mathscr{E}^{\ast}_{\overline
      {Y}/D_{Y}}.
  \end{equation}
  For more details see \cite{Jannsen:DcHD} \S1.
\end{proof}

\nnpar{Poincar\'e duality for non proper varieties.}
A direct consequence of theorem \ref{thm:21} is the following

\begin{corollary}[Poincar\'e duality] 
Let $X$ be an equidimensional complex algebraic manifold of 
dimension $d$. Then, there is a natural isomorphism
\begin{displaymath}
H^{n}_{\mathcal{D}}(X,\mathbb{R}(p))\longrightarrow
{}'H^{n-2d}_{\mathcal{D}}(X,\mathbb{R}(p-d))=H_{2d-n}^
{\mathcal{D}}(X,\mathbb{R}(d-p)).
\end{displaymath}
\hfill $\square$
\end{corollary}

\begin{corollary}
The definition of Deligne-Beilinson homology given here 
agrees with the definition given in \cite{Beilinson:hr} 
and \cite{Jannsen:DcHD}. 
\hfill $\square$
\end{corollary}

\nnpar{Direct images.} Let $f:X\longrightarrow Y$ be a proper
equidimensional morphism between complex
algebraic manifolds, of relative dimension $e$. Then, using Poincar\'e
duality and the covariance of Deligne-Beilinson homology, for each
$n,p$, we obtain an induced morphism  
\begin{displaymath}
f_{!}:H^{n}_{\mathcal{D}}(X,\mathbb{R}(p))
\longrightarrow H^{n-2e}_{\mathcal{D}}(Y,\mathbb{R}(p-e)).
\end{displaymath}

\nnpar{Deligne-Beilinson homology of arbitrary varieties.}  
In order to define Deligne-Beilinson homology for a possibly 
singular variety, we will use simplicial resolutions. For 
details on simplicial resolutions and cohomological descent 
the reader is referred to \cite{Deligne:THIII} or \cite{SGA4}.
 
Let $X$ be a variety over $\mathbb{C}$. Then, there is a
smooth simplicial scheme $X_{\cdot}$ with an augmentation 
$\pi:X_{\cdot}\longrightarrow X$, which is a proper hypercovering,
and hence satisfies cohomological descent. Furthermore, there
is a proper smooth simplicial scheme $\overline{X}_{\cdot}$
together with an open immersion $X_{\cdot}\longrightarrow
\overline{X}_{\cdot}$ such that the complement $D_{\cdot}=
\overline{X}_{\cdot}\setminus X_{\cdot} $ is a normal crossing 
divisor. By the covariance of currents, we obtain a simplicial 
graded complex $\mathcal{D}^{\ast}({}'E_{\overline{X}_{\cdot}/
D_{\cdot}},p)$. Let ${}'\mathcal{N}$ be the canonical 
normalization functor transforming a simplicial group into a 
cochain complex, i.e., this functor transforms a simplicial 
group into a chain complex and then reverses the signs of the 
grading.
 
\begin{definition}
The \emph{real Deligne-Beilinson homology groups of $X$} are
defined by
\begin{displaymath}
{}'H^{\ast}_{\mathcal{D}}(X,\mathbb{R}(p))=H^{\ast}(s({}'
\mathcal{N}\mathcal{D}({}'E_{\overline{X}_{\cdot}/D_{\cdot}},p))).   
\end{displaymath}
As usual, we will also write
\begin{displaymath}
H_{n}^{\mathcal{D}}(X,\mathbb{R}(p))={}'H^{-n}_{\mathcal{D}}
(X,\mathbb{R}(-p)).
\end{displaymath}
\end{definition}

For a proof of the next result see \cite{Jannsen:DcHD}. It is based on
the same principles as the proof of theorem \ref{thm:20}.

\begin{theorem} \label{thm:25}
  Real Deligne-Beilinson homology is well defined. That is, it does not
  depend on the 
  choice of a compactification of a proper hypercovering. Moreover it
  is covariant for proper 
  morphisms between complex varieties.   
\end{theorem}

\begin{remark}
If $Y$ is a normal crossing divisor of a proper smooth variety $W$, 
the proof of theorem \ref{thm:21} shows that we do not need to use 
a simplicial resolution of $Y$ to compute the Deligne-Beilinson 
homology of $Y$. Indeed, the complex ${}'\mathscr{E}^{\ast}_{Y}$ 
is quasi-isomorphic to the one obtained by a simplicial resolution. 
Therefore, we have 
\begin{displaymath}     
{}'H^{\ast}_{\mathcal{D}}(Y,\mathbb{R}(p))=H^{\ast}(\mathcal{D}
({}'E_{Y},p)).
\end{displaymath}
More generally, if $Y$ is a normal crossing divisor of a smooth
variety $W$, and $\overline{W}$ is a smooth compactification of
$W$ such that $Z=\overline{W}\setminus W$ and $Y\cup Z$ are normal 
crossing divisors, we write
\begin{displaymath}
{}'\mathscr{E}^{\ast}_{(Y\cup Z)/Z}={}'\mathscr{E}^{\ast}_
{Y\cup Z}/{}'\mathscr{E}^{\ast}_{Z}.
\end{displaymath}
Then, we can use the graded complex $\mathcal{D}^{\ast}({}'E_
{(Y\cup Z)/Z},p)$ to compute the real Deligne-Beilinson 
homology of $Y$. 
\end{remark}

\nnpar{Poincar\'e duality for general varieties.}
We will recall how to construct the long exact sequence in 
homology associated to a closed subvariety, and the Poincar\'e 
duality isomorphism between homology and cohomology with 
support. For details on the proof we refer to \cite{Jannsen:DcHD}.

Let $X$ be a variety over $\mathbb{C}$, and $Y$ a closed subvariety 
of $X$. Then, there is a smooth simplicial scheme $X_{\cdot}$ 
with an augmentation $\pi:X_{\cdot}\longrightarrow X$, which 
is a proper hypercovering, and hence satisfies cohomological 
descent. Furthermore, there is a proper smooth simplicial scheme 
$\overline{X}_{\cdot}$ together with an open immersion $X_{\cdot}
\longrightarrow\overline{X}_{\cdot}$ such that $D_{\cdot}=
\overline{X}_{\cdot}\setminus X_{\cdot}$, $Y_{\cdot}=\pi^{-1}
(Y)$, and $Z_{\cdot}=D_{\cdot}\cup Y_{\cdot}$ are normal crossing 
divisors. We observe that $Y_{\cdot}\longrightarrow Y$ also 
satisfies cohomological descent. For $p\in\mathbb{Z}$, we now
obtain a short exact sequence
\begin{align}
0\longrightarrow s({}'\mathcal{N}\mathcal{D}^{\ast}({}'E_
{Z_{\cdot}/D_{\cdot}},p))&\longrightarrow s({}'\mathcal{N}
\mathcal{D}^{\ast}({}'E_{\overline{X}_{\cdot}/D_{\cdot}},p))
\notag \\
&\longrightarrow s({}'\mathcal{N}\mathcal{D}^{\ast}({}'E_
{\overline{X}_{\cdot}/Z_{\cdot}},p))\longrightarrow 0.
\label{eq:27}
\end{align}
From this exact sequence we obtain a long exact sequence of 
homology groups
\begin{displaymath}
\dots\longrightarrow{}'H^{\ast}_{\mathcal{D}}(Y,\mathbb{R}(p))
\longrightarrow{}'H^{\ast}_{\mathcal{D}}(X,\mathbb{R}(p))
\longrightarrow{}'H^{\ast}_{\mathcal{D}}(X\setminus Y,\mathbb
{R}(p))\longrightarrow\dots
\end{displaymath}

\begin{theorem}[Poincar\'e duality]
\label{thm:23} 
Let $X$ be an equidimensional complex algebraic manifold of 
dimension $d$, and $Y$ a closed subvariety of $X$. Then, there 
is a natural isomorphism
\begin{displaymath}
H_{n}^{\mathcal{D}}(Y,\mathbb{R}(p))\longrightarrow H^{2d-n}_
{\mathcal{D},Y}(X,\mathbb{R}(d-p)).
\end{displaymath}
\end{theorem}
\begin{proof}
Let $\overline{X}$ be a smooth compactification of $X$ with
$D=\overline{X}\setminus X$ a normal crossing divisor. Then, 
we may assume that the simplicial resolution $\overline{X}_
{\cdot}$ used to construct the exact sequence \eqref{eq:27} 
satisfies the conditions that $\overline{X}_{0}$ is a proper 
modification of $\overline{X}$, which is an isomorphism over 
the complement of the adherence of $Y$ and furthermore that $\overline
X_{\cdot}\longrightarrow \overline X$ is a proper hypercovering.  We
define $D_{\cdot}$, 
$Y_{\cdot}$, and $Z_{\cdot}$ as before. Then $\overline{X}_{0}
\setminus Z_{0}$ can be identified with $X\setminus Y$. Therefore,
the cohomology group $H^{\ast}_{\mathcal{D},Y}(X,\mathbb{R}(p))$ 
is computed as the cohomology of the simple complex associated 
to the morphism of graded complexes
\begin{displaymath}
\mathcal{D}^{\ast}(E_{\overline{X}}(\log D),p)\longrightarrow
\mathcal{D}^{\ast}(E_{\overline{X}_{0}}(\log Z_{0}),p).
\end{displaymath}
For any variety $S$ occurring in the simplicial resolution
$\overline{X}_{\cdot}$, we put
\begin{displaymath}
\overline{\mathscr{D}}^{\ast}_{S}={}'\mathscr{E}^{\ast}_{S}
[-2d](-d), 
\end{displaymath}
i.e., all the dimensions are considered relatively to $X$. 
The condition of $\overline
X_{\cdot}\longrightarrow \overline X$ being a proper hypercovering
implies the existence of a 
quasi-isomorphism 
\begin{displaymath}
  s({}'\mathcal{N}\mathcal{D}^{\ast}(\overline{D}_{\overline{X}_
{\cdot}/D_{\cdot}},p))\longrightarrow
\mathcal{D}^{\ast}(\overline{D}_{\overline{X}/D},p). 
\end{displaymath}
Then, there is a commutative diagram   
\begin{equation}
\label{eq:28}
\xymatrix{
\mathcal{D}^{\ast}(E_{\overline{X}}(\log D),p)\ar[r]\ar[d]_
{\sim}&\mathcal{D}^{\ast}(E_{\overline{X}_{0}}(\log D_{0}),p)
\ar[r]\ar[d]_{\sim}&\mathcal{D}^{\ast}(E_{\overline{X}_{0}}
(\log Z_{0}),p)\ar[d]_{\sim} \\
\mathcal{D}^{\ast}(\overline{D}_{\overline{X}/D},p)&\mathcal
{D}^{\ast}(\overline{D}_{\overline{X}_{0}/D_{0}},p)\ar[l]\ar[r]
\ar[d]&\mathcal{D}^{\ast}(\overline{D}_{\overline{X}_{0}/Z_{0}},
p)\ar[d]_{\sim} \\ 
&s({}'\mathcal{N}\mathcal{D}^{\ast}(\overline{D}_{\overline{X}_
{\cdot}/D_{\cdot}},p))\ar[r]\ar[ul]_{\sim}&s({}'\mathcal{N}
\mathcal{D}^{\ast}(\overline{D}_{\overline{X}_{\cdot}/Z_{\cdot}},
p))\,, 
}
\end{equation}
where the arrows marked with $\sim$ are quasi-isomorphisms. This
shows that the cohomology groups $H^{\ast}_{\mathcal{D},Y}(X,
\mathbb{R}(p))$ can be computed as the cohomology of the simple 
complex associated to the morphism of graded complexes 
\begin{displaymath}
s({}'\mathcal{N}\mathcal{D}^{\ast}(\overline{D}_{\overline{X}_
{\cdot}/D_{\cdot}},p))\longrightarrow s({}'\mathcal{N}\mathcal
{D}^{\ast}(\overline{D}_{\overline{X}_{\cdot}/Z_{\cdot}},p)).
\end{displaymath}
The result now follows from the exact sequence \eqref{eq:27} and
Poincar\'e duality for smooth varieties.
\end{proof}

\nnpar{Direct images.} We are finally able to define direct images for
Deligne-Beilinson cohomology with support. Let $f:X\longrightarrow Y$
be a proper equidimensional morphism between smooth complex varieties,
of relative dimension $e$. Let $Z\subset X$ and $Z'\subset Y$ be
closed subsets such that $f(Z)\subset Z'$. Then, for each
$n,p$, there is a morphism 
\begin{displaymath}
f_{!}:H^{n}_{\mathcal{D},Z}(X,\mathbb{R}(p))
\longrightarrow H^{n-2e}_{\mathcal{D},Z'}(Y,\mathbb{R}(p-e)).
\end{displaymath}
induced by Poincar\'e
duality and the covariance of Deligne-Beilinson homo\-logy.

There are two special cases where the direct image can be easily
described in terms of differential forms with logarithmic
singularities. The proof is left to the reader.

\begin{proposition}\label{prop:33}
  \begin{enumerate}
  \item Let $X$ be a proper smooth complex variety and let $Z\subset
    Y$ be closed subsets. Let $\widetilde X_{Z}$ be an embedded
    resolution of singularities of $Z$ with normal crossing divisor
    $D_{Z}$ and let $\widetilde X_{Y}$ be an embedded
    resolution of singularities of $Y$ with normal crossing divisor
    $D_{Y}$ that dominates $\widetilde X_{Z}$. Then the morphism
    $H^{n}_{\mathcal{D},Z}(X,\mathbb{R}(p))
    \longrightarrow H^{n}_{\mathcal{D},Y}(X,\mathbb{R}(p))$ is induced by
    the natural morphism
    \begin{multline*}
      s(\mathcal{D}^{\ast}(E_{X},p)\longrightarrow 
      \mathcal{D}^{\ast}(E_{\widetilde
        {X}_{Z}}(\log D_{Z}),p))
      \longrightarrow \\
      s(\mathcal{D}^{\ast}(E_{X},p)\longrightarrow 
      \mathcal{D}^{\ast}(E_{\widetilde
        {X}_{Y}}(\log D_{Y}),p)).
    \end{multline*}
  \item Let $f:X\longrightarrow Y$ be a proper smooth morphism between
  proper varieties, of relative dimension $e$. Let $Z\subset Y$ be a
  closed subset. Let $\widetilde Y$ be an embedded resolution of
  singularities of $Z$ with normal crossing divisor $D$. We write 
  \begin{displaymath}
    \widetilde X=X\underset{Y}{\times}\widetilde Y,\quad
    B=X\underset{Y}{\times}\widetilde D.
  \end{displaymath}
  By the smoothness of $f$, $\widetilde X$ is an embedded resolution
  of singularities of $f^{-1}(Z)$ with normal crossing divisor
  $B$. Then integration along the fiber induces a morphism  
  \begin{displaymath}
    f_{!}:\mathcal{D}^{\ast}(E_{\widetilde X}(\log B),p)\longrightarrow 
    \mathcal{D}^{\ast}(E_{\widetilde Y}(\log D),p-e)[2e]
  \end{displaymath}
  and the morphism
    $H^{n}_{\mathcal{D},f^{-1}(Z)}(X,\mathbb{R}(p))
    \longrightarrow H^{n-2e}_{\mathcal{D},Z}(Y,\mathbb{R}(p-e))$ is
    induced by 
    the  morphism
    \begin{multline*}
      f_{!}:s(\mathcal{D}^{\ast}(E_{X},p)\longrightarrow 
      \mathcal{D}^{\ast}(E_{\widetilde
        {X}}(\log B),p))
      \longrightarrow \\
      s(\mathcal{D}^{\ast}(E_{Y},p-e)\longrightarrow 
      \mathcal{D}^{\ast}(E_{\widetilde
        {Y}}(\log D),p-e))[-2e].
    \end{multline*}
  \end{enumerate}

\end{proposition}
\subsection{Classes of cycles and line bundles}
\label{sec:classes-cycles-line}

\nnpar{The class of a cycle in cohomology with support.} 
We are now in a position to describe the class of a cycle in 
real Deligne-Beilinson cohomology with support. To do this, 
let $X$ be an equidimensional complex algebraic manifold of 
dimension $d$. Since the class of a cycle in a quasi-projective 
variety is the restriction of the class of a cycle in any 
compactification, we may assume that $X$ is projective. Then, 
let $y$ be a $p$-codimensional, i.e., $(d-p)$-dimensional 
cycle of $X$, and $Y$ the support of $y$. Letting $Y=\bigcup_
{j}Y_{j}$ be the decomposition of $Y$ into irreducible 
components, we have $y=\sum_{j}n_{j}Y_{j}$ with certain 
multiplicities $n_{j}$. If $\widetilde{Y}_{\cdot}$ denotes 
a proper hypercovering of $Y$ such that $\widetilde{Y}_{0}$ 
is a resolution of singularities of $Y$, then $\widetilde
{Y}_{0}$ is a disjoint union of irreducible components 
$\widetilde{Y}_{0}=\bigcup_{j}\widetilde{Y}_{0,j}$ with
$\widetilde{Y}_{0,j}$ a resolution of singularities of 
$Y_{j}$. If $n_{y}$ denotes the locally constant function 
with value $n_{j}$ on the component $\widetilde{Y}_{0,j}$, 
then $n_{y}$ defines an element $[n_{y}]$ of
\begin{displaymath}
\mathcal{D}^{-2d+2p}({}'E_{\widetilde{Y}_{0}},-d+p)=
{}'E^{-2d+2p}_{\widetilde{Y}_{0},\mathbb{R}}(-d+p)
\cap{}'E^{-d+p,-d+p}_{\widetilde{Y}_{0}}
\end{displaymath}
given by 
\begin{displaymath} 
[n_{y}](\eta)=\frac{1}{(2\pi i)^{d-p}}\int_{\widetilde{Y}_{0}} 
n_{y}\cdot\eta\,.   
\end{displaymath}
The element $[n_{y}]$ is a cycle of the complex $s(\mathcal{D}^
{\ast}({}'E_{\widetilde{Y}_{\cdot}},-d+p))$, and hence defines 
a cohomology class in 
\begin{align*}
H_{2d-2p}^{\mathcal{D}}(Y,\mathbb{R}(d-p))&={}'H^{-2d+2p}_
{\mathcal{D}}(Y,\mathbb{R}(-d+p)) \\
&=H^{-2d+2p}(s(\mathcal{D}({}'E_{\widetilde{Y}_{\cdot}},-d+p))).
\end{align*}
The image of this class under the Poincar\'e duality isomorphism
\begin{displaymath}
H_{2d-2p}^{\mathcal{D}}(Y,\mathbb{R}(d-p))\overset{\cong}
{\longrightarrow}H^{2p}_{\mathcal{D},Y}(X,\mathbb{R}(p))
\end{displaymath}
is merely the class of the cycle $y$ in $H^{2p}_{\mathcal{D},
Y}(X,\mathbb{R}(p))$. Observing that the real mixed Hodge structure 
of $H^{2p}_{Y}(X,\mathbb{R})$ is pure of type $(p,p)$ and that
$H^{2p-1}_{Y}(X,\mathbb{R})=0$, the exact sequence \eqref{eq:25} 
implies that the natural morphism   
\begin{equation}
\label{eq:33}
r_{p}:H^{2p}_{\mathcal{D},Y}(X,\mathbb{R}(p))\longrightarrow 
H^{2p}_{Y}(X,\mathbb{R}(p))
\end{equation}
is an isomorphism. Therefore, the class of $y$ in $H^{2p}_
{\mathcal{D},Y}(X,\mathbb{R}(p))$ is determined by its image 
in the Betti cohomology group $H^{2p}_{Y}(X,\mathbb{R}(p))$.  
Since the latter group is determined as the cohomology of the 
simple complex
\begin{displaymath}
s(D^{\ast}_{X}(p)\longrightarrow D^{\ast}_{X\setminus Y}(p)),
\end{displaymath}
the class of the cycle $y$ can be represented by the pair
$(\delta_{y},0)$. 

\nnpar{Differential forms representing the class of a cycle}. 
By the above considerations we have a representative for the 
class of a cycle in terms of currents. We now give a criterion, 
when a pair of differential forms represents this class. For 
this let 
\begin{displaymath}
(\omega,g)\in s^{2p}(\mathcal{D}_{\log}(X,p)\longrightarrow
\mathcal{D}_{\log}(X\setminus Y,p));
\end{displaymath}
by \cite{Burgos:Gftp}, 3.8.2, the form $g$ is locally integrable 
on $X$. If $U=X\setminus Y$, we write $[g]_{U}$ for the current 
associated to $g$ as a smooth form on $U$; we write $[g]_{X}$ 
for the current associated to $g$ as a locally integrable form 
on $X$. As usual, we put $\dd^{c}=(4\pi i)^{-1}(\partial-\bar
{\partial})$.   
\begin{proposition}
\label{prop:30} 
Let $X$ be a complex algebraic manifold, and $y$ a $p$-codimensional 
cycle on $X$ with support $Y$. Let $(\omega ,g)$ be a cycle in
\begin{displaymath}
  s^{2p}(\mathcal{D}_{\log}(X,p)\longrightarrow
\mathcal{D}_{\log}(X\setminus Y,p)).
\end{displaymath}
 Then, we have the following statements:
\begin{enumerate}
\item[(i)]   
The class of the cycle $(\omega,g)$ in $H^{2p}_{\mathcal{D},Y}
(X,\mathbb{R}(p))$ is equal to the class of $y$, if and only if 
\begin{equation} 
\label{eq:35}
-2\partial\bar{\partial}[g]_{X}=[\omega]-\delta_{y}. 
\end{equation}
\item[(ii)]  
Assume that $y=\sum_{j}n_{j}Y_{j}$ with irreducible subvarieties 
$Y_{j}$ and certain multiplicities $n_{j}$. Then, the cycle 
$(\omega,g)$ represents the class of $y$, if and only if the
equality
\begin{equation}
\label{eq:36} 
-\lim_{\varepsilon\to 0}\int_{\partial B_{\varepsilon}(Y)}\alpha 
\dd^{c}g=\frac{(2\pi i)^{p-1}}{2}\sum_{j}n_{j}\int_{Y_{j}}\alpha 
\end{equation}
holds for any differential form $\alpha$; here $B_{\varepsilon}(Y)$ 
is an $\varepsilon$-neighborhood of $Y$ such that the orientation 
of $\partial B_{\varepsilon}(Y)$ is induced from the orientation of
$B_{\varepsilon}(Y)$. 
\end{enumerate}
\end{proposition}
\begin{proof}
(i) First, assume that the cycle $(\omega,g)$ represents the class 
of $y$ in the group $H^{2p}_{\mathcal{D},Y}(X,\mathbb{R}(p))$. 
Let $\overline{X}$ be a smooth compactification of $X$, and 
$\overline{y}$ a cycle extending $y$ with support $\overline{Y}$. 
By \cite{Burgos:Gftp}, 4.8, there exists a pair $(\omega_{1}',
g_{1}')$ representing the class of $\overline{y}$ in the simple complex
$s^{2p}(\mathcal{D}_{\log}(\overline{X},p)\longrightarrow
\mathcal{D}_{\log}(\overline{X}\setminus\overline{Y},p))$ such 
that $-2\partial\bar{\partial}[g_{1}']_{\overline{X}}=[\omega_
{1}']-\delta_{\overline{y}}$. We denote by $\omega_{1}$, resp. 
$g_{1}$ the restriction of $\omega_{1}'$, resp. $g_{1}'$ to 
$X$. By the functoriality of the class of a cycle, the pairs 
$(\omega_{1},g_{1})$ and $(\omega,g)$ represent the same class. 
Therefore, there are elements $a\in\mathcal{D}^{2p-1}_{\log}
(X,p)$ and $b\in\mathcal{D}^{2p-2}_{\log}(X\setminus Y,p)$
such that     
\begin{displaymath}     
(\dd_{\mathcal{D}}a,a-\dd_{\mathcal{D}}b)=(\omega,g)-(\omega_{1},   
g_{1}).   
\end{displaymath} 
This shows 
\begin{displaymath} 
\dd_{\mathcal{D}}[g]_{X}=\dd_{\mathcal{D}}[g_{1}]_{X}+\dd_
{\mathcal{D}}[a]=[\omega_{1}]-\delta_{y}+[\omega]-[\omega_{1}]=
[\omega]-\delta_{y}. 
\end{displaymath}
Since $\dd_{\mathcal{D}}=-2\partial\bar{\partial}$ in degree 
$2p-1$, we find that the pair $(\omega,g)$ satisfies equation
\eqref{eq:35}.

Conversely, assume that $(\omega,g)$ satisfies equation \eqref{eq:35}. 
By \cite{Burgos:Gftp}, 3.8.3, we know that $r_{p}([g]_{X})=[r_{p}
(g)]_{X}$, where $r_{p}$ is as in proposition \ref{prop:19} (iii). 
Since we have $-2\partial\bar{\partial}[g]_{X}=\dd r_{p}([g]_{X})$, 
equation \eqref{eq:35} implies 
\begin{equation}
\label{eq:30}
\dd[r_{p}(g)]_{X}=[r_{p}(\omega)]-\delta_{y}.
\end{equation}
Hence, we obtain
\begin{displaymath}
\dd([r_{p}(g)]_{X},0)=([r_{p}(\omega)],[r_{p}(g)]_{U})-
(\delta_{y},0).
\end{displaymath}
Since the natural morphism \eqref{eq:33} is an isomorphism, we 
therefore obtain that $(\omega,g)$ represents the class of $y$.
  
(ii) Let $\alpha$ be a differential form and put $n=\deg\alpha$. 
Equation \eqref{eq:30} is equivalent to the equation
\begin{displaymath}
(-1)^{n+1}\frac{1}{(2\pi i)^{n}}\int_{X}\dd\alpha\land r_{p}(g)=
\frac{1}{(2\pi i)^{n}}\int_{X}\alpha\land\omega-\frac{1}{(2\pi i)^
{n-p}}\sum_{j}n_{j}\int_{Y_{j}}\alpha.   
\end{displaymath}
Using the fact that $(\omega,g)$ is a cycle and that $\alpha$ has even
degree, the above equation is equivalent to
\begin{displaymath}
\int_{X}\dd(\alpha\land r_{p}(g))=(2\pi i)^{p}\sum_{j}n_{j}
\int_{Y_{j}}\alpha.
\end{displaymath}
The result now follows from Stokes theorem and the fact that 
\begin{displaymath}
\dd^{c}g=\frac{1}{4\pi i}\,r_{p}(g), 
\end{displaymath}
as shown in proposition \ref{prop:19} (iii).
\end{proof}

\nnpar{The first Chern form of a line bundle.} 
The morphism ${\rm c}_{1}$ in \eqref{eq:29} normalizes the first 
Chern class of a line bundle and therefore of all Chern classes 
of vector bundles. We recall now the usual way of how to obtain 
representatives for the first Chern class of a line bundle. Let 
$X$ be a quasi-projective complex algebraic manifold, $\overline
{X}$ a smooth compactification of $X$, and $L$ a line bundle over 
$X$. Then, $L$ can be extended to a line bundle $\overline{L}$ 
over $\overline{X}$. Let $\overline{h}$ be a smooth hermitian 
metric on $\overline{L}$, and $h$ the restriction of $\overline
{h}$ to $L$. If $s$ is a non-vanishing, rational section of $L$, 
we put $\Vert s\Vert^{2}=h(s,s)$, $y=\dv(s)$, $Y=\supp y$, $U=
X\setminus Y$, and
\begin{align}
g_{s}&=-\frac{1}{2}\log(\Vert s\Vert^{2}), \\
\omega_{s}&=-2\partial\bar{\partial}g_{s}.
\end{align}
Moreover, we put $\overline{y}=\dv(s)$ with $s$ viewed as a 
section on $\overline{X}$, and $\overline{Y}=\supp\overline{y}$.  
We call
\begin{displaymath}
{\rm c}_{1}(L,h)=\omega_{s}=-2\partial\bar{\partial}g_{s}=\partial
\bar{\partial}\log(\Vert s\Vert^{2})
\end{displaymath}
the \emph{first Chern form of $(L,h)$}. We note that this definition 
of the first Chern form differs by a factor $(2\pi i)^{-1}$ from the 
corresponding definition in cohomology with real coefficients because
of the $2\pi i$-twist in Deligne-Beilinson cohomology. 

\begin{proposition} 
\label{prop:26}
With the above assumptions the following statements hold: The 
form $\omega_{s}$ belongs to $\mathcal{D}^{2}_{\log}(X,1)$, 
and the form $g_{s}$ belongs to $\mathcal{D}^{1}_{\log}(U,1)$. 
The pair $(\omega_{s},g_{s})$ is a cycle of the simple complex 
\begin{displaymath}
s^{2}(\mathcal{D}_{\log}(X,1)\longrightarrow \mathcal{D}_
{\log}(U,1)).
\end{displaymath}
Moreover, this pair represents the class of $\dv(s)$ in the 
cohomology group $H^{2}_{\mathcal{D},Y}(X,\mathbb{R}(1))$.
\end{proposition}    
\begin{proof} 
Let $s'$ be a non-vanishing, regular section of $L$ in an open 
subset $V$ of $X$. Then, $\omega_{s}=\partial\bar{\partial}\log
(\Vert s'\Vert^{2})$, which shows that $\omega_{s}$ is smooth 
on the whole of $\overline{X}$. Furthermore, since $g_{s}$ is 
a real function and $-2\partial\bar{\partial}$ is an imaginary 
operator, we obtain     
\begin{displaymath}
\omega_{s}\in E^{2}_{\overline{X},\mathbb{R}}(1)\cap E^{1,1}_
{\overline{X}}=\mathcal{D}^{2}(E_{\overline{X}},1)\subseteq
\mathcal{D}^{2}_{\log}(X,1).
\end{displaymath}
Let $\widetilde{X}$ be an embedded resolution of singularities 
of $\overline{Y}$, and $E$ the pre-image of $\overline{Y}$ in
$\widetilde{X}$. If $E$ is locally described by the equation
$z_{1}\cdots z_{m}=0$, the section $s$ can locally be written 
as $z_{1}^{\alpha_{1}}\cdots z_{m}^{\alpha_{m}}\cdot s'$, where 
$s'$ is a suitable non-vanishing, regular section on $\widetilde
{X}$. This shows  
\begin{displaymath}     
g_{s}\in E^{0}_{\log,\mathbb{R}}(U)=\mathcal{D}^{1}_{\log}(U,1).
\end{displaymath}
The claim that $(\omega_{s},g_{s})$ is a cycle follows from 
the definitions. The Poincar\'e-Lelong formula 
\begin{displaymath}
-2\partial\bar{\partial}[g_{s}]_{X}=[\omega_{s}]-\delta_{\dv(s)} 
\end{displaymath}
together with proposition \ref{prop:30} (i) finally proves that
the pair $(\omega_{s},g_{s})$ represents the class of $\dv(s)$.
\end{proof}

\subsection{Real varieties}
\label{sec:real-varieties}

\nnpar{Real Deligne-Beilinson cohomology of real varieties.} 
Recall that a real variety $X_{\mathbb{R}}$ is a pair $(X_{\mathbb
{C}},F_{\infty})$ with $X_{\mathbb{C}}$ a variety over $\mathbb{C}$
and $F_{\infty}$ an antilinear involution of $X_{\mathbb{C}}$. 
Analogously, a real algebraic manifold is a pair $(X_{\mathbb{C}},
F_{\infty})$ with $X_{\mathbb{C}}$ a complex algebraic manifold.  
 
\begin{notation}
\label{def:19}
Given a real variety $X_{\mathbb{R}}$, and a sheaf $V$ of complex 
vector spaces with a real structure, we will denote by $\sigma$ 
the involution given by
\begin{displaymath}
\omega\longmapsto\overline{F_{\infty}^{\ast}(\omega)}. 
\end{displaymath}
We will use the same notation for any subsheaf of abelian groups 
of $V$, which is invariant under complex conjugation.
\end{notation}

Let  $W_{\mathbb{R}}$ now be a proper real algebraic manifold, 
and $D_{\mathbb{R}}$ a normal crossing divisor in $W_{\mathbb{R}}$
defined over $\mathbb{R}$. We put $X_{\mathbb{R}}=W_{\mathbb{R}}
\setminus D_{\mathbb{R}}$. The antilinear involution $\sigma$ of
$E^{\ast}_{W_{\mathbb{C}}}(\log D_{\mathbb{C}})$ respects the real 
structure and the Hodge filtration. It thus induces an involution 
of $\mathcal{D}^{\ast}(E_{W_{\mathbb{C}}}(\log D_{\mathbb{C}}),p)$, 
which we denote again by $\sigma$. The real Deligne-Beilinson 
cohomology groups of $X_{\mathbb{R}}$ are defined as the cohomology 
of the complex of fixed elements $\mathcal{D}^{\ast}(E_{W_{\mathbb{C}}}
(\log D_{\mathbb {C}}),p)^{\sigma}$, i.e.,
\begin{displaymath}
H^{\ast}_{\mathcal{D}}(X_{\mathbb{R}},\mathbb{R}(p))=H^{\ast}
(\mathcal{D}(E_{W_{\mathbb{C}}}(\log D_{\mathbb{C}}),p)^{\sigma}).
\end{displaymath}
An analogous definition can be given for the real Deligne-Beilinson 
cohomology groups with support. From the corresponding result for 
complex varieties, we obtain

\begin{theorem}
Real Deligne-Beilinson cohomology is a Gillet cohomology for 
regular real schemes.
\hfill $\square$
\end{theorem}

\nnpar{A Gillet complex for real Deligne-Beilinson cohomology.} 
Here we denote by $C$ the site of regular schemes in $\ZAR
(\Spec(\mathbb{R}))$. A scheme $X$ in $C$ defines a real
algebraic manifold $X_{\mathbb{R}}=(X_{\mathbb{C}},F_{\infty})$.

\begin{definition} 
\label{def:18b}
For any integers $n,p$, let $\mathcal{D}^{n}_{\log}(p)$ denote 
the presheaf (in fact, the sheaf) over $C$, which assigns to $X$ the
group  
\begin{displaymath}
\mathcal{D}^{n}_{\log}(X,p)=\mathcal{D}^{n}(E_{\log}(X_
{\mathbb{R}}(\mathbb{C})),p)^{\sigma}
\end{displaymath}
with $\sigma$ as in notation \ref{def:19}. For any scheme $X$ 
in $C$, we will denote the induced presheaf of graded complexes 
of real vector spaces on $X$ by $\mathcal{D}_{\log,X}=\mathcal
{D}^{\ast}_{\log,X}(\ast)$.
\end{definition}

\begin{proposition} 
\label{prop:18}
Let $X$ be a scheme in $C$. For any integers $n,p$, the presheaf 
$\mathcal{D}^{n}_{\log,X}(p)$ is a totally acyclic sheaf. 
\end{proposition}
\begin{proof}
Since we are working with complexes of vector spaces, the functor
$(\cdot)^{\sigma}$ is exact. Therefore, the result follows from
proposition \ref{prop:24}.
\end{proof}

\begin{theorem} 
\label{thm:26b}
The graded complex of sheaves of abelian groups $\mathcal{D}_
{\log}$ is a Gillet complex for regular schemes over $\mathbb{R}$, 
which computes real Deligne-Beilinson cohomology. Moreover, the 
pair $(\mathcal{D}_{\log},\bullet)$ is a graded commutative and
pseudo-associative algebra for real Deligne-Beilinson cohomology.
\end{theorem}
\begin{proof}
This is a consequence of theorem \ref{thm:26} and the fact that 
all operations are compatible with the involution $\sigma$.
\end{proof}

\newpage
\section{Examples of $\mathcal{D}_{\log}$-arithmetic Chow groups}

In this section we will use the abstract theory of arithmetic 
Chow groups to define contravariant and covariant arithmetic Chow 
groups. The former were introduced in \cite{Burgos:CDB}; they have 
a ring structure after tensoring with $\mathbb{Q}$ and agree with 
the arithmetic Chow groups defined by Gillet and Soul\'e in
\cite{GilletSoule:ait} for arithmetic varieties with projective 
generic fiber. The latter were introduced in \cite{Burgos:acr}; 
they are covariant for arbitrary proper morphisms and have a 
module structure over the contravariant Chow groups.

In all the examples in this paper the complex $\mathcal{D}_{\log}$
will play the role of the Gillet complex, i.e., $\Gi=\mathcal{D}_
{\log}$.

\subsection{Contravariant $\mathcal{D}_{\log}$-arithmetic Chow rings}
\label{sec:arithm-chow-groups}

\nnpar{Definition and exact sequences.}
Let $A$ be an arithmetic ring. A natural example of a $\mathcal{D}_
{\log}$-complex is $\mathcal{D}_{\log}$ itself by the identity
morphism. Then, all the properties of $\mathcal{D}_{\log}$ as a 
Gillet complex (multiplicativity, functoriality) imply the same 
properties for $\mathcal{D}_{\log}$ as a $\mathcal{D}_{\log}$-complex.

For any arithmetic variety $X$ over $A$, we put $\Gi=\mathcal{D}_
{\log}$ and $\cc=\mathcal{D}_{\log,X}$ in definition \ref{def:20}.  
Thus the pair $(X,\mathcal{D}_{\log})$ is a $\mathcal{D}_{\log}
$-arithmetic variety over $A$. By means of definition \ref{def:gacg},
we obtain the arithmetic Chow groups $\cha^{p}(X,\mathcal{D}_{\log})$. 
We recover the properties of these groups by applying the theory 
developed in section \ref{sec:AC}; these properties have already 
been stated in \cite{Burgos:CDB}. We start by writing down the exact
sequences of theorem \ref{thm:15}. For this we recall
\begin{align*}
\widetilde{\mathcal{D}}^{2p-1}_{\log}(X,p)&=\mathcal{D}^{2p-1}_
{\log}(X,p)\left/\Img(\dd_{\mathcal{D}})\right. \\
&=\left\{\omega\in E^{p-1,p-1}_{\log,\mathbb{R}}(X_{\infty})
(p-1)\mid F^{\ast}_{\infty}(\omega)=\bar{\omega}\right\}\left/
(\Img\partial+\Img\bar{\partial})\right..
\end{align*}

\begin{theorem}
For an arithmetic variety $X$ over $A$, there are exact sequences:
\begin{align}
&\CH^{p-1,p}(X)\overset{\rho}{\longrightarrow}\widetilde 
{\mathcal{D}}^{2p-1}_{\log}(X,p)\overset{\amap}{\longrightarrow} 
\cha^{p}(X,\mathcal{D}_{\log})\overset{\zeta}{\longrightarrow}
\CH^{p}(X)\longrightarrow 0, 
\label{eq:13} \\[3mm]
&\CH^{p-1,p}(X)\overset{\rho}{\longrightarrow}H^{2p-1}_{\mathcal
{D}}(X_{\mathbb{R}},\mathbb{R}(p))\overset{\amap}{\longrightarrow}
\cha^{p}(X,\mathcal{D}_{\log})\overset{(\zeta,-\omega)}
{\longrightarrow} 
\notag \\ 
&\phantom{CH^{p-1,p}(X)\overset{\rho}{\longrightarrow}}\CH^{p}(X) 
\oplus{\rm Z}\mathcal{D}^{2p}_{\log}(X,p)\overset{\cl+h}
{\longrightarrow}H^{2p}_{\mathcal{D}}(X_{\mathbb{R}},\mathbb{R}(p)) 
\longrightarrow 0, 
\label{eq:14} \\[3mm]
&\CH^{p-1,p}(X)\overset{\rho}{\longrightarrow}H^{2p-1}_
{\mathcal{D}}(X_{\mathbb{R}},\mathbb{R}(p))\overset{\amap}
{\longrightarrow}\cha^{p}(X,\mathcal{D}_{\log})_{0}\overset
{\zeta}{\longrightarrow}\CH^{p}(X)_{0}\longrightarrow 0.
\label{eq:24}
\end{align}
\hfill $\square$
\end{theorem}

\nnpar{Green forms for a cycle.} 
We now translate the result of proposition \ref{prop:30} into 
the language of Green objects.  

\begin{proposition}
\label{prop:28} 
Let $X$ be an arithmetic variety over $A$, and $y$ a $p$-codimen\-sional 
cycle on $X$. A pair $\mathfrak{g}=(\omega,\widetilde{g})\in\widehat
{H}^{2p}_{\mathcal{D}_{\log},\mathcal{Z}^{p}}(X,p)$ is a Green object
for the class of $y$, if and only if 
\begin{equation}
\label{prop:31}
-2\partial\bar{\partial}[g]_{X}=[\omega]-\delta_{y}. 
\end{equation}
\hfill $\square$
\end{proposition}  

\noindent
To simplify notations, in the above proposition, we wrote $[g]_{X}$ 
instead of $[g]_{X_{\infty}}$, and $\delta_{y}$ instead of $\delta_
{y_{\mathbb{R}}}$. We will call Green objects with values in
$\mathcal{D}_{\log}$ \emph{Green forms}.  

\nnpar{The arithmetic cycle associated to a rational function.} 
We make the Green object and the arithmetic cycle associated
to a rational function explicit.
 
\begin{proposition}
Let $f\in k^{\ast}(X)$, $y=\dv(f)$, $Y=\supp y$, and $U=X\setminus 
Y$. Then, the class $\cl_{\mathcal{D}}(f)$ of $f$ in the group
$H^{1}_{\mathcal{D}_{\log}}(U,1)=H^{1}_{\mathcal{D}}(U_{\mathbb{R}},
\mathbb{R}(1))$ is represented by the function 
\begin{displaymath}
\frac{1}{2}\log(f\bar{f})\in\mathcal{D}^{1}_{\log}(U,1); 
\end{displaymath}
we note that in order to ease notation, we have written $\cl_
{\mathcal{D}}(f)$ instead of $\cl_{\mathcal{D}_{\log}}(f)$.
Furthermore, we have
\begin{displaymath}
\diva(f)=\left(\dv(f),\left(0,-\frac{1}{2}\log(f\bar{f})\right)
\right).
\end{displaymath}
\end{proposition}
\begin{proof} 
Since we are interested in the cohomology of $U$, we may assume 
that $Y$ is a normal crossing divisor; we denote the inclusion by $j:U 
\longrightarrow X$. Then, $f$ is a global section 
of $\mathcal{O}^{\times}_{\alg}$ over $U$.

To prove the first assertion, we have to show that the morphism 
\begin{displaymath}
H^{0}(U,\mathcal{O}^{\times}_{\alg})\longrightarrow H^{1}_
{\mathcal{D}}(U_{\mathbb{R}},\mathbb{R}(1))\,,
\end{displaymath}
given by mapping $f$ to $\frac{1}{2}\log(f\bar{f})$, is compatible 
with the map ${\rm c}_{1}$ of (\ref{eq:15}). The key ingredient 
for this compatibility is provided by the commutative diagram of 
sheaves on $U$  
\begin{equation}
\label{eq:43}
\xymatrix{ 
s(\mathbb{Z}(1)\rightarrow\mathcal{O}_{\alg})\ar[r]^-{e}\ar[d]
&\mathcal{O}^{\times}_{\alg}[-1]\ar[d]^{l} \\
s(\mathscr{E}^{0}_{U,\mathbb{R}}(1)\rightarrow\mathscr{E}^{0}_{U})
\ar[r]^-{\pi}&\mathscr{E}^{0}_{U,\mathbb{R}}[-1]\,,}
\end{equation}
where $e(a,b)=\exp(b)$, $l(f)=\frac{1}{2}\log(f\bar{f})$, and 
$\pi(a,b)=\frac{1}{2}(b+\bar{b})$. Denoting by $G$ the complex 
of sheaves on $U$ given by
\begin{displaymath}
\mathcal{O}^{\times}_{\alg}\overset{\dd\log}{\longrightarrow}
\Omega^{1}_{U}\longrightarrow\Omega^{2}_{U}\longrightarrow\dots,
\end{displaymath}
the commutative diagram (\ref{eq:43}) induces a commutative
diagram of sheaves on $X$
\begin{displaymath}
\xymatrix{ 
s(F^{1}\Omega^{\ast}_{X}(\log Y)\rightarrow Rj_{\ast}s(\mathbb
{Z}(1)\rightarrow\Omega^{\ast}_{U}))\ar[r]^-{e}\ar[d]&s(F^{1}
\Omega^{\ast}_{X}(\log Y)\rightarrow Rj_{\ast}G)\ar[d]^{l} \\
s(F^{1}\mathscr{E}^{\ast}_{X}(\log Y)\rightarrow Rj_{\ast}
s(\mathscr{E}^{\ast}_{U,\mathbb{R}}(1)\rightarrow\mathscr{E}^
{\ast}_{U}))\ar[r]^-{\pi}&s(F^{1}\mathscr{E}^{\ast}_{X}(\log Y)
\rightarrow Rj_{\ast}\mathscr{E}^{\ast}_{U,\mathbb{R}})\,.}
\end{displaymath}
This commutative diagram shows that the class $\cl_{\mathcal{D}}
(f)={\rm c}_{1}(f)$ is represented by the pair $(\dd\log f,\frac
{1}{2}\log(f\bar{f}))$ in the complex $s(F^{1}\mathscr{E}^{\ast}_
{X}(\log Y)\rightarrow Rj_{\ast}\mathscr{E}^{0}_{U,\mathbb{R}})$. 
This implies that the class $\cl_{\mathcal{D}}(f)$ is represented
by $\frac{1}{2}\log(f\bar{f})$ in the complex $\mathcal{D}^{\ast}_
{\log}(U,1)$.

The second assertion follows from the first and the fact that
the map $\bmap$ is given by $\bmap(b)=(0,-b)$ (see definition
\ref{def:8}).
\end{proof}

\nnpar{Multiplicative properties.}
By theorem \ref{thm:26b}, the complex $\mathcal{D}_{\log}$
satisfies all the properties required to apply theorem
\ref{thm:10}. In particular, we obtain

\begin{theorem}
The direct sum
\begin{displaymath}
\cha^{\ast}(X,\mathcal{D}_{\log})_{\mathbb{Q}}=\bigoplus_{p}
\cha^{p}(X,\mathcal{D}_{\log})_{\mathbb{Q}}
\end{displaymath}
has the structure of a commutative and associative $\mathbb{Q}$-algebra 
with unit. This structure is compatible with the maps $\omega $, $\cl$ 
and $\zeta $; furthermore, it is compatible with the algebra structures 
of the Deligne complex, the Deligne-Beilinson cohomology and the Chow   
ring. 
\hfill $\square$
\end{theorem}

Moreover, we observe that theorem \ref{thm:17} and corollary \ref{cor:4}
also apply to the present situation. We leave it to the reader to write 
down the corresponding statements.

\nnpar{Formulas for the $*$-product.}  
Using partitions of unity we will give alternative formulas for the 
$*$-product of Green objects. These intermediate formulas will be 
useful for the explicit computations carried out in our applications. 

Let $y$ be a $p$-codimensional, resp. $z$ a $q$-codimensional cycle
of $X$, and $Y=\supp y_{\mathbb{R}}$, resp. $Z=\supp z_{\mathbb{R}}$. 
Let $\mathfrak{g}_{y}=(\omega_{y},\widetilde{g}_{y})\in\widehat{H}^
{2p}_{\mathcal{D}_{\log},Y}(X,p)$, resp. $\mathfrak{g}_{z}=(\omega_
{z},\widetilde{g}_{z})\in\widehat{H}^{2q}_{\mathcal{D}_{\log},Z}
(X,q)$ be Green objects for $y$, resp. $z$.

\begin{proposition}
\label{prop:explicit}
\begin{enumerate}
\item[(i)]
\label{prop:explicit1} 
If $\mathfrak{g}_{y}=\amap(\widetilde{x})$ for some $\widetilde{x}
\in\widetilde{\mathcal{D}}^{2p-1}_{\log}(X,p)$, we have the
equality
\begin{displaymath}
\mathfrak{g}_{y}*\mathfrak{g}_{z}=\amap(\widetilde{x\land\omega_{z}})
\end{displaymath}
in the group $\widehat{H}^{2p+2q}_{\mathcal{D}_{\log},Y\cap Z}(X,p+q)$.
\item[(ii)] 
\label{prop:explicit2} 
If $Y=Z$, we have the equality
\begin{displaymath}
\mathfrak{g}_{y}*\mathfrak{g}_{z}=(\omega_{y}\land\omega_{z},
\widetilde{g_{y}\land\omega_{z}})
\end{displaymath}
in the group $\widehat{H}^{2p+2q}_{\mathcal{D}_{\log},Y}(X,p+q)$.
\end{enumerate}
\end{proposition} 
\begin{proof} 
The first claim is an immediate consequence of proposition 
\ref{prop:restrlem} (iii) and the explicit formulas for the 
product in Deligne algebras given after remark \ref{remark:fulea}. 
The second claim is covered by the explicit formula for the 
$*$-product given in proposition \ref{pro:prodtrun}, again
taking into account the formulas for the product in Deligne 
algebras given after remark \ref{remark:fulea}.
\end{proof}

The cases $Y\subsetneq Z$ and $Z\subsetneq Y$ can be treated 
similarly. But the most interesting case is when neither $Y$ 
nor $Z$ are contained in $Y\cap Z$. In this case, we put $U=
X\setminus Y$ and $V=X\setminus Z$. Then, the Mayer-Vietoris 
sequence 
\begin{displaymath}
0\to\mathcal{D}^{n}_{\log}(U\cup V,p)\overset{i}{\to}\mathcal
{D}^{n}_{\log}(U,p)\oplus\mathcal{D}^{n}_{\log}(V,p)\overset{j}
{\to}\mathcal{D}^{n}_{\log}(U\cap V,p)\to 0,
\end{displaymath}
where $i(\eta)=(\eta,\eta)$ and $j(\omega,\eta)=-\omega+
\eta$, gives rise to the kernel-simple quasi-isomorphism
\begin{displaymath}
\iota:\mathcal{D}_{\log}(U\cup V,p)\longrightarrow s(-j),
\end{displaymath} 
which was an essential tool in the proof of theorem \ref{thm:ap-com}. 
We will now construct a section of $j$. Adapting the argument given
in \cite{Burgos:Gftp}, we can find a resolution of singularities 
$\pi:\widetilde{X}\longrightarrow X_{\mathbb{R}}$ of $Y\cup Z$, 
which factors through embedded resolutions of $Y$, $Z$, $Y\cap Z$. 
In particular, we can assume that
\begin{displaymath}
\pi^{-1}(Y),\,\pi^{-1}(Z),\,\pi^{-1}(Y\cap Z)
\end{displaymath}
are also normal crossing divisors. We denote by $\widehat{Y}$ the
normal crossing divisor formed by the components of $\pi^{-1}(Y)$
which are not contained in $\pi^{-1}(Y\cap Z)$; analogously we 
denote by $\widehat{Z}$ the normal crossing divisor formed by the
components of $\pi^{-1}(Z)$ which are not contained in $\pi^{-1}
(Y\cap Z)$. Then, $\widehat{Y}$ and $\widehat{Z}$ are closed subsets 
of $\widetilde{X}$ which do not meet. Therefore, there exist two 
smooth, $F_{\infty}$-invariant functions $\sigma_{_{YZ}}$ and 
$\sigma_{_{ZY}}$ satisfying $0\leq\sigma_{_{YZ}},\sigma_{_{ZY}}
\leq 1$, $\sigma_{_{YZ}}+\sigma_{_{ZY}}=1$, $\sigma_{_{YZ}}=1$ 
in a neighborhood of $\widehat{Y}$, and $\sigma_{_{ZY}}=1$ in 
a neighborhood of $\widehat{Z}$. Let now $\omega\in E^{n}_{\log}
(U\cap V)$. Since $\sigma_{_{YZ}}$ is zero in a 
neighborhood of $\widehat{Z}$, we find that $\sigma_{_{YZ}}
\omega\in E^{n}_{\log}(U)$; similarly, we get $\sigma_{_{ZY}}
\omega\in E^{n}_{\log}(V)$. Moreover, one easily checks that 
$j(-\sigma_{_{YZ}}\omega,\sigma_{_{ZY}}\omega)=\omega$. Therefore,
the assignment $\omega\mapsto(-\sigma_{_{YZ}}\omega,\sigma_{_{ZY}}
\omega)$ gives rise to a section of $j$. Consequently, the map
$\mathcal{D}^{n}_{\log}(U\cap V,p)\longrightarrow\mathcal{D}^{n}_
{\log}(U,p)\oplus\mathcal{D}^{n}_{\log}(V,p)$ given by $c\mapsto 
(-\sigma_{_{YZ}}c,\sigma_{_{ZY}}c)$ determines a section of $j$.
We are now in a position to apply proposition \ref{prop:quisinv} 
(ii).
  
\begin{lemma}
\label{lem:quisinvdb}
The map $s(-j)\longrightarrow\mathcal{D}_{\log}(U\cup V,p)$ 
given by 
\begin{displaymath}
((a,b),c)\longmapsto\sigma_{_{ZY}}a+\sigma_{_{YZ}}b+\sigma_
{_{YZ}}\dd_{\mathcal{D}}c-\dd_{\mathcal{D}}(\sigma_{_{YZ}}c)
\end{displaymath}
is a morphism of complexes. It is a left inverse of the 
kernel-simple quasi-isomorphism $\iota$.
\hfill $\square$
\end{lemma}

\begin{theorem}
\label{thm:partition} 
Let $\mathfrak{g}_{y}=(\omega_{y},\widetilde{g}_{y})$, $\mathfrak
{g}_{z}=(\omega_{z},\widetilde{g}_{z})$, and $\sigma_{_{YZ}}$, 
$\sigma_{_{ZY}}$ be as above. In the group $\widehat{H}^{2p+2q}
_{\mathcal{D}_{\log},Y\cap Z}(X,p+q)$, we then have the identity 
\begin{align*}     
\mathfrak{g}_{y}*\mathfrak{g}_{z}=\left(\omega_{y}\land\omega_{z},
(-2\sigma_{_{ZY}}g_{y}\land\partial\bar{\partial}g_{z}-2\partial\bar
{\partial}(\sigma_{_{YZ}}g_{y})\land g_{z})^{\widetilde{\phantom{=}}}
\right)\,.  
\end{align*}
\end{theorem}
\begin{proof} 
By theorem \ref{thm:wstp}, we obtain
\begin{displaymath}
\mathfrak{g}_{y}*\mathfrak{g}_{z}=\left(\omega_{y}\bullet\omega_{z},
((g_{y}\bullet\omega_{z},\omega_{y}\bullet g_{z}),-g_{y}\bullet g_
{z})^{\widetilde{\phantom{=}}}\right) 
\end{displaymath}
in the group $\widehat{H}^{2p+2q}(\mathcal{D}_{\log}(X),s(-j),p+q)$. 
By means of lemma \ref{lem:quisinvdb}, the latter element corresponds 
to   
\begin{displaymath}
\left(\omega_{y}\bullet\omega_{z},(\sigma_{_{ZY}}g_{y}\bullet\dd_
{\mathcal{D}}g_{z}+\dd_{\mathcal{D}}(\sigma_{_{YZ}}g_{y})\bullet 
g_{z})^{\widetilde{\phantom{=}}}\right)\in\widehat{H}^{2p+2q}_
{\mathcal{D}_{\log},Y\cap Z}(X,p+q).
\end{displaymath}
The stated formula follows now from the explicit description
of $\bullet$ and $\dd_{\mathcal{D}}$ in these degrees.
\end{proof}

\nnpar{Inverse images.} 
Let $f:X\longrightarrow Y$ be a morphism of arithmetic varieties
over $A$. Since $\mathcal{D}_{\log}$ is a sheaf in the big Zariski 
site of smooth schemes over $\mathbb{R}$, there exists a contravariant
$f$-morphism  
\begin{displaymath}
f^{\#}:\mathcal{D}_{\log,Y}\longrightarrow f_{\ast}\mathcal{D}_
{\log,X}.
\end{displaymath}
Thus, the complex $\mathcal{D}_{\log}$ satisfies all the properties 
required to apply theorems \ref{thm:6} and \ref{thm:4.30}. In 
particular, we obtain

\begin{theorem}
Let $f:X\longrightarrow Y$ be a morphism of arithmetic varieties 
over $A$. Then, there is a pull-back morphism  
\begin{displaymath}
f^{\ast}:\cha^{p}(Y,\mathcal{D}_{\log})\longrightarrow\cha^{p}
(X,\mathcal{D}_{\log}),
\end{displaymath}
which is compatible with the pull-back of differential forms via
the morphism $\omega$, and with the pull-back
of algebraic cycles via the morphism $\zeta$.

If $g:Y\longrightarrow Z$ is another such morphism, then the equality  
$(g\circ f)^{\ast}=f^{\ast}\circ g^{\ast}$ holds. 

Moreover, the induced map
\begin{displaymath}
f^{\ast}:\cha^{p}(Y,\mathcal{D}_{\log})_{\mathbb{Q}}\longrightarrow 
\cha^{p}(X,\mathcal{D}_{\log})_{\mathbb{Q}}    
\end{displaymath}
is an algebra morphism.
\hfill $\square$
\end{theorem}

\nnpar{Direct images of differential forms with logarithmic
singularities.} 
Differential forms with logarithmic singularities at infinity are 
not well suited to define covariant morphisms for $\mathcal{D}_
{\log}$-complexes, as the following example shows.

\begin{example}
Let $f:C\longrightarrow C'$ be a morphism of smooth complex
projective curves, and let $\Sigma$ be the pre-image of the 
singular values of the map $f$. Assume that there is a point 
$P\in\Sigma$ such that the morphism $f$ is given by the local 
expression $w=f(z)=z^{2}$ in an analytic neighborhood of $P$.
While the differential form $\dd z\land\dd\bar{z}$ is smooth 
in a neighborhood of $P$, its push-forward $f_{\ast}(\dd z
\land\dd\bar{z})$ is locally of the form $\dd w\land\dd\bar
{w}/(w\bar{w})^{1/2}$, which is not smooth and does not have 
logarithmic singularities. This shows that, even if the 
restriction
\begin{displaymath}
f:C\setminus\Sigma\longrightarrow C'\setminus f(\Sigma ) 
\end{displaymath}
is smooth, it does not induce a morphism between $E^{\ast}_
{\log}(C\setminus\Sigma)$ and $E^{\ast}_{\log}(C'\setminus 
f(\Sigma))$.
\end{example}

Nevertheless, if $f:X\longrightarrow Y$ is a smooth proper
equidimensional morphism of projective varieties over $\mathbb
{R}$ of relative dimension $e$, we can define a covariant 
$f$-pseudo-morphism as follows. For every open subset $U$ 
of $Y$ we define
\begin{displaymath}
\mathcal{F}^{n}(U,p)=\lim_{\substack{\longrightarrow \\ 
\overline{Y},D}}\,\mathcal{D}^{n}\left(E_{X\underset{Y}
{\times}\overline{Y}}(\log X\underset{Y}{\times}D),p\right),
\end{displaymath}
where the limit is taken over all compactifications $\overline{Y}$ 
of $U$ with $D=\overline{Y}\setminus U$ a normal crossing divisor
such that $\overline{Y}$ dominates $Y$. Since $f$ is smooth, we
note that $X\underset{Y}{\times}D$ is a normal crossing divisor
of the smooth variety $X\underset{Y}{\times}\overline{Y}$. Let  
\begin{displaymath}
u:\mathcal{F}\longrightarrow f_{\ast}\mathcal{D}_{\log,X}
\end{displaymath}
be the natural morphism, and let 
\begin{displaymath}
v:\mathcal{F}\longrightarrow\mathcal{D}_{\log,Y}(-e)[-2e]
\end{displaymath}
be the morphism induced by integration along the fiber, i.e., 
\begin{displaymath}
v(\omega)=f_{!}\omega=\frac{1}{(2\pi i)^{e}}\int_{f}\omega\,.
\end{displaymath}
By proposition \ref{prop:33} the morphisms $u$ and $v$ 
determine a covariant $f$-pseudo-morphism. Using definition 
\ref{def:22} together with remarks \ref{rem:4}, \ref{rem:5}, 
and proposition \ref{prop:20} we obtain

\begin{theorem} 
\label{thm:properpushlog}
Let $f:X\longrightarrow Y$ be a smooth proper morphism of 
proper varieties over $\mathbb{R}$ of relative dimension 
$e$. Then, there is a push-forward morphism
\begin{displaymath}
f_{\#}:\widehat{H}^{2p}_{\mathcal{D}_{\log},\mathcal{Z}^{p}}
(X,p)\longrightarrow\widehat{H}^{2p-2e}_{\mathcal{D}_{\log},
\mathcal{Z}^{p-e}}(Y,p-e),
\end{displaymath}
which is compatible with the push-forward of differential 
forms via the morphism $\omega$, and 
with the direct image into relative Deligne-Beilinson 
cohomology via the morphism $\cl$.

If $g:Y\longrightarrow Z$ is another such morphism, then 
the equality $(g\circ f)_{\#}=g_{\#}\circ f_{\#}$ holds. 

Moreover, if $\alpha\in\widehat{H}^{2p}_{\mathcal{D}_{\log},
\mathcal{Z}^{p}}(Y,p)$ and $\beta\in\widehat{H}^{2q}_{\mathcal
{D}_{\log},\mathcal{Z}^{q}}(X,q)$, we have the formula
\begin{displaymath} 
f_{\#}(f^{\#}(\alpha)\cdot\beta)=\alpha\cdot f_{\#}(\beta)\in
\widehat{H}^{2p+2q-2e}_{\mathcal{D}_{\log},\mathcal{Z}^{p+q-e}}
(Y,p+q-e).
\end{displaymath} 
\hfill $\square$ 
\end{theorem}

\begin{remark}
In order to obtain more general push-forward morphisms one has to
use different complexes. For instance, using a real Deligne-Beilinson 
complex made with differential forms satisfying logarithmic growth 
conditions instead of logarithmic singularities at infinity, we 
expect that one can avoid the projectivity condition for the 
varieties under consideration. Another option is the use of the 
real Deligne-Beilinson complex with currents as carried out in the 
next section; then, one obtains direct images for arbitrary proper 
morphisms, at the price of losing some multiplicativity properties.
\end{remark}

\nnpar{Direct images of contravariant arithmetic Chow rings.}
Using the technique of covariant pseudo-morphisms as before, we 
can define a push-forward morphism associated to a morphism between 
arithmetic varieties, which are generically projective and smooth. 
As we will see below, this suffices to define arithmetic degrees. The 
following result is a consequence of theorem \ref{thm:properpushlog} 
and theorem \ref{thm:18}.

\begin{theorem}
Let $f:X\longrightarrow Y$ be a proper morphism of arithmetic 
varieties over $A$ of relative dimension $e$. Assume that $f_
{\mathbb{R}}:X_{\mathbb{R}}\longrightarrow Y_{\mathbb{R}}$ is
a smooth proper morphism of proper varieties. Then, there is 
a push-forward morphism
\begin{displaymath}
f_{*}:\cha^{p}(X,\mathcal{D}_{\log})\longrightarrow\cha^{p-e}(Y,
\mathcal{D}_{\log}),
\end{displaymath}
which is compatible with the push-forward of differential forms 
via the morphism $\omega$, and with the
push-forward of algebraic cycles via the morphism $\zeta $.

If $g:Y\longrightarrow Z$ is another such morphism, then the
equality $(g\circ f)_{*}=g_{*}\circ f_{*}$ holds. 

Moreover, if $\alpha\in\cha^{p}(Y,\mathcal{D}_{\log})$ and 
$\beta\in\cha^{q}(X,\mathcal{D}_{\log})$, we have the formula
\begin{displaymath}
f_{*}(f^{*}(\alpha)\cdot\beta)=\alpha\cdot f_{*}(\beta)\in
\cha^{p+q-e}(Y,\mathcal{D}_{\log})_{\mathbb{Q}}.
\end{displaymath}
\hfill $\square$
\end{theorem}

\nnpar{Homotopy invariance.}
The fact that the definition of contravariant arithmetic Chow 
groups for non proper varieties has good Hodge theoretical 
properties is reflected in the homotopy invariance of the groups 
$\cha^{p}(X,\mathcal{D}_{\log})_{0}$.

\begin{proposition}
Let $X$ be an arithmetic variety over $A$, and $\pi:M
\longrightarrow X$ a geometric vector bundle. Then, the 
induced morphism
\begin{displaymath}  
\pi^{*}:\cha^{p}(X,\mathcal{D}_{\log})_{0}\longrightarrow
\cha^{p}(M,\mathcal{D}_{\log})_{0} 
\end{displaymath}
is an isomorphism.
\end{proposition}
\begin{proof}
For a proof we refer to \cite{Burgos:CDB}, 7.5.
\end{proof}

This theorem has the following variant used by H. Hu in his 
Ph.D. thesis to construct an arithmetic intersection pairing 
based on the deformation to the normal cone technique. With 
the hypothesis of the above theorem, we write
\begin{displaymath}
\mathcal{D}^{\ast}_{\log}(M,p)_{\verti}=\pi^{\ast}\left(
\mathcal{D}^{\ast}_{\log}(X,p)\right)\subseteq\mathcal{D}^
{\ast}_{\log}(M,p),
\end{displaymath}
and put
\begin{displaymath}
\cha^{p}(M,\mathcal{D}_{\log})_{\verti}=\omega ^{-1}\left(
\mathcal{D}^{2p}_{\log}(M,p)_{\verti}\right).
\end{displaymath}
Then, we have

\begin{proposition}
Let $X$ be an arithmetic variety over $A$, and $\pi:M
\longrightarrow X$ a geometric vector bundle. Then, the   
induced morphism
\begin{displaymath}  
\pi^{*}:\cha^{p}(X,\mathcal{D}_{\log})\longrightarrow
\cha^{p}(M,\mathcal{D}_{\log})_{\verti} 
\end{displaymath}
is an isomorphism. 
\hfill $\square$
\end{proposition}

\nnpar{Comparison with the arithmetic Chow groups defined 
by Gillet and Soul\'e.} 
Let $X$ be a $d$-dimensional arithmetic variety over $A$. We 
denote by $\cha^{p}(X)$ the arithmetic Chow groups defined by 
Gillet and Soul\'e. Let $y$ be a $p$-codimensional cycle of 
$X$, $Y=\supp y_{\mathbb{R}}$, $U=X_{\mathbb{R}}\setminus Y$, 
and 
\begin{displaymath}
\mathfrak{g}_{y}=(\omega_{y},\widetilde{g}_{y})\in\widehat{H}^
{2p}\left(\mathcal{D}_{\log}(X_{\mathbb{R}},p),\mathcal{D}_
{\log}(U,p)\right)
\end{displaymath}
a Green object for the cycle $y$. 

\begin{lemma}
The current $2(2\pi i)^{d-p+1}[g_{y}]_{X}$ is a Green current 
in the sense of \cite{GilletSoule:ait} for the cycle $y$. 
\hfill $\square$
\end{lemma}
\begin{proof} 
Since the definitions of the current associated to a differential 
form and the current associated to a cycle used in this paper 
differ from those in \cite{GilletSoule:ait} by a normalization  
factor, we will write the proof explicitly in terms of integrals.

We may assume that the arithmetic ring $A$ is $\mathbb{R}$ and 
that $y$ is a prime cycle, hence $Y$ is irreducible. We denote
by $\widetilde{Y}$ a resolution of singularities of $Y_{\infty}$,
and by $\imath:\widetilde{Y}\longrightarrow X_{\infty}$ the 
induced map. 

A current $g\in D^{p-1,p-1}(X_{\infty})$ is a Green current for 
the cycle $y$, if the equality
\begin{displaymath}
F_{\infty}^{\ast}(g)=(-1)^{p-1}g
\end{displaymath}
holds, and if there exists a smooth differential form $\omega$ 
satisfying
\begin{equation}
\label{eq:green_current1}
g({\rm dd^{c}}\eta)=\int_{X_{\infty}}\eta\land\omega-\int_
{\widetilde{Y}}\imath^{\ast}\eta
\end{equation}
for any test form $\eta$.

We put $g=2(2\pi i)^{d-p+1}[g_{y}]_{X}$, and fix a test form 
$\eta$. Then, the definition of the current associated to a 
differential form and the relation $(4\pi i){\rm dd^{c}}=
-2\partial\bar{\partial}$ leads to
\begin{displaymath}
g({\rm dd^{c}}\eta)=\frac{1}{(2\pi i)^{p}}\int_{X_{\infty}}
(-2\partial\bar{\partial}\eta)\land g_{y}.
\end{displaymath}
Since $\sigma(g_{y})=\overline{F_{\infty}^{\ast}(g_{y})}=g_{y}$ 
and $\overline{g_{y}}=(-1)^{p-1}g_{y}$, we conclude that $F_
{\infty}^{\ast}(g)=(-1)^{p-1}g$. Therefore, the first condition 
for a Green current is satisfied. Writing out explicitly all 
normalization factors in proposition \ref{prop:28}, we find 
\begin{equation}
\label{eq:green_form1}
\frac{1}{(2\pi i)^{d}}\int_{X_{\infty}}(-2\partial\bar{\partial} 
\eta)\land g_{y}=\frac{1}{(2\pi i)^{d}}\int_{X_{\infty}}\eta\land 
\omega_{y}-\frac{1}{(2\pi i)^{d-p}}\int_{\widetilde{Y}}\imath^
{\ast}\eta.      
\end{equation}
Equation \eqref{eq:green_form1} is now easily seen to imply that 
$g$ satisfies equation \eqref{eq:green_current1} with $\omega=
(2\pi i)^{-p}\omega_{y}$, which concludes the proof of the lemma.
\end{proof}

\begin{theorem}
\label{thm:24}
The assignment $[y,(\omega_{y},\widetilde{g}_{y})]\mapsto
[y,2(2\pi i)^{d-p+1}[g_{y}]_{X}]$ induces a well-defined 
map
\begin{displaymath}
\Psi:\cha^{p}(X,\mathcal{D}_{\log})\longrightarrow\cha^{p}(X),
\end{displaymath}
which is compatible with products and pull-backs. Moreover, 
if $X_{\mathbb{R}}$ is projective, then it is an isomorphism,
which is compatible with push-forwards.
\end{theorem}
\begin{proof}
In the projective case, the proof is given in \cite{Burgos:Gftp}
and \cite{Burgos:CDB}. On the other hand, we note that the
well-definedness of the map and its compatibility with products 
and pull-backs carries over to the quasi-projective case.
\end{proof}

\nnpar{Arithmetic Picard group.} 
Let $X$ be a projective arithmetic variety over $A$, and $L$ a 
line bundle on $X$ equipped with a smooth hermitian metric $h$
on the induced line bundle $L_{\infty}$ over $X_{\infty}$, which 
is invariant under $F_{\infty}$. As usual, we write $\overline{L}=
(L,\Vert\cdot\Vert)$, and refer to it as a \emph{hermitian line 
bundle}. Given a rational section $s$ of $L$, we write $\Vert s
\Vert^{2}=h(s,s)$ for the point-wise norm of the induced section 
of $L_{\infty}$. We say that two hermitian line bundles $\overline
{L}_{1}$ and $\overline{L}_{2}$ are isometric, if $\overline{L}_
{1}\otimes\overline{L}_{2}^{-1}\cong(\mathcal{O}_{X},\vert\cdot
\vert)$, where $\vert\cdot\vert$ is the standard absolute value. 
The \emph{arithmetic Picard group $\pica(X)$} is the group of 
isometry classes of hermitian line bundles with group structure 
given by the tensor product. One easily proves that there is an 
isomorphism 
\begin{displaymath}
\ca_{1}:\pica(X)\longrightarrow\cha^{1}(X,\mathcal{D}_{\log}),
\end{displaymath}
given by sending the class of $\overline{L}$ to the class 
$[\dv(s),(\omega_{s},\widetilde{g}_{s})]$, where $s$, $\omega_
{s}$, and $g_{s}$ are as in proposition \ref{prop:26}. We call 
the element
\begin{displaymath}
\ca_{1}(\overline{L})\in\cha^{1}(X,\mathcal{D}_{\log})
\end{displaymath}
the \emph{first arithmetic Chern class of $\overline{L}$}. 

\nnpar{Arithmetic degree map.}
Let $K$ be a number field, and $\mathcal{O}_{K}$ its ring of 
integers. According to \cite{GilletSoule:ait}, $\mathcal{O}_{K}$ 
can be viewed as an arithmetic ring in a canonical way. Putting
$S=\Spec(\mathcal{O}_{K})$, theorem \ref{thm:24} provides an 
isomorphism 
\begin{displaymath}
\Psi:\cha^{1}(S,\mathcal{D}_{\log})\overset{\cong}{\longrightarrow}
\cha^{1}(S).
\end{displaymath}
In particular, the computations of \cite{GilletSoule:ait}, 3.4.3,
carry over to this case with some minor changes. Since $S_{\infty}$
consists of a finite number of points, we have
\begin{displaymath}
\widehat{H}^{2}_{\mathcal{D}_{\log},\mathcal{Z}^{1}}(S,1)= 
\left(\sum\limits_{\sigma\in\Sigma}\mathbb{R}(0)\right)^{+}=
\mathbb{R}^{r_{1}+r_{2}}, 
\end{displaymath}
where $^{+}$ means the invariants under complex conjugation 
acting on the set of the $r_{1}$ real, resp. $2r_{2}$ complex 
immersions of $K$. The exact sequence \eqref{eq:13} therefore
specializes to the exact sequence
\begin{displaymath}
1\longrightarrow\mu(K)\longrightarrow\mathcal{O}^{\ast}_{K}\overset
{\rho}{\longrightarrow}\mathbb{R}^{r_{1}+r_{2}}\overset{\amap}
{\longrightarrow}\cha^{1}(S,\mathcal{D}_{\log})\longrightarrow 
\Cl(\mathcal{O}_{K})\longrightarrow 0;
\end{displaymath}
here $\mu(K)$ is the group of roots of unity of $K$, $\rho$ the 
Dirichlet regulator map, and $\Cl(\mathcal{O}_{K})$ the ideal class 
group of $K$. We observe that, with the normalizations used in this 
paper, the map $\rho$ equals the Dirichlet regulator map and there 
is no factor of $(-2)$ as in \cite{GilletSoule:ait}. The discrepancy 
by the factor $2$ is explained by theorem \ref{thm:24}; the discrepancy 
about the sign is due to the fact that our map $\amap$ is minus the 
corresponding map in \cite{GilletSoule:ait}. 

Due to the product formula for the valuations of $K$, we have 
as in \cite{GilletSoule:ait} a well-defined \emph{arithmetic 
degree map}
\begin{align}
\label{eq:ar-deg}
\dega:\cha^{1}(S,\mathcal{D}_{\log})\longrightarrow\mathbb{R},
\end{align}
induced by the assignment
\begin{displaymath}
\left(\sum_{\mathfrak{p_{j}}\in S}n_{j}\mathfrak{p}_{j}, 
\sum_{\sigma\in\Sigma}(0,\widetilde{g}_{\sigma})\right)
\mapsto\sum_{\mathfrak{p_{j}}\in S}n_{j}\log\big|\mathcal
{O}_{K}\big/\mathfrak{p}_{j}\big|+\sum_{\sigma\in\Sigma} 
g_{\sigma}.
\end{displaymath}
In particular, this is a group homomorphism, which is an 
isomorphism in the case $K=\mathbb{Q}$; it is common to 
identify $\cha^{1}(\Spec(\mathbb{Z}),\mathcal{D}_{\log})$ 
with $\mathbb{R}$. We note that in spite of the many different
normalizations used in this paper, the arithmetic degree
map defined above is compatible with the arithmetic degree 
map defined in \cite{GilletSoule:ait} under the isomorphism
$\Psi$ of theorem \ref{thm:24}.


\subsection{Covariant $\mathcal{D}_{\log}$-arithmetic Chow groups}
\label{sec:gi-complex-from}

In this section we will use currents to define arithmetic 
Chow groups which are covariant for arbitrary proper morphisms. 
These groups do not have a ring structure, but they are modules 
over the contravariant arithmetic Chow ring defined in the 
previous section. We will grade the cycles by their codimension, 
and hence we will use the cohomological complexes of currents 
$\mathscr{D}^{\ast}_{X}$; this forces us to restrict ourselves 
to equidimensional varieties $X$. Note however that we could 
also define covariant arithmetic Chow groups indexed by the 
dimension as in \cite{Burgos:acr} avoiding this restriction.

\nnpar{Currents with good Hodge properties.}
We start by noting that, if $Y$ is a normal crossing divisor 
on a proper equidimensional complex algebraic manifold $X$, 
then the complex $\mathscr{D}^{\ast}_{X/Y}$ has the same 
cohomological properties as the complex $\mathscr{E}^{\ast}_
{X}(\log Y)$. Nevertheless, we cannot use this complex of 
currents, since it does not form a presheaf. For instance, 
if $Y$ is a closed subvariety of $X$, and $\pi:\widetilde{X}
\longrightarrow X$ an embedded resolution of singularities 
of $Y$ with exceptional divisor $D$, there does not exist 
a natural morphism $\mathscr{D}^{\ast}_{X}\longrightarrow
\mathscr{D}^{\ast}_{\widetilde{X}/D}$, in general. For this 
reason the theory of covariant arithmetic Chow groups is not 
fully satisfactory.  

\nnpar{Currents with support on a subvariety.} 
Instead, we will use another complex of currents for quasi-projective 
subvarieties, which forms a presheaf, but has worse Hodge theoretical 
properties. To do this, let $X$ be a proper equidimensional complex
algebraic manifold, and $\mathscr{D}^{\ast}_{X}$ the complex of 
currents on $X$ as in section \ref{sec:deligne-homology}. 

\begin{definition}
\label{def:15}
Let $Y\subseteq X$ be a closed subvariety. Then, the \emph{group 
of currents of degree $n$ on $X$ with support on $Y$} is defined 
by
\begin{displaymath}
\mathscr{D}^{n}_{Y_{\infty}}=\left\{T\in\mathscr{D}^{n}_{X}\mid 
\supp T\subseteq Y\right\}; 
\end{displaymath}
we write
\begin{displaymath}
\mathscr{D}^{n}_{X/Y_{\infty }}=\mathscr{D}^{n}_{X}\left/\mathscr
{D}_{Y_{\infty }}^{n}\right..
\end{displaymath}
\end{definition}

\noindent
We note that $\mathscr{D}^{\ast}_{Y_{\infty}}$ and $\mathscr{D}^
{\ast}_{X/Y_{\infty}}$ are complexes, since $\supp\dd T\subseteq
\supp T$.

\nnpar{Invariance under birational morphisms.}
The main reason for using the complex $\mathscr{D}^{\ast}_{X/Y_
{\infty}}$ is explained by the following result obtained by J.P. Poly  
(for a proof, see \cite{Poly:shcsesa}).

\begin{theorem}
\label{thm:19}
Let $f:X'\longrightarrow X$ be a proper morphism of equidimensional
complex algebraic manifolds, and $Y\subsetneq X$ a closed subvariety 
satisfying $Y'=f^{-1}(Y)$. If  
\begin{displaymath}
f|_{X'\setminus Y'}:X'\setminus Y'\longrightarrow X\setminus Y
\end{displaymath}
is an isomorphism, then the induced morphism 
\begin{displaymath}
f_{\ast}:\mathscr{D}^{\ast}_{X'/Y'_{\infty}}\longrightarrow 
\mathscr{D}^{\ast}_{X/Y_{\infty}}
\end{displaymath}
is an isomorphism.
\hfill $\square$
\end{theorem}

\nnpar{Cohomological properties.}
Let $X$ be a proper equidimensional complex algebraic manifold, 
$Y\subsetneq X$ a closed subvariety, and $U=X\setminus Y$. 
Since the complex $\overline{\mathscr{D}}^{\ast}_{Y}$ of currents on $Y$ 
is a subcomplex of $\mathscr{D}^{\ast}_{Y_{\infty}}$, there is 
an induced morphism $\mathscr{D}^{\ast}_{X/Y}\longrightarrow
\mathscr{D}^{\ast}_{X/Y_{\infty}}$; we note that these two 
complexes do not agree. As usual, we put  
\begin{displaymath}
D^{\ast}_{X/Y_{\infty}}=\Gamma(X,\mathscr{D}^{\ast}_{X/Y_
{\infty}}),
\end{displaymath}
and observe that the complex $D_{X/Y_{\infty}}=(D^{\ast}_
{X/Y_{\infty},\mathbb{R}},\dd)$ is a Dolbeault complex.

\begin{proposition}\label{prop:36}
For any integers $p,q$, the assignment which sends an open subset $U$
of $X$
 to $D^{p,q}_{X/Y_{\infty}}$ (where $Y=X\setminus U$),
is a totally acyclic sheaf.  
\end{proposition}
\begin{proof}  
Let $U$ and $V$ be open subsets of $X$. 
We put $Y=X\setminus U$ and $Z=X\setminus V$. 
Since the complexes of currents depend only on 
$U\cup V$ and not on $X$, we may assume that $Y=Y_{0}\cup W$, $Z=
Z_{0}\cup W$ such that $Y_{0}\cap Z_{0}=\emptyset$ with closed 
subsets $Y_{0}$, $Z_{0}$, and $W$. 
 We have to show that the
sequence
\begin{displaymath}
0\longrightarrow D^{p,q}_{(Y\cap Z)_{\infty}}\longrightarrow D^{p,q}_
{Y_{\infty}}\oplus D^{p,q}_{Z_{\infty}}\longrightarrow D^{p,q}_{(Y\cup 
Z)_{\infty}}\longrightarrow 0
\end{displaymath}
is exact.
By definition, one obviously has
\begin{displaymath}
D^{p,q}_{(Y\cap Z)_{\infty}}=D^{p,q}_{Y_{\infty}}\cap D^{p,q}_{Z_{\infty}}.
\end{displaymath}
Letting $\{\sigma_{_{Y}},\sigma_{_{Z}}\}$ be a partition of unity 
subordinate to the open covering $\{X\setminus Z_{0},X\setminus 
Y_{0}\}$, and writing any current $T\in D^{p,q}_{(Y\cup Z)_{\infty}}$
as $T=\sigma_{_{Y}}T+\sigma_ 
{_{Z}}T$, we see on the other hand that
\begin{displaymath}
D^{p,q}_{(Y\cup Z)_{\infty}}=D^{p,q}_{Y_{\infty}}+D^{p,q}_{Z_{\infty}}.
\end{displaymath}
This proves the proposition.
\end{proof}

Let now $\pi:\widetilde{X}\longrightarrow X$ be an embedded
resolution of singularities of $Y$ with $D=\pi^{-1}(Y)$ a 
normal crossing divisor, i.e., $\widetilde{X}$ can be viewed 
as a smooth compactification $\overline{U}$ of $U$ with $D=
\overline{U}\setminus U$ a normal crossing divisor. By means 
of theorem \ref{thm:21}, we obtain a morphism of complexes
\begin{displaymath}
E^{\ast}_{\overline{U}}(\log D)\longrightarrow D^{\ast}_
{\widetilde{X}/D}\longrightarrow D^{\ast}_{\widetilde{X}/D_
{\infty}},
\end{displaymath}
which leads, by theorem \ref{thm:19}, to a morphism of presheaves of
Dolbeault  
complexes
\begin{displaymath}
E_{\log}(U)^{\circ}=(E^{\ast}_{\log,\mathbb{R}}(U)^{\circ},\dd)\longrightarrow 
(D^{\ast}_{X/Y_{\infty},\mathbb{R}},\dd)=D_{X/Y_{\infty}},
\end{displaymath}
which, in turn, induces a morphism of sheaves of Dolbeault complexes 
\begin{equation}
\label{eq:1}
E_{\log}(U)=(E^{\ast}_{\log,\mathbb{R}}(U),\dd)\longrightarrow 
(D^{\ast}_{X/Y_{\infty},\mathbb{R}},\dd)=D_{X/Y_{\infty}}.
\end{equation}
By another result of J.P. Poly (see \cite{Poly:shcsesa}) this 
morphism is a quasi-isomor\-phism. Nevertheless, it is not a filtered
quasi-isomorphism with
respect to the Hodge filtration. The Hodge filtration of the 
Dolbeault complex $D_{X/Y_{\infty}}$ is related to the formal 
Hodge filtration studied by Ogus in \cite{Ogus:fHf}.

\begin{remark}
We have now obtained a complex of currents on a quasi-projective
variety which does not depend on the choice of a compactification. 
The price we have to pay is that the Hodge filtration of this 
complex is not the desired one. It would be useful to have a 
complex of currents on quasi-projective varieties which is 
independent on the compactification and compatible with the 
right Hodge filtration.
\end{remark}

\nnpar{A $\mathcal{D}_{\log}$-complex constructed by means of currents.}
Let $C$ denote the site of regular schemes in $\ZAR(\Spec(\mathbb{R}))$. 
A scheme $X$ in $C$ defines a real algebraic manifold $X_{\mathbb{R}}=
(X_{\mathbb{C}},F_{\infty})$.

\begin{definition}
Let $X$ be an equidimensional scheme in $C$, $U$ an open subset of 
$X$, $Y=X\setminus U$, and $\pi:\widetilde{X}\longrightarrow X$ an 
embedded resolution of singularities of $Y$ with $D=\pi^{-1}(Y)$ a
normal crossing divisor. For any integers $n,p$, let $\mathcal{D}^
{n}_{\D,X}(p)=\mathcal{D}^{n}_{\D,X}(\cdot,p)$ denote the presheaf 
in the Zariski topology of $X$, which assigns to $U$ the group
\begin{displaymath}
\mathcal{D}^{n}_{\D,X}(U,p)=\mathcal{D}^{n}\left(D_{\widetilde{X}_
{\mathbb{C}}/D_{\mathbb{C}\,\infty}},p\right)^{\sigma},
\end{displaymath}
with $\sigma$ as in notation \ref{def:19}.
\end{definition}

\begin{proposition}
Let $X$ be an equidimensional scheme in $C$. For any integers $n,p$, 
the presheaf $\mathcal{D}^{n}_{\D,X}(p)$ is a totally acyclic sheaf. 
Moreover, it has a natural structure of a $\mathcal{D}_{\log}$-complex.  
\end{proposition}
\begin{proof}  
The first statement is consequence of proposition \ref{prop:36} and
the exactness of the functors $\mathcal{D}^{n}(\cdot,p)$. The
structure as a $\mathcal{D}_{\log}$-complex is given by the morphism
\eqref{eq:1}. 
\end{proof}

\nnpar{Multiplicative properties.} 
\begin{proposition}
Let $X$ be an equidimensional scheme in $C$ of dimension $d$. Then,
the sheaf $\mathcal{D}_{\D,X}=\mathcal{D}^{\ast}_{\D,X}(\ast)$ is a
$\mathcal{D}_{\log}$-module over $\mathcal{D}_{\log,X}=\mathcal{D}^
{\ast}_{\log,X}(\ast)$.
\end{proposition}
\begin{proof}
Let $U$ be an open subset of $X$, $Y=X\setminus U$, and $\pi:
\widetilde{X}\longrightarrow X$ an embedded resolution of 
singularities of $Y$ with $D=\pi^{-1}(Y)$ a normal crossing 
divisor. The space of currents $D^{n}_{\widetilde{X}_{\mathbb
{C}}/D_{\mathbb{C}\,\infty}}$ is the topological dual of the 
space of differential forms of degree $2d-n$ on $\widetilde
{X}_{\mathbb{C}}$, which are flat along $D_{\mathbb{C}}$ (for 
a proof in the case $n=2d$, see \cite{Malgrange:Idf}). Since   
the product of a form with logarithmic singularities along 
$D_{\mathbb{C}}$ with a form which is flat along $D_{\mathbb
{C}}$, is again flat along $D_{\mathbb{C}}$ according to
\cite{Tougeron:Ifd}, IV.4.2, there is a product
\begin{displaymath}  
E^{n}_{\log}(U_{\mathbb{C}})\otimes D^{m}_{\widetilde{X}_
{\mathbb{C}}/D_{\mathbb{C}\,\infty}}\overset{\land}{\longrightarrow } 
D^{m+n}_{\widetilde{X}_{\mathbb{C}}/D_{\mathbb{C}\,\infty}}
\end{displaymath}
given by $\varphi\land T(\omega)=T(\omega\land\varphi)$. This
pairing turns $D_{\widetilde{X}_{\mathbb{C}}/D_{\mathbb{C}\,
\infty}}=(D^{\ast}_{\widetilde{X}_{\mathbb{C}}/D_{\mathbb{C}
\,\infty},\mathbb{R}},\dd)$ into a Dolbeault module over the
Dolbeault algebra $E_{\log}(U)=(E^{\ast}_{\log,\mathbb{R}}(U),
\dd)$. An application of proposition \ref{prop:27} now shows 
that $\mathcal{D}_{\D,X}$ becomes a $\mathcal{D}_{\log}$-module 
over $\mathcal{D}_{\log,X}$.   
\end{proof}

\nnpar{Functorial properties.} 
The following result is a direct consequence of the fact that 
the pull-back of a flat differential form is again flat.

\begin{proposition}
Let $f:X\longrightarrow Y$ be a proper morphism of equidimensional 
schemes in $C$ of relative dimension $e$. Then, the push-forward 
of currents induces a covariant $f$-morphism of $\mathcal{D}_{\log}
$-complexes
\begin{displaymath}
f_{\#}:f_{\ast}\mathcal{D}_{\D,X}\longrightarrow\mathcal{D}_{\D,Y}
(-e)[-2e].
\end{displaymath}
Moreover, if $f^{\#}$ denotes the pull-back of differential forms, 
and $\bullet_{X}$ (resp. $\bullet_{Y}$) is the pairing between  
$\mathcal{D}_{\log,X}$ and $\mathcal{D}_{\D,X}$ (resp. $\mathcal{D}_
{\log,Y}$ and $\mathcal{D}_{\D,Y}$), then $(f^{\#},f_{\#},f_{\#},
\bullet_{X},\bullet_{Y})$ is a projection five-tuple.  
\end{proposition}

\nnpar{Covariant $\mathcal{D}_{\log}$-arithmetic Chow groups.}
Let $A$ be an arithmetic ring, and $X$ an arithmetic variety over
$A$ such that $X_{\mathbb{R}}$ is equidimensional. Then, the pair 
$(X,\mathcal{D}_{\D})$ is a $\mathcal{D}_{\log}$-arithmetic variety. 
Therefore, we can apply the results of section \ref{sec:AC}; in
particular, we can define the arithmetic Chow groups $\cha^{\ast}
(X,\mathcal{D}_{\D})$.

The main properties of these groups are
summarized in the subsequent theorem, which is a consequence of 
section \ref{sec:AC} and the properties of currents discussed 
above. 

\begin{remark}
Note that, when $X(\mathbb{C})$ is not compact, the cohomology groups
of the complex $\mathcal{D}_{\D,X}$ 
are not the Deligne-Beilinson cohomology groups of $X$. Nevertheless, it
is possible to prove that the complex
$\mathcal{D}_{\D,X}$ satisfies the weak purity property. This fact is
reflected by the exact sequences given in the next theorem. 
\end{remark}

\begin{theorem}
\label{thm:logD}
With the above notations, we have the following statements:
\begin{enumerate}
\item[(i)] There is an exact sequence
\begin{displaymath}
\CH^{p-1,p}(X)\overset{\rho}{\longrightarrow}\widetilde 
{\mathcal{D}}^{2p-1}_{\D}(X,p)\overset{\amap}{\longrightarrow} 
\cha^{p}(X,\mathcal{D}_{\D})\overset{\zeta}{\longrightarrow}
\CH^{p}(X)\longrightarrow 0. 
\end{displaymath}
Moreover, if $X(\mathbb{C})$ is projective, then there is an exact sequence
\begin{align}
&\CH^{p-1,p}(X)\overset{\rho}{\longrightarrow}H^{2p-1}_{\mathcal
{D}}(X_{\mathbb{R}},\mathbb{R}(p))\overset{\amap}{\longrightarrow}
\cha^{p}(X,\mathcal{D}_{\D})\overset{(\zeta,-\omega)}
{\longrightarrow} 
\notag \\ 
&\phantom{CH^{p-1,p}(X)\overset{\rho}{\longrightarrow}}\CH^{p}(X) 
\oplus{\rm Z}\mathcal{D}^{2p}_{\D}(X,p)\overset{\cl+h}
{\longrightarrow}H^{2p}_{\mathcal{D}}(X_{\mathbb{R}},\mathbb{R}(p)) 
\longrightarrow 0.
\end{align}
\item[(ii)]
For any arithmetic variety $X$ over $A$ such that $X_{\mathbb{R}}$ 
is equidimensional, there is a covariant morphism of $\mathcal{D}_
{\log}$-arithmetic varieties   
\begin{displaymath}
(X,\mathcal{D}_{\log,X})\longrightarrow(X,\mathcal{D}_{\D,X}), 
\end{displaymath}
which induces a morphism of arithmetic Chow groups
\begin{displaymath}
\cha^{p}(X,\mathcal{D}_{\log})\longrightarrow \cha^{p}
(X,\mathcal{D}_{\D}).
\end{displaymath}
When $X(\mathbb{C})$ is compact this morphism is injective. 
Moreover, if $X(\mathbb{C})$ has dimension zero, this morphism is an
isomorphism.  
\item[(iii)]
For any proper morphism $f:X\longrightarrow Y$ of arithmetic 
varieties over $A$ of relative dimension $e$, there is a 
covariant morphism of $\mathcal{D}_{\log}$-arithmetic varieties
of relative dimension $e$     
\begin{displaymath}
(X,\mathcal{D}_{\D,X})\longrightarrow(Y,\mathcal{D}_{\D,Y}), 
\end{displaymath}
which induces a morphism of arithmetic Chow groups
\begin{displaymath}
f_{\ast}:\cha^{p}(X,\mathcal{D}_{\D})\longrightarrow\cha^{p-e}
(Y,\mathcal{D}_{\D}).    
\end{displaymath}
If $g:Y\longrightarrow Z$ is another such morphism, the equality
$(g\circ f)_{\ast}=g_{\ast}\circ f_{\ast}$ holds. Moreover, if 
$f_{\mathbb{R}}:X_{\mathbb{R}}\longrightarrow Y_{\mathbb{R}}$ is
a smooth proper morphism of proper varieties, then $f_{\ast}$ is 
compatible with the direct image of contravariant arithmetic Chow 
groups.    
\item[(iv)]
The group $\cha^{\ast}(X,\mathcal{D}_{\D})$ is a module over $\cha^
{\ast}(X,\mathcal{D}_{\log})$ with an associative action. Moreover,
this action satisfies the projection formula for proper morphisms.
\end{enumerate}
\hfill $\square$
\end{theorem}

\nnpar{Canonical class of a cycle.}
An interesting property of covariant arithmetic Chow groups 
is that any algebraic cycle has its canonical class. One may
think of the contravariant arithmetic Chow groups as the 
operational Chow groups where the characteristic classes of 
vector bundles live, and the covariant arithmetic Chow groups 
as the groups where the algebraic cycles live.

If $y$ is a $p$-codimensional cycle with $Y=\supp y$, and $U=
X\setminus Y$, the pair $(\delta_{y},0)$ represents the class 
of the cycle $y$ in the group 
\begin{displaymath}
H^{2p}_{\mathcal{D}_{\D},Y}(X,p)=H^{2p}(\mathcal{D}_{\D}(X,p),
\mathcal{D}_{\D}(U,p)).
\end{displaymath}
Therefore, it is a Green object for $y$. 

\begin{definition}
\label{def:cancyc}
Let $X$ be an equidimensional arithmetic variety over $A$, 
and $y$ a $p$-codimensional cycle on $X$. Then, we denote 
by $\widehat{y}\in\cha^{p}(X,\mathcal{D}_{\D})$ the class 
of the arithmetic cycle $(y,(\delta_{y},0))$, i.e.,
\begin{displaymath}
\widehat{y}=[y,(\delta_{y},0)]\in\cha^{p}(X,\mathcal{D}_{\D}).
\end{displaymath}
\end{definition}

\subsection{Height of a cycle} 

Let $K$ be a number field, $\mathcal{O}_{K}$ its ring of integers, 
and $S=\Spec(\mathcal{O}_{K})$. Let $X$ be a $d$-dimensional 
projective arithmetic variety over $\mathcal{O}_{K}$ with structural 
morphism $\pi:X\longrightarrow S$. Since $S_{\infty}$ consists of 
a finite number of points, we have $E^{\ast}_{\log,\mathbb{R}}
(S_{\mathbb{C}})=D^{\ast}_{S_{\mathbb{C}},\mathbb{R}}$. Therefore,
we have 
\begin{displaymath} 
\cha^{\ast}(S,\mathcal{D}_{\D})=\cha^{\ast}(S,\mathcal{D}_{\log})
\cong\cha^{\ast}(S), 
\end{displaymath}
where the isomorphism is provided by the map $\Psi$ given in 
theorem \ref{thm:24}.

The key ingredient in the construction of the height of a cycle 
in \cite{BostGilletSoule:HpvpGf} is a biadditive pairing, the
so-called height pairing,
\begin{equation}
\label{eq:44}
(\cdot\mid\cdot):\cha^{p}(X)\otimes{\rm Z}^{q}(X)\longrightarrow
\cha^{p+q-d}(S)_{\mathbb{Q}}.
\end{equation}
We refer the reader to \cite{BostGilletSoule:HpvpGf} for more
details about this pairing. Here, we want to interpret the height 
pairing \eqref{eq:44} in terms of the contravariant and covariant  
arithmetic Chow groups developed in this chapter.

\begin{definition}
\label{def:first-height-pairing}
With the above notations and conventions, we define for $\alpha
\in\cha^{p}(X,\mathcal{D}_{\log})$ and $z\in{\rm Z}^{q}(X)$
\begin{equation} 
\label{eq:can-height}
(\alpha\mid z)=\pi_{\ast}(\alpha\cdot\widehat{z})\in\cha^{p+q-d}
(S,D_{\log})_{\mathbb{Q}}.
\end{equation}
If $p+q=d+1$, we call the real number 
\begin{displaymath}
\textrm{ht}_{\alpha}(z)=\dega(\alpha\mid z)
\end{displaymath}
the \emph{height of $z$ with respect to $\alpha$}.
\end{definition}

\begin{proposition} 
There is a commutative diagram 
\begin{displaymath}
\begin{CD}
\cha^{p}(X,\mathcal{D}_{\log})\otimes{\rm Z}^{q}(X)
@>(\cdot\mid\cdot)>>
\cha^{p+q-d}(S,\mathcal{D}_{\log})_{\mathbb{Q}} \\
@V\Psi\otimes\Id VV@VV\Psi V \\
\cha^{p}(X)\otimes{\rm Z}^{q}(X)
@>(\cdot\mid\cdot)>>  
\cha^{p+q-d}(S)_{\mathbb{Q}},
\end{CD}
\end{displaymath}
where the map $\Psi$ is given by theorem \ref{thm:24}, and the 
horizontal maps are provided by the pairings \eqref{eq:can-height}
and \eqref{eq:44}, respectively. 
\end{proposition}
\begin{proof}
We have to show that the height pairing given in definition
\ref{def:first-height-pairing} translates into the height pairing
of \cite{BostGilletSoule:HpvpGf} using the isomorphism $\Psi$; 
in particular, we will see that our formula \eqref{eq:can-height} 
translates into the formulas (2.3.1) and (2.3.2) therein.

If $p+q<d$, we have $\pi_{\ast}(\alpha\cdot\widehat{z})=0$, since
$X$ has dimension $d$ over $S$. This trivially proves the claimed
commutativity.

In order to treat the case $p+q=d$, we choose a representative
$(y,\mathfrak{g}_{y})$ of $\alpha$ such that $y_{K}$ intersects 
$z_{K}$ properly and such that $\mathfrak{g}_{y}$ is a Green 
object for the cycle $y$. By the $\cha^{\ast}(X,\mathcal{D}_
{\log})$-module structure of $\cha^{\ast}(X,\mathcal{D}_{\D})$, 
we now compute
\begin{align*}
&(y,\mathfrak{g}_{y})\cdot(z,(\delta_{z},0))= \\
&\left(\left(y_{K}\cdot z_{K},\mathfrak{g}_{y}\ast(\delta_{z},0)
\right),[y \cdot z]_{\fin}\right)\in\za^{p+q}(X_{K},\mathcal{D}_
{\D})_{\mathbb{Q}}\oplus\CH^{p+q}_{\fin}(X)_{\mathbb{Q}}\,,
\end{align*}
which has class $\alpha\cdot\widehat{z}$ in $\cha^{p+q}(X,\mathcal
{D}_{\D})_{\mathbb{Q}}$. Therefore, we obtain in this case
\begin{align*}
&(\alpha\mid z)=\pi_{\ast}(\alpha\cdot\widehat{z})= \\
&[\pi_{\ast}(y_{K}\cdot z_{K}),(0,0)]\in\cha^{0}(S,\mathcal{D}_
{\D})_{\mathbb{Q}}\cong\CH^{0}(S)_{\mathbb{Q}}\,.
\end{align*}
Taking into account formula (2.3.1) of \cite{BostGilletSoule:HpvpGf},
this immediately implies
\begin{displaymath}
\Psi(\alpha\mid z)=[\pi_{\ast}(y_{K}\cdot z_{K}),0]=(\Psi(\alpha)
\mid z),
\end{displaymath}
which is the claimed commutativity.

In order to treat the case $p+q=d+1$, we let $\mathfrak{g}_{y}=
(\omega_{y},\widetilde{g_{y}})$. Since $y_{K}\cap z_{K}=\emptyset$,
we find that $g_{y}\land\delta_{z}$ is a well-defined closed current 
in $\mathcal{D}_{\D}^{2d+1}(X,d+1)$. Using the explicit description 
of the $\ast$-product, we obtain 
\begin{displaymath}
\mathfrak{g}_{y}\ast(\delta_{z},0)=\left(\omega_{y}\land\delta_{z}, 
\widetilde{g_{y}\land\delta_{z}}\right)\,,
\end{displaymath}
which implies 
\begin{align*}
&(\alpha\mid z)=\pi_{\ast}(\alpha\cdot\widehat{z})= \\
&\left[\pi_{\ast}\left([y\cdot z]_{\fin}\right),\left(0,\widetilde
{\pi_{\#}(g_{y}\land\delta_{z})}\right)\right]\in\cha^{1}(S,
\mathcal{D}_{\D})_{\mathbb{Q}}=\cha^{1}(S,\mathcal{D}_{\log})_
{\mathbb{Q}}\,.
\end{align*}
By our normalization of the current associated to a cycle, we
observe
\begin{displaymath}
\pi_{\#}(g_{y}\land\delta_{z})=\frac{1}{(2\pi i)^{d-q}}\int_
{z_{\infty}}g_{y}\,;
\end{displaymath}
this shows
\begin{displaymath}
\Psi(\alpha\mid z)=\left[\pi_{\ast}\left([y\cdot z]_{\fin}\right),
\frac{2}{(2\pi i)^{d-q}}\int_{z_{\infty}}g_{y}\right]\,.
\end{displaymath}
On the other hand, since $\Psi(\alpha)=[y,2(2\pi i)^{d-p+1}
[g_{y}]_{X}]$, using formula (2.3.2) of
\cite{BostGilletSoule:HpvpGf} and respecting our definition
of the current $[g_{y}]_{X}$, we compute 
\begin{align*}
(\Psi(\alpha)\mid z)&=\left[\pi_{\ast}\left([y\cdot z]_{\fin}
\right),\frac{1}{(2\pi i)^{d}}\int_{z_{\infty}}2(2\pi i)^{d-p+1}
g_{y}\right] \\
&=\left[\pi_{\ast}\left([y\cdot z]_{\fin}\right),\frac{2}{(2\pi 
i)^{p-1}}\int_{z_{\infty}}g_{y}\right]\,,
\end{align*}
which proves the claimed commutativity in the case $p+q=d+1$.
\end{proof}

We recall that there is an alternative formula for the height 
pairing in \cite{BostGilletSoule:HpvpGf} given in terms of a 
choice of a Green object for the cycle $z$. Since this formula 
is useful for explicit computations, we describe it with the 
notations and normalizations used in this paper.

\begin{lemma}
Let $\alpha\in\cha^{p}(X,\mathcal{D}_{\log})$ and $z\in{\rm Z}^
{q}(X)$. Then, the height pairing $(\alpha\mid z)$ is given by 
the formula
\begin{displaymath}
\label{def:bgs-height}
(\alpha\mid z)=\pi_{\ast}\left(\alpha\cdot[z,\mathfrak{g}_{z}]
\right)+\amap\left([g_{z}]_{X}(\omega(\alpha))\right), 
\end{displaymath}
where $\mathfrak{g}_{z}=(\omega_{z},\widetilde{g_{z}})\in
\widehat{H}^{2q}_{\mathcal{D}_{\log},\mathcal{Z}^{q}}(X,q)$ 
is an arbitrarily chosen Green object for the class of $z$ 
and, by abuse of notation, the quantity $[g_{z}]_{X}(\omega
(\alpha))$ is the real number given by 
\begin{align*}
[g_{z}]_{X}(\omega(\alpha))&= 
\begin{cases} 
[g_{z}]_{X}(\omega(\alpha)),&\textrm{if }p+q=d+1, \\ 
0,&\textrm{if }p+q\neq d+1.
\end{cases}
\end{align*}
\end{lemma}
\begin{proof}
Using theorem \ref{thm:logD} (i), we note that the element 
\begin{displaymath}
[z,\mathfrak{g}_{z}]=[z,(\omega_{z},\widetilde{g_{z}})]
\in\cha^{q}(X,\mathcal{D}_{\log})
\end{displaymath}
maps to the element
\begin{displaymath}
\left[z,\left([\omega_{z}],\widetilde{[g_{z}]_{X}}\right)
\right]\in\cha^{q}(X,\mathcal{D}_{\D}).
\end{displaymath}
By means of the $\cha^{\ast}(X,\mathcal{D}_{\log})$-module 
structure of $\cha^{\ast}(X,\mathcal{D}_{\D})$, we now compute
\begin{align*}
\pi_{\ast}\left(\alpha\cdot[z,\mathfrak{g}_{z}]\right)&+
\amap\left([g_{z}]_{X}(\omega(\alpha))\right) \\
&=\pi_{\ast}\left(\alpha\cdot\left[z,\left([\omega_{z}],
\widetilde{[g_{z}]_{X}}\right)\right]\right)+\amap\left
(\pi_{\#}\left(\widetilde{[\omega(\alpha)\land g_{z}]_{X}}
\right)\right) \\
&=\pi_{\ast}\left(\alpha\cdot\left[z,\left([\omega_{z}],
\widetilde{[g_{z}]_{X}}\right)\right]\right)+\pi_{\ast} 
\left(\amap\left(\widetilde{[\omega(\alpha)\land g_{z}]_
{X}}\right)\right) \\
&=\pi_{\ast}\left(\alpha\cdot\left[z,\left([\omega_{z}],
\widetilde{[g_{z}]_{X}}\right)\right]\right)+\pi_{\ast} 
\left(\alpha\cdot\amap\left(\widetilde{[g_{z}]_{X}}\right) 
\right) \\
&=\pi_{\ast}\left(\alpha\cdot\left(\left[z,\left([\omega_
{z}],\widetilde{[g_{z}]_{X}}\right)\right]+\amap\left
(\widetilde{[g_{z}]_{X}}\right)\right)\right)\,.
\end{align*}
Putting $\gamma=[g_{z}]_{X}$ and observing
\begin{align*}
\amap\left(\widetilde{[g_{z}]_{X}}\right)&=\amap(\widetilde
{\gamma})=[0,(-\dd_{\mathcal{D}}\gamma,-\widetilde{\gamma})] \\
&=\left[0,\left(-[\omega_{z}]+\delta_{z},-\widetilde{[g_{z}]_
{X}}\right)\right] \\
&=-\left[z,\left([\omega_{z}],\widetilde{[g_{z}]_{X}}\right)
\right]+\widehat{z}\,,
\end{align*}
we find
\begin{displaymath}
\pi_{\ast}\left(\alpha\cdot[z,\mathfrak{g}_{z}]\right)+
\amap\left([g_{z}]_{X}(\omega(\alpha))\right)=\pi_{\ast}
(\alpha\cdot\widehat{z})=(\alpha\mid z).
\end{displaymath}
\end{proof}

\newpage
\section{Arithmetic Chow rings with pre-log-log\\forms}
\label{sec:loglog}

\subsection{Pre-log-log forms}

\nnpar{Notations.}
Let $X$ be a complex algebraic manifold of dimension $d$ and
$D$ a normal crossing divisor of $X$. Write $U=X\setminus D$,
and let $j:U\longrightarrow X$ be the inclusion.

Let $V$ be an open coordinate subset of $X$ with coordinates
$z_{1},\dots,z_{d}$; we put $r_{i}=|z_{i}|$. We say that $V$
\emph{is adapted to $D$}, if the divisor $D$ is locally given
by the equation $z_{1}\cdots z_{k}=0$. We assume that the 
coordinate neighborhood $V$ is small enough; more precisely, 
we will assume that all the coordinates satisfy $r_{i}<1/e^{e}$, 
which implies that $\log 1/r_{i}>e$ and $\log(\log 1/r_{i})>1$.

If $f$ and $g$ are two functions with non-negative real values,
we will write $f\prec g$, if there exists a constant $C>0$
such that $f(x)\le C\cdot g(x)$ for all $x$ in the domain of
definition under consideration.

\nnpar{log-log growth forms.}
\begin{definition}
\label{def:loglog}
We say that a smooth complex function $f$ on $X\setminus D$ 
has \emph{log-log growth along $D$}, if we have
\begin{equation}
\label{eq:loglog1}
|f(z_{1},\dots,z_{d})|\prec\prod_{i=1}^{k}\log(\log(1/r_{i}))^
{M}
\end{equation}
for any coordinate subset $V$ adapted to $D$ and some positive
integer $M$. The \emph{sheaf of differential forms on $X$ with
log-log growth along $D$} is the subalgebra of $j_{\ast}\mathscr
{E}^{\ast}_{U}$ generated, in each coordinate neighborhood $V$
adapted to $D$, by the functions with log-log growth along $D$
and the differentials
\begin{alignat*}{2}
&\frac{\dd z_{i}}{z_{i}\log(1/r_{i})},\,\frac{\dd\bar{z}_{i}}
{\bar{z}_{i}\log(1/r_{i})},&\qquad\text{for }i&=1,\dots,k, \\
&\dd z_{i},\,\dd\bar{z}_{i},&\qquad\text{for }i&=k+1,\dots,d.
\end{alignat*}
A differential form with log-log growth along $D$ will be called
a \emph{log-log growth form}.
\end{definition}

\nnpar{Dolbeault algebra of pre-log-log forms.}
\begin{definition}
A log-log growth form $\omega$ such that $\partial\omega$,
$\bar{\partial}\omega$ and $\partial\bar{\partial}\omega$
are also log-log growth forms is called a \emph{pre-log-log
form}. The sheaf of pre-log-log forms is the subalgebra of
$j_{\ast}\mathscr{E}^{\ast}_{U}$ generated by the pre-log-log
forms. We will denote this complex by $\mathscr{E}^{\ast}_{X}
\langle\langle D\rangle\rangle_{\w}$.
\end{definition}

The sheaf $\mathscr{E}^{\ast}_{X}\langle\langle D\rangle\rangle_
{\w}$, together with its real structure, its bigrading, and the
usual differential operators $\partial$, $\bar{\partial}$ is 
easily shown to be a sheaf of Dolbeault algebras. Moreover, 
it is the maximal subsheaf of Dolbeault algebras of the sheaf 
of differential forms with log-log growth.

\nnpar{Comparison with good forms.}
We start by recalling the notion of good forms from
\cite{Mumford:Hptncc}.

\begin{definition}
\label{def:14}
Let $X$, $D$, $U$, and  $j$ be as above. A smooth function on
$X\setminus D$ has \emph{Poincar\'e growth along $D$}, if it
is bounded in a neighborhood of each point of $D$. The sheaf
of differential forms with \emph{Poincar\'e growth along $D$}
is the subalgebra of $j_{\ast}\mathscr{E}^{\ast}_{U}$ generated
by the functions with Poincar\'e growth along $D$ and the
differentials
\begin{alignat*}{2}
&\frac{\dd z_{i}}{z_{i}\log(1/r_{i})},\,\frac{\dd\bar{z}_{i}}
{\bar{z}_{i}\log(1/r_{i})},&\qquad\text{for }i&=1,\dots,k, \\
&\dd z_{i},\,\dd\bar{z}_{i},&\qquad\text{for }i&=k+1,\dots,d
\end{alignat*}
for any coordinate subset $V$ adapted to $D$. A differential
form $\eta$ is said to be \emph{good}, if $\eta$ and $\dd\eta$
have Poincar\'e growth.
\end{definition}

From the definition, it is clear that a form with Poincar\'e
growth has log-log growth, but a log-log growth form, in 
general, does not have Poincar\'e growth. In contrast, a good 
form does not need to be a pre-log-log form. Note however 
that a differential form which is good, of pure bidegree 
and $\partial\bar{\partial}$-closed, is a pre-log-log form. 
In particular, a closed form of pure bidegree and Poincar\'e 
growth is also a pre-log-log form.

\nnpar{Inverse images.}
The pre-log-log forms are functorial with respect to inverse
images. More precisely, we have the following result.

\begin{proposition}
\label{prop:25}
Let $f:X\longrightarrow Y$ be a morphism of complex algebraic
manifolds, let $D_{X}$, $D_{Y}$ be normal crossing divisors on
$X$, $Y$ respectively, satisfying $f^{-1}(D_{Y})\subseteq D_{X}$.
If $\eta$ is a section of $\mathscr{E}^{\ast}_{Y}\langle\langle
D_{Y}\rangle\rangle_{\w}$, then $f^{\ast}\eta$ is a section of
$\mathscr{E}^{\ast}_{X}\langle\langle D_{X}\rangle\rangle_{\w}$.
\end{proposition}
\begin{proof}
Let $P$ be a point of $X$ and $Q=f(P)$. Let $V_{X}$ (resp.
$V_{Y}$) be a coordinate neighborhood of $P$ (resp. $Q$) 
adapted to $D_{X}$ (resp. $D_{Y}$) with coordinates $z_{1},
\dots,z_{d}$ (resp. $w_{1},\dots,w_{n}$) such that $f(V_{X})
\subseteq V_{Y}$. Assume that $D_{X}$ has equation $z_{1}\cdots 
z_{k}=0$ and $D_{Y}$ has equation $w_{1}\cdots w_{l}=0$. In 
these coordinate neighborhoods we can write $f=(f_{1},\dots,
f_{n})$. The condition $f^{-1}(D_{Y})\subseteq D_{X}$ means 
that the divisor of $f_{i}$ is contained in $D_{X}$ for $i=1,
\dots,l$. Therefore, we can write 
\begin{displaymath}
f_{i}=u_{i}\prod_{j=1}^{k}z_{j}^{\alpha_{i,j}},
\end{displaymath}
where the $\alpha_{i,j}$ are non-negative integers and the 
$u_{i}$ are units.

We start by showing that the pre-image of a log-log growth 
form has log-log growth. If $h$ is a log-log growth function 
on $Y\setminus D_{Y}$, then in order to show that $h\circ f$ 
has also log-log growth on $X\setminus D_{X}$, one uses the 
inequality
\begin{displaymath}
\log(a+b)<2(\log(a)+\log(b))\le 4\log(a)\log(b)
\end{displaymath}
for $a,b\ge e$. Next we consider $f^{\ast}(\dd w_{i}/(w_{i}\log
(1/|w_{i}|)))$; we have
\begin{align*}
f^{\ast}\left(\frac{\dd w_{i}}{w_{i}\log(1/|w_{i}|)}\right)&=
\frac{\dd f_{i}}{f_{i}\log(1/|f_{i}|)} \\
&=\frac{\sum_{j=1}^{k}\alpha_{i,j}\dd z_{j}/z_{j}+\dd u_{i}/u_{i}}
{\sum_{j=1}^{k}\alpha_{i,j}\log(1/r_{j})+\log(1/|u_{i}|)}.
\end{align*}
Since, for any $j_{0}$, the functions
\begin{displaymath}
\frac{\log(1/r_{j_{0}})}{\sum_{j=1}^{k}\alpha_{i,j}\log(1/r_{j})+
\log(1/u_{i})}\, , \,\,\frac{\log(1/u_{i})}{\sum_{j=1}^{k}\alpha_
{i,j}\log(1/r_{j})+\log(1/u_{i})}
\end{displaymath}
have log-log growth, we obtain that $f^{\ast}(\dd w_{i}/(w_{i}
\log(1/|w_{i}|)))$ is a log-log growth form. The analogous result
holds true for the complex conjugate form.

Finally, in order to show that the pre-image of a pre-log-log
form is a pre-log-log form one uses the compatibility between
inverse images and the operators $\partial$ and $\bar{\partial}$.
\end{proof}

\nnpar{Integrability.}
Despite the fact that pre-log-log forms and good forms are not
exactly the same, they share many properties. For instance the
following result is the analogue of propositions 1.1 and 1.2 in
\cite{Mumford:Hptncc}.

\begin{proposition}
\label{prop:23}
\begin{enumerate}[(i)]
\item[(i)]
Any log-log growth form is locally integrable.
\item[(ii)]
If $\eta$ is a pre-log-log form, and $[\eta]_{X}$ is the
associated current, then
\begin{displaymath}
[\dd\eta]_{X}=\dd[\eta]_{X}.
\end{displaymath}
The same holds true for the differential operators $\partial$,
$\bar{\partial}$ and $\partial\bar{\partial}$.
\end{enumerate}
\end{proposition}
\begin{proof}
Recall that $d$ is the dimension of $X$. For the first statement it 
is enough to show that, if $\eta$ is a $(d,d)$-form on $X$ with 
log-log growth along $D$ and $V$ is an open neighborhood adapted 
to $D$, then we have
\begin{displaymath}
\left|\int_{V}\eta\right|<\infty.
\end{displaymath}
Let us denote by $\Delta^{\ast}_{1/e^{e}}\subseteq V$ the
punctured disc of radius $1/e^{e}$. We now have the estimates
\begin{align*}
\left|\int_{V}\eta\right|&<C_{1}\prod_{i=1}^{k}\int_{\Delta^{\ast}_
{1/e^{e}}}\left|\log(\log(1/r_{i}))^{M}\frac{\dd z_{i}\wedge\dd\bar
{z}_{i}}{z_{i}\bar{z}_{i}\log(1/r_{i})^{2}}\right| \\
&<C_{2}\prod_{i=1}^{k}\int_{0}^{1/e^{e}}\log(\log(1/r_{i}))^{M}
\frac{\dd r_{i}}{r_{i}\log(1/r_{i})^{2}} \\
&<C_{3}\prod_{i=1}^{k}\int_{0}^{1/e^{e}}\frac{\dd r_{i}}{r_{i}
\log(1/r_{i})^{1+\varepsilon}}<\infty
\end{align*}
for some positive real constants $C_{1}$, $C_{2}$, $C_{3}$ and
$\varepsilon$, and a positive integer $M$. The proof of the 
second statement is analogous to the proof of proposition 1.2. 
in \cite{Mumford:Hptncc}.
\end{proof}

\subsection{Pre-log forms}
\label{sec:log-growth-cond}
We will now define a complex which contains pre-log-log forms
as well as the differential forms with logarithmic singularities
introduced in section \ref{sec:real-deligne-beil}.

\nnpar{Log growth forms.}
Let $X$, $D$, $U$ and $j$ be as in the previous section.

\begin{definition}
\label{def:log}
We say that a smooth complex function $f$ on $U$ has \emph{log
growth along $D$}, if we have
\begin{equation}
\label{eq:log1}
|f(z_{1},\dots,z_{d})|\prec\prod_{i=1}^{k}\log(1/r_{i})^{M}
\end{equation}
for any coordinate subset $V$ adapted to $D$ and some positive
integer $M$. The \emph{sheaf of differential forms on $X$ with
log growth along $D$} is the subalgebra of $j_{\ast}\mathscr
{E}^{\ast}_{V}$ generated, in each coordinate neighborhood $V$
adapted to $D$, by the functions with log growth along $D$ and
the differentials
\begin{alignat*}{2}
&\frac{\dd z_{i}}{z_{i}},\,\frac{\dd\bar{z}_{i}}{\bar{z}_{i}},
&\qquad\text{for }i&=1,\dots,k, \\
&\dd z_{i},\,\dd\bar{z}_{i},&\qquad\text{for }i&=k+1,\dots,d.
\end{alignat*}
A differential form with log growth along $D$ will be called a
\emph{log growth form}.
\end{definition}

\nnpar{Dolbeault algebra of pre-log forms.}
\begin{definition}
A log growth form $\omega$ such that $\partial\omega$, $\bar
{\partial}\omega$ and $\partial\bar{\partial}\omega$ are also
log growth forms is called a \emph{pre-log form}. The sheaf
of pre-log forms is the subalgebra of $j_{\ast}\mathscr{E}^
{\ast}_{U}$ generated by the pre-log forms. We will denote
this complex by $\mathscr{E}^{\ast}_{X}\langle D\rangle_{\w}$.
\end{definition}

The sheaf $\mathscr{E}^{\ast}_{X}\langle D\rangle_{\w} $, together
with its real structure, its bigrading and the usual differential
operators $\partial$, $\bar{\partial}$ is easily shown to be a
sheaf of Dolbeault algebras. Moreover, it is the maximal subsheaf
of Dolbeault algebras of the sheaf of differential forms with log
growth.

\subsection{Mixed forms}
For the general situation which interests us, we need a
combination of the concepts of pre-log-log and pre-log forms.

\nnpar{Mixed growth forms.}
Let $X$, $D$, $U$ and $j$ be as in the previous section. Let
$D_{1}$ and $D_{2}$ be normal crossing divisors, which may have
common components, and such that $D=D_{1}\cup D_{2}$. We denote
by $D_{2}'$ the union of the components of $D_{2}$ which are
not contained in $D_{1}$. We will say that the open coordinate
subset $V$ is adapted to $D_{1}$ and $D_{2}$, if $D_{1}$ has
equation $z_{1}\cdots z_{k}=0$, $D_{2}'$ has equation $z_{k+1}
\cdots z_{l}=0$ and that $r_{i}=|z_{i}|<1/e^{e}$ for $i=1,\dots,d$.

\begin{definition}
We define the \emph{sheaf of forms with log growth along $D_{1}$
and log-log growth along $D_{2}$} to be the subalgebra of $j_{\ast}
\mathscr{E}^{\ast}_{U}$ generated by differential forms with log
growth along $D_{1}$ and log-log growth along $D_{2}$.

A differential form with log growth along $D_{1}$ and log-log
growth along $D_{2}$ will be called a \emph{mixed growth form},
if the divisors $D_{1}$ and $D_{2}$ are clear from the context.
\end{definition}

\nnpar{Dolbeault algebra of mixed forms.}
\begin{definition}
Let $X$, $D=D_{1}\cup D_{2}$, $U$ and $j$ be as before. A
mixed growth form $\omega$ such that $\partial\omega$, $\bar
{\partial}\omega$ and $\partial\bar{\partial}\omega$ are also
mixed growth forms is called a \emph{mixed form}. The sheaf of
mixed forms is the subalgebra of $j_{\ast}\mathscr{E}^{\ast}_
{U}$ generated by the mixed forms. We will denote this complex
by $\mathscr{E}^{\ast}_{X}\langle D_{1}\langle D_{2}\rangle
\rangle_{\w}$.
\end{definition}

The sheaf $\mathscr{E}^{\ast}_{X}\langle D_{1}\langle D_{2}\rangle
\rangle_{\w}$ together with its real structure, its bigrading and
the usual differential operators $\partial$, $\bar{\partial}$ is
easily checked to be a sheaf of Dolbeault algebras. Observe that
we have by definition
\begin{displaymath}
\mathscr{E}^{\ast}_{X}\langle D_{1}\langle D_{2}\rangle\rangle_{\w}=
\mathscr{E}^{\ast}_{X}\langle D_{1}\langle D_{2}'\rangle\rangle_{\w}.
\end{displaymath}

\nnpar{Inverse images.}
We can generalize proposition \ref{prop:25} as follows. The
proof is similar to the proof of proposition \ref{prop:25},
and is therefore left to the reader.

\begin{proposition}
\label{prop:invimagewll}
Let $f:X\longrightarrow Y$ be a morphism of complex algebraic
manifolds. Let $D_{1}$, $D_{2}$ and $E_{1}$, $E_{2}$ be normal
crossing divisors on $X$, $Y$, respectively, such that $D_{1}
\cup D_{2}$ and $E_{1}\cup E_{2}$ are also normal crossing divisors.
Furthermore, assume that $f^{-1}(E_{1})\subseteq D_{1}$ and $f^
{-1}(E_{2})\subseteq D_{1}\cup D_{2}$. If $\eta$ is a section of
$\mathscr{E}^{\ast}_{Y}\langle E_{1}\langle E_{2}\rangle\rangle_
{\w}$, then $f^{\ast}\eta$ is a section of $\mathscr{E}^{\ast}_
{X}\langle D_{1}\langle D_{2}\rangle\rangle_{\w}$.
\hfill $\square$
\end{proposition}

\nnpar{Integrability.}
Let $y$ be a $p$-codimensional cycle of $X$, $Y=\supp y$, and
$U=X\setminus Y$. Let $\pi:\widetilde{X}\longrightarrow X$ be
an embedded resolution of singularities of $Y$ with normal
crossing divisors $D_{Y}=\pi^{-1}(Y)$ and $\widetilde{D}=\pi^
{-1}(D)$.

\begin{lemma}
\label{lemm:3}
Assume that $g\in\Gamma(\widetilde{X},\mathscr{E}^{n}_{\widetilde
{X}}\langle D_{Y}\langle\widetilde{D}\rangle\rangle_{\w})$. Then,
the following statements hold:
\begin{enumerate}
\item[(i)]
If $n<2p$, then $g$ is locally integrable on the whole of $X$.
We denote by $[g]_{X}$ the current associated to $g$.
\item[(ii)]
If $n<2p-1$, then $\dd [g]_{X}=[\dd g]_{X}$.
\end{enumerate}
\end{lemma}
\begin{proof}
Recall again that $d$ is the dimension of $X$. To prove the first
statement, we have to show that for any differential form $\alpha$
of degree $2d-n>2d-2p$ and any compact set $K\subseteq X$, the
integral
\begin{displaymath}
\int_{K}\alpha\wedge g
\end{displaymath}
is convergent. The restriction of $\pi^{\ast}\alpha$ to $D_{Y}$
vanishes, since the map $D_{Y}\longrightarrow X$ factors through
$Y$. Therefore, if $z_{1}\cdots z_{k}=0$ is a local equation for
$D_{Y}$, $\pi^{\ast}\alpha$ can be decomposed as a sum of terms,
each of which contains a factor $\dd z_{i}$, $\dd \bar{z}_{i}$,
$z_{i}$, or $\bar{z}_{i}$ for $i=1,\dots,k$. The result now follows
from the local expression describing logarithmic growth. The second
statement is proven analogously.
\end{proof}

\nnpar{Good Metrics.}
We recall Mumford's notion of a good metric (see \cite{Mumford:Hptncc})
in the case of line bundles.

\begin{definition}
Let $X$ be a complex algebraic manifold, $D$ a normal crossing
divisor, and put $U=X\setminus D$. Let $L$ be a line bundle on $X$,
and $L_{\mid U}$ the restriction of $L$ to $U$. A smooth hermitian
metric $h$ on $L_{\mid U}$ is said to be \emph{good along $D$}, if
we have for all $x\in D$, all neighborhoods $V$ of $x$ adapted to
$D$ (with coordinates $z_{1},\dots,z_{d}$), and all non-vanishing,
rational sections $s$ of $L$ (writing $\|s\|^{2}=h(s,s)$):
\begin{enumerate}
\item[(i)]
$\|s\|,\|s\|^{-1}\le C\left|\prod_{i=1}^{k}\log(r_{i})\right|^{N}$
for some $C>0$ and some $N\in\mathbb{N}$;
\item[(ii)]
the $1$-form $\partial\log\|s\|$ is good on $V$.
\end{enumerate}
A line bundle $L$ equipped with a good hermitian metric will be
called a \emph{good hermitian line bundle}. The pair $(L,h)$ will
be denoted by $\overline{L}$.
\end{definition}

Let $X$, $D$, $U$ be as in the preceding definition, and let
$\overline{L}$ be a good hermitian line bundle on $X$. For
a non-vanishing, rational section $s$ of $L$, let $Y$ denote 
the support of $\dv(s)$, and put $V=X\setminus Y$. Let $\pi:
\widetilde{X}\longrightarrow X$ be an embedded resolution of 
singularities of $Y$ with the property that
\begin{displaymath}
D_{Y}=\pi^{-1}(Y),\,\widetilde{D}=\pi^{-1}(D),\,\pi^{-1}
(Y\cup D)
\end{displaymath}
are normal crossing divisors.

\begin{proposition}
\label{prop:metrgo}
With the above notations we have
\begin{displaymath}
\log\|s\|\in\Gamma(\widetilde{X},\mathscr{E}^{0}_{\widetilde{X}}
\langle D_{Y}\langle\widetilde{D}\rangle\rangle_{\w})\,,\,\,
\partial\bar{\partial}\log\|s\|\in\Gamma(\widetilde{X},\mathscr
{E}^{1,1}_{\widetilde{X}}\langle\langle\widetilde{D}\rangle\rangle_
{\w}).
\end{displaymath}
\end{proposition}
\begin{proof}
First we show that $\omega_{s}=\partial\bar{\partial}\log\|s\|$
is a pre-log-log form. If $s'$ is another non-vanishing, rational
section of $L$ in an open subset $V$ of $X$, we note that $\omega_
{s}|_{V}=\omega_{s'}|_{V}$. This shows that $\partial\bar{\partial}
\log\|s\|$ gives rise to a differential form $\omega$ on $X$, which
is independent of the choice of $s$. By the very definition of a good
metric, $\omega$ has Poincar\'e growth along $D\cap V$. Since this
is true for any open covering of $X$, $\omega$ has Poincar\'e growth
along $D$. Therefore, it is a log-log growth form. Since $\omega$ is
closed, $\partial\omega$, $\bar{\partial}\omega$, and $\partial\bar
{\partial}\omega$ have also log-log growth. This proves that $\omega$
is a pre-log-log form along $D$.

We now show that $\log\|s\|$ is a pre-log form along $D_{Y}$ and
a pre-log-log form along $\widetilde{D}$. In a neighborhood of
a point of $D_{Y}\setminus\widetilde{D}$, the function $\log\|s\|$
has logarithmic singularities along $D_{Y}$, as shown in proposition
\ref{prop:26}. Therefore, $\log\|s\|$ is a pre-log form. On the
other hand, in a neighborhood of a point of $\widetilde{D}\setminus
D_{Y}$, the function $\|s\|$ has log growth along $\widetilde{D}$,
since $h$ is a good metric; hence, $\log\|s\|$ has log-log growth.
Moreover, the forms $\partial\log\|s\|$, $\bar{\partial}\log\|s\|$,
and $\partial\bar{\partial}\log\|s\|$ have Poincar\'e growth. This
shows that $\log\|s\|$ is a pre-log-log form along $\widetilde{D}
\setminus D_{Y}$. Finally, in a neighborhood of a point of $D_{Y}
\cap\widetilde{D}$, we can write $\log\|s\|$ as the sum of a form
with logarithmic singularities and a pre-log-log form. Therefore,
$\log\|s\|$ is a pre-log form there.
\end{proof}

\subsection{A $\mathcal{D}_{\log}$-complex with pre-log-log forms}
\label{subsection:cpllf}
In this section we will construct a $\mathcal{D}_{\log}$-complex
using pre-log-log forms along a fixed normal crossing divisor.

\nnpar{Varieties with a fixed normal crossing divisor.}
Let $X$ be a complex algebraic manifold of dimension $d$, and 
$D$ a normal crossing divisor. We will denote the pair $(X,D)$ 
by $\underline{X}$. If $V\subseteq X$ is an open subset, we 
will write $\underline{V}=(V,D\cap V)$.

In the sequel all operations have to be adapted with respect 
to the pair $\underline{X}$. For instance, if $Y\subseteq X$ 
is a closed algebraic subset and $V=X\setminus Y$, we mean by 
an embedded resolution of singularities of $Y$ adapted to $D$ 
an embedded resolution of singularities $\pi:\widetilde{X}
\longrightarrow X$ of $Y$ with the property that
\begin{displaymath}
D_{Y}=\pi^{-1}(Y),\,\widetilde{D}=\pi^{-1}(D),\,\pi^{-1}
(Y\cup D)
\end{displaymath}
are normal crossing divisors. Using Hironaka's theorem on the 
resolution of singularities \cite{Hironaka:rs}, one can show 
that such an embedded resolution of singularities exists; for 
a more detailed description of this fact we refer to theorem 
\ref{thm:thm:hironaka_r_s}, below.

Analogously, a normal crossing compactification of $\underline{X}$
will be a smooth compactification $\overline{X}$ such that the
adherence $\overline{D}$ of $D$ and the subsets $B_{\overline{X}}=
\overline{X}\setminus X$ and $B_{\overline{X}}\cup
\overline{D}$ are normal
crossing divisors.

\nnpar{Pre-log along infinity.}
Given a diagram of normal crossing compactifications of 
$\underline{X}$
\begin{displaymath}
\xymatrix{
\overline{X}'\ar[r]^{\varphi}&\overline{X} \\
&X\ar[ul]\ar[u],}
\end{displaymath}
with divisors $B_{\overline{X}'}$ and $B_{\overline{X}}$ at
infinity, respectively, proposition \ref{prop:invimagewll}
gives rise to an induced morphism
\begin{displaymath}
\varphi^{\ast }:\mathscr{E}^{\ast}_{\overline{X}}\langle
B_{\overline{X}}\langle\overline{D}\rangle\rangle_{\w}
\longrightarrow\mathscr{E}^{\ast}_{\overline{X}'}\langle
B_{\overline{X}'}\langle\overline{D}'\rangle\rangle_{\w}.
\end{displaymath}
In order to have a complex which is independent of the choice
of a particular compactification, as in section
\ref{sec:real-deligne-beil}, we take the limit over all possible
compactifications.
For $\underline{X}=(X,D)$ as before, we then denote 
\begin{displaymath}
E^{\ast}_{\wlg}(\underline{X})^{\circ}=\lim_{\longrightarrow}\Gamma
(\overline{X},\mathscr{E}^{\ast}_{\overline{X}}\langle
B_{\overline{X}}\langle\overline{D}\rangle\rangle_{\w}),
\end{displaymath}
where the limit is taken over all normal crossing compactifications
$\overline{X}$ of $\underline{X}$.

The assignment which sends an open subset $U$ of $X$ to
$E^{\ast}_{\wlg}(\underline{U})^{\circ}$ is a presheaf in the Zariski
topology; we denote by $E^{\ast}_{\wlg,\underline X}$ its associated
sheaf.

\begin{definition}
Let $\underline{X}=(X,D)$ be as above. Then, we define the
\emph{complex $E^{\ast}_{\wlg}(\underline{X})$ of differential
forms on $X$, pre-log along infinity and pre-log-log along
$D$} as the complex of global sections of $E^{\ast}_{\wlg,\underline
X}$, i.e.
\begin{displaymath}
E^{\ast}_{\wlg}(\underline{X})=
\Gamma (X,E^{\ast}_{\wlg,\underline X}).
\end{displaymath}
\end{definition}

\begin{remark} Since in the definition of the presheaf
  $E^{\ast}_{\wlg}(\underline{X})^{\circ}$ we are using growth
  conditions it might be possible that this presheaf
  is already a sheaf. 
\end{remark}

\nnpar{A $\mathcal{D}_{\log}$-complex.}
Let $X$ be a smooth real variety and $D$ a normal crossing divisor
defined over $\mathbb{R}$; as before, we write $\underline{X}=(X,D)$.
For any $U\subseteq X$, the complex $E^{\ast}_{\wlg}(\underline{U})$
is a Dolbeault algebra with respect to the wedge product.

\begin{definition}
For any Zariski open subset $U\subseteq X$, we put
\begin{displaymath}
\mathcal{D}^{\ast}_{\wlg,\underline{X}}(U,p)=(\mathcal{D}^{\ast}_
{\wlg,\underline{X}}(U,p),\dd_{\mathcal{D}})=(\mathcal{D}^{\ast}
(E_{\wlg}(\underline{U}_{\mathbb{C}}),p)^{\sigma},\dd_{\mathcal{D}}),
\end{displaymath}
where $\sigma$ is as in notation \ref{def:19}.
\end{definition}

\begin{theorem}
\label{thm:26wll}
The complex $\mathcal{D}_{\wlg,\underline{X}}$ is a $\mathcal{D}_
{\log}$-complex on $X$. Moreover, it is a pseudo-associative and
commutative $\mathcal{D}_{\log}$-algebra.
\end{theorem}
\begin{proof}
As in proposition \ref{prop:24}, we obtain that the presheaf 
$\mathcal{D}^{n}_{\wlg,\underline{X}}(\cdot,p)$ satisfies the 
Mayer-Vietoris principle for any pair of integers $p,n$. Therefore, 
it is a totally acyclic sheaf. Clearly, there are morphisms of 
sheaves of algebras 
\begin{displaymath}
\mathcal{D}^{\ast}_{\log,X}\longrightarrow\mathcal{D}^{\ast}_
{\wlg,\underline{X}}.
\end{displaymath}
The claim now follows.
\end{proof}

\begin{remark}
\label{rem:2}
Because of technical reasons we are not able to prove or 
disprove a filtered Poincar\'e lemma for the complex $E^{\ast}_
{\wlg,\underline{X}}$. Therefore, we do not know the exact 
cohomology of the complex $\mathcal{D}^{\ast}_{\wlg,\underline
{X}}$. Nevertheless, if $X$ is projective, we note that the 
composition
\begin{displaymath}
\mathcal{D}^{\ast}(E_{X},p)\longrightarrow\mathcal{D}^{\ast}
(E_{X}\langle\langle D\rangle\rangle_{\w},p)\longrightarrow
\mathcal{D}^{\ast}(D_{X},p)
\end{displaymath}
is a quasi-isomorphism. Therefore, the least we can say is 
that the cohomology of $\mathcal{D}^{\ast}(E_{X}\langle\langle 
D\rangle\rangle_{\w},p)$ has the usual real Deligne-Beilinson 
cohomology as a direct summand. As has been mentioned before, 
this problem can be solved by imposing growth conditions on 
\emph{all} derivatives of the differential forms under 
consideration.
\end{remark}

\subsection{Properties of Green objects with values in $\mathcal{D}_
{\wlg}$.}
We start by noting that theorem \ref{thm:26wll} together with 
section \ref{sec:go} provides us with a theory of Green objects 
with values in $\mathcal{D}_{\wlg}$.

\nnpar{Mixed forms representing the class of a cycle.}
In the case of pre-log-log forms the analogue of proposition
\ref{prop:30} is less precise, since we do not know the 
cohomology of the complex of pre-log-log forms.


\begin{proposition}
\label{prop:302}
Let $X$ be a smooth real variety, and $y$ a $p$-codimensional 
cycle on $X$ with support $Y$. Let $(\omega,g)$ be a cycle in
\begin{displaymath}
  s^{2p}(\mathcal{D}_{\wlg,\underline X}(X,p)\longrightarrow
\mathcal{D}_{\wlg,\underline X}(X\setminus Y,p)).
\end{displaymath}
 Then, we have the following 
statements:
\begin{enumerate}
\item[(i)]
If the class of the cycle $(\omega,g)$ in $H^{2p}_{\mathcal{D}_
{\wlg},Y}(X,p)$ is equal to the class of $y$, then
\begin{align}
\label{eq:eqgreenwlog}
-2\partial\bar{\partial}[g]_{X}=[\omega]-\delta_{y}.
\end{align}
\item[(ii)]
Assume that $y=\sum_{j}n_{j}Y_{j}$ with irreducible subvarieties
$Y_{j}$ and certain multiplicities $n_{j}$. If the cycle 
$(\omega,g)$ represents the class of $y$, then the equality
\begin{equation}
\label{eq:dcwlg}
-\lim_{\varepsilon\rightarrow 0}\int_{\partial B_{\varepsilon}
(Y)}\alpha\dd^{{\rm c}}g=\frac{(2\pi i)^{p-1}}{2}\sum_{j}n_{j}
\int_{Y_{j}}\alpha
\end{equation}
holds for any differential form $\alpha$; here $B_{\varepsilon}
(Y)$ is an $\varepsilon$-neighborhood of $Y$ such that the 
orientation of $\partial B_{\varepsilon}(Y)$ is induced from 
the orientation of $B_{\varepsilon}(Y)$.
\item[(iii)]
Let $g_y$ be a a differential form on $X\setminus ( Y\cup
  D)$ such that there exists an embedded resolution of singularities
  $\pi:\widetilde X\to X$
  of $Y$ in $\underline{X}$ with
  $D_Y=\pi^{-1}(Y)$ and $\widetilde D=\pi^{-1}(D)$, such that in any
  coordinate neighborhood adapted to $D_Y$, we have
  \begin{align}\label{eq:basicgreen}
    \pi^*(g_{y})= \sum_{i =1}^k \alpha_i\log (1/r_{i}) + \beta,
  \end{align}
  where $\alpha_i$ are smooth forms on $\widetilde X$ and $\beta $ is
  the pull-back of a pre-log-log form on $\underline{X}$.  If the pair 
  $(\omega _{y},g_{y})=(-2\partial\bar \partial
  g_{y},g_{y})$ is a cycle in 
  \begin{displaymath}
    s^{2p}(\mathcal{D}_{\wlg,\underline X}(X,p)\longrightarrow
    \mathcal{D}_{\wlg,\underline X}(X\setminus Y,p))
  \end{displaymath}
  and $g_y$ satisfies one of the equations \eqref{eq:eqgreenwlog} or
  \eqref{eq:dcwlg}, then the pair $(\omega _{y},\widetilde {g}_{y})$
  is a pre-log-log Green object for the 
  cycle $y$.
\end{enumerate}
\end{proposition}
\begin{proof}
(i) Let $(\omega',g')$ be a cycle representing the class of $y$ 
in $H^{2p}_{\mathcal{D},Y}(X,\mathbb{R}(p))$. By proposition 
\ref{prop:30}, we then have
\begin{displaymath}
-2\partial\bar{\partial}[g']_{X}=[\omega']-\delta_{y}.
\end{displaymath}
By our assumption, the pair $(\omega,g)$ and the image of 
$(\omega',g')$ represent the same class in $H^{2p}_{\mathcal{D}_
{\wlg},Y}(X,p)$. Therefore, there are elements $a\in\mathcal{D}^
{2p-1}_{\wlg,\underline{X}}(X,p)$ and $b\in\mathcal{D}^{2p-2}_
{\wlg,\underline{X}}(X\setminus Y,p)$ such that
\begin{displaymath}
(\dd_{\mathcal{D}}a,a-\dd_{\mathcal{D}}b)=(\omega,g)-(\omega',g').
\end{displaymath}
The result now follows as in proposition \ref{prop:30} (i)
using proposition \ref{prop:23} and lemma \ref{lemm:3}.

(ii) The second statement follows as in proposition \ref{prop:30}
(ii).

(iii) Since pre-log-log forms have 
no residues we may deduce this statement from proposition \ref{prop:30}.
\end{proof}

\nnpar{The first Chern form of a line bundle.}
Let $X$ be a quasi-projective smooth real variety, $D$ a normal
crossing divisor defined over $\mathbb{R}$, and $\overline{L}=
(L,h)$ a hermitian line bundle over $X$, where the hermitian
metric $h$ on $L\otimes\mathbb{C}$ is good along $D$ and satisfies
$F^{\ast}_{\infty}h=h$. If $s$ is a non-vanishing, rational section
of $L$, we put $\|s\|^{2}=h(s,s)$, $y=\dv(s)$, $Y=\supp y$, $U=
X\setminus Y$, and (see section \ref{sec:classes-cycles-line})
\begin{align*}
g_{s}&=-\frac{1}{2}\log(\|s\|^{2}), \\
\omega_{s}&=-2\partial\bar{\partial}g_{s}.
\end{align*}
Moreover, we put $\underline{X}=(X,D)$, and $\underline{U}=
(U,U\cap D)$. We call
\begin{displaymath}
{\rm c}_{1}(\overline{L})=\omega_{s}=-2\partial\bar{\partial}
g_{s}=\partial\bar{\partial}\log(\Vert s\Vert^{2})
\end{displaymath}
the \emph{first Chern form of $\overline{L}$}.

\begin{proposition}
\label{prop:29}
With the above assumptions the following statements hold: The form
$\omega_{s}$ belongs to $\mathcal{D}^{2}_{\wlg,\underline{X}}(X,1)$,
the form $g_{s}$ belongs to $\mathcal{D}^{1}_{\wlg,\underline{X}}
(U,1)$. The pair $(\omega_{s},g_{s})$ is a cycle of the simple 
complex
\begin{displaymath}
s^{2}(\mathcal{D}_{\wlg,\underline{X}}(X,1)\longrightarrow\mathcal
{D}_{\wlg,\underline{X}}(U,1)).
\end{displaymath}
Moreover, this pair represents the class of $\dv(s)$ in the
cohomology group $H^{2}_{\mathcal{D}_{\wlg},Y}(X,1)$.
\end{proposition}
\begin{proof}
We first show that the form $\omega_{s}$ belongs to $\mathcal{D}^
{2}_{\wlg,\underline{X}}(X,1)$. By proposition \ref{prop:metrgo}, 
the form $\omega_{s}$ is a pre-log-log form. The invariance of $h$ 
with respect to $F_{\infty}$ now shows that $\omega_{s}$ belongs 
to $\mathcal{D}^{2}_{\wlg,\underline{X}}(X,1)$. Analogously, again
using proposition \ref{prop:metrgo}, we obtain that $g_{s}\in
\mathcal{D}^{1}_{\wlg,\underline{X}}(U,1)$.

Let now $h'$ be a hermitian metric on $L$, which is invariant 
under $F_{\infty}$ and smooth on the whole of $X$. Let us write
\begin{align*}
g_{s}'&=-\frac{1}{2}\log(\|s\|'^{2}), \\
\omega_{s}'&=-2\partial\bar{\partial}g_{s}'.
\end{align*}
By proposition \ref{prop:26}, the class of $\dv(s)$ in $H^{2}_
{\mathcal{D}_{\log},Y}(X,1)$ is represented by the pair $(\omega_
{s}',g_{s}')$. On the other hand, we have
\begin{displaymath}
g_{s}-g_{s}'\in\mathcal{D}^{1}_{\wlg,\underline{X}}(X,1),
\end{displaymath}
and therefore
\begin{displaymath}
\dd_{\mathcal{D}}(g_{s}-g_{s}',0)=(\omega_{s},g_{s})-(\omega_{s}',
g_{s}').
\end{displaymath}
Thus both pairs represent the same cohomology class.
\end{proof}

\nnpar{Formulas for the $*$-product.}
The formulas of proposition \ref{prop:explicit} still hold true
in this more general context. Moreover, if we consider embedded
resolutions adapted to $\underline{X}$, it is not difficult to
show that theorem \ref{thm:partition} also remains true.

\nnpar{Inverse images.}
\begin{proposition}
\label{prop:invimageloglog}
Let $f:X\longrightarrow Y$ be a morphism of smooth real varieties,
let $D_{X}$, $D_{Y}$ be normal crossing divisors on $X$, $Y$
respectively, satisfying $f^{-1}(D_{Y})\subseteq D_{X}$. Put
$\underline{X}=(X,D_{X})$ and $\underline{Y}=(Y,D_{Y})$. Then,
there exists a contravariant $f$-morphism
\begin{displaymath}
f^{\#}:\mathcal{D}_{\wlg,\underline{Y}}\longrightarrow f_{\ast}
\mathcal{D}_{\wlg,\underline{X}}.
\end{displaymath}
\end{proposition}
\begin{proof}
By proposition \ref{prop:invimagewll}, the pull-back of differential
forms induces a morphism of the corresponding Dolbeault algebras 
of mixed forms. This morphism is compatible with the involution 
$\sigma$. Thus, this morphism gives rise to an induced morphism 
between the corresponding Deligne algebras.
\end{proof}

\nnpar{Push-forward.}
We will only state the most basic property concerning direct 
images, which is necessary to define arithmetic degrees. Note 
however that we expect log-log forms to be useful in the study 
of non smooth, proper, surjective morphisms. By proposition 
\ref{prop:23}, we have

\begin{proposition}
\label{prop:pushforward-point}
Let $\underline{X}=(X,D)$ be a proper, smooth real variety with
fixed normal crossing divisor $D$. Let $f:X\longrightarrow\Spec
(\mathbb{R})$ denote the structural morphism. Then, there exists
a covariant $f$-morphism
\begin{displaymath}
f_{\#}:f_{\ast}\mathcal{D}_{\wlg,\underline{X}}\longrightarrow
\mathcal{D}_{\log,\Spec(\mathbb{R})}.
\end{displaymath}
\hfill $\square$
\end{proposition}

\noindent
In particular, if $X$ has dimension $d$, we obtain a well-defined
morphism
\begin{displaymath}
f_{\#}:\widehat{H}^{2d+2}_{\mathcal{D}_{\wlg},\mathcal{Z}^{d+1}}
(X,d)\longrightarrow\widehat{H}^{2}_{\mathcal{D}_{\log},\mathcal
{Z}^{1}}(\Spec(\mathbb{R}),1)=\mathbb{R}.
\end{displaymath}
Note that, by dimension reasons, we have $\mathcal{Z}^{d+1}=
\emptyset$, and
\begin{displaymath}
\widehat{H}^{2d+2}_{\mathcal{D}_{\wlg},\mathcal{Z}^{d+1}}(X,d)=
H^{2d+1}(\mathcal{D}_{\wlg,\underline{X}}(X,d+1)).
\end{displaymath}
Thus, every element of $\widehat{H}^{2d+2}_{\mathcal{D}_{\wlg},
\mathcal{Z}^{d+1}}(X,d)$ is represented by a pair $\mathfrak{g}=
(0,\widetilde{g})$. The morphism $f_{\#}$ mentioned above, is
then given by
\begin{displaymath}
\mathfrak{g}=(0,\widetilde{g})\longmapsto\left(0,\frac{1}{(2\pi i)^
{d}}\int_{X}g\right).
\end{displaymath}

\begin{notation}
For $\mathfrak{g}=(0,\widetilde{g})\in\widehat{H}^{2d+2}_{\mathcal
{D}_{\wlg},\mathcal{Z}^{d+1}}(X,d)$, we will write
\begin{displaymath}
\frac{1}{(2\pi i)^{d}}\int_{X}\mathfrak{g}=\frac{1}{(2\pi i)^{d}}
\int_{X}g.
\end{displaymath}
\end{notation}

\subsection{Push-forward of a $*$-product.}

\nnpar{Embedded resolution of singularities.}  
We now come to a more detailed description of embedded resolutions
of singularities of closed algebraic subsets adapted to a fixed 
normal crossing divisor, which has been announced in section 
\ref{subsection:cpllf}. This description is a consequence of the
following precise version of Hironaka's theorem on the resolution 
of singularities, which is contained in \cite{Hironaka:rs}. These
considerations provide an important tool for explicit computations 
of $*$-products.

\begin{theorem}[Hironaka]
\label{thm:thm:hironaka_r_s}
Let $X$ be a smooth variety over a field of characteristic zero. 
Let $D$ be a normal crossing divisor of $X$, and $Y$ an irreducible 
reduced subscheme of $X$. Then, there is a finite sequence of 
smooth varieties $\widetilde{X}_{k}$ and subschemes $D_{k},Y_{k},
W_{k}$ ($k=0,\dots,N$) satisfying:
\begin{enumerate}
\item[(i)]
$\widetilde{X}_{0}=X$, $D_{0}=D$, $Y_{0}=Y$;
\item[(ii)]
for $k=0,\dots,N$, the subscheme $W_{k}$ is contained in $Y_{k}$, 
the subscheme $Y_{k}$ is normally flat along $W_{k}$, the subscheme
$D_{k}$ is a normal crossing divisor, and the pair $(D_{k},W_{k})$ 
has only normal crossings; 
\item[(iii)]
$\widetilde{X}_{k+1}$ is the blow-up of $\widetilde{X}_{k}$ along 
$W_{k}$, $D_{k+1}$ is the pre-image of $D_{k}\cup W_{k}$ by this 
blow-up, and $Y_{k+1}$ is the strict transform of $Y_{k}$;
\item[(iv)]
$W_{k}\neq Y_{k}$ for $k=0,\dots,N-1$, $W_{N-1}=Y_{N-1}$, and 
$W_{N}=\emptyset$.
\end{enumerate}
These conditions imply that $Y_{N-1}$ is smooth and that $\widetilde
{X}_{N}$ is an embedded resolution of singularities of $Y$. Moreover, 
the class of embedded resolutions of singularities of $Y$, which can 
be obtained by this method is cofinal among the class of all embedded 
resolutions of singularities of $Y$.   
\hfill $\square$
\end{theorem}

We apply theorem \ref{thm:thm:hironaka_r_s} to the case when $X$ 
is a smooth real variety of dimension $d$, $D$ a fixed normal 
crossing divisor of $X$, and $Y$ an irreducible reduced subscheme 
of codimension $p$ in $X$ intersecting $D$ properly. Theorem
\ref{thm:thm:hironaka_r_s} now provides an embedded resolution
of singularities $\pi:\widetilde{X}_{N}\longrightarrow X$ of 
$Y$ adapted to $D$. We note that $D_{N}=\pi^{-1}(Y\cup D)$, and 
that $\pi^{-1}(Y)$ and $\pi^{-1}(D)$ have either no component 
or at least one component $E$ in common. In the latter case, the 
component $E$ appears in some intermediate step of the resolution 
of singularities, so we may assume that $E$ arises as the 
exceptional divisor from $W_{k}$ by the map $\pi_{k}:\widetilde
{X}_{N}\longrightarrow\widetilde{X}_{k}$ obtained by composing 
the corresponding blow-ups. Since $Y$ is not contained in $D$, 
we have $\codim W_{k}>\codim Y$. For later purposes, we denote 
by $\rho_{k}:\widetilde{X}_{k}\longrightarrow X$ the map obtained 
by composing the corresponding remaining blow-ups, i.e., $\pi=
\rho_{k}\circ\pi_{k}$.

We now describe the above situation in terms of local coordinates.
We let $U$ denote a coordinate neighborhood of $x\in\widetilde
{X}_{N}$ with local coordinates $z_{1},\dots,z_{d}$ adapted to 
$D_{N}$. This means in particular that there are subsets $S$, 
resp. $T$ of $\{1,\dots,d\}$ such that the normal crossing 
divisors $\pi^{-1}(Y)$, resp. $\pi^{-1}(D)$ are locally given 
by the equations
\begin{displaymath}
\prod_{i\in S}z_{i}=0\,,\,\,{\rm resp.}\,\,\prod_{i\in T}z_{i}=0\,.
\end{displaymath}
The case when $\pi^{-1}(Y)$ and $\pi^{-1}(D)$ have no component
in common is then characterized by $S\cap T=\emptyset$, whereas 
in the other case we have $S\cap T\neq\emptyset$; in the latter
case, we assume that $i\in S\cap T$, i.e., the common component 
$E$ is given by the equation $z_{i}=0$. Furthermore, we denote 
by $V$ a coordinate neighborhood of $\pi_{k}(x)\in\widetilde
{X}_{k}$ with local coordinates $t_{1},\dots,t_{d}$ adapted to 
the pair $(D_{k},W_{k})$; this means that $D_{k}\cap V$ is a 
union of coordinate hyperplanes and that $W_{k}$ is contained
in the intersection of at least $p+1$ coordinate hyperplanes. 
For simplicity, we will assume that $W_{k}$ is contained in 
the subset defined by $t_{1}=\dots=t_{p+1}=0$. After shrinking 
$U$, if necessary, we may assume that $\pi_{k}(U)\subseteq V$; 
furthermore, we may assume that $\pi_{k}(z_{1},\dots,z_{d})=
(t_{1},\dots,t_{d})$. The condition $\pi_{k}(E)\subseteq W_{k}$
implies that these local coordinates satisfy
\begin{equation}
\label{eq:W_k_equation}
(t_{1},\dots,t_{d})=(z_{i}^{n_{1}}u_{1},\dots,z_{i}^{n_{p+1}}
u_{p+1},\ast,\dots,\ast),
\end{equation}
where $n_{1},\dots,n_{p+1}$ are positive integers and $u_{1},
\dots,u_{p+1}$ are holomorphic functions, whose divisor of 
zeros does not contain $E$.

\nnpar{Basic pre-log-log-Green forms.}  
In many cases, we can derive a formula for the push-forward 
of a $*$-product of top degree, which is similar to the 
push-forward of the $*$-product defined by Gillet and Soul\'e. 
A basic ingredient to derive such formula is the concept of 
basic pre-log-log Green forms, which we define below. These 
forms are the analogues of the basic Green forms introduced 
in \cite{Burgos:CDc}.

\begin{definition}
\label{def:basic-pre-log-log}
Let $X$ be a smooth real variety with fixed normal crossing 
divisor $D$ as above, $y$ a $p$-codimensional cycle of $X$, 
and $Y=\supp y$. A \emph{basic pre-log-log Green form for
$y$} is a differential form $g_{y}$ on $X$ satisfying:
\begin{enumerate}
\item[(i)] 
$(-2\partial\bar{\partial}g_{y},\widetilde{g_{y}})$ is a 
Green object with values in $\mathcal{D}_{\wlg}$ for the 
cycle $y$;
\item[(ii)]
there exists an embedded resolution of singularities $\pi:
\widetilde{X}\longrightarrow X$ of $Y$ adapted to $D$ such 
that in any coordinate neighborhood with coordinates $z_
{1},\dots,z_{d}$ adapted to $\pi^{-1}(Y)$ the equality
\begin{equation}
\label{eqbpllgf}
g_{y}=\sum_{i}\alpha_{i}\log(1/r_{i})+\beta
\end{equation}
holds, where $\alpha_{i}$ are smooth differential forms 
on $\widetilde{X}$ and $\beta$ is a pre-log-log form on
$\widetilde{X}$.
\end{enumerate}
\end{definition}

\noindent
We note that in equation \eqref{eqbpllgf}, the quantity $g_{y}$ 
should be replaced by $\pi^{\ast}(g_{y})$. We permit ourselves 
this slight abuse of notation here and in the subsequent 
considerations in order to make our notations less heavy. Note also
that in view of proposition \ref{prop:302}, if condition (ii) is
satisfied then condition (i) holds if any of the of the equations
\eqref{eq:eqgreenwlog} or \eqref{eq:dcwlg} is satisfied.

By the cofinality property in the stated above version of 
Hironaka's theorem, we may assume that the embedded resolution 
of singularities appearing in the preceding definition of a 
basic pre-log-log Green form is obtained by means of theorem
\ref{thm:thm:hironaka_r_s}, i.e., we may assume $\widetilde{X}=
\widetilde{X}_{N}$ and choose the local coordinates according
to the discussion following theorem \ref{thm:thm:hironaka_r_s}.
In particular, we can then rewrite equation \eqref{eqbpllgf} 
as follows
\begin{displaymath}
g_{y}=\sum_{i\in S}\alpha_{i}\log(1/r_{i})+\beta.
\end{displaymath}
We note that $g_{y}$ is a basic Green form in the sense of 
\cite{Burgos:CDc}, if $\beta$ is a smooth differential form. We 
observe that the function $-\log\|s\|$ is a basic pre-log-log 
Green form for $\dv(s)$, if $\overline{L}=(L,\|\cdot\|)$ is a 
good hermitian line bundle and $s$ a non-vanishing, rational 
section of $L$.

\nnpar{Push-forward of a $*$-product.}
\begin{theorem}
\label{STAR-PRODUCT}
Let $\underline{X}=(X,D)$ be a proper, smooth real variety of
dimension $d$ with fixed normal crossing divisor $D$. Let $y$ be 
an irreducible $p$-codimensional cycle of $X$, $Y=\supp y$, and
$g_{y}$ a basic pre-log-log Green form for $y$. Let $z$ be a
$q$-codimensional cycle of $X$, $Z=\supp z$, and $\mathfrak{g}_
{z}=(\omega_{z},\widetilde{g_{z}})$ a Green object with values 
in $\mathcal{D}_{\wlg}$ for the cycle $z$. Assume that $p+q=d+1$, 
that $Y$ and $Z$ intersect properly, i.e., $Y\cap Z=\emptyset$, 
and that $Y$ intersects $D$ properly. Then, the following formula 
holds
\begin{align*}
\frac{1}{(2\pi i)^{d}}\int_{X}\mathfrak{g}_{y}\ast\mathfrak{g}_
{z}=\frac{1}{(2\pi i)^{q-1}}\int_{Y}g_{z}+\frac{1}{(2\pi i)^{d}}
\int_{X}g_{y}\wedge\omega_{z}\,.
\end{align*}
\end{theorem}
\begin{proof}
Let $\sigma_{_{YZ}}$ and $\sigma_{_{ZY}}$ be as in lemma 
\ref{lem:quisinvdb}. Recalling $-2\partial\bar{\partial}=
(4\pi i){\rm dd^{c}}$ in conjunction with theorem 
\ref{thm:partition}, which is also valid for pre-log-log 
forms, we obtain
\begin{displaymath}
\mathfrak{g}_{y}\ast\mathfrak{g}_{z}=\left(\omega_{y}\wedge
\omega_{z},4\pi i\left({\rm dd^{c}}(\sigma_{_{YZ}}g_{y})\wedge
g_{z}+\sigma_{_{ZY}}g_{y}\wedge{\rm dd^{c}}g_{z}\right)^{\sim}
\right).
\end{displaymath}
We have to investigate the integral
\begin{displaymath}
\int_{X}\mathfrak{g}_{y}\ast\mathfrak{g}_{z}=4\pi i\int_{X}
\left({\rm dd^{c}}(\sigma_{_{YZ}}g_{y})\wedge g_{z}+\sigma_
{_{ZY}}g_{y}\wedge{\rm dd^{c}}g_{z}\right).
\end{displaymath}
In order to perform these calculations, we put
\begin{displaymath}
X_{\varepsilon}=X\setminus(B_{\varepsilon}(D)\cup B_{\varepsilon}
(Y)\cup B_{\varepsilon}(Z)),
\end{displaymath}
where $B_{\varepsilon}(\cdot)$ denotes an $\varepsilon$-neighborhood
of the quantity in question. We observe that on $X_{\varepsilon}$ we
can split up the integral
\begin{displaymath}
\int_{X_{\varepsilon}}\left({\rm dd^{c}}(\sigma_{_{YZ}}g_{y})\wedge
g_{z}+\sigma_{_{ZY}}g_{y}\wedge{\rm dd^{c}}g_{z}\right)
\end{displaymath}
into two summands. For the first summand we obtain
\begin{align*}
\int_{X_{\varepsilon}}{\rm dd^{c}}(\sigma_{_{YZ}}g_{y})&\wedge g_{z} \\
&=\int_{X_{\varepsilon}}\left({\rm dd^{c}}g_{y}\wedge g_{z}-{\rm dd^{c}}
(\sigma_{_{ZY}}g_{y})\wedge g_{z}\right) \\
&=\int_{X_{\varepsilon}}({\rm dd^{c}}g_{y}\wedge g_{z}-\dd\left(\dd^
{{\rm c}}(\sigma_{_{ZY}}g_{y})\wedge g_{z}\right)+\dd^{{\rm c}}(\sigma_
{_{ZY}}g_{y})\wedge\dd g_{z}).
\end{align*}
For the second summand we find
\begin{align*}
\int_{X_{\varepsilon}}\sigma_{_{ZY}}g_{y}\wedge{\rm dd^{c}}g_{z}&=
\int_{X_{\varepsilon}}\left(g_{y}\wedge{\rm dd^{c}}g_{z}-\sigma_
{_{YZ}}g_{y}\wedge{\rm dd^{c}}g_{z}\right) \\
&=\int_{X_{\varepsilon}}\left(g_{y}\wedge{\rm dd^{c}}g_{z}-\dd
\left(\sigma_{_{YZ}}g_{y}\wedge\dd^{{\rm c}}g_{z}\right)+\dd
(\sigma_{_{YZ}}g_{y})\wedge\dd^{{\rm c}}g_{z}\right) \\
&=\int_{X_{\varepsilon}}\left(g_{y}\wedge{\rm dd^{c}}g_{z}-\dd
\left(\sigma_{_{YZ}}g_{y}\wedge\dd^{{\rm c}}g_{z}\right)+\dd^
{{\rm c}}(\sigma_{_{YZ}}g_{y})\wedge\dd g_{z}\right).
\end{align*}
Summing up all the terms in question now yields
\begin{align}
\label{eq:keyformulall}
&\int_{X_{\varepsilon}}\left({\rm dd^{c}}(\sigma_{_{YZ}}g_{y})
\wedge g_{z}+\sigma_{_{ZY}}g_{y}\wedge{\rm dd^{c}}g_{z}\right)
\notag \\
&\qquad=\int_{X_{\varepsilon}}\left(g_{y}\wedge{\rm dd^{c}}g_{z}+
{\rm dd^{c}}g_{y}\wedge g_{z}+\dd^{{\rm c}}g_{y}\wedge\dd g_{z}\right)
\notag \\
&\qquad-\int_{X_{\varepsilon}}\dd\left(\dd^{{\rm c}}(\sigma_{_{ZY}}
g_{y})\wedge g_{z}+\sigma_{_{YZ}}g_{y}\wedge\dd^{{\rm c}}g_{z}\right)
\notag \\
&\qquad=\int_{X_{\varepsilon}}g_{y}\wedge{\rm dd^{c}}g_{z}+\int_
{X_{\varepsilon}}\dd\left(\dd^{{\rm c}}(\sigma_{_{YZ}}g_{y})\wedge 
g_{z}-\sigma_{_{YZ}}g_{y}\wedge\dd^{{\rm c}}g_{z}\right).
\end{align}
The claim of theorem \ref{STAR-PRODUCT} will now follow as a
consequence of the subsequent two lemmas.
\end{proof}

\begin{lemma}
\label{lemm:integrability_s_p} 
With the notations and assumptions of theorem \ref{STAR-PRODUCT}, 
the differential form $g_{y}\wedge{\rm dd^{c}}g_{z}$ is locally 
integrable on $X$. 
\end{lemma}
\begin{proof} 
In order to show that the differential form $g_{y}\wedge{\rm dd^
{c}}g_{z}$ is locally integrable on $X$, it suffices to show that 
it is locally integrable on an embedded resolution of singularities 
$\pi:\widetilde{X}_{N}\longrightarrow X$ of $Y\cup Z$ adapted to
$D$ as described in theorem \ref{thm:thm:hironaka_r_s}. We let $U$ 
be a coordinate neighborhood of $x\in\widetilde{X}_{N}$ with local 
coordinates $z_{1},\dots,z_{d}$ adapted to $D_{N}$ as described in
the discussion following theorem \ref{thm:thm:hironaka_r_s}. Then,
we have 
\begin{displaymath}
g_{y}=\sum_{i\in S}\alpha_{i}\log(1/r_{i})+\beta, 
\end{displaymath}
where $\alpha_{i}$ and $\beta$ are as in definition 
\ref{def:basic-pre-log-log}. Since the differential form $\beta
\wedge{\rm dd^{c}}g_{z}$ has log-log growth along $\pi^{-1}(D)$, 
it is locally integrable on $\widetilde{X}_{N}$.

In order to study the local integrability of the forms $\alpha_
{i}\log(1/r_{i})\wedge{\rm dd^{c}}g_{z}$ ($i\in S$), we have to 
distinguish the following two cases.

\emph{Case 1:} $i\not\in T$. This is the easy case. Namely, 
if we write 
\begin{displaymath}
\alpha_{i}\log(1/r_{i})\wedge{\rm dd^{c}}g_{z}=f(z_{1},
\dots,z_{d})\dd z_{1}\wedge\dots\wedge\dd\bar{z}_{d},
\end{displaymath}
then $f$ satisfies the estimate
\begin{displaymath}
|f(z_{1},\dots,z_{d})|\prec\log(1/r_{i})\prod_{j\in T}
\frac{\log(\log(1/r_{j}))^{M}}{r_{j}^{2}\log(1/r_{j})^
{2}}
\end{displaymath}
for some positive integer $M$. Therefore, $\alpha_{i}\log
(1/r_{i})\wedge{\rm dd^{c}}g_{z}$ is locally integrable 
on $\widetilde{X}_{N}$.

\emph{Case 2:} $i\in T$. Let $E$ be the divisor given by
the equation $z_{i}=0$; then $E$ is a common component of 
$\pi^{-1}(Y)$ and $\pi^{-1}(D)$. Letting $\widetilde{X}_{k}$ 
with the local coordinates $t_{1},\dots,t_{d}$ be as in 
the discussion following theorem \ref{thm:thm:hironaka_r_s}, 
the component $E$ is mapped to the subset $W_{k}$, which is 
contained in the subset defined by $t_{1}=\dots=t_{p+1}=0$. 
In these coordinates, we now write
\begin{displaymath}
\rho_{k}^{\ast}({\rm dd^{c}}g_{z})=\sum_{K',L'}k_{K',L'}
(t_{1},\dots,t_{d})\dd t_{K'}\wedge\dd\bar{t}_{L'}\,, 
\end{displaymath}
where $K',L'\subseteq\{1,\dots,d\}$ are strictly ordered 
subsets of cardinality $q$. We fix $K',L'$, and put $\eta=
k_{K',L'}(t_{1},\dots,t_{d})\dd t_{K'}\wedge\dd\bar{t}_
{L'}$. Since $q=d-p+1$, the subset $K'$ contains at least 
two elements  
\begin{equation}
\label{eq:elements_a_b}
a,b\in\{1,\dots,p+1\},
\end{equation}
which implies that $\eta$ contains the differentials $\dd
t_{a}$ and $\dd t_{b}$.

We are left to show that $\alpha_{i}\log(1/r_{i})\wedge
\pi_{k}^{\ast}(\eta)$ is locally integrable. To this end,
we write
\begin{align*}
\alpha_{i}\log(1/r_{i})&=\sum_{I,J}g_{I,J}(z_{1},\dots,
z_{d})\dd z_{I}\wedge\dd\bar{z}_{J}\,, \\
\pi_{k}^{\ast}(\eta)&=\sum_{K,L}h_{K,L}(z_{1},\dots,z_{d})
\dd z_{K}\wedge\dd\bar{z}_{L}\,,
\end{align*}
where $I,J\subseteq\{1,\dots,d\}$ are strictly ordered subsets
of cardinality $p-1$, and $K,L\subseteq\{1,\dots,d\}$ are 
strictly ordered subsets of cardinality $q$. In this way, we 
obtain the $(d,d)$-form
\begin{displaymath}
\alpha_{i}\log(1/r_{i})\wedge\pi_{k}^{\ast}(\eta)=\sum_
{\substack{I,J,K,L\\I\cap K=\emptyset\\J\cap L=\emptyset}}
f_{I,J,K,L}(z_{1},\dots,z_{d})\dd z_{1}\wedge\dots\wedge
\dd\bar{z}_{d} 
\end{displaymath}
with $f_{I,J,K,L}=g_{I,J}\cdot h_{K,L}$. We now show that each
function $f_{I,J,K,L}$ is locally integrable on $\widetilde{X}_
{N}$. If $i\in I\cap J$, the local integrability of $f_{I,J,K,L}$
follows as in case 1. If $i\in I\cap L$ or $i\in J\cap K$, the 
local integrability of $f_{I,J,K,L}$ follows from the estimate
\begin{displaymath}
|f_{I,J,K,L}(z_{1},\dots,z_{d})|\prec\frac{\log(\log(1/r_{i}))^
{M}}{r_{i}}\prod_{\substack{j\in T\\j\neq i}}\frac{\log(\log(1/r_{j}))^
{M}}{r_{j}^{2}\log(1/r_{j})^{2}}\,. 
\end{displaymath}
We are left to investigate the last sub-case $i\in K\cap L$. Letting 
$a,b$ be as in \eqref{eq:elements_a_b}, and recalling the relation 
$t_{a}=z_{i}^{n_{a}}u_{a}$ from \eqref{eq:W_k_equation}, we first
observe 
\begin{align}
\label{ukuehn}
&\pi_{k}^{\ast}\left(\frac{\log(\log(1/|t_{a}|))^{M}}{t_{a}\log
(1/|t_{a}|)}\,\dd t_{a}\right)=\notag \\
&\phantom{xxxxxxx}\frac{(\log(n_{a}\log(1/r_{i})+\log(1/|u_{a}|)))
^{M}}{u_{a}(n_{a}\log(1/r_{i})+\log(1/|u_{a}|))}\,\sum_{j=1}^{d}
\frac{\partial u_{a}}{\partial z_{j}}\dd z_{j}+\notag \\
&\phantom{xxxxxxx}\frac{n_{a}(\log(n_{a}\log(1/r_{i})+\log(1/|u_
{a}|)))^{M}}{z_{i}(n_{a}\log(1/r_{i})+\log(1/|u_{a}|))}\,\dd z_{i},
\end{align}
and a similar formula for the index $b$. If $j\neq i$ in the above 
sum, we obtain the estimate
\begin{displaymath}
\left|\frac{(\log(n_{a}\log(1/r_{i})+\log(1/|u_{a}|)))^{M}}{u_{a}
(n_{a}\log(1/r_{i})+\log(1/|u_{a}|))}\cdot\frac{\partial u_{a}}
{\partial z_{j}}\right|\prec\frac{\displaystyle{\prod_{k\in T}}
\log(\log(1/r_{k}))^{M'}}{r_{j}(\log(1/r_{i})+\log(1/r_{j}))}
\end{displaymath}
with some positive integer $M'$. Since the divisor of zeros of 
$u_{a}$ does not contain the component $E$, we obtain for $j=i$
the estimate
\begin{displaymath}
\left|\frac{(\log(n_{a}\log(1/r_{i})+\log(1/|u_{a}|)))^{M}}{u_{a}
(n_{a}\log(1/r_{i})+\log(1/|u_{a}|))}\cdot\frac{\partial u_{a}}
{\partial z_{i}}\right|\prec\frac{\displaystyle{\prod_{k\in T}}
\log(\log(1/r_{k}))^{M'}}{\log(1/r_{i})+\log(1/r_{j})}\,.
\end{displaymath}
Moreover, for the last summand, we note the estimate
\begin{displaymath}
\left|\frac{n_{a}(\log(n_{a}\log(1/r_{i})+\log(1/|u_{a}|)))^{M}}
{z_{i}(n_{a}\log(1/r_{i})+\log(1/|u_{a}|))}\right|\prec\frac
{\displaystyle{\prod_{k\in T}}\log(\log(1/r_{k}))^{M'}}{r_{i}
\log(1/r_{i})}\,.
\end{displaymath}
Now, we note that the differential form $f_{I,J,K,L}(z_{1},\dots,
z_{d})\dd z_{1}\wedge\dots\wedge\dd\bar{z}_{d}$ is built up in 
particular by the differential forms $\pi_{k}^{\ast}(\dots\dd t_
{a})$ and $\pi_{k}^{\ast}(\dots\dd t_{b})$ (cf. \eqref{ukuehn}).
We then observe that the summand of the differential $\pi_{k}^
{\ast}(\dots\dd t_{a})$ containing $\dd z_{i}$ is multiplied with 
the corresponding summands of $\pi_{k}^{\ast}(\dots\dd t_{b})$ 
containing the differentials $\dd z_{j}$ for $j\neq i$. Hence, 
the function $f_{I,J,K,L}$ can be decomposed as a sum of functions 
$F_{j}$ ($j\neq i$), each of which satisfying the estimate   
\begin{align*}
|F_{j}(z_{1},\dots,z_{d})|\prec&\frac{\log(1/r_{i})}{r_{i}^
{2}\log(1/r_{i})^{2}}\cdot\frac{(\log(\log(1/r_{i}))\log
(\log(1/r_{j})))^{M'}}{r_{j}(\log(1/r_{i})+\log(1/r_{j}))}
\cdot \\
&\cdot\frac{1}{r_{j}\log(1/r_{j})}\cdot\prod_{\substack{k
\in T\\k\neq i,j}}\frac{\log(\log(1/r_{k}))^{M'}}{r_{k}^{2}
\log(1/r_{k})^{2}}\,. 
\end{align*}
Finally, using the inequality between the arithmetic and the 
geometric mean, we obtain
\begin{displaymath}
|F_{j}(z_{1},\dots,z_{d})|\prec\frac{1}{r_{i}^{2}r_{j}^{2}
(\log(1/r_{i})\log(1/r_{j}))^{1+\varepsilon }}\cdot\prod_
{\substack{k\in T\\k\neq i,j}}\frac{\log(\log(1/r_{k}))^{M'}}
{r_{k}^{2}\log(1/r_{k})^{2}}
\end{displaymath}
for some $\varepsilon>0$. Adding up, this proves the local
integrability of $f_{I,J,K,L}$, and thus concludes the proof 
of the local integrability of the differential form $\alpha_
{i}\log(1/r_{i})\wedge{\rm dd^{c}}g_{z}$ on $\widetilde{X}_
{N}$.
\end{proof}

\begin{lemma}
\label{lemm:residue_s_p}
With the notations and assumptions of theorem \ref{STAR-PRODUCT}, 
the equality
\begin{displaymath}
\lim_{\varepsilon\rightarrow 0}\int_{X_{\varepsilon}}\dd\left
(\dd^{{\rm c}}(\sigma_{_{YZ}}g_{y})\wedge g_{z}-\sigma_{_{YZ}}
g_{y}\wedge\dd^{{\rm c}}g_{z}\right)=\frac{(2\pi i)^{p-1}}{2}
\int_{Y}g_{z}
\end{displaymath}
holds true.
\end{lemma}
\begin{proof} 
For any $\varepsilon,\rho>0$, we put
\begin{align*}
X_{\varepsilon,\rho}&=X\setminus\left(B_{\varepsilon}(D)\cup 
B_{\rho}(Y)\cup B_{\varepsilon}(Z)\right), \\
Y_{\varepsilon}&=Y\setminus\left(B_{\varepsilon}(D)\cup 
B_{\varepsilon}(Z)\right).
\end{align*}
Since $\sigma_{_{YZ}}$ vanishes in a neighborhood of $Z$, we
obtain by means of Stokes theorem 
\begin{align*}
\int_{X_{\varepsilon,\rho}}\dd(\dd^{{\rm c}}(\sigma_{_{YZ}}g_{y})
\wedge g_{z}&-\sigma_{_{YZ}}g_{y}\wedge\dd^{{\rm c}}g_{z}) \\
&=\int_{\partial(B_{\varepsilon}(D)\cup B_{\rho}(Y))}\left(\dd^
{{\rm c}}(\sigma_{_{YZ}}g_{y})\wedge g_{z}-\sigma_{_{YZ}}g_{y}
\wedge\dd^{{\rm c}}g_{z}\right).
\end{align*}
Since we derive from equation \eqref{eq:keyformulall} in 
conjunction with lemma \ref{lemm:integrability_s_p} that 
the limit
\begin{displaymath}
\lim_{\varepsilon,\rho\rightarrow 0}\int_{\partial(B_{\varepsilon}
(D)\cup B_{\rho}(Y))}\left(\dd^{{\rm c}}(\sigma_{_{YZ}}g_{y})\wedge 
g_{z}-\sigma_{_{YZ}}g_{y}\wedge\dd^{{\rm c}}g_{z}\right)
\end{displaymath}
exists, we can compute it as an iterated limit by taking first the
limit as $\rho$ tends to zero, and then the limit as $\varepsilon$ 
goes to zero. Since $g_{y}$ is a basic Green form in the sense of
\cite{Burgos:CDc} for the cycle $y$ on $X\setminus\left(B_{\varepsilon}
(D)\cup B_{\varepsilon}(Z)\right)$, proposition \ref{prop:30} (ii) 
shows
\begin{align*}
&\lim_{\rho\rightarrow 0}\int_{\partial(B_{\varepsilon}(D)\cup
B_{\rho}(Y))}\left(\dd^{{\rm c}}(\sigma_{_{YZ}}g_{y})\wedge g_
{z}-\sigma_{_{YZ}}g_{y}\wedge\dd^{{\rm c}}g_{z}\right)= \\
&\int_{\partial B_{\varepsilon}(D)}\left(\dd^{{\rm c}}(\sigma_
{_{YZ}}g_{y})\wedge g_{z}-\sigma_{_{YZ}}g_{y}\wedge\dd^{{\rm c}}
g_{z}\right)+\frac{(2\pi i)^{p-1}}{2}\int_{Y_{\varepsilon }}g_{z}.
\end{align*}
By the functoriality and local integrability of the pre-log-log 
forms, we have 
\begin{displaymath}
\lim_{\varepsilon\rightarrow 0}\int_{Y_{\varepsilon}}g_{z}=
\int_{Y}g_{z}.
\end{displaymath}
Therefore, we are left to show that 
\begin{displaymath}
\lim_{\varepsilon\rightarrow 0}\int_{\partial B_{\varepsilon}(D)}
\left(\dd^{{\rm c}}(\sigma_{_{YZ}}g_{y})\wedge g_{z}-\sigma_{_{YZ}}
g_{y}\wedge\dd^{{\rm c}}g_{z}\right)=0.
\end{displaymath}
To perform this calculation, we will again work on an embedded 
resolution of singularities $\pi:\widetilde{X}_{N}\longrightarrow X$ 
of $Y\cup Z$ adapted to $D$ together with a coordinate neighborhood
$U$ of $x\in\widetilde{X}_{N}$ with local coordinates $z_{1},\dots,
z_{d}$ adapted to $D_{N}$. Following the discussion after theorem
\ref{thm:thm:hironaka_r_s}, we may assume without loss of generality 
that the situation is local, i.e.,
\begin{displaymath}
\pi^{-1}(D)=\bigcup_{i\in T}E_{i}\,,
\end{displaymath}
where $E_{i}$ is given by the equation $z_{i}=0$ ($i\in T$); the 
choice of local coordinates determines differentiable maps
\begin{displaymath}
p_{i}:\partial B_{\varepsilon}(E_{i})\longrightarrow E_{i}\quad
(i\in T).
\end{displaymath}
From this we obtain (recalling our abuse of notation)
\begin{align*}
&\lim_{\varepsilon\rightarrow 0}\int_{\partial B_{\varepsilon}(D)}
\left(\dd^{{\rm c}}(\sigma_{_{YZ}}g_{y})\wedge g_{z}-\sigma_{_{YZ}}
g_{y}\wedge\dd^{{\rm c}}g_{z}\right)= \\
&\sum_{i\in T}\lim_{\varepsilon\rightarrow 0}\int_{\partial B_
{\varepsilon}(E_{i})}\left(\dd^{{\rm c}}(\sigma_{_{YZ}} g_{y})
\wedge g_{z}-\sigma_{_{YZ}}g_{y}\wedge\dd^{{\rm c}}g_{z}\right)= \\
&\sum_{i\in T}\lim_{\varepsilon\rightarrow 0}\int_{E_{i}}\int_
{p_{i}}\left(\dd^{{\rm c}}(\sigma_{_{YZ}}g_{y})\wedge g_{z}-\sigma_
{_{YZ}}g_{y}\wedge\dd^{{\rm c}}g_{z}\right).
\end{align*}
We will now fix $i_{0}\in T$, and show that  
\begin{equation}
\label{eq:bound_r_s_p}
\int_{p_{i_{0}}}\left|\dd^{{\rm c}}(\sigma_{_{YZ}}g_{y})\wedge 
g_{z}-\sigma_{_{YZ}}g_{y}\wedge\dd^{{\rm c}}g_{z}\right|\prec 
h_{1}(\varepsilon)\cdot h_{2},
\end{equation}
where $\lim_{\varepsilon\rightarrow 0}h_{1}(\varepsilon)=0$ and 
$h_{2}$ is locally integrable on $E_{i_{0}}$; we emphasize that 
equality \eqref{eq:bound_r_s_p} is an inequality between forms 
meaning that the inequality holds true for the absolute values 
of all the components involved. To do this, we observe that the 
differential form 
\begin{displaymath}
\dd^{{\rm c}}(\sigma_{_{YZ}}g_{y})\wedge g_{z}-\sigma_{_{YZ}}
g_{y}\wedge\dd^{{\rm c}}g_{z}
\end{displaymath}
is the sum of two forms, one of which is of type $(d,d-1)$, and the 
other of type $(d-1,d)$. We will prove the lemma for the 
term of type $(d,d-1)$ noting that the argument for the term of
type $(d-1,d)$ is analogous. 

Since $g_{y}$ is a basic pre-log-log Green form, we have (as in
the proof of lemma \ref{lemm:integrability_s_p}) 
\begin{displaymath}
g_{y}=\sum_{i\in S}\alpha_{i}\log(1/r_{i})+\beta.
\end{displaymath}
Therefore, we have to establish the bound \eqref{eq:bound_r_s_p}
for each summand of $g_{y}$ in the above decomposition. For the
summand $\beta$ the bound \eqref{eq:bound_r_s_p} follows easily,
since $\beta$ is a pre-log-log form. In order to treat the 
remaining summands, we write
\begin{align*}
\sigma_{_{YZ}}g_{y}\wedge\partial g_{z}&=f(z_{1},\dots,z_{d})
\dd z_{i_{0}}\wedge\prod_{\substack{j=1\\j\neq i_{0}}}^{d}\dd z_{j}
\wedge\dd\bar{z}_{j}+\dots\,, \\
\partial(\sigma_{_{YZ}}g_{y})\wedge g_{z}&=g(z_{1},\dots,z_{d})
\dd z_{i_{0}}\wedge\prod_{\substack{j=1\\j\neq i_{0}}}^{d}\dd z_{j}
\wedge\dd\bar{z}_{j}+\dots\,.
\end{align*}
As in the previous lemma, we have to distinguish two cases.

\emph{Case 1:} $i\notin T$. Then, we easily find the estimates
\begin{align*}
&|f(z_{1},\dots,z_{d})|\prec\log(1/r_{i})\cdot\frac{\log(\log
(1/r_{i_{0}}))^{M}}{r_{i_{0}}\log(1/r_{i_{0}})}\cdot\prod_{\substack{j
\in T\\j\neq i_{0}}}\frac{\log(\log(1/r_{j}))^{M}}{r_{j}^
{2}\log(1/r_{j})^{2}}\,, \\
&|g(z_{1},\dots,z_{d})|\prec\frac{1}{r_{i}}\cdot\frac{\log
(\log(1/r_{i_{0}}))^{M}}{r_{i_{0}}\log(1/r_{i_{0}})}\cdot
\prod_{\substack{j\in T\\j\neq i_{0}}}\frac{\log(\log(1/r_{j}))^
{M}}{r_{j}^{2}\log(1/r_{j})^{2}}\,,
\end{align*}
from which the bound \eqref{eq:bound_r_s_p} follows.

\emph{Case 2:} $i\in T$. In this case, the functions $f$ and 
$g$ will satisfy the same type of estimates; thus, we will 
discuss here only the estimate for the function $f$. If $i=
i_{0}$, by arguing as in the second part of the discussion
of case 2 in the proof of lemma \ref{lemm:integrability_s_p},
we find that $f$ can be decomposed as a sum of functions 
$F_{j}$ ($j\neq i$), each of which satisfying the estimate
\begin{align*} 
|F_{j}(z_{1},\dots,z_{d})|\prec&\frac{\log(1/r_{i_{0}})}
{r_{i_{0}}\log(1/r_{i_{0}})}\cdot\frac{(\log(\log(1/r_{i_
{0}}))\log(\log(1/r_{j})))^{M'}}{r_{j}(\log(1/r_{i_{0}})+
\log(1/r_{j}))}\cdot \\
&\cdot\frac{1}{r_{j}\log(1/r_{j})}\cdot\prod_{\substack{k\in T\\
k\neq i_{0},j}}\frac{\log(\log(1/r_{k}))^{M'}}{r_{k}^{2}
\log(1/r_{k})^{2}} \\
\prec&\frac{1}{r_{i_{0}}r_{j}^{2}\log(1/r_{i_{0}})^{\varepsilon}
\log(1/r_{j})^{1+\varepsilon}}\cdot\prod_{\substack{k\in T\\k\neq 
i_{0},j}}\frac{\log(\log(1/r_{k}))^{M'}}{r_{k}^{2}\log(1/r_
{k})^{2}}\,;
\end{align*}
here $M'$ and $\varepsilon$ are as in the proof of lemma 
\ref{lemm:integrability_s_p}. If $i\neq i_{0}$, the function 
$f$ can be decomposed as a sum of functions $F_{j}$ ($j\neq 
i$), which are seen to satisfy for $j=i_{0}$ the estimate 
\begin{displaymath}
|F_{i_{0}}(z_{1},\dots,z_{d})|\prec\frac{1}{r_{i_{0}}r_{i}^{2}
\log(1/r_{i_{0}})^{\varepsilon}\log(1/r_{i})^{1+\varepsilon}}
\cdot\prod_{\substack{k\in T\\k\neq i_{0},i}}\frac{\log(\log(1/r_
{k}))^{M'}}{r_{k}^{2}\log(1/r_{k})^{2}}\,,
\end{displaymath}
and for $j\neq i_{0}$ the estimate
\begin{multline*}
|F_{j}(z_{1},\dots,z_{d})|\prec \\ 
\frac{1}{r_{i_{0}}r_{i}^{2}r_{j}^{2}\log(1/r_{i_{0}})\log
(1/r_{i})^{1+\varepsilon}\log(1/r_{j})^{1+\varepsilon}}\cdot
\prod_{\substack{k\in T\\k\neq i_{0},i,j}}\frac{\log(\log(1/r_{k}))^
{M'}}{r_{k}^{2}\log(1/r_{k})^{2}}\,.
\end{multline*}
All in all, this leads to the desired estimate \eqref{eq:bound_r_s_p}
and concludes the proof of the lemma. 
\end{proof}

\begin{remark}
Clearly,  the theorem also holds when $D=\emptyset$, and in 
particular for Green objects with values in $\mathcal{D}_{\log}$.
In this case we recover the classical formula for the $*$-product 
of Green currents used in \cite{GilletSoule:ait} and \cite{Soule:lag} 
(see also \cite{Burgos:Gftp}, \cite{Burgos:CDB}). If $X$ is a 
Riemann surface, then formula \eqref{eq:keyformulall} implies the 
formulas for the generalized arithmetic intersection numbers given 
in \cite{Kuehn:gainc}.
\end{remark}

\subsection{Arithmetic Chow rings with pre-log-log forms}
\label{sec:aasDa}

\nnpar{Definition of arithmetic Chow groups with pre-log-log
forms.} Let $(A,\Sigma,F_{\infty})$ be an arithmetic ring and
$X$ an arithmetic variety over $A$. Let $D$ be a fixed normal
crossing divisor of $X_{\Sigma}$, which is stable under $F_
{\infty}$. Following the notation of the previous section we
will denote the pair $(X_{\mathbb{R}},D)$ by $\underline{X}$ .
The natural inclusion $\mathcal{D}_{\log}\longrightarrow\mathcal
{D}_{\wlg}$ induces a $\mathcal{D}_{\log}$-complex structure
in $\mathcal{D}_{\wlg}$. Then, $(X,\mathcal{D}_{\wlg})$ is a
$\mathcal{D}_{\log}$-arithmetic variety. Here, we use 
the convention that whenever $\underline{X}$ is clear from the
context, we write $\mathcal{D}_{\wlg}$ instead of $\mathcal{D}_
{\wlg,\underline{X}}$ Therefore, applying the theory developed
in chapter \ref{sec:AC}, we can define the arithmetic Chow groups
$\cha^{\ast}(X,\mathcal{D}_{\wlg})$.

The rest of this section will be devoted to studying the basic
properties of the arithmetic Chow groups $\cha^{\ast}(X,\mathcal
{D}_{\wlg})$. In particular, we will see that they satisfy
properties similar to ones of the classical arithmetic Chow
groups stated in \cite{GilletSoule:ait} and most of the properties
stated in \cite{MaillotRoessler:cdl}.

\nnpar{Exact sequences.}
We will use the notations of section \ref{sec:arithm-chow-groups}
applied to the sheaf $\mathcal{D}_{\wlg}$. We start by writing the
exact sequences given in theorem \ref{thm:15}.

\begin{theorem}
\label{thm:16}
Let $X$ be an arithmetic variety over $A$, and $\underline{X}=
(X_{\mathbb{R}},D)$ as above. Then, we have the following exact
sequences:
\begin{align}
&\CH^{p-1,p}(X)\stackrel{\rho}{\longrightarrow}\widetilde
{\mathcal{D}}^{2p-1}_{\wlg}(X,p)^{\pure}\stackrel{\amap}
{\longrightarrow}\cha^{p}(X,\mathcal{D}_{\wlg})\stackrel
{\zeta}{\longrightarrow}\CH^{p}(X)\longrightarrow 0,\notag
\\[3mm]
&\CH^{p-1,p}(X)\stackrel{\rho}{\longrightarrow}H^{2p-1}_
{\mathcal{D}_{\wlg}}(X,p)^{\pure}\stackrel{\amap}{\longrightarrow}
\cha^{p}(X,\mathcal{D}_{\wlg})\stackrel{(\zeta,-\omega)}
{\longrightarrow}\notag \\
&\phantom{CH^{p-1,p}(X)\stackrel{\rho}{\longrightarrow}}
\CH^{p}(X)\oplus{\rm Z}\mathcal{D}^{2p}_{\wlg}(X,p)\stackrel
{\cl+h}{\longrightarrow}H^{2p}_{\mathcal{D}_{\wlg}}(X,p)
\longrightarrow 0,\notag \\[3mm]
&\CH^{p-1,p}(X)\stackrel{\rho}{\longrightarrow}H^{2p-1}_
{\mathcal{D}_{\wlg}}(X,p)^{\pure}\stackrel{\amap}{\longrightarrow}
\cha^{p}(X,\mathcal{D}_{\wlg})_0\stackrel{\zeta}{\longrightarrow}
\CH^{p}(X)_{0}\longrightarrow 0.\notag
\end{align}
\end{theorem}
\begin{proof}
This is a direct consequence of theorem \ref{thm:15}.
\end{proof}

\noindent
Since $\mathcal{D}_{\wlg}$ is a $\mathcal{D}_{\log}$-complex,
there is a morphism
\begin{displaymath}
H^{\ast}_{\mathcal{D}}(X_{\mathbb{R}},p)\longrightarrow
H^{\ast}_{\mathcal{D}_{\wlg}}(X,p).
\end{displaymath}
Thus, we obtain from theorem \ref{thm:16}

\begin{corollary}
\label{cor:5}
There is a complex of abelian groups
\begin{displaymath}
H^{2p-1}_{\mathcal{D}}(X_{\mathbb{R}},p)\stackrel{\amap}
{\longrightarrow}\cha^{p}(X,\mathcal{D}_{\wlg})\stackrel
{(\zeta,-\omega)}{\longrightarrow}\CH^{p}(X)\oplus{\rm Z}
\mathcal{D}^{2p}_{\wlg}(X,p).
\end{displaymath}
\end{corollary}

\begin{remark}
\label{rem:3}
Since we have not computed the cohomology of the complex $\mathcal
{D}_{\wlg}$, these exact sequences do not give all the possible
information. We do not even know whether the corresponding cohomology
groups satisfy the weak purity condition.

As we have already pointed out, by imposing bounds on all derivatives
of the functions, one can define a smaller complex, the complex of
log-log forms, for which we can prove the Poincar\'e Lemma. Using
this complex one obtains arithmetic Chow groups that satisfy the
exact sequences of theorem \ref{thm:16}, but with real Deligne-Beilinson
cohomology (see remark \ref{rem:2}). This will be developed in a
forthcoming paper (\cite{BurgosKramerKuehn:accavb}).
\end{remark}

\nnpar{Green forms for a cycle.}
We now translate the result of proposition \ref{prop:302} into
the language of Green objects.

\begin{proposition}
\label{prop:282}
Let $X$ be an arithmetic variety over $A$, and $y$ a $p$-codimen\-sional
cycle on $X$. If a pair $\mathfrak{g}=(\omega,\widetilde{g})\in\widehat
{H}^{2p}_{\mathcal{D}_{\wlg},\mathcal{Z}^{p}}(X,p)$ is a Green object
for the class of $y$, then
\begin{equation}
\label{prop:312}
-2\partial\bar{\partial}[g]_{X}=[\omega]-\delta_{y}.
\end{equation}
\hfill $\square$
\end{proposition}

\noindent
We recall that in the above proposition $[g]_{X}$ stands for
$[g]_{X_{\infty}}$, and $\delta_{y}$ for $\delta_{y_{\mathbb
{R}}}$.

\nnpar{Multiplicative properties.}
By theorem  \ref{thm:26wll}, the complex $\mathcal{D}_{\wlg}$
satisfies all the properties required to apply theorem
\ref{thm:10}.

\begin{theorem}
The abelian group
\begin{displaymath}
\cha^{\ast}(X,\mathcal{D}_{\wlg})_{\mathbb{Q}}=\bigoplus_{p\ge 0}
\cha^{p}(X,\mathcal{D}_{\wlg})\otimes\mathbb{Q}
\end{displaymath}
is a commutative and associative $\mathbb{Q}$-algebra with unit.
\hfill $\square$
\end{theorem}

\nnpar{Inverse images.}
Let $f: X\longrightarrow Y$ be a morphism  of arithmetic varieties
over $A$. By proposition \ref{prop:invimagewll}, the complex
$\mathcal{D}_{\wlg}$ satisfies all the properties required to
apply theorem \ref{thm:6}.

\begin{theorem}
\label{thm:wll-inverse-image}
Let $f:X\longrightarrow Y$ be a morphism of arithmetic varieties
over $A$, let $D_{X}$, $D_{Y}$ be a normal crossing divisor on
$X_{\mathbb{R}}$, $Y_{\mathbb{R}}$, respectively, satisfying $f^
{-1}(D_{Y})\subseteq D_{X}$. Put $\underline{X}=(X_{\mathbb{R}},
D_{X})$ and $\underline{Y}=(Y_{\mathbb{R}},D_{Y})$. Then, there
is an inverse image morphism
\begin{displaymath}
f^{\ast}:\cha^{\ast}(Y,\mathcal{D}_{\wlg,\underline{Y}})
\longrightarrow\cha^{\ast}(X,\mathcal{D}_{\wlg,\underline{X}}).
\end{displaymath}
Moreover, this morphism is a morphism of rings after tensoring
with $\mathbb{Q}$.
\hfill $\square$
\end{theorem}

\nnpar{Push-forward.}
We will state only the consequence of proposition
\ref{prop:pushforward-point}.

\begin{theorem}
\label{thm:pushforward-point}
Let $X$ be a $d$-dimensional projective arithmetic variety
over $A$, $\underline{X}=(X_{\mathbb{R}},D)$ as above, and
$\pi:X\longrightarrow\Spec(A)$ the structural morphism. Then,
there is a direct image morphism
\begin{displaymath}
\pi_{\ast}:\cha^{d+1}(X,\mathcal{D}_{\wlg})\longrightarrow
\cha^{1}(\Spec(A)).
\end{displaymath}
\end{theorem}
\begin{proof}
Since $\Spec(A)_{\mathbb{R}}$ is a finite collection of points,
the complex $\mathcal{D}_{\log,\Spec(A)}$ is constant. By
proposition \ref{prop:pushforward-point}, there is a covariant
morphism of arithmetic varieties
\begin{displaymath}
\pi:(X,\mathcal{D}_{\wlg})\longrightarrow(\Spec(A),\mathcal{D}_
{\log}).
\end{displaymath}
The claim now follows from theorem \ref{thm:18}.
\end{proof}

\nnpar{Relationship with other arithmetic Chow groups.}
Let $X$ be an arithmetic variety over $A$, and $\underline{X}=
(X_{\mathbb{R}},D)$ as above. The structural morphism
\begin{displaymath}
\mathcal{D}_{\log}\longrightarrow\mathcal{D}_{\wlg}
\end{displaymath}
induces a morphism of arithmetic Chow groups
\begin{displaymath}
\cha^{\ast}(X,\mathcal{D}_{\log})\longrightarrow\cha^{\ast}
(X,\mathcal{D}_{\wlg}),
\end{displaymath}
which is compatible with inverse images, intersection products
and arithmetic degrees. If $X$ is projective, the isomorphism
between $\cha^{\ast}(X,\mathcal{D}_{\log})$ and the arithmetic
Chow groups defined by Gillet and Soul\'e (denoted by $\cha^
{\ast}(X)$) induces a morphism
\begin{equation}
\label{eq:34}
\cha^{\ast}(X)\longrightarrow\cha^{\ast}(X,\mathcal{D}_{\wlg}),
\end{equation}
which is also compatible with inverse images, intersection
products and arithmetic degrees. Since we have not proven the
purity of the cohomology of $\mathcal{D}_{\wlg}$, we cannot
state that the above morphism is an isomorphism, when $D$
is trivial, except when the dimension of $X_{F}$ is zero (see
remarks \ref{rem:2} and \ref{rem:3}). By contrast, if we impose growth 

\nnpar{Arithmetic Picard group.}
The theory of line bundles equipped with smooth hermitian metrics
developed in section \ref{sec:arithm-chow-groups} can be generalized
to the case of good metrics. For this let $X$ be a projective arithmetic
variety over $A$, $\underline{X}=(X_{\mathbb{R}},D)$ as above, and
$L$ a line bundle on $X$ equipped with a hermitian metric $h$ on the
induced line bundle $L_{\infty}$ over $X_{\infty}$, which is good along
$D$ and invariant under $F_{\infty}$. As usual, we write $\overline{L}=
(L,\|\cdot\|)$, and refer to it as a \emph{good hermitian line bundle}.
Given a rational section $s$ of $L$, we write $\|s\|^{2}=h(s,s)$ for
the point-wise norm of the induced section of $L_{\infty}$. We say
that two good hermitian line bundles $\overline{L}_{1}$ and $\overline
{L}_{2}$ are isometric, if $\overline{L}_{1}\otimes\overline{L}_{2}^
{-1}\cong(\mathcal{O}_{X},|\cdot|)$, where $|\cdot|$ is the standard
absolute value. The \emph{arithmetic Picard  group $\pica(\underline
{X})$} is the group of isometry classes of good hermitian line bundles
with group structure given by the tensor product.

\begin{proposition}
Let $X$ be a projective arithmetic variety over $A$, $\underline
{X}=(X_{\mathbb{R}},D)$ as above, and $\overline{L}=(L,\|\cdot\|)$
a good hermitian line bundle on $X$. Then, there is a map
\begin{displaymath}
\ca_{1}:\pica(\underline{X})\longrightarrow\cha^{1}(X,\mathcal
{D}_{\wlg}),
\end{displaymath}
given by sending the class of $\overline{L}$ to the class $[\dv(s),
(\omega_{s},\widetilde{g}_{s})]\in\cha^{1}(X,\mathcal{D}_{\wlg})$,
where $s$, $\omega_{s}$, and $g_{s}$ are as in proposition
\ref{prop:29}.
\end{proposition}
\begin{proof}
We have to show that the map $\ca_{1}$ is well defined. By
proposition \ref{prop:29}, the pair $(\omega_{s},g_{s})$
represents a Green object for the class of $\dv(s)$. If $s'$
is another non-vanishing, rational section, then $s'=f\cdot s$
with a rational function $f$, and we have
\begin{displaymath}
(\dv(s'),(\omega_{s'},\widetilde{g}_{s'}))=(\dv(s),(\omega_{s},
\widetilde{g}_{s}))+\diva(f).
\end{displaymath}
Hence, $(\dv(s),(\omega_{s},\widetilde{g}_{s}))$ and $(\dv(s'),
(\omega_{s'},\widetilde{g}_{s'}))$ represent the same class
in the arithmetic Chow group $\cha^{1}(X,\mathcal{D}_{\wlg})$.
Therefore, the map $\ca_{1}$ is well defined.
\end{proof}

\begin{definition}
Let $X$ be a projective arithmetic variety over $A$, $\underline{X}=
(X_{\mathbb{R}},D)$ as above, and $\overline{L}$ a good hermitian
line bundle on $X$. Then $\ca_{1}(\overline{L})\in\cha^{1}(X,\mathcal
{D}_{\wlg})$ is called the \emph{first arithmetic Chern class of
$\overline{L}$}.
\end{definition}

\nnpar{Arithmetic degree.}
Let $K$ be a number field, $\mathcal{O}_{K}$ its ring of integers, and
$S=\Spec(\mathcal{O}_{K})$. From section \ref{sec:arithm-chow-groups},
we recall the arithmetic degree map
\begin{displaymath}
\dega:\cha^{1}(S,\mathcal{D}_{\log})\longrightarrow\mathbb{R}.
\end{displaymath}

\begin{definition}
Let $X$ be a $d$-dimensional projective arithmetic variety over
$\mathcal{O}_{K}$ with structural morphism $\pi:X\longrightarrow S$,
$\underline{X}= (X_{\mathbb{R}},D)$ as above, and $\overline{L}$
a good hermitian line bundle on $X$. The real number
\begin{displaymath}
\overline{L}^{d+1}=\dega\left(\pi_{\ast}\ca_{1}(\overline{L})^
{d+1}\right)
\end{displaymath}
is called the \emph{arithmetic degree of $\overline{L}$}, or
the \emph{arithmetic self-intersection number of $\overline{L}$}.
\end{definition}

\begin{remark}
Logarithmically singular line bundles on arithmetic surfaces
as defined in \cite{Kuehn:gainc} are good hermitian line
bundles. It is straightforward to show that the intersection
number at the infinite places given in \cite{Kuehn:gainc}
(see \cite{Kuehn:gainc}, lemma 3.9), agrees with the explicit
formula for the $*$-product given in theorem \ref{STAR-PRODUCT}.
We also note that by \cite{Kuehn:gainc}, proposition 7.4, the
arithmetic self-intersection number of a good hermitian line
bundle calculated by either of the formulas in \cite{Bost:Lfg},
\cite{Moriwaki:ipacdGc}, \cite{Kuehn:gainc}, or theorem
\ref{STAR-PRODUCT} agree.
\end{remark}

\nnpar{Height pairings.} Let $K$ be a number field, $\mathcal{O}_
{K}$ its ring of integers, and $S=\Spec(\mathcal{O}_{K})$. Let 
$X$ be a $d$-dimensional projective arithmetic variety over 
$\mathcal{O}_{K}$ with structural morphism $\pi:X\longrightarrow 
S$. Furthermore, fix a normal crossing divisor $D$ of $X_{\mathbb
{R}}$.

We want to generalize the height pairing given in definition 
\ref{def:first-height-pairing} to the arithmetic Chow groups 
with values in $\mathcal{D}_{\wlg}$. This generalization will 
in particular include the logarithmic heights for points considered 
by Faltings in \cite{Faltings:EaVZ}, \cite{Faltings86:ftavnf}. 
Since the height of a cycle, whose generic part is supported in 
$D$, might be infinite, one cannot expect that the height pairing 
\eqref{eq:can-height} unconditionally generalizes to a height 
pairing between the arithmetic Chow groups $\cha^{p}(X,\mathcal
{D}_{\wlg})$ and the \emph{whole} group of cycles ${\rm Z}^{q}(X)$. 
Therefore, we restrict ourselves to considering a subgroup of cycles, 
for which a height pairing with respect to $\cha^{p}(X,\mathcal
{D}_{\wlg})$ can always be defined.

Putting $U=X_{\mathbb{R}}\setminus D$, we write ${\rm Z}^{q}_{U}
(X)$ for the group of the $q$-codimensional cycles $z$ of $X$ such 
that $z_{\mathbb{R}}$ intersects $D$ properly. A case of particular 
interest is given when the normal crossing divisor $D$ is defined 
over $K$; we then introduce $U_{K}=X_{K}\setminus D_{K}$, and
observe that there is a natural injective map
\begin{displaymath}
{\rm Z}^{q}(U_{K})\longrightarrow{\rm Z}^{q}_{U}(X).
\end{displaymath}

\begin{definition}
\label{defi:sing-height}
We want to define a bilinear pairing
\begin{equation}
\label{eq:sing-height}
(\cdot\mid\cdot):\,\cha^{p}(X,\mathcal{D}_{\wlg})\otimes{\rm Z}^{q}_
{U}(X)\longrightarrow\cha^{p+q-d}(S,\mathcal{D}_{\log})_{\mathbb{Q}}.
\end{equation}
Let $z\in{\rm Z}^{q}_{U}(X)$ be an irreducible, reduced cycle and 
$\alpha\in\cha^{p}(X,\mathcal{D}_{\wlg})$. We represent $\alpha$ 
by the class of an arithmetic cycle $(y,\mathfrak{g}_y)$, where 
$y$ is a $p$-codimensional cycle such that $y_{K}$ intersects $z_
{K}$ properly and where $\mathfrak{g}_{y}=(\omega_{y},\widetilde
{g_{y}})$ is a pre-log-log Green object for $y$. If $p+q<d$, we 
put $(\alpha\mid z)=0$. If $p+q=d$, we define
\begin{displaymath}
(\alpha\mid z)=[\pi_{\ast}(y_{K}\cdot z_{K}),(0,0)]\in\cha^{0}
(S,\mathcal{D}_{\log})_{\mathbb{Q}}=\CH^{0}(S)_{\mathbb{Q}}.
\end{displaymath}
Finally, if $p+q=d+1$, we define
\begin{displaymath}
(\alpha\mid z)=\left[\pi_{\ast}([y\cdot z]_{\fin}),\left(0,
\widetilde{\pi_{\#}(g_{y}\wedge\delta_{z})}\right)\right]\in
\cha^{1}(S,\mathcal{D}_{\log})_{\mathbb{Q}}.
\end{displaymath}
Here, the quantity $\widetilde{\pi_{\#}(g_{y}\wedge\delta_{z})}$ 
has to be understood as follows: Let $Z=\supp z_{\mathbb{R}}$ 
and $\imath:\widetilde{Z}\longrightarrow Z$ be a resolution of 
singularities of $Z$ adapted to $D$. Since $y_{K}\cap z_{K}=
\emptyset$, the functoriality of pre-log-log forms shows that
$\imath^{\ast}(g_{y})$ is a pre-log-log form on $\widetilde{Z}$, 
hence it is locally integrable on $\widetilde{Z}$, and we have
\begin{displaymath}
\widetilde{\pi_{\#}(g_{y}\wedge\delta_{z})}=\frac{1}{(2\pi i)^
{p-1}}\int_{\widetilde{Z}}\imath^{\ast}(g_{y}).
\end{displaymath}
The pairing \eqref{eq:sing-height} is now obtained by linearly 
extending the above definitions.
\end{definition}

Note that this definition a priori depends on the choice of 
the representative $(y,\mathfrak{g}_{y})$ of $\alpha$.

\begin{proposition}
\label{prop:log-sing-height-pairing}
With the above notations and assumptions, let $\alpha\in\cha^{p}
(X,\mathcal{D}_{\wlg})$, $z\in{\rm Z}_{U}^{q} (X)$, $Z=\supp z_
{\mathbb{R}}$, and $p+q=d+1$. If we choose a basic pre-log-log
Green form $g_{z}$ for $z$ and put $\mathfrak{g}_{z}=(-2\partial
\bar{\partial}g_{z},\widetilde{g_z})$, the height pairing 
\eqref{eq:sing-height} satisfies the equality 
\begin{equation}
\label{eq:log-sing-heightpairing}
(\alpha\mid z)=\pi_{\ast}\left(\alpha\cdot[z,\mathfrak{g}_{z}]
\right)+\amap\left(\pi_{\#}\left(\widetilde{[\omega(\alpha)
\wedge g_{z}]_{X}}\right)\right)
\end{equation}
in $\cha^{1}(S,\mathcal{D}_{\log})_{\mathbb{Q}}$. 
\end{proposition}
\begin{proof}  
By lemma \ref{lemm:integrability_s_p}, the pre-log form $\omega
(\alpha)\wedge g_{z}$ is integrable. Therefore, the right-hand 
side of the claimed formula is well defined. We now represent 
$\alpha$ by the arithmetic cycle $(y,\mathfrak{g}_y)$, where 
$y$ is a $p$-codimensional cycle such that $y_{K}$ intersects 
$z_{K}$ properly and where $\mathfrak{g}_{y}=(\omega(\alpha),
\widetilde{g_{y}})$ is a pre-log-log Green object for the cycle 
$y$. Then, the first summand of the right-hand side of equation
\eqref{eq:log-sing-heightpairing} equals 
\begin{displaymath}
\pi_{\ast}\left(\alpha\cdot[z,\mathfrak{g}_{z}]\right)=\left
[\pi_{\ast}([y\cdot z]_{\fin}),\left(0,\frac{1}{(2\pi i)^{d}}
\int_{X}\mathfrak{g}_{y}\ast\mathfrak{g}_{z}\right)\right].
\end{displaymath}
Using theorem \ref{STAR-PRODUCT} (with the roles of $y$ and $z$ 
interchanged), i.e., the equality
\begin{displaymath}
\frac{1}{(2\pi i)^{d}}\int_{X}\mathfrak{g}_{y}\ast\mathfrak{g}_
{z}=\frac{1}{(2\pi i)^{p-1}}\int_{Z}g_{y}+\frac{1}{(2\pi i)^{d}}
\int_{X}\omega(\alpha)\wedge g_{z}\,,
\end{displaymath}
we find
\begin{align*}
\pi_{\ast}\left(\alpha\cdot[z,\mathfrak{g}_{z}]\right)+&
\amap\left(\pi_{\#}\left(\widetilde{[\omega(\alpha)\wedge 
g_{z}]_{X}}\right)\right) \\
=&\left[\pi_{\ast}([y\cdot z]_{\fin}),\left(0,\widetilde{\pi_{\#}
(g_{y}\wedge\delta_{z})}+\pi_{\#}\left(\widetilde{[\omega(\alpha)
\wedge g_{z}]_{X}}\right)\right)\right] \\
&+\amap\left(\pi_{\#}\left(\widetilde{[\omega(\alpha)\wedge 
g_{z}]_{X}}\right)\right) \\
=&\left[\pi_{\ast}([y\cdot z]_{\fin}),\left(0,\widetilde{\pi_{\#}
(g_{y}\wedge\delta_{z})}\right)\right] \\
=&(\alpha\mid z),
\end{align*}
which proves the claim.
\end{proof}

In the following theorem, we will show that the pairing 
\eqref{eq:sing-height} is well defined.

\begin{theorem}
\label{thm:log-sing-height-pairing}
The bilinear pairing \eqref{eq:sing-height} is well defined.
Furthermore, there is a commutative diagram
\begin{displaymath}
\begin{CD}
\cha^{p}(X,\mathcal{D}_{\log})\otimes{\rm Z}^{q}_{U}(X)
@>(\cdot\mid\cdot)>>\cha^{p+q-d}(S,\mathcal{D}_{\log})_
{\mathbb{Q}} \\
@VVV@VV\Id V \\
\cha^{p}(X,\mathcal{D}_{\wlg})\otimes{\rm Z}^{q}_{U}(X) 
@>(\cdot\mid\cdot)>>\cha^{p+q-d}(S,\mathcal{D}_{\log})_
{\mathbb{Q}},
\end{CD}
\end{displaymath}
where horizontal maps are given by the pairings \eqref{eq:can-height} 
and \eqref{eq:sing-height}, respectively.
\end{theorem}
\begin{proof}
We have to show that the height pairing $(\alpha\mid z)$ does not
depend on the choice of a representative $(y,\mathfrak{g}_y)$ for
$\alpha$. Since the height pairing \eqref{eq:sing-height} coincides
with the height pairing \eqref{eq:can-height} for $p$, $q$ satisfying
$p+q\le d$, it suffices to treat the case $p+q=d+1$. In order to prove
the well-definedness in the latter case, we recall that the height 
pairing satisfies formula \eqref{eq:log-sing-heightpairing}. By the 
well-definedness of the arithmetic intersection product and the map 
$\omega$, the right-hand side of \eqref{eq:log-sing-heightpairing} 
turns out to be independent of the representative $(y,\mathfrak{g}_
{y})$,which proves the well-definedness of $(\alpha\mid z)$.

The commutativity of the  diagram follows from lemma 
\ref{def:bgs-height} observing the formula
\begin{displaymath}
\amap\left([g_{z}]_{X}(\omega(\alpha))\right)=\amap\left(\pi_{\#}
\left(\widetilde{[\omega(\alpha)\wedge g_{z}]_{X}}\right)\right).
\end{displaymath}
\end{proof}

\nnpar{Faltings heights.}  
We keep the notations of the previous section. In particular, we
let $X$, $D$, $U$ be as before and $p$, $q$ integers satisfying
$p+q=d+1$. The height pairing given by definition \ref{defi:sing-height}
is of particular interest, when $\alpha=\ca_{1}(\overline{L})^{p}$ 
for some good hermitian line bundle $\overline{L}$ on $X$. If $z\in
{\rm Z}^{q}_{U}(X)$, we call the real number 
\begin{displaymath}
\operatorname{ht}_{\overline{L}}(z)=\dega\left(\ca_{1}(\overline{L})^
{p}\mid z\right)
\end{displaymath}
the \emph{Faltings height of $z$ (with respect to $\overline{L}$)}. 
This name is justified by the following result due to Faltings.

\begin{proposition}
In addition to the data fixed above, assume that $L$ is an ample 
line bundle. Then, for any $c\in\mathbb{R}$, the cardinality of 
the set
\begin{displaymath}
\left\{z\in{\rm Z}^{d}_{U}(X)\,\big|\,\operatorname{ht}_{\overline
{L}}(z)<c\right\}
\end{displaymath}
is finite.
\end{proposition}
\begin{proof}
Since the metric of a good hermitian line bundle is logarithmically
singular as defined by Faltings, the proof is standard (see, e.g.,
\cite{Silverman86:ttohf}).
\end{proof}

\subsection{Application to products of modular curves}
As an example, we show in this final section how the arithmetic
intersection theory with pre-log-log forms developed in this
chapter can be applied to compute the arithmetic self-intersection
number of the line bundle of modular forms equipped with the natural
invariant metric on the product of two modular curves. We point
out that because of the singularities of the Petersson metric,
these numbers are only well defined in the setting of our newly
developed extended arithmetic intersection theory. 
Related but more elaborate results for the case of  Hilbert modular
surfaces have recently
been obtained in \cite{BruinierBurgosKuehn:bpaihs}.
Further calculations  of other
naturally metrized automorphic line bundles have been carried out in
\cite{BruinierKuehn:iaghd}.
The general theory of arithmetic characteristic classes of automorphic
vector bundles of arbitrary rank is developed in
\cite{BurgosKramerKuehn:accavb}.

\nnpar{Modular curves.}
\label{sss:modularcurves}
Let $\mathbb{H}$ denote the upper half plane with complex coordinate
$z=x+iy$, and
\begin{displaymath}
X(1)={\rm SL}_{2}(\mathbb{Z})\backslash\mathbb{H}\cup\{S_{\infty}\}
\end{displaymath}
the modular curve with the cusp $S_{\infty}$. For more details on
the subsequent facts we refer to \cite{Kuehn:gainc}.

Let $\mathcal{X}(1)=\mathbb{P}^{1}_{\mathbb{Z}}$ be the regular
model for the modular curve $X(1)$. With $s_{\infty}$ denoting the
Zariski closure of (the normal crossing divisor) $S_{\infty}\subset
X(1)$ and $k$ a positive integer satisfying $12|k$, we define the
{\em line bundle of modular forms of weight $k$} by $\mathcal{M}_
{k}=\mathcal{O}(s_{\infty})^{\otimes k/12}$. The line bundle
$\mathcal{M}_{k}$ is equipped with the Petersson metric $\|\cdot
\|$, which is a good hermitian metric along $S_{\infty}$; hence,
$\overline{\mathcal{M}}_{k}=(\mathcal{M}_{k},\|\cdot\|)$ is a
good hermitian line bundle, and we have
\begin{displaymath}
\ca_{1}(\overline{\mathcal{M}}_{k})\in\cha^{1}(\mathcal{X}(1),
\mathcal{D}_{\wlg}).
\end{displaymath}
For the first Chern form ${\rm c}_{1}(\overline{\mathcal{M}}_{k})$
of $\overline{\mathcal{M}}_{k}$, we recall the formulas
\begin{displaymath}
{\rm c}_{1}(\overline{\mathcal{M}}_{k})=k\cdot\frac{\dd z\wedge
\dd\bar{z}}{(z-\bar{z})^{2}}=-\frac{k}{2i}\cdot\frac{\dd x\wedge
\dd y}{y^{2}}\,,\quad\int_{X(1)}{\rm c}_{1}(\overline{\mathcal{M}}_
{k})=-\frac{\pi k}{6i}\,.
\end{displaymath}
Fixing a $(1,1)$-form $\omega$ satisfying $\int_{X(1)}\omega=-2\pi
i$, the main result of \cite{Kuehn:gainc} shows
\begin{equation}
\label{shit1}
\ca_{1}(\overline{\mathcal{M}}_{k})^{2}=[0,\overline{\mathcal{M}}_
{k}^{2}\cdot\amap(\omega)]\in\cha^{2}(\mathcal{X}(1),\mathcal{D}_
{\wlg})_\mathbb{Q},
\end{equation}
where the arithmetic self-intersection number is given by
\begin{displaymath}
\overline{\mathcal{M}}_{k}^{2}=k^{2}\cdot\left(\frac{1}{2}\zeta_
{\mathbb{Q}}(-1)+\zeta_{\mathbb{Q}}'(-1)\right)
\end{displaymath}
with the Riemann zeta function $\zeta_{\mathbb{Q}}(s)$.

\nnpar{Products of modular curves.}
In this example we consider the arithmetic threefold $\mathcal{H}=
\mathcal{X}(1)\times_{\mathbb{Z}}\mathcal{X}(1)$; we let $p_{1}$,
resp. $p_{2}$ denote the projection onto the first, resp. second
factor. The divisor
\begin{displaymath}
D=p_{1}^{\ast}\,\mathcal{X}(1)\times p_{2}^{\ast}\,s_{\infty}+p_
{1}^{\ast}\,s_{\infty}\times p_{2}^{\ast}\,\mathcal{X}(1)
\end{displaymath}
induces a normal crossing divisor $D_{\mathbb{R}}$ on $\mathcal{H}_
{\mathbb{R}}$. For $k,l\in\mathbb{N}$, $12|k$, $12|l$, we define
the hermitian line bundle
\begin{displaymath}
\overline{\mathcal{L}}(k,l)=p_{1}^{\ast}\,\overline{\mathcal{M}}_{k}
\otimes p_{2}^{\ast}\,\overline{\mathcal{M}}_{l}.
\end{displaymath}
It can be easily checked that the hermitian metric of $\overline
{\mathcal{L}}(k,l)$ is good along $D_{\mathbb{R}}$.

\begin{theorem}
Let $\mathcal{H}$ and $\overline{\mathcal{L}}(k,l)$ be as above.
Then, the arithmetic self-intersection number of $\overline
{\mathcal{L}}(k,l)$ is well defined, and given by the formula
\begin{align*}
\overline{\mathcal{L}}(k,l)^{3}=\frac{k^{2}\cdot l+l^{2}\cdot k}
{4}\left(\frac{1}{2}\zeta_{\mathbb{Q}}(-1)+\zeta_{\mathbb{Q}}'(-1)
\right).
\end{align*}
\end{theorem}
\begin{proof}
By theorem \ref{thm:wll-inverse-image}, we have
\begin{displaymath}
\ca_{1}(\overline{\mathcal{L}}(k,l))\in\cha^{1}(\mathcal{H},
\mathcal{D}_{\wlg});
\end{displaymath}
this can also be checked directly. Letting $\pi:\mathcal{H}
\longrightarrow\Spec(\mathbb{Z})$ denote the structural morphism,
the arithmetic self-intersection number of $\overline{\mathcal{L}}
(k,l)$ is given by
\begin{displaymath}
\overline{\mathcal{L}}(k,l)^{3}=\widehat{\deg}\left(\pi_{\ast}
\ca_{1}(\overline{\mathcal{L}}(k,l))^{3}\right),
\end{displaymath}
which is well defined by theorem \ref{thm:pushforward-point}.
Using the properties of arithmetic Chow groups, we derive the
following identities in $\cha^{3}(\mathcal{H},\mathcal{D}_
{\wlg})_{\mathbb{Q}}$:
\begin{align*}
\ca_{1}(\overline{\mathcal{L}}(k,l))^{3}&=\left(p_{1}^{\ast}
\,\ca_{1}(\overline{\mathcal{M}}_{k})+p_{2}^{\ast}\,\ca_{1}
(\overline{\mathcal{M}}_{l})\right)^{3} \\
&=\sum_{j=0}^{3}{3\choose j}p_{1}^{\ast}\,\ca_{1}(\overline
{\mathcal{M}}_{k})^{j}\cdot p_{2}^{\ast}\,\ca_{1}(\overline
{\mathcal{M}}_{l})^{3-j} \\
&=\sum_{j=1}^{2}{3\choose j}p_{1}^{\ast}\,\ca_{1}(\overline
{\mathcal{M}}_{k})^{j}\cdot p_{2}^{\ast}\,\ca_{1}(\overline
{\mathcal{M}}_{l})^{3-j}\,.
\end{align*}
By means of formula (\ref{shit1}), we find
\begin{align*}
\ca_{1}(\overline{\mathcal{L}}(k,l))^{3}&={3\choose 2}\cdot
\bigg(\left[0,\overline{\mathcal{M}}_{k}^{2}\cdot\amap(p_{1}^
{\ast}\,\omega)\right]\cdot p_{2}^{\ast}\,\ca_{1}(\overline
{\mathcal{M}}_{l}) \\
&+p_{1}^{\ast}\,\ca_{1}(\overline{\mathcal{M}}_{k})\cdot
\left[0,\overline{\mathcal{M}}_{l}^{2}\cdot\amap(p_{2}^
{\ast}\,\omega)\right]\bigg).
\end{align*}
Applying proposition \ref{prop:explicit}, we obtain
\begin{align*}
\ca_{1}(\overline{\mathcal{L}}(k,l))^{3}&=3\cdot\bigg(\left
[0,\overline{\mathcal{M}}_{k}^{2}\cdot\amap\left(p_{1}^{\ast}
\,\omega\wedge p_{2}^{\ast}\,{\rm c}_{1}(\overline{\mathcal{M}}_
{l})\right)\right] \\
&+\left[0,\overline{\mathcal{M}}_{l}^{2}\cdot\amap\left(p_{1}^
{\ast}\,{\rm c}_{1}(\overline{\mathcal{M}}_{k})\wedge p_{2}^{\ast}
\,\omega\right)\right]\bigg).
\end{align*}
Taking the push-forward of the latter quantity by $\pi$ and
applying the arithmetic degree map, we finally obtain
\begin{align*}
\overline{\mathcal{L}}(k,l)^{3}&=\frac{3}{(2\pi i)^{2}}\left
(\frac{1}{2}\zeta_{\mathbb{Q}}(-1)+\zeta_{\mathbb{Q}}'(-1)
\right)\cdot \\
&\cdot\left(k^{2}\int_{X(1)}(-\omega)\cdot\int_{X(1)}{\rm c}_{1}
(\overline{\mathcal{M}}_{l})+l^{2}\int_{X(1)}{\rm c}_{1}(\overline
{\mathcal{M}}_{k})\cdot\int_{X(1)}(-\omega)\right),
\end{align*}
from which the claimed formula follows.
\end{proof}

\nnpar{Faltings heights of Hecke correspondences.}
Let $N$ be a positive integer, and $M_{N}$ the set of integral
$(2\times 2)$-matrices of determinant $N$. Recall that the group
${\rm SL}_{2}(\mathbb{Z})$ acts form the right on the set $M_{N}$
and that a complete set of representatives for this action is
given by the set
\begin{displaymath}
R_{N}=\left\{\gamma=\begin{pmatrix}a&b\\0&d\end{pmatrix}\,\bigg|\,
a,b,d\in\mathbb{Z};\,ad=N;\,d>0;\,0\le b<d\right\}.
\end{displaymath}
The cardinality $\sigma(N)$ of $R_{N}$ is given by
\begin{displaymath}
\sigma(N)=\sum_{d|N}d\,,
\end{displaymath}
where the sum is taken over all positive divisors $d$ of $N$.
Furthermore, we put
\begin{displaymath}
T_{N}=\left\{(z_{1},z_{2})\in\mathbb{H}\times\mathbb{H}\,\big|
\,\exists\gamma\in M_{N}:\,z_{1}=\gamma z_{2}\right\}\,.
\end{displaymath}
This defines a divisor on $X(1)\times X(1)$, whose normalization
is a finite sum of modular curves $X_{0}(m)$ associated to the
congruence subgroups $\Gamma_{0}(m)$ for $m$ dividing $N$. These
curves are also referred to as the graphs of the Hecke correspondences.
Denoting the discriminant by $\Delta(z)$ and the $j$-function by $j(z)$,
the Hilbert modular form
\begin{displaymath}
s_{N}(z_{1},z_{2})=\Delta(z_{1})^{\sigma(N)}\Delta(z_{2})^{\sigma(N)}
\prod_{\gamma\in R_{N}}\left(j(\gamma z_{1})-j(z_{2})\right)
\end{displaymath}
has divisor $T_{N}$. Since the Fourier coefficients of $s_{N}(z_{1},
z_{2})$ are integral, it defines a section of $\mathcal{L}(12\sigma
(N),12\sigma(N))$; we put
\begin{displaymath}
\mathcal{T}_{N}=\dv(s_{N})\subseteq\mathcal{H}.
\end{displaymath}

\begin{theorem}
The Faltings height of the Hecke correspondence $\mathcal{T}_{N}$
with respect to $\overline{\mathcal{L}}(k,k)$ is given by the
formula
\begin{align*}
\operatorname{ht}_{\overline{\mathcal{L}}(k,k)}(\mathcal{T}_{N})&=
(2k)^{2}\bigg(\sigma(N)\left(\frac{1}{2}\zeta_{\mathbb{Q}}(-1)+
\zeta_{\mathbb{Q}}'(-1)\right) \\
&+\sum_{d|N}\frac{d\log(d)}{24}-\frac{\sigma(N)\log(N)}{48}\bigg)\,.
\end{align*}
\end{theorem}
\begin{proof}
Using Rohrlich's modular Jensen's formula, we obtain for any fixed
$z_{0}\in\mathbb{H}$
\begin{align*}
&\int_{X(1)}\log\|s_{N}(z,z_{0})\|\frac{\dd x\wedge \dd y}
{4\pi y^{2}}= \\
&-12\sigma(N)\left(\frac{1}{2}\zeta_{\mathbb{Q}}(-1)+\zeta_{\mathbb
{Q}}'(-1)\right)+\frac{1}{12}\sum_{\gamma\in R_{N}}\log\left(\frac
{\|\Delta(z)\|}{\|\Delta(\gamma z)\|}\right).
\end{align*}
Since we have
\begin{displaymath}
\prod_{\gamma\in R_{N}}\frac{\Delta(z)}{\Delta(\gamma z)}=1\,,
\quad{\rm Im}\left(\frac{az+b}{d}\right)=\frac{N}{d^{2}}\cdot
{\rm Im}z\quad\left(\begin{pmatrix}a&b\\0&d\end{pmatrix}\in
R_{N}\right)\,,
\end{displaymath}
we obtain for the sum on the right hand side
\begin{align*}
\frac{1}{12}\sum_{\gamma\in R_{N}}\log\left(\frac{\|\Delta(z)\|}
{\|\Delta(\gamma z)\|}\right)&=\frac{1}{12}\log\bigg|\prod_{\gamma
\in R_{N}}\frac{\Delta(z)}{\Delta(\gamma z)}\bigg|+\frac{1}{2}
\sum_{\gamma\in R_{N}}\log\left(\frac{d^{2}}{N}\right) \\
&=\sum_{d|N}d\log(d)-\frac{1}{2}\sigma(N)\log(N).
\end{align*}
Writing $z_{1}=x_{1}+iy_{1}$, $z_{2}=x_{2}+iy_{2}$, and observing
\begin{displaymath}
{\rm c}_{1}(\overline{\mathcal{L}}(k,k))^{2}=-\frac{k^{2}}{2}\cdot
\frac{\dd x_{1}\wedge\dd y_{1}}{y_{1}^{2}}\cdot\frac{\dd x_{2}\wedge
\dd y_{2}}{y_{2}^{2}}\,,
\end{displaymath}
we find by means of Rohrlich's formula
\begin{align*}
&\frac{1}{(2\pi i)^{2}}\int_{X(1)\times X(1)}\log\|s_{N}(z_{1},
z_{2})\|\cdot{\rm c}_{1}(\overline{\mathcal{L}}(k,k))^{2}= \\
&k^{2}\bigg(-2\sigma(N)\left(\frac{1}{2}\zeta_{\mathbb{Q}}(-1)+
\zeta_{\mathbb{Q}}'(-1)\right)+ \\
&\frac{1}{6}\sum_{d|N}d\log(d)-\frac{1}{12}\sigma(N)\log(N)
\bigg)\,.
\end{align*}
Since $s_{N}$ defines a section of $\mathcal{L}(12\sigma(N),12\sigma
(N))$ and $\log\|s_{N}(z_{1},z_{2})\|$ is a basic pre-log-log form, 
we obtain by proposition \ref{prop:log-sing-height-pairing} and the 
previous considerations

\begin{align*}
\operatorname{ht}_{\overline{\mathcal{L}}(k,k)}(\mathcal{T}_
{N})=&\dega\left(\ca_{1}(\overline{\mathcal{L}}(k,k))^{2}
\mid\mathcal{T}_{N}\right) \\[3mm]
=&\dega\bigg(\pi_{\ast}\bigg(\ca_{1}(\overline{\mathcal{L}}
(k,k))^{2}\cdot \\[2mm]
&\cdot\left[\mathcal{T}_{N},\left({\rm c}_{1}(\overline
{\mathcal{L}}(12\sigma(N),12\sigma(N))),\widetilde{-\log
\|s_{N}(z_{1},z_{2})\|}\right)\right]\bigg) \\[2mm]
&+\amap\left(\pi_{\#}\left(\left[-\log\|s_{N}(z_{1},z_{2})\|
\wedge{\rm c}_{1}(\overline{\mathcal{L}}(k,k))^{2}\right]_
{X(1)\times X(1)}\right)\right)\bigg)
\end{align*}
\begin{align*}
\phantom{\operatorname{ht}_{\overline{\mathcal{L}}(k,k)}
(\mathcal{T}_{N})}
=&\dega\left(\pi_{\ast}\left(\ca_{1}(\overline{\mathcal{L}}
(k,k))^{2}\cdot\ca_{1}(\overline{\mathcal{L}}(12\sigma(N),
12\sigma(N)))\right)\right) \\[2mm]
&+\frac{1}{(2\pi i)^{2}}\int_{X(1)\times X(1)}\log\|s_{N}
(z_{1},z_{2})\|\cdot{\rm c}_{1}(\overline{\mathcal{L}}(k,k))^
{2} \\[2mm]
=&\,6k^{2}\sigma(N)\left(\frac{1}{2}\zeta_{\mathbb{Q}}(-1)+
\zeta_{\mathbb{Q}}'(-1)\right) \\[2mm]
&+k^{2}\bigg(-2\sigma(N)\left(\frac{1}{2}\zeta_{\mathbb{Q}}
(-1)+\zeta_{\mathbb{Q}}'(-1)\right) \\[2mm]
&+\frac{1}{6}\sum_{d|N}d\log(d)-\frac{1}{12}\sigma(N)\log(N)
\bigg)\,.
\end{align*}
\end{proof}

\begin{remark}
  We point out that the same formula has been obtained independently
  by P. Autissier \cite{Autissier:hhc}.  His calculations, which are
  based on results of U. K\"uhn \cite{Kuehn:gainc}  and use  J.-B.
  Bost's theory of $L^{2}_{1}$-singular metrics on arithmetic surfaces
  \cite{Bost:Lfg}, are compatible with definition
  \eqref{eq:sing-height}.
We also note that the arithmetic
  intersection numbers obtained via Bost's $L^{2}_{1}$-theory applied
  to logarithmically singular metrics equal those in
  \cite{Kuehn:gainc} and coincide with those obtained
by applying the theory developed here in the
  2-dimensional setting.
\end{remark}

\nnpar{Arithmetic generating series.}
We put
\begin{align*}
\widehat{T(N)}=
\begin{cases}
0\,,&\quad\textrm{if }N<0, \\
[\mathcal{T}_{N},\mathfrak{g}_{N}]\,,&\quad\textrm{if }N>0,
\end{cases}
\end{align*}
where $\mathfrak{g}_{N}(z)=\left(2\partial\bar{\partial}\log
\|s_{N}(z_{1},z_{2})\|,\widetilde{-\log\|s_{N}(z_{1},z_{2})
\|}\right)$. With the divisor $D\subseteq\mathcal{H}$, we
additionally put
\begin{displaymath}
\widehat{T(0)}=-\frac{1}{24}\ca_{1}\left(\overline{\mathcal
{L}}(12,12)\right)+\frac{1}{8\pi y}[D,\mathfrak{g}_{D}(z)]
\in\cha^{1}(\mathcal{H},\mathcal{D}_{\wlg})_{\mathbb{R}},
\end{displaymath}
where $\mathfrak{g}_{D}(z)=\left(2\partial\bar{\partial}\log
\|\Delta(z_{1})\Delta(z_{2})\|,\widetilde{-\log\|\Delta(z_{1})
\Delta(z_{2})\|}\right)$.
With these ingredients we define the following generating
series of arithmetic cycles
\begin{displaymath}
\widehat{\phi}_{\rm{can}}(q)=\sum_{N\in\mathbb{Z}}\widehat
{T(N)}\,q^{N}\in\cha^{1}(\mathcal{H},\mathcal{D}_{\wlg})_
{\mathbb{R}}\otimes\mathbb{C}((q))\,;
\end{displaymath}
here $q=e^{2\pi iz}$.

\begin{theorem}
\label{thm:naivekudla}
The generating series $\widehat{\phi}_{{\rm can}}$ is a
non-holomorphic modular form with values in $\cha^{1}
(\mathcal{H},\mathcal{D}_{\wlg})_{\mathbb{R}}$ of weight
$2$ for ${\rm SL}_{2}(\mathbb{Z})$.
\end{theorem}
\begin{proof} We observe that
\begin{align*}
\widehat{T(N)}=\ca_{1}\left(\overline{\mathcal{L}}(12,12)\right)
\cdot
\begin{cases}
0\,,&\textrm{if }N<0, \\
\left(-\frac{1}{24}+\frac{1}{8\pi y}\right)\,,&\textrm{if }N=0, \\
\sigma(N)\,,&\textrm{if }N>0.
\end{cases}
\end{align*}
Recall now that there is the non-holomorphic Eisenstein series
$E_{2}(z)$ of weight $2$ obtained by means of Hecke's convergence
trick from the non-holomorphic Eisenstein series $E_{2}(z,s)$ at
$s=1$; it has Fourier expansion
\begin{displaymath}
E_{2}(z,1)=-\frac{1}{24}+\frac{1}{8\pi y}+\sum_{N>0}\sigma(N)q^{N}.
\end{displaymath}
Therefore, we obtain
\begin{displaymath}
\widehat{\phi}_{{\rm can}}(q)=\ca_{1}\left(\overline{\mathcal
{L}}(12,12)\right)\otimes E_{2}(z)\in\cha^{1}(\mathcal{H},
\mathcal{D}_{\wlg})_{\mathbb{R}}\otimes\mathcal{A}_{2}({\rm
SL}_{2}(\mathbb{Z})),
\end{displaymath}
where $\mathcal{A}_{2}({\rm SL}_{2}(\mathbb{Z}))$ denotes the
space of non-holomorphic modular forms of weight $2$ with respect
to ${\rm SL}_{2}(\mathbb{Z})$.
\end{proof}

By theorem \ref{thm:naivekudla} any linear functional on $\cha^{1}
(\mathcal{H},\mathcal{D}_{\wlg})_{\mathbb{R}}$ applied to $\widehat
{\phi}_{{\rm can}}$ gives rise to a modular form of weight $2$. As
possible applications we present an integral version of Borcherds's
generating series
\begin{displaymath}
\zeta(\widehat{\phi}_{{\rm can}})=-\frac{1}{24}{\rm c}_{1}\left(
\mathcal{L}(12,12)\right)+\frac{1}{8\pi y}D+\sum_{N>0}\mathcal
{T}_{N}q^{N}\in\CH(\mathcal{H})_{\mathbb{R}}\otimes\mathcal{A}_
{2}({\rm SL}_{2}(\mathbb{Z})),
\end{displaymath}
or the Hirzebruch-Zagier formula, which states that for any
$(1,1)$-form $\eta$, we have
\begin{displaymath}
[\eta]_{\mathcal{H}_{\infty}}\left(\cl(\zeta(\widehat{\phi}_
{{\rm can}}))\right)=\sum_{N\ge 0}\left(\frac{1}{2\pi i}
\int_{T_{N}}\eta\right)\,q^{N}\in\mathcal{A}_{2}({\rm SL}_
{2}(\mathbb{Z})).
\end{displaymath}
In particular, we emphasize that choosing $\eta={\rm c}_{1}\left
(\overline{\mathcal{L}}(12,12)\right)$ yields a multiple of
the Eisenstein series $E_{2}(z,1)$.

We conclude by noting that according to Kudla's conjectures there
should  also exist a generating series related to the derivative
of $E_{2}(z,s)$ at $s=1$.

\bibliographystyle{amsplain}
\newcommand{\noopsort}[1]{} \newcommand{\printfirst}[2]{#1}
  \newcommand{\singleletter}[1]{#1} \newcommand{\switchargs}[2]{#2#1}
\providecommand{\bysame}{\leavevmode\hbox to3em{\hrulefill}\thinspace}
\providecommand{\MR}{\relax\ifhmode\unskip\space\fi MR }
\providecommand{\MRhref}[2]{%
  \href{http://www.ams.org/mathscinet-getitem?mr=#1}{#2}
}
\providecommand{\href}[2]{#2}

\end{document}